\documentclass[10pt,draft
]{article}

\usepackage[all]{xy}

\usepackage[latin1]{inputenc}
\usepackage{amsfonts}
\usepackage{amsmath}
\usepackage{amssymb}
\usepackage{amscd}
\usepackage{amsthm}
\usepackage{indentfirst}
\usepackage[hmargin=2.98cm
,vmargin=4.0cm
]{geometry}

\usepackage{epsfig}

\newtheorem{thm}{Theorem}[section]
\newtheorem{lemma}[thm]{Lemma}
\newtheorem{prop}[thm]{Proposition}
\newtheorem{coroll}[thm]{Corollary}

\theoremstyle{definition}

\newtheorem{rem}[thm]{Remark}

\newtheorem*{acknow}{Acknowledgements}

\newtheorem*{prf}{Proof}

\newcommand{\R}{{\mathbb{R}}}

\newcommand{\T}{{\mathbb{T}}}
\newcommand{\Z}{{\mathbb{Z}}}

\newcommand{\C}{{\mathbb{C}}}

\newcommand{\cA}{{\mathcal{A}}}

\newcommand{\cD}{{\mathcal{D}}}
\newcommand{\cE}{{\mathcal{E}}}
\newcommand{\cF}{{\mathcal{F}}}

\newcommand{\cJ}{{\mathcal{J}}}

\newcommand{\cL}{{\mathcal{L}}}
\newcommand{\cM}{{\mathcal{M}}}

\newcommand{\cO}{{\mathcal{O}}}
\newcommand{\cP}{{\mathcal{P}}}

\newcommand{\cR}{{\mathcal{R}}}
\newcommand{\cS}{{\mathcal{S}}}
\newcommand{\cT}{{\mathcal{T}}}
\newcommand{\cU}{{\mathcal{U}}}
\newcommand{\cV}{{\mathcal{V}}}

\newcommand{\cX}{{\mathcal{X}}}
\newcommand{\cY}{{\mathcal{Y}}}

\newcommand{\mft}{{\mathfrak{t}}}

\newcommand{\mfo}{{\mathfrak{o}}}

\newcommand{\fc}{{:\ }}

\newcommand{\ol}{\overline}
\newcommand{\wt}{\widetilde}
\newcommand{\wh}{\widehat}

\newcommand{\tb}{\textbf}

\DeclareMathOperator{\RE}{Re}

\DeclareMathOperator{\Hom}{Hom}

\DeclareMathOperator{\Crit}{Crit}

\DeclareMathOperator{\im}{im}
\DeclareMathOperator{\id}{id}

\DeclareMathOperator{\PSS}{PSS}

\DeclareMathOperator{\pr}{pr}

\DeclareMathOperator{\pt}{pt}

\DeclareMathOperator{\ev}{ev}

\DeclareMathOperator{\ind}{ind}

\DeclareMathOperator{\const}{const}

\DeclareMathOperator{\inward}{inward}

\DeclareMathOperator{\Lie}{Lie}

\DeclareMathOperator{\Aut}{Aut}

\DeclareMathOperator{\coker}{coker}

\DeclareMathOperator{\Pin}{Pin}

\DeclareMathOperator{\ddd}{d}

\DeclareMathOperator{\std}{std}

\title{The Lagrangian Floer--quantum--PSS package and canonical orientations in Floer theory}
\author{Frol Zapolsky \footnote{Department of Mathematics, Faculty of Natural Sciences, University of Haifa, Israel. \texttt{frol.zapolsky@gmail.com}} \footnote{Partially supported by grant number 1281 from the GIF, the German--Israeli Foundation for Scientific Research and Development, and by grant number 1825/14 from the Israel Science Foundation. Portions of this manuscript were written while I was a postdoc at Ludwig--Maximilian--Universit\"at, where I was supported by the DFG grant number DFG/CI45/5-1.} }

\begin{document}

\renewcommand{\labelenumi}{(\roman{enumi})}

\maketitle

\begin{abstract}
The purpose of this paper is to extend the construction of the PSS-type isomorphism between the Floer homology and the quantum homology of a monotone Lagrangian submanifold $L$ of a symplectic manifold $M$, provided that the minimal Maslov number of $L$ is at least two, to \emph{arbitrary coefficients}. We provide a proof, again over arbitrary coefficients, that this isomorphism respects the natural algebraic structures on both sides, such as the quantum product and the quantum module action. This isomorphism serves as the technical foundation for the construction of Lagrangian spectral invariants in \cite{Leclercq_Zapolsky_Spectral_invts_monotone_Lags}. Our constructions work when the second Stiefel--Whitney class of $L$ vanishes on the image of the boundary homomorphism $\pi_3(M,L) \to \pi_2(L)$, a condition strictly weaker than being relatively Pin; in particular we do not require $L$ to be orientable. The constructions are done using canonical orientations, and require no further choices such as relative Pin-structures. Such structures do however play a significant role when endowing the various complexes and homologies with structures of modules over Novikov rings, and in calculations.
\end{abstract}

\tableofcontents

\section{Introduction and main result}\label{s:intro}

\subsection{Introduction}

The Piunikhin--Salamon--Schwarz isomorphism between the Hamiltonian Floer homology and the quantum homology of a symplectic manifold $M$ was constructed in \cite{PSS}. Nowadays it has become a standard tool of symplectic topology and is known under the abbreviated name of the PSS isomorphism. The analogs of the PSS isomorphism for Lagrangian Floer homology were constructed at varying levels of generality and rigor in \cite{Katic_Milinkovic_PSS_Lagr_intersections, Cornea_Lalonde_Cluster_homology, Biran_Cornea_Quantum_structures_Lagr_submfds, Albers_Lagr_PSS_comparison_morphisms_HF, Leclercq_spectral_invariants_Lagr_FH, Biran_Cornea_Rigidity_uniruling, HLL11}. Applications of these isomorphisms include spectral invariants, see \cite{Oh_Construction_sp_invts_Ham_paths_closed_symp_mfds, Leclercq_spectral_invariants_Lagr_FH, Leclercq_Zapolsky_Spectral_invts_monotone_Lags} and references therein, as well as problems such as homological Lagrangian monodromy \cite{HLL11}.

Our contribution in the present paper is as follows:
\begin{itemize}
\item We extend the construction of the PSS isomorphism for a monotone Lagrangian submanifold $L \subset M$ with minimal Maslov number at least two, appearing in \cite{Biran_Cornea_Quantum_structures_Lagr_submfds, Biran_Cornea_Rigidity_uniruling}, to arbitrary coefficients, under a certain assumption on the second Stiefel--Whitney class of $L$.
\item We show that the natural algebraic structures on Lagrangian Floer and quantum homology, such as products and the quantum module action, are intertwined by this isomorphism. This has been known to experts (see \cite{Biran_Cornea_Quantum_structures_Lagr_submfds, Biran_Cornea_Rigidity_uniruling}), and here we provide details of the relevant construction, again with arbitrary coefficients.
\item We use and develop the approach of canonical orientations, appearing for instance in \cite{Welschinger_Open_strings_Lag_conductors_Floer_fctr, Seidel_The_Book, Abouzaid_Symp_H_Viterbo_thm}, to tackle the issue of signs, rather than using the usual approach of coherent orientations; this allows us to generalize the construction of the Lagrangian Floer and quantum complexes, as well as of the PSS isomorphism, and in particular our Lagrangian is not assumed to be orientable or to carry a relative Pin-structure.
\end{itemize}



The canonical Floer and quantum complexes that we construct in this paper distinguish homotopy classes of cappings. In applications it is often more desirable to have smaller complexes in which cappings are only distinguished by their area or homology class. Relative Pin-structures allow us to do that. They are also used in order to endow the canonical complexes and their homologies with module structures over Novikov rings.

The existing literature does not explicitly cover one particular case appearing when proving that the quantum boundary operator squares to zero, namely the case analogous to bubbling off of a holomorphic disk of Maslov index $2$ in Floer homology. Even though it is fairly straightforward, at least with coefficients in $\Z_2$ \cite{Paul_Octav_Private_comm_July_2015}, we include a treatment of this case with arbitrary coefficients in \S \ref{ss:boundary_op_squares_zero_bubbling_QH} for the sake of completeness.

We also provide a description of various auxiliary algebraic structures, which have already appeared in the literature, such as augmentations and duality \cite{Biran_Cornea_Quantum_structures_Lagr_submfds, Biran_Cornea_Rigidity_uniruling}, and spectral sequences \cite{Oh_HF_spectral_sequences_Maslov_class, Biran_Cornea_Quantum_structures_Lagr_submfds, Biran_Cornea_Rigidity_uniruling, Biran_Lag_non_intersections}. These require careful formulation due to the presence of orientations.


Lastly, we compute the canonical quantum complexes for the standard monotone $\R P^n \subset \C P^n$, $n \geq 1$, as well as certain monotone Lagrangian tori in $\C P^2$ and $\C P^1 \times \C P^1$.

\subsection{Main result}

Let us fix a closed \footnote{With straightforward modifications, everything said here can be formulated for open symplectic manifolds which are convex at infinity \cite{Eliashberg_Gromov_Convex_symp_mfds}.} connected symplectic manifold $(M,\omega)$ of dimension $2n$ and a closed connected Lagrangian submanifold $L$. There are two natural homomorphisms associated to $L$:
$$\omega \fc \pi_2(M,L) \to \R\,, \quad \text{the \tb{symplectic area}, and}$$
$$\mu \fc \pi_2(M,L) \to \Z\,, \quad\text{the \tb{Maslov class}.}$$
We say that $L$ is \tb{monotone} if there is a positive constant $\tau$ such that $\omega = \tau \mu$. We denote by $N_L$ the positive generator of $\mu(\pi_2(M,L)) \subset \Z$ if this subgroup is nonzero, otherwise we put $N_L = \infty$. This is called the \tb{minimal Maslov number} of $L$.

\medskip

\tb{Throughout this paper we assume that $L$ is monotone and that its minimal Maslov number is at least $2$.}

\medskip

Since we are dealing with Lagrangian Floer theory with arbitrary coefficients, we need to impose a condition on the second Stiefel--Whitney class of $L$.

\medskip

\tb{Assumption (O): the class $w_2(TL)$ vanishes on the image of the boundary homomorphism $\pi_3(M,L) \to \pi_2(L)$.}

\medskip

\noindent We impose assumption \tb{(O)} throughout.

\begin{rem}
If only coefficients in a ring of characteristic $2$ are required, this assumption is not needed. See also \S\ref{ss:arbitrary_rings_loc_coeffs}.
\end{rem}

\begin{rem}
Assumption \tb{(O)} is implied by and is strictly weaker than being relatively $\Pin^\pm$, see \S\ref{ss:relative_Pin_structs_coh_ors_disks}.
\end{rem}

Given a commutative ground ring $R$, we construct the Lagrangian Floer complex
$$(CF_*(H:L),\partial_{H,J})$$
of $L$ over $R$, relative to a Hamiltonian perturbation $H$, where $J$ is an $\omega$-compatible almost complex structure on $M$, and the Lagrangian quantum complex
$$(QC_*(\cD:L),\partial_\cD)\,,$$
where $\cD=(f,\rho,I)$ is a quantum datum, $(f,\rho)$ being a Morse--Smale pair on $L$, and $I$ another compatible almost complex structure. We prove that the homologies $HF_*(L)$ and $QH_*(L)$ of these complexes are independent of the auxiliary data. These carry the structure of unital associative algebras over $R$. Our main result in this paper is
\begin{thm}\label{thm:main_result_existence_PSS_iso}
There is a canonical PSS isomorphism
$$\PSS \fc HF_*(L) \to QH_*(L)\,,$$
which is a unital algebra isomorphism.
\end{thm}
\noindent We wish to emphasize that the Lagrangian Floer and quantum complexes, and the PSS isomorphism that we construct, are a generalization of existing constructions, and that the main points here are the use of canonical orientations, the generalization of the construction of the PSS isomorphism to arbitrary coefficients, and a proof of the fact that the natural algebraic structures are intertwined by it.

This theorem is proved as part of the constructions of \S\ref{s:PSS}. In addition the PSS isomorphism respects the so-called quantum module structures on both sides. Namely, $HF_*(L)$ is a superalgebra over the Hamiltonian Floer homology $HF_*(M)$ of $M$, and analogously $QH_*(L)$ is a superalgebra over the quantum homology $QH_*(M)$. The PSS isomorphism we construct intertwines the two structures.

\subsection{Overview of the construction}\label{ss:overview_construction}

Let us briefly review these constructions. We have the path space
$$\Omega_L = \{\gamma \fc [0,1] \to M\,|\, \gamma(0),\gamma(1) \in L\,,[\gamma] = 0 \in \pi_1(M,L)\}$$
and its covering space
$$\wt\Omega_L = \{[\gamma,\wh\gamma] \,|\, \gamma \in \Omega_L\,,\wh\gamma \text{ a capping of }\gamma\}\,,$$
where a capping of $\gamma$ is a map from the closed half-disk to $M$ where the diameter maps to $\gamma$ while the boundary maps to $L$. Two cappings $\wh\gamma,\wh\gamma'$ of the same path $\gamma$ are equivalent if their concatenation $\wh\gamma \sharp - \wh\gamma'$ is nullhomotopic relative to $L$, and $[\gamma,\wh\gamma]$ denotes the equivalence class of a capping. Given a Hamiltonian $H \fc M \times [0,1] \to \R$ its action functional is
$$\cA_{H:L} \fc \wt\Omega_L \to \R\,,\quad \cA_{H:L}([\gamma,\wh\gamma]) = \int_0^1 H_t(\gamma(t))\,dt - \int \wh\gamma^*\omega\,.$$
The critical points $\Crit \cA_{H:L}$ are those points $[\gamma,\wh\gamma]$ for which $\gamma$ is a Hamiltonian arc of $H$, that is $\dot \gamma = X_H\circ\gamma$. We call $H$ nondegenerate if for every such critical point the linearized map
$$\phi_{H,*} \fc T_{\gamma(0)}M \to T_{\gamma(1)}M$$
maps $T_{\gamma(0)}L$ to a subspace transverse to $T_{\gamma(1)}L$. For such $H$ the underlying module of its Floer complex is defined as
$$CF_*(H:L) = \bigoplus_{\wt\gamma \in \Crit \cA_{H:L}}C(\wt\gamma)\,,$$
where $C(\wt\gamma)$ is a canonical rank $1$ free $R$-module \footnote{For $R = \Z$ it is a ``fake $\Z$,'' as J.-Y. Welschinger would call it.} associated to $\wt\gamma$. It is generated by the two possible orientations of the determinant line bundle of a certain natural family of Fredholm operators, which are formal linearizations of the Floer PDE defined on cappings $\wh\gamma$ in class $\wt\gamma$. This module is graded by the Conley--Zehnder index $m_{H:L}$.

\begin{rem}
We wish to remark that until quite recently there was no complete treatment of the topic of determinant lines of Fredholm operators, except the paper \cite{Knudsen_Mumford_Projectivity_moduli_space_stable_curves_I_prelims_det_Div}, which however is very abstract, and in which the useful properties were not formulated. To the best of our knowledge the first complete exposition including such properties, only appeared in 2013 in Zinger's paper \cite{Zinger_Det_line_bundle_Fredholm_ops}, on which the present paper relies very heavily. The somewhat unfortunate reality of this topic is that the word ``canonical'' is very much used, often without specifying which one of the possible canonical choices is selected. Zinger's paper remedies the situation. See \S\ref{s:determinants}.
\end{rem}

Given a sufficiently generic compatible almost complex structure $J$ the Floer boundary operator
$$\partial_{H,J} \fc CF_*(H:L) \to CF_{*-1}(H:L)$$
is defined via its matrix elements which are homomorphisms
$$\sum_{[u] \in \cM(H,J;\wt\gamma_-,\wt\gamma_+)} C(u) \fc C(\wt\gamma_-) \to C(\wt\gamma_+)\,,$$
where for $\wt\gamma_\pm \in \Crit \cA_{H:L}$ of index difference $1$ we have the moduli space of solutions $\cM(H,J;\wt\gamma_-,\wt\gamma_+)$ of the Floer PDE, which is a finite set. Here
$$C(u) \fc C(\wt\gamma_-) \to C(\wt\gamma_+)$$
is an isomorphism defined roughly as follows: linear gluing of Fredholm operators allows us to produce, starting from a formal linearization of the Floer PDE on a capping $\wh\gamma_-$, a formal linearization on a capping $\wh\gamma_+$ in the class $\wt\gamma_+$. It then follows from the existence of canonical isomorphisms of determinant lines for Fredholm operators that there is a bijection between the orientations of the linearized operator $D_u$ of $u$ and isomorphisms $C(\wt\gamma_-) \simeq C(\wt\gamma_+)$. The isomorphism $C(u)$ is the one corresponding to the orientation of $D_u$ given by the positive direction of the natural $\R$-action on its kernel.

We then prove that $\partial_{H,J}$ squares to zero. This includes a fairly standard argument involving the compactification of the $2$-dimensional solution spaces of the Floer PDE by broken trajectories, as well as bubbling analysis in case $N_L = 2$. We denote by
$$HF_*(H,J:L)$$
the resulting Lagrangian Floer homology. Usual continuation maps yield canonical isomorphisms
$$\Phi_{H,J}^{H',J'} \fc HF_*(H,J:L) \to HF_*(H',J': L)\,,$$
which satisfy the cocycle identity and therefore allow us to define the abstract Floer homology $HF_*(L)$. This is endowed with the structure of an associative unital algebra using moduli spaces of solutions of the Floer PDE on the disk with three boundary punctures. The Hamiltonian Floer homology $HF_*(M)$ is made to act on $HF_*(L)$ using the disk with two boundary and one interior puncture.

The quantum complex of a quantum datum $\cD = (f,\rho,I)$ as above has
$$QC_*(\cD:L) = \bigoplus_{\substack{q \in \Crit f\\ A \in \pi_2(M,L,q)}}C(q,A)$$
as the underlying $R$-module, where $C(q,A)$ is a certain rank $1$ free $R$-module, again generated by the orientations of a certain natural family of Fredholm operators. The boundary operator
$$\partial_\cD \fc QC_*(\cD:L) \to QC_{*-1}(\cD:L)$$
has quite an involved definition via the pearly spaces \cite{Biran_Cornea_Quantum_structures_Lagr_submfds, Biran_Cornea_Rigidity_uniruling}, and therefore we will not review it here. The quantum homology
$$QH_*(\cD:L)$$
is the homology of the chain complex $(QC_*(\cD:L),\partial_\cD)$. We include a treatment of the particular case $N_L = 2$ when ``bubbling'' may arise, see \S \ref{ss:boundary_op_squares_zero_bubbling_QH}. We construct a product and a superalgebra structure over $QH_*(M)$.

Finally the PSS morphism
$$\PSS_{H,J}^\cD \fc CF_*(H:L) \to QC_*(\cD:L)$$
is defined via mixed Floer-pearly moduli spaces. We show that it induces an isomorphism on homology and that it respects continuation maps, which implies that it induces an isomorphism
$$HF_*(L) \to QH_*(\cD:L)\,.$$
We also construct an opposite isomorphism
$$QH_*(\cD:L) \to HF_*(L)$$
using an analogous strategy. Composing the two isomorphisms for different quantum data we obtain ``continuation maps'' for quantum homology
$$\Phi_\cD^{\cD'} \fc QH_*(\cD:L) \to QH_*(\cD':L)\,,$$
which are isomorphisms by construction, and which satisfy a cocycle identity. This allows us to define the abstract quantum homology $QH_*(L)$. PSS isomorphisms are shown to respect the algebraic structures on both sides, such as the quantum product and the quantum module action, which means that $QH_*(L)$ inherits a product and a superalgebra structure over $QH_*(M)$, and in particular the product is unital and associative.

\subsection{Relation with previous results and constructions}\label{ss:relation_prev_constructions}

Lagrangian Floer homology was, of course, defined by Floer \cite{Floer_Morse_thry_Lagr_intersections, Floer_unregularized_grad_flow_symp_action, Floer_Witten_cx_inf_dim_Morse_thry} for weakly exact Lagrangians and extended to the monotone case by Oh \cite{Oh_FH_Lagr_intersections_hol_disks_I}, both over $\Z_2$. For the pants product on Lagrangian Floer homology see \cite{Oh_sympl_topology_action_fcnl_II} and references therein. The superalgebra structure of $HF_*(L)$ over $HF_*(M)$ is known to experts and is a generalization of the so-called ``closed-open'' map, see also \cite{Albers_Lagr_PSS_comparison_morphisms_HF}. Lagrangian quantum homology was constructed by Biran--Cornea \cite{Biran_Cornea_Quantum_structures_Lagr_submfds, Biran_Cornea_Rigidity_uniruling} using the pearly complex which was described by Oh \cite{Oh_Relative_Floer_quantum_cohomology}. In \cite{Biran_Cornea_Quantum_structures_Lagr_submfds, Biran_Cornea_Rigidity_uniruling} the authors construct a product and a quantum module action on the Lagrangian quantum homology. Lagrangian Floer homology with arbitrary coefficients is defined, in the most general case, in \cite{FOOO_Lagr_intersection_Floer_thry_anomaly_obstr_I, FOOO_Lagr_intersection_Floer_thry_anomaly_obstr_II}. Hu--Lalonde treat the monotone case in \cite{Hu_Lalonde_Relative_Seidel_morphism_Albers_map}. Also see Seidel's book \cite{Seidel_The_Book} for the exact case; the approach in it is closest to what is done in this paper. Quantum homology with orientations is defined in \cite{Biran_Cornea_Lagr_top_enumerative_invts} for $L$ being oriented and Spin. The use of (relative) (S)Pin structures to produce coherent orientations appears in \cite{FOOO_Lagr_intersection_Floer_thry_anomaly_obstr_I, FOOO_Lagr_intersection_Floer_thry_anomaly_obstr_II}, as well as in \cite{Solomon_Intersection_thry_moduli_space_holo_curves_Lag_boundary_conds, Hu_Lalonde_Relative_Seidel_morphism_Albers_map, Biran_Cornea_Lagr_top_enumerative_invts, Wehrheim_Woodward_Orientations_pseudoholo_quilts}.

The Lagrangian PSS isomorphism appears already in \cite{Cornea_Lalonde_Cluster_homology}, in the context of cluster homology. The papers \cite{Biran_Cornea_Quantum_structures_Lagr_submfds, Biran_Cornea_Rigidity_uniruling} define this isomorphism between Lagrangian quantum and Floer homology, in a way that is used in the present paper. Certain particular cases of the Lagrangian PSS morphism for monotone Lagrangians were handled by Albers \cite{Albers_Lagr_PSS_comparison_morphisms_HF}. The case of the zero section of a cotangent bundle appears in \cite{Katic_Milinkovic_PSS_Lagr_intersections}. For conormal bundles see \cite{Duretic_PSS_isos_spectral_invts_conormal_bundle}. For weakly exact Lagrangians a description of the PSS isomorphism appears in \cite{Leclercq_spectral_invariants_Lagr_FH} over $\Z_2$ and over arbitrary rings in \cite{HLL11}. The fact that the Lagrangian PSS morphism respects the natural product structures was proved in \cite{Katic_Milinkovic_Simcevic_cohomology_rings_iso_HF_Morse} for the zero section of a cotangent bundle. In general, the fact that the algebraic structures on Lagrangian Floer and quantum homology are canonically isomorphic has been known to experts, although a proof does not seem to be written anywhere.

Coherent orientations were introduced into symplectic topology by Floer and Hofer \cite{Floer_Hofer_Coherent_orientations}. This approach is mainstream now. The use of canonical orientations was inspired by conversations with J.-Y. Welschinger and by his paper \cite{Welschinger_Open_strings_Lag_conductors_Floer_fctr}, as well as by Seidel's book \cite{Seidel_The_Book}. In particular our spaces $C(\wt\gamma)$ and $C(q,A)$ appearing in \S\ref{ss:overview_construction} are similar to the orientation spaces of \cite[Chapter 11]{Seidel_The_Book}. See also Abouzaid's expository paper \cite{Abouzaid_Symp_H_Viterbo_thm}. Coherent orientations for the PSS isomorphism in the cotangent bundle case were constructed by Kati\'c--Milinkovi\'c \cite{Katic_Milinkovic_coherent_oris_mixed_moduli_spaces}. Coherent orientations in the context of cluster complexes can be found in \cite{Charest_Source_spaces_perturbations_cluster_complexes}; the approach of that paper resembles what we do here to an extent.

Lagrangian spectral invariants using Floer theory were defined in \cite{Oh_sympl_topology_action_fcnl_I} for the zero section of a cotangent bundle, in \cite{Leclercq_spectral_invariants_Lagr_FH} for weakly exact Lagrangians, and most recently in \cite{Leclercq_Zapolsky_Spectral_invts_monotone_Lags} for monotone Lagrangians using the technical results of the present paper.

\subsection{Overview of the paper}

Since our approach is to use canonical orientations, and the existing literature on this topic is quite scarce, we decided to include an extensive treatment of Floer and quantum homology, which also serves to establish notation. The methods developed while constructing these serve as a foundation for the construction of the PSS isomorphism and for the proofs of its various properties.

It is best to read the paper linearly. In order to avoid redundancies, each section builds upon the previous ones. Since the subject matter is quite technical, the exposition is terse, therefore short summaries are given at the beginning of each section.

In \S\ref{s:determinants} we review the basic notion of the determinant line of a Fredholm operator and state the properties necessary for the rest of the paper.

In \S\ref{s:HF} we define the canonical Floer complex $(CF_*(H:L),\partial_{H,J})$ associated to a regular Floer datum $(H,J)$. We define the boundary operator and prove that it squares to zero, when bubbling of Maslov $2$ disks is absent. Various algebraic structures are defined, such as the quantum product, the unit, and the quantum module structure, and their properties are proved, such as the associativity of the product and the like. The abstract Floer homology $HF_*(L)$ is defined.

In \S\ref{s:QH} we define the quantum complex $(QC_*(\cD:L),\partial_\cD)$, using the canonical orientation approach. We similarly define the boundary operator, prove that it squares to zero in absence of bubbling, and define the quantum product and quantum module action.

In \S\ref{s:PSS} we construct the canonical PSS maps between the Floer and quantum homology of $L$, and prove that they are in fact isomorphisms. As a result of their properties, PSS isomorphisms are used to define ``continuation maps'' in quantum homology, which we use instead of the more direct approach of Biran--Cornea \cite{Biran_Cornea_Quantum_structures_Lagr_submfds, Biran_Cornea_Rigidity_uniruling}. Thus we obtain the abstract quantum homology $QH_*(L)$. We prove that the PSS isomorphisms intertwine the algebraic structures on both sides, such as products and quantum module actions.

In \S\ref{s:boundary_op_squares_zero_bubbling} we prove that the Lagrangian Floer and quantum boundary operators square to zero when there is bubbling present.

In \S\ref{s:quotient_cxs} we describe the construction of quotient complexes, which are the more familiar objects in Floer theory. We start with a summary of relative Pin structures and how to use them to define a system of coherent orientations on formal linearized operators corresponding to disks with boundary on $L$. We obtain a simplification of the usual process of constructing such orientations via framings. We describe the construction of quotient complexes in Hamiltonian Floer homology and the quantum homology of $M$, which require no choices at all, and then proceed with the Lagrangian case, which is much less trivial and for which coherent orientations are needed.

In \S\ref{s:examples_computations} we compute the canonical quantum complexes for $\R P^n \subset \C P^n$, as well as for three monotone tori, the Clifford torus and the Chekanov torus in $\C P^2$, and the exotic torus in $\C P^1 \times \C P^1$.

\begin{acknow}
I wish to thank David Blanc, Strom Borman, Fran\c cois Charette, Boris Chorny, Tobias Ekholm, Misha Entov, Oli Fabert, Penka Georgieva, Luis Haug, Vladimir Hinich, Marco Mazzucchelli, Dusa McDuff, Cedric Membrez, Will Merry, Fabien Ng\^o, Tony Rieser, Dietmar Salamon, Matthias Schwarz, and Egor Shelukhin for numerous stimulating discussions and general interest, Paul Biran and Octav Cornea for reading a portion of the manuscript and providing valuable comments which allowed me to write the abstract and the introduction more clearly, and Leonid Polterovich for useful suggestions. The approach chosen in this project was largely inspired by a conversation with J.-Y. Welschinger at a pierogi restaurant in \L \'od\'z; he taught me canonical orientations in Morse theory, and answered many questions about generalizations to Floer homology. Jacqui Espina kindly agreed to listen to a number of preliminary results and took notes, which later helped me write this manuscript. Special thanks go to Judy Kupferman who listened to me ramble about this paper over innumerable cups of coffee, and to my collaborator R\'emi Leclercq for reading a portion of the text, and for his support and patience during the long time that it took me to complete this project, which is necessary for our \cite{Leclercq_Zapolsky_Spectral_invts_monotone_Lags}.
\end{acknow}

\section{Determinant lines for Fredholm operators}\label{s:determinants}

In this section we collect the necessary preliminaries concerning Fredholm operators, their determinant lines and properties. We refer the reader to the wonderful paper by Zinger \cite{Zinger_Det_line_bundle_Fredholm_ops}, which constructs determinant lines of Fredholm operators and proves their various properties and classification in complete detail.

Given two real Banach spaces $X,Y$, a Fredholm operator between them is a bounded linear operator
$$D \fc X \to Y$$
with closed range, whose kernel and cokernel are both finite-dimensional. We let $\cF(X,Y)$ be the set of Fredholm operators. It is an open subspace of the space of bounded linear operators between $X$ and $Y$ relative to the norm topology. The index of $D$ is the integer
$$\ind D = \dim \ker D - \dim \coker D\,.$$
The index is a continuous function on $\cF(X,Y)$ and consequently it is locally constant.

A \tb{$\Z_2$-graded line}, or a graded line for short, is a one-dimensional real vector space together with a grading, which is an element of $\Z_2 = \{\ol 0, \ol 1\}$. For a graded line $L$ we let $\deg L \in \Z_2$ be its grading. For a finite-dimensional real vector space $V$ we denote by
$$\ddd(V)$$
the graded line whose underlying line is the top exterior power $\bigwedge^{\text{max}}V$ and whose grading is $\dim V \pmod 2$. The \tb{determinant line} of a Fredholm operator $D \in \cF(X,Y)$ is the graded line
$$\ddd(D) = \ddd(\ker D) \otimes \big(\ddd(\coker D)\big)^\vee\quad \text{with grading} \quad \deg \ddd(D) = \ind D \pmod 2\,.$$
We let
$$\ddd_{X,Y} = \biguplus_{D \in \cF(X,Y)}\ddd(D)\,.$$
For now this is just a set with fiberwise linear structure with respect to the obvious projection
$$\ddd_{X,Y} \to \cF(X,Y)\,.$$
We wish to put a topology on $\ddd_{X,Y}$ so that it becomes a line bundle over $\cF(X,Y)$. However, not any topology would do. The reason for this is that in applications there are a plethora of natural operations on Fredholm operators, and we wish our topology to be compatible with these operations.

Abstractly, we have the collection of spaces (yet without topology) $\ddd_{X,Y}$, indexed by pairs of real Banach spaces $(X,Y)$. Therefore \emph{a priori} we can topologize each one of these spaces to make them into real line bundles. As we just pointed out, not every such system of topologies will have the desired properties. Zinger \cite{Zinger_Det_line_bundle_Fredholm_ops} lists a number of very natural properties which should be satisfied by any useful such system of topologies, and proves the following fundamental result.
\begin{thm}\label{thm:exist_systems_topologies_det_lines}
There exist systems of topologies on the collection of spaces $\ddd_{X,Y}$, indexed by pairs of Banach spaces $(X,Y)$, satisfying the desired natural properties. \qed
\end{thm}
\noindent It would take us too far afield to list all the properties formulated by Zinger. Instead we refer the reader to his paper, which also contains a precise formulation of Theorem \ref{thm:exist_systems_topologies_det_lines}.

In this paper we fix once and for all a single compatible system of topologies whose existence is guaranteed by Theorem \ref{thm:exist_systems_topologies_det_lines}, for instance the system explicitly described in \cite[Section 4.2]{Zinger_Det_line_bundle_Fredholm_ops}. We will now list the properties which will be needed in the applications presented in the paper. Note that in this section we only list properties pertaining to abstract Fredholm operators. Later in the paper we will specialize to real Cauchy--Riemann operators (see \S \ref{s:HF}, \ref{s:QH}); these have additional properties specific to them, and we will introduce them at appropriate points in the text.

\begin{rem}
Systems of topologies whose existence is asserted in Theorem \ref{thm:exist_systems_topologies_det_lines} are not unique. In fact, Zinger classifies all such systems in \cite[Theorem 2]{Zinger_Det_line_bundle_Fredholm_ops}. From this classification it is apparent that there are infinitely many such systems. However, this has no bearing on results of computations, since they are based solely on the properties of these topologies, which ultimately boil down to canonical constructions involving solely finite-dimensional vector spaces. Since all the compatible systems of topologies satisfy the same properties, computations will always yield the same result, whatever system of topologies is chosen.
\end{rem}

\subsection{Properties of determinant lines}\label{ss:pties_det_lines}

\paragraph{Naturality with respect to isomorphisms} If $f \fc V \to W$ is an isomorphism of finite-dimensional real vector spaces, we let
$$\ddd(f) \fc \ddd(V) \to \ddd(W)$$
be the induced isomorphism of determinant lines. If $\phi \fc X \to X'$, $\psi \fc Y \to Y'$ are Banach space isomorphisms, there is an induced homeomorphism
$$\cF(X,Y) \to \cF(X',Y')\,,\quad D \mapsto \psi \circ D \circ \phi^{-1}\,.$$
This lifts to an isomorphism of line bundles
$$\ddd_{X,Y} \to \ddd_{X',Y'}\,,$$
as follows. For $D \in \cF(X,Y)$ let $D' = \psi \circ D \circ \phi^{-1}$. The map between $\ddd(D)$ and $\ddd(D')$ is given by
$$\ddd(\phi|_{\ker D}) \otimes \big(\ddd(\ol \psi)^{-1}\big)^\vee \fc \ddd(\ker D) \otimes \big(\ddd(\coker D)\big)^\vee \to \ddd(\ker D') \otimes \big(\ddd(\coker D')\big)^\vee\,,$$
where $\ol \psi \fc \coker D \to \coker D'$ is the isomorphism induced by $\psi$.

\begin{rem}[Banach bundles] Given a system of topologies on the spaces $\ddd_{X,Y}$ satisfying the naturality property, we can topologize the determinant line of a Fredholm morphism between two Banach bundles. In more detail, let $p \fc \cX \to B$ and $q \fc \cY \to B$ be two locally trivial Banach bundles over a space $B$, with fibers $X$ and $Y$, respectively. A \tb{Fredholm morphism} between these is a fiber-preserving fiberwise linear continuous map $\cD \fc \cX \to \cY$ such that for each $b \in B$, $\cD_b \in \cF(\cX_b,\cY_b)$. We have the associated space
$$\ddd(\cD) = \biguplus_{b\in B}\ddd(\cD_b)\,.$$
Local trivializations of $\cX$ and $\cY$ over $U \subset B$ conjugate $\cD|_U$ to a map $D_U \fc U \to \cF(X,Y)$, which leads to the line bundle $\ddd(D_U) := (D_U)^* \ddd_{X,Y}$. We can push the topology on $\ddd(D_U)$ to a topology on the space $\ddd(\cD)|_U$. The induced topology on $\ddd(\cD)$ is then well-defined thanks to the naturality property, and it makes $\ddd(\cD)$ into a line bundle over $B$.
\end{rem}

\paragraph{Exact triples} An \tb{exact triple of Fredholm operators}, or an exact Fredholm triple for short, is a commutative diagram with exact rows
$$\xymatrix{0 \ar[r] & X' \ar[r] \ar[d]^{D'} & X \ar[r] \ar[d]^{D} & X'' \ar[r] \ar[d]^{D''} & 0 \\
0 \ar[r] & Y' \ar[r] & Y \ar[r] & Y'' \ar[r] & 0}$$
where $X',X,X'',Y',Y,Y''$ are Banach spaces, $D',D,D''$ are Fredholm operators. We will denote an exact triple by $\mft = (D',D,D'')$ with the Banach spaces and maps between them being implicit. To each such exact Fredholm triple there corresponds an isomorphism
$$\Psi_\mft \fc \ddd(D') \otimes \ddd(D'') \to \ddd(D)\,,$$
called the exact triple isomorphism. These isomorphisms lift to line bundle isomorphisms over the spaces of exact triples with fixed Banach spaces.
\begin{rem}
We will often omit the Banach spaces and write an exact triple in an abbreviated form as
$$0 \to D' \to D \to D'' \to 0\,.$$
This should cause no confusion.
\end{rem}

\begin{rem}
Oftentimes in applications exact triples come in the form of a family of exact triples where instead of fixed Banach spaces we have Banach bundles and Fredholm morphisms between them. It is then apparent that the exact triple isomorphisms depend continuously on the base space of the bundle.
\end{rem}

\paragraph{Normalization}\label{par:normalization_pty} Given a short exact sequence of finite-dimensional real vector spaces
$$0 \to V' \xrightarrow{\iota} V \xrightarrow{\pi} V'' \to 0\,,$$
there is a naturally induced isomorphism
\begin{equation}\label{eqn:iso_det_lines_exact_seq_vector_spaces}
\ddd(V') \otimes \ddd(V'') \to \ddd(V)\,,
\end{equation}
defined as follows. Pick ordered bases $v_1',\dots,v_k' \in V'$ and $v_{k+1}'',\dots, v_{k+l}'' \in V''$, where $k = \dim V'$, $l = \dim V''$, define $v_i = \iota(v_i')$ for $i \leq k$, and let $v_i \in V$ be such that $\pi(v_i) = v_i''$ for $i > k$. The isomorphism then sends
$$\textstyle \bigwedge_{i=1}^k v_i' \otimes \bigwedge_{i=1}^lv_{k+i}'' \mapsto \bigwedge_{i=1}^{k+l}v_i\,.$$
It can be checked that this is independent of the chosen bases.

Note that if $0$ denotes the zero vector space, we have canonically $\ddd(0) \equiv \R$. If $D$ is a surjective Fredholm operator, then there is a canonical isomorphism
\begin{equation}\label{eqn:identify_det_Fredholm_op_det_its_kernel}
\ddd(\ker D) \to \ddd(D) = \ddd(\ker D) \otimes \R^\vee \,,\quad \sigma \mapsto \sigma \otimes 1^\vee\,.
\end{equation}
We will sometimes tacitly identify the determinant line of a surjective operator with the determinant line of its kernel.

If $\mft=(D',D,D'')$ is an exact Fredholm triple of \emph{surjective} operators, then there is an induced short exact sequence of kernels:
$$0 \to \ker D' \to \ker D \to \ker D'' \to 0\,.$$
The isomorphism
$$\ddd(\ker D') \otimes \ddd(\ker D'') \to \ddd(\ker D)$$
induced from this short exact sequence coincides with the isomorphism $\Psi_\mft$ if we use the identification \eqref{eqn:identify_det_Fredholm_op_det_its_kernel}.

Also, if $V$ is a finite-dimensional space, we have the Fredholm operator $0_V \fc V \to 0$ and an obvious isomorphism
$$\ddd(0_V) = \ddd(V)\,.$$

\paragraph{The exact squares property}\label{par:exact_squares_pty} Given two graded lines $L_1,L_2$ we define the \tb{interchange isomorphism}
$$R \fc L_1 \otimes L_2 \to L_2 \otimes L_1\,,\quad v_1 \otimes v_2 \mapsto (-1)^{\deg L_1\cdot \deg L_2} v_2 \otimes v_1\,.$$

An \tb{exact square} of vector spaces is by definition a short exact sequence of short exact sequences of vector spaces, that is a commutative diagram with exact rows and columns consisting of nine vectors spaces and maps between them, plus bounding zero vector spaces. An \tb{exact square of Fredholm operators} is a commutative diagram consisting of two layers of exact squares of Banach spaces and nine Fredholm operators between those layers. We denote such a square schematically as follows, where the bounding zeroes are implicit and are omitted for the sake of economy of space:
$$\xymatrix{D_{\text{LT}} \ar[r] \ar[d] & D_{\text{CT}} \ar[r] \ar[d] & D_{\text{RT}} \ar[d] \\
D_{\text{LM}} \ar[r] \ar[d] & D_{\text{CM}} \ar[r] \ar[d] & D_{\text{RM}} \ar[d] \\
D_{\text{LB}} \ar[r] & D_{\text{CB}} \ar[r] & D_{\text{RB}}}$$
The various exact triple isomorphisms form the following commutative diagram:
$$\xymatrix{\ddd(D_{\text{LT}}) \otimes \ddd(D_{\text{RT}}) \otimes \ddd(D_{\text{LB}}) \otimes \ddd(D_{\text{RB}}) \ar[r]^-{\Psi_{\text{T}} \otimes \Psi_{\text{B}}} \ar[d]^{(\Psi_{\text{L}} \otimes \Psi_{\text{R}}) \circ (\id \otimes R \otimes \id)} & \ddd(D_{\text{CT}}) \otimes \ddd(D_{\text{CB}}) \ar[d]^{\Psi_{\text{C}}}\\
\ddd(D_{\text{LM}}) \otimes \ddd(D_{\text{RM}}) \ar[r]^-{\Psi_{\text{M}}} & \ddd(D_{\text{CM}}) }$$
where $\Psi_{\text{L}}$, $\Psi_{\text{C}}$, $\Psi_{\text{R}}$, $\Psi_{\text{T}}$, $\Psi_{\text{M}}$, $\Psi_{\text{B}}$ denote the exact triple isomorphisms corresponding to the left, center, right, top, middle, and bottom exact triples appearing in the diagram, respectively. A parametrized version of this commutative diagram exists when the given exact squares involve Banach bundles and Fredholm morphisms.

\paragraph{Direct sum isomorphisms} \label{par:direct_sum_isos} Given Fredholm operators $D_i \in \cF(X_i,Y_i)$, $i=1,2$, there is the direct sum operator $D_1 \oplus D_2 \fc X_1 \oplus X_2 \to Y_1 \oplus Y_2$ and an obvious exact triple
$$0 \to D_1 \to D_1 \oplus D_2 \to D_2 \to 0\,.$$
This exact triple gives rise to an isomorphism
$$\ddd(D_1) \otimes \ddd(D_2) \to \ddd(D_1 \oplus D_2)\,.$$
This isomorphism pervades the present paper and is of great importance. We refer to it as the \tb{direct sum isomorphism}.

Since direct sum isomorphisms are a particular case of exact triple isomorphisms, the exact squares property of the latter implies two properties of the former, namely supercommutativity and associativity. Using the exact Fredholm square
$$\xymatrix{0 \ar[r] \ar[d] & D_1 \ar[r] \ar[d] & D_1 \ar[d] \\ D_2 \ar[r] \ar[d] & D_1 \oplus D_2 \ar[r] \ar[d] & D_1 \ar[d] \\ D_2 \ar[r] & D_2 \ar[r] & 0}$$
we obtain the following commutative diagram:
$$\xymatrix{\ddd(D_1) \otimes \ddd(D_2) \ar@{=}[r] \ar[d]_{R} & \ddd(D_1) \otimes \ddd(D_2) \ar[d]^{\oplus} \\ \ddd(D_2) \otimes \ddd(D_1) \ar[r]^{\oplus} & \ddd(D_1 \oplus D_2)}$$
where $\oplus$ denotes the direct sum isomorphism. This means that the composition of direct sum isomorphisms
$$\ddd(D_1) \otimes \ddd(D_2) \to \ddd(D_1 \oplus D_2) \to \ddd(D_2) \otimes \ddd(D_1)$$
coincides with the interchange isomorphism $R$. This is the supercommutativity property of the direct sum isomorphisms.

Next, if we have a third operator $D_3 \in \cF(X_3,Y_3)$, then we have the exact square
$$\xymatrix{D_1 \ar[r] \ar[d] & D_1 \oplus D_2 \ar[r] \ar[d] & D_2 \ar[d] \\
D_1 \ar[r] \ar[d] & D_1 \oplus D_2 \oplus D_3 \ar[r] \ar[d] & D_2 \oplus D_3 \ar[d]\\
0 \ar[r] & D_3 \ar[r] & D_3}$$
which yields the commutative diagram
$$\xymatrix{\ddd(D_1) \otimes \ddd(D_2) \otimes \ddd(D_3) \ar[r]^{\oplus \otimes \id} \ar[d]^{\id \otimes \oplus} & \ddd(D_1 \oplus D_2) \otimes \ddd(D_3) \ar[d]^{\oplus} \\
\ddd(D_1) \otimes \ddd(D_2 \oplus D_3) \ar[r]^{\oplus} & \ddd(D_1 \oplus D_2 \oplus D_3)}$$
which means that the direct sum isomorphisms are associative.

Again, all these properties hold in parametric versions as well.

\paragraph{Final remark on terminology} In this paper we will be interested in \emph{orientations} of a Fredholm operator $D$, which are elements of the two-point set $(\ddd(D) - \{0\})/\R_{>0}$, and therefore we will often say that a diagram of real lines and isomorphisms \tb{commutes} to mean that it commutes up to multiplication by a positive real number. Also we will say that a given isomorphism between real lines is \tb{canonical}, even if it is only canonically defined up to a positive multiple.

\section{Floer homology}\label{s:HF}

In this section we construct the canonical chain complexes computing Lagrangian and Hamiltonian Floer homology. We define the various algebraic structures on these complexes, such as products, the module action of the Hamiltonian Floer homology on the Lagrangian Floer homology, and prove various relations between these operations. We develop these structures in a TQFT-like framework, which allows us to describe all of them in a transparent and unified manner.

In \S\ref{ss:punctures_RS_gluing} we define punctured Riemann surfaces, their gluing, and related concepts. \S\ref{ss:CROs_dets_gluing} we define real-linear Cauchy--Riemann operators, their determinant lines and gluing of such operators defined on surfaces undergoing gluing. In \S\ref{ss:b_smooth_maps_pregluing} we define the technical notion of b-smooth maps, which are maps on punctured Riemann surfaces satisfying the property that they extend to suitable compactifications of the surfaces; the main property of b-smooth maps is that solutions of the Floer PDE are b-smooth. In \S\ref{ss:CROs_from_b_smooth_maps} we describe the Cauchy--Riemann operators arising from b-smooth maps as a formal linearization of the Floer operator, and how pregluing b-smooth maps relates to gluing of the corresponding formal linearizations. In \S\ref{ss:Floer_PDE} we define the Floer PDE on a punctured Riemann surface, as well as on a family of such surfaces. In \S\ref{ss:moduli_spaces_sols_Floer_PDE} we define solution spaces and moduli spaces of solutions of the Floer PDE, and describe their compactness properties. In \S\ref{ss:orientatiosn} we describe the canonical orientations of the linearized operators corresponding to solutions of the Floer PDE of the lowest dimension as well as the orientations of compactified moduli spaces induced on them by the canonical orientations corresponding to the boundary points. In \S\ref{ss:operations} we define the matrix elements of operations in Floer homology, using the canonical orientations defined in \S\ref{ss:orientatiosn}, and prove that the various matrix elements are subject to identities, a fact whose proof uses the induced orientations. In \S\ref{ss:HF} we define the Floer complexes and homology, various algebraic operations on them, and prove their properties. \S\ref{ss:arbitrary_rings_loc_coeffs} deals with the case of arbitrary rings and twisting by local systems. In \S\ref{ss:duality_HF} we treat duality in Floer homology and define the augmentation map as the dual of the unit.

We refer the reader to the papers \cite{Floer_Morse_thry_Lagr_intersections, Floer_unregularized_grad_flow_symp_action, Floer_Witten_cx_inf_dim_Morse_thry, Oh_FH_Lagr_intersections_hol_disks_I}, Schwarz's thesis \cite{Schwarz_PhD_thesis}, Seidel's book \cite{Seidel_The_Book}, and references therein for the analytical results used here. We do not provide precise references for all the results, mainly because they are more or less standard by now. The material presented here is largely borrowed from the wonderful book \cite{Seidel_The_Book}, especially Part II, which has also been extremely influential on the style of exposition chosen here.

\subsection{Punctured Riemann surfaces and their gluing}\label{ss:punctures_RS_gluing}

Fix a compact connected Riemann surface $\wh\Sigma$ with (possibly empty) boundary, and a finite subset $\Theta \subset \wh\Sigma$. The elements of $\Theta$ are called \tb{punctures}. We let $\Sigma = \wh\Sigma - \Theta$; the surface $\Sigma$ is called a \tb{punctured Riemann surface}. The set of punctures $\Theta$ is decomposed into two disjoint subsets, $\Theta = \Theta^+ \cup \Theta^-$ of \tb{positive} and \tb{negative} punctures. A puncture $\theta$ is called \tb{boundary} if $\theta \in \partial\wh\Sigma$, otherwise it is \tb{interior}.

Set $\R^\pm = \{s \in \R\,|\, \pm s \geq 0\}$. Throughout we use the following standard surfaces with boundary and corners:
$$S = \R \times [0,1]\,,\quad C = \R \times S^1\,,$$
the \tb{standard strip} and the \tb{standard cylinder}, and
$$S^\pm = \R^\pm \times [0,1]\,,\quad C^\pm = \R^\pm \times S^1$$
the \tb{standard half-strips} and \tb{half-cylinders}. The strips $S,S^\pm$ are given the conformal structures coming from the obvious embeddings into $\C$ while $C,C^\pm$ are given conformal structures by viewing $S^1 = \R/\Z$. We let $(s,t)$ be the standing notation for the standard conformal coordinates on $S^\pm,C^\pm$.

Let $\theta \in\Theta^\pm$ be a puncture. An \tb{end associated to} $\theta$ is a proper conformal embedding
$$\epsilon_\theta \fc C^\pm \to \Sigma \quad \text{if } \theta \text{ is interior,}$$
$$\epsilon_\theta \fc S^\pm \to \Sigma \quad \text{if } \theta \text{ is boundary,}$$
in which case we also require that $\epsilon_\theta^{-1}(\partial \Sigma) = \R^\pm \times \{0,1\}$, where $\lim_{s \to \pm\infty}\epsilon_\theta(s,t) = \theta$ in $\wh\Sigma$.

A \tb{choice of ends} for $\Sigma$ is a family $\{\epsilon_\theta\}_{\theta\in\Theta}$ of ends associated to all the punctures of $\Sigma$, with pairwise disjoint images.

\paragraph{Gluing punctured Riemann surfaces} \label{par:gluing_punctured_Riem_surf} Let $\cT$ be a finite tree with vertex set $\cV=\cV(\cT)$ and edge set $\cE=\cE(\cT) \subset \cV \times \cV$. Assume each vertex $v \in \cV$ is labeled by a punctured Riemann surface $\Sigma_v$ with puncture set $\Theta_v$, and that each edge $e = (v,v')$ is labeled by a pair of punctures $(\theta,\theta') \in \Theta_v^+ \times \Theta_{v'}^-$ of the same type (both boundary or both interior) and a positive real number $R_e$ called a gluing length. Fix a choice of ends for each $\Sigma_v$. We can define the \tb{glued surface} $\Sigma_\cT$ corresponding to these data as follows. Take the disjoint union
$$\biguplus_{v \in \cV}\Sigma_v\,,$$
and for each edge $e$ labeled by $(\theta,\theta') \in \Theta_v^+ \times \Theta_{v'}^-$, remove from it the subset
$$\epsilon_\theta([R_e,\infty)\times[0,1]) \cup \epsilon_{\theta'}((-\infty,-R_e]\times [0,1])\quad \text{if } \theta,\theta' \text{ are boundary, or}$$
$$\epsilon_\theta([R_e,\infty)\times S^1) \cup \epsilon_{\theta'}((-\infty,-R_e]\times S^1)\quad \text{ if they are interior,}$$
where $\epsilon_\theta,\epsilon_{\theta'}$ are the ends associated to $\theta,\theta'$. On the resulting subset of $\biguplus_v\Sigma_v$ make the identification
$$\epsilon_\theta(s,t)\simeq \epsilon_{\theta'}(-R_e+s,t)$$
for $e=(\theta,\theta')$ and $s\in (0,R_e)$.
The glued surface $\Sigma_\cT$ inherits a conformal structure from the $\Sigma_v$. It also inherits punctures and a choice of ends from the $\Sigma_v$, namely all the punctures that did not appear in labels of the edges of $\cT$, and the ends associated to them.

Gluing of punctured Riemann surfaces satisfies an associativity property. To formulate it, let $\cF = \cT_1\uplus \dots \uplus \cT_k$ be a subforest of $\cT$, that is a graph obtained from $\cT$ by deleting some of the edges. We can then form the glued surfaces $\Sigma_{\cT_1},\dots,\Sigma_{\cT_k}$ according to the procedure just described. The quotient tree $\ol\cT = \cT/\cF$ has vertex set $\ol\cV = \{\cT_i\}_i$. An edge of $\ol\cT$ is an edge of $\cT$ connecting a pair of the subtrees $\cT_i$. Label the vertices of $\ol\cT$ by the surfaces $\Sigma_{\cT_i}$. Note that all the $\Sigma_{\cT_i}$ have a choice of ends and puncture sets coming from gluing. Labels of edges of $\cT$ not appearing in $\cF$ define in a natural way labels of edges of $\ol\cT$. Therefore we can form the glued surface $\Sigma_{\ol\cT}$. The associativity of gluing is expressed by means of an obvious canonical identification
$$\Sigma_\cT  = \Sigma_{\ol\cT}\,,$$
which preserves the conformal structures, the punctures, and the choice of ends.

\subsection{Cauchy--Riemann operators, their determinant lines and gluing}\label{ss:CROs_dets_gluing}

Let $\Sigma$ be a punctured Riemann surface and endow it with a set of ends. Let $(E,F) \to (\Sigma,\partial\Sigma)$ be a Hermitian bundle pair, that is $E$ is a vector bundle endowed with a symplectic form $\omega$ and a compatible almost complex structure $J$, and $F \subset E|_{\partial\Sigma}$ is a Lagrangian subbundle. Choose a connection $\nabla$ on $E$. A \tb{real Cauchy--Riemann operator} is an operator of the form
$$\ol\partial_\nabla = \nabla^{0,1} \fc C^\infty (\Sigma,\partial\Sigma; E,F) \to C^\infty(\Sigma,\Omega^{0,1}_\Sigma \otimes E)\,,$$
that is the complex-antilinear part of $\nabla$, which is defined as follows:
$$\nabla^{0,1}\xi = \xi + J\nabla_{j\cdot}\xi\,,$$
where $j$ denotes the conformal structure on $\Sigma$. When $\Sigma$ is compact, such an operator can be extended to suitable Sobolev completions, where it becomes a Fredholm operator. In order to have an analogous statement for noncompact $\Sigma$, we need to have sufficient control on the behavior of our data at infinity.

A \tb{limiting datum} at a boundary puncture $\theta$ is a quintuple
$$(E_\theta,F_\theta,\omega_\theta,J_\theta,\nabla_\theta)\,,$$
where $(E_\theta,F_\theta) \to ([0,1],\{0,1\})$ is a Hermitian bundle pair with symplectic form $\omega_\theta$ and complex structure $J_\theta$, and $\nabla_\theta$ is a symplectic connection on $E_\theta$, meaning its parallel transport maps along $[0,1]$ are symplectic. If $\theta$ is interior, a limiting datum is a quadruple $(E_\theta,\omega_\theta,J_\theta,\nabla_\theta)$, where $E_\theta \to S^1$ is a Hermitian bundle, while the rest of the symbols carry the same meaning as in the boundary case.

Assume we have fixed a choice of limiting data at all the punctures of $\Sigma$. Let $\pi \fc S^\pm \to [0,1]$ or $\pi \fc C^\pm \to S^1$ be the projection onto the $t$ variable and assume that for every puncture $\theta$ we have fixed identifications
\begin{equation}\label{eqn:identification_Hermitian_bundle_with_limiting_bundle_at_puncture}
\epsilon_\theta^*E \simeq \pi^*E_\theta
\end{equation}
with respect to which all the data on $\Sigma$ are asymptotic in suitable topologies to the limiting data at $\theta$, that is $\omega \to \omega_\theta$, $J \to J_\theta$, $F \to F_\theta$, and $\nabla \to \nabla_\theta$. In this case we call $\ol\partial_\nabla$ \tb{admissible}. Note for future use that these identifications can be deformed so that $\omega,J,F,\nabla$ all become constant on the ends. This will be useful for gluing in the next subsection.

We call the connection $\nabla_\theta$ \tb{nondegenerate} if, in case $\theta$ is boundary, the parallel transport map along $[0,1]$ maps $F_{\theta,0}$ to a subspace of $E_{\theta,1}$ transverse to $F_{\theta,1}$, and in case $\theta$ is interior, the parallel transport map around $S^1$ does not have $1$ as an eigenvalue. If all the connections $\nabla_\theta$ are nondegenerate, we call the operator $\ol\partial_\nabla$ nondegenerate.

Using the measure induced on $\Sigma$ by the choice of ends, one can define Sobolev completions
$$W^{1,p}(\Sigma,\partial\Sigma; E,F)\,,\;L^p(\Sigma,\Omega^{0,1}_\Sigma\otimes E)$$
of the corresponding spaces of smooth sections for $p > 2$; these are Banach spaces. An admissible nondegenerate Cauchy--Riemann operator extends to a Fredholm operator
$$\ol\partial_\nabla \fc W^{1,p}(\Sigma,\partial\Sigma;E,F) \to L^p(\Sigma,\Omega^{0,1}_\Sigma\otimes E)\,.$$

\paragraph{Gluing Cauchy--Riemann operators and their determinant lines}\label{par:gluing_CROs_and_det_lines}
Let now $\cT$ be a gluing tree as in \S \ref{par:gluing_punctured_Riem_surf}, whose vertex set $\cV$ is labeled by punctured Riemann surfaces $\{\Sigma_v\}_{v \in \cV}$ with puncture sets $\Theta_v$ and a choice of ends for each $\Sigma_v$. Assume each edge $e = (v,v')$ of $\cT$ is labeled by a positive gluing length $R_e$ and by a pair $(\theta,\theta') \in \Theta_v^+ \times \Theta_{v'}^-$, where $\theta,\theta'$ are of the same type (both boundary or interior). As described in \S \ref{par:gluing_punctured_Riem_surf}, we can form the glued surface $\Sigma_\cT$.

For each $v \in \cV$ let $(E_v,F_v) \to (\Sigma_v,\partial\Sigma_v)$ be a Hermitian bundle pair with Hermitian structure $(\omega_v,J_v)$ and let $\nabla_v$ be a connection on $E_v$. Then we have the corresponding Cauchy--Riemann operator $D_v:=\ol\partial_{\nabla_v}$. Assume now that all the operators $D_v$ are admissible and nondegenerate, and assume that if a pair of punctures $(\theta,\theta')$ label an edge, then they have identical limiting data. Assume also that the identifications \eqref{eqn:identification_Hermitian_bundle_with_limiting_bundle_at_puncture} are such that the data $\omega_v,J_v,F_v,\nabla_v$ are constant at each puncture undergoing gluing. We can then glue the bundle pairs, the Hermitian structures, and the connections in an obvious manner to form a bundle pair
$$(E_\cT,F_\cT) \to (\Sigma_\cT,\partial\Sigma_\cT)$$
with Hermitian structure $(\omega_\cT,J_\cT)$ and a connection $\nabla_\cT$. The corresponding operator $D_\cT:= \ol\partial_{\nabla_\cT}$ is admissible and nondegenerate.

We can use cutoff functions to patch sections of $(E_v,F_v)$ together to form sections of $(E_\cT,F_\cT)$. Using the orthogonal projection \footnote{Even though we only defined the Cauchy--Riemann operators for $p > 2$, they can also be defined for $p = 2$, in which case the Sobolev spaces involved are Hilbert spaces. It is a standard fact that the kernel and cokernel of such an operator are independent of $p$. Therefore we can take elements of the kernels, patch them together, get an element in $W^{1,2}$ and use the inner product to project. The resulting section belongs to $W^{1,p}$ for all $p$. The same applies to cokernels.} onto $\ker D_\cT$ we get a linear map
$$\bigoplus_v\ker D_v \to \ker D_\cT\,,$$
which for large enough gluing lengths $R_e$ becomes an isomorphism. Similarly, we get a linear map
$$\bigoplus_v \coker D_v \to \coker D_\cT\,,$$
which is also an isomorphism for large gluing lengths. Therefore we have that
$$\ind D_\cT = \sum_v \ind D_v$$
and we obtain an isomorphism
$$\ddd \Big( \bigoplus_v D_v\Big) \simeq \ddd(D_\cT)\,,$$
independent of the choices. We refer to it as the \tb{gluing isomorphism} throughout.

This gluing of Cauchy--Riemann operators and the resulting gluing isomorphism of the determinant lines are associative, in the following sense. If $\cF=\biguplus_i \cT_i \subset \cT$ is a subforest as in \S\ref{par:gluing_punctured_Riem_surf}, we can form the glued operators $\cD_{\cT_i}$ over $\Sigma_{\cT_i}$, which are admissible and nondegenerate. These have matching limiting data at punctures labeling edges of $\ol\cT = \cT/\cF$, and therefore can also be glued. The resulting operator $D_{\ol\cT}$ can be canonically identified with $D_{\cT}$. The corresponding diagram of isomorphisms of determinant lines commutes:
$$\xymatrix{\ddd \big( \bigoplus_{v\in\cV} D_v\big) \ar [r] \ar [d]& \ddd \big(\bigoplus_i D_{\cT_i} \big) \ar [d]\\
\ddd(D_\cT) \ar@{=}[r]& \ddd(D_{\ol\cT})}$$
where all the arrows except $=$ are the gluing isomorphisms.

\paragraph{Gluing isomorphisms and direct sum isomorphisms commute} We will now formulate a crucial property satisfied by gluing isomorphisms, which is the foundation of many computations leading to various properties of algebraic operations in Floer theory. Recall the direct sum isomorphism defined in \S\ref{par:direct_sum_isos}. If we have operators $D_i$, $i=1,2,3$, where $D_1,D_2$ can be glued, that is they are defined on Hermitian bundle pairs with matching limiting data at a gluing puncture, and we let $D_1 \sharp D_2$ be the resulting glued operator, then the following diagram commutes:
$$\xymatrix{\ddd(D_1 \oplus D_2 \oplus D_3) \ar[r]^{\oplus} \ar[d]^{\sharp} & \ddd(D_1 \oplus D_2) \otimes \ddd(D_3) \ar[d]^{\sharp \otimes \id}\\ \ddd(D_1 \sharp D_2 \oplus D_3) \ar[r]^{\oplus}& \ddd(D_1 \sharp D_2) \otimes \ddd(D_3)}$$

\subsection{B-smooth maps and their pregluing}\label{ss:b_smooth_maps_pregluing}

If $H \fc [0,1] \times M \to \R$ is a time-dependent Hamiltonian, the associated Hamiltonian vector field is defined by
$$\omega(X_H^t,\cdot) = -dH_t\,.$$
We call a smooth curve $\gamma \fc [0,1] \to M$ a \tb{Hamiltonian orbit} of $H$ with $L$ \tb{as the boundary condition} if $\gamma(0),\gamma(1) \in L$ and
$$\dot\gamma(t) = X^t_H(\gamma(t))\,.$$
Similarly if $H \fc S^1 \times M \to \R$ is time-periodic, a smooth loop $x \fc S^1 \to M$ is a \tb{periodic Hamiltonian orbit} of $H$ if
$$\dot x(t) = X_H^t(x(t))\,.$$

Let $\ol S{}^+ = [0,\infty] \times [0,1]$, and similarly define $\ol S{}^-$ and $\ol C{}^\pm$. We endow $\ol S{}^+$ with the unique structure of a smooth manifold with corners by declaring the map
$$\ol S{}^+ \to [0,1] \times [0,1]\,,\quad (s,t) \mapsto \Big(\frac{s}{\sqrt{1+s^2}},t\Big)$$
to be a diffeomorphism, and we do the same with $\ol S{}^-,\ol C{}^\pm$. If $\Sigma$ is a punctured Riemann surface endowed with a choice of ends for its punctures, we let $\ol \Sigma$ be the smooth manifold with corners obtained from $\Sigma$ by gluing $\ol S{}^\pm, \ol C{}^\pm$ along the ends $\epsilon_\theta$. Intuitively this amounts to compactifying $\Sigma$ by adding a copy of the interval $[0,1]$ for each boundary puncture and a copy of $S^1$ for each interior puncture.

A smooth map $f \fc \Sigma \to M$ is called \tb{b-smooth} \footnote{This concept is borrowed from \cite{Schwarz_PhD_thesis}, where it appears under a different name.} if it has a smooth extension to the whole of $\ol\Sigma$. This in particular means that $f$ extends continuously to the compactified surface $\ol\Sigma$ and that its derivatives (with respect to the $s$ variable on the ends) decay sufficiently rapidly at the punctures. We let
$$C^\infty_b(\Sigma,M)$$
be the set of b-smooth maps and
$$C^\infty_b(\Sigma,\partial\Sigma;M,L)$$
be the subset mapping $\partial\Sigma\to L$.

For $u\in C^\infty_b(\Sigma,\partial\Sigma;M,L)$ and a puncture $\theta$ of $\Sigma$ we let
$$u_\theta = \lim_{|s| \to \infty}u(\epsilon_\theta(s,\cdot))$$
be the limiting curve at $\theta$, defined either on $[0,1]$ or on $S^1$. Assume that we have chosen a time-dependent Hamiltonian $H^\theta$ for every puncture $\theta$, with the additional condition that $H^\theta$ is time-periodic if $\theta$ is interior, and a Hamiltonian orbit $y_\theta$ of $H^\theta$, which has $L$ as the boundary condition if $\theta$ is boundary or is a loop if $\theta$ is interior. We let
$$C^\infty_b(\Sigma,\partial\Sigma;M,L;\{y_\theta\}_\theta) = \{u \in C^\infty_b(\Sigma,\partial\Sigma;M,L)\,|\,u_\theta = y_\theta \text{ for all }\theta\}\,.$$
For $u \in C^\infty_b(\Sigma,\partial\Sigma;M,L;\{y_\theta\}_\theta)$, a fixed puncture $\theta$, and a Riemannian metric $\rho$ on $M$, which in case $\theta$ is boundary satisfies the additional condition that $L$ be totally geodesic with respect to it, we can express $u$ near $\theta$ using the exponential map of $\rho$. To illustrate, assume that $\theta$ is a positive interior puncture. Then there is a section $U \in C^\infty(C^+,\pi^*y_\theta^*TM)$ ($\pi \fc C^+ \to S^1$ the projection) such that
$$u(\epsilon_\theta(s,t)) = \exp_{\rho,y_\theta(t)}(U(s,t))$$
for all sufficiently large $s$.

\paragraph{Pregluing b-smooth maps}\label{par:pregluing_b_smooth_maps}

Let again $\cT$ be a gluing tree as in \S\ref{par:gluing_punctured_Riem_surf} with vertices labeled by $\Sigma_i$. Assume that for each $i$ and for each puncture of $\Sigma_i$ we have chosen a time-dependent Hamiltonian, which is time-periodic in case the puncture is interior, and a Hamiltonian orbit thereof, which has $L$ as the boundary condition if the puncture is boundary, and is a loop if the puncture is interior. Assume that if two punctures label an edge of $\cT$, then the corresponding Hamiltonians and Hamiltonian orbits coincide. Let $u_i \in C^\infty_b(\Sigma_i,\partial\Sigma_i; M,L)$ have these Hamiltonian orbits as asymptotics.

Using the expression of $u_i$ via exponential maps as above at each puncture undergoing gluing, we can piece together the corresponding vector fields with the help of cutoff functions to get a b-smooth map $u \in C^\infty_b(\Sigma,\partial\Sigma;M,L)$ defined on the surface $\Sigma$ obtained from the $\Sigma_i$ by gluing according to $\cT$, where we take the gluing lengths to be large enough. Note that $u$ has asymptotics dictated by the orbits corresponding to the punctures of $\Sigma$. We refer to a b-smooth map obtained in such a fashion from the maps $u_i$ as the result of \tb{pregluing} the $u_i$.

\subsection{Cauchy--Riemann operators associated to b-smooth maps}\label{ss:CROs_from_b_smooth_maps}

We will now describe how to construct an admissible Cauchy--Riemann operator starting from a b-smooth map asymptotic to Hamiltonian orbits. Let $\Sigma$ be a punctured Riemann surface endowed with a choice of ends. Assume that to each puncture $\theta$ we associate a Hamiltonian $H^\theta$ and an orbit $y_\theta$ of $H^\theta$. Fix
$$u \in C^\infty_b(\Sigma,\partial\Sigma;M,L;\{y_\theta\}_\theta)\,.$$
Let $E_u = u^*TM$, $F_u = (u|_{\partial\Sigma})^*TL$ and $\omega_u = u^*\omega$. Assume that for each puncture $\theta$ we have a family of almost complex structures $J^\theta$ on $M$ compatible with $\omega$, such that $J^\theta$ is parametrized by $[0,1]$ if $\theta$ is boundary, and by $S^1$ if it's interior. We let $J_u$ be any compatible almost complex structure on $E_u$ satisfying
$$J_u(\epsilon_\theta(s,t)) = J^\theta_t(u(\epsilon_\theta(s,t)))$$
for all $\theta$ and $(s,t)$.

For each puncture $\theta$ let $E_\theta = y_\theta^*TM$, $\omega_\theta = y_\theta^*\omega$, $J_\theta(t) = J^\theta_t(y_\theta(t))$, and let $\nabla_\theta$ be the symplectic connection whose parallel transport maps along $t$ are given by the linearized flow of $H^\theta$ along $y_\theta$. Finally, if $\theta$ is boundary, let $F_\theta = (y_\theta|_{\{0,1\}})^*TL$. Since $u$ is b-smooth, it extends to a smooth map $\ol u \fc \ol\Sigma \to M$, see \S\ref{ss:b_smooth_maps_pregluing}. Choose a connection $\nabla_u$ on $\ol u^*TM$ which is a smooth extension of the set of connections $\nabla_\theta$ over the limiting curves $u_\theta=y_\theta$.

We therefore have all the necessary data, $(E_u,F_u,\omega_u,J_u,\nabla_u)$, in order to define the associated admissible Cauchy--Riemann operator
$$D_u = \ol\partial_{\nabla_u}\,.$$

Let us say that a Hamiltonian orbit $\gamma \fc [0,1] \to M$ of $H$ with endpoints on $L$ is \tb{nondegenerate} if $\phi_{H*,\gamma(0)}^1(T_{\gamma(0)}L)$ is transverse to $T_{\gamma(1)}L$, and an orbit $x \fc S^1 \to M$ is nondegenerate if $\phi_{H*,x(0)}^1$ has no eigenvalue equal to $1$.

If all the orbits $y_\theta$ are nondegenerate, so are the connections $\nabla_\theta$ and therefore $D_u$ is nondegenerate. Thus it extends as a Fredholm operator to the Sobolev completions:
$$D_u \fc W^{1,p}(\Sigma,\partial\Sigma; E_u,F_u) \to L^p(\Sigma,\Omega_\Sigma^{0,1}\otimes E_u)\,.$$
We refer to $D_u$ as a \tb{formal linearization} of $D_u$. See Remark \ref{rem:relating_linearized_and_formal_linearized_ops} for a relation between this operator and the linearization of Floer's PDE.

\paragraph{Gluing Cauchy--Riemann operators and pregluing b-smooth maps} Now we describe how gluing Cauchy--Riemann operators on surfaces relates to pregluing b-smooth maps. We keep the notation of \S\ref{par:pregluing_b_smooth_maps}: $\cT$ is a gluing tree with vertices $\Sigma_i$, to each puncture of every $\Sigma_i$ there is associated a Hamiltonian and a nondegenerate orbit thereof, so that for every pair of punctures appearing as a label of an edge the corresponding Hamiltonians and orbits coincide. Finally, let $u_i \fc \Sigma_i \to M$ be b-smooth maps asymptotic to the chosen Hamiltonian orbits. We also have a preglued b-smooth map $u$ defined on the glued surface $\Sigma$.

We then have Cauchy--Riemann operators $D_i=D_{u_i}$ constructed as above, the operator $D$ glued from the $D_i$ as described in \S\ref{par:gluing_CROs_and_det_lines}, and the operator $D_u$, where the almost complex structure $J_u$ and the connection $\nabla_u$ coincide with the data $J_{u_i}$ and $\nabla_{u_i}$ outside the parts of $\Sigma_i$ participating in the gluing process. We can make an identification of the bundle $E$, obtained from gluing the bundles $E_{u_i}$, and the bundle $E_u$. The operator $D$ acting on $E$ and the operator $D_u$ acting on $E_u$ can now be deformed into one another relative to this identification, keeping the limiting data intact. This deformation induces an isomorphism
$$\ddd(D) \simeq \ddd(D_u)\,,$$
independent of the choices. Details are left to the reader. We refer to this isomorphism as the \tb{deformation isomorphism} below.

\subsection{The Floer PDE}\label{ss:Floer_PDE}

A \tb{Floer datum} associated to a puncture $\theta$ of $\Sigma$ is a pair $(H,J)$ where $H$ is a smooth time-dependent Hamiltonian on $M$, while $J$ is a smooth time-dependent family of $\omega$-compatible almost complex structures. Both $H,J$ are required to be time-periodic in case $\theta$ is interior. We call a Floer datum $(H,J)$ associated to a boundary puncture \tb{nondegenerate} if all the Hamiltonian orbits of $H$ with boundary on $L$ are nondegenerate. If the datum $(H,J)$ is associated to an interior puncture, it is called nondegenerate if all the periodic Hamiltonian orbits of $H$ are nondegenerate.

Fix a punctured Riemann surface $\Sigma$, a choice of ends $\{\epsilon_\theta\}_{\theta\in\Theta}$ for it, and a Floer datum $(H^\theta,J^\theta)$ associated to each puncture $\theta$. A \tb{perturbation datum} on $\Sigma$ is a pair $(K,I)$ where $K$ is a smooth $1$-form on $\Sigma$ with values in $C^\infty(M)$, satisfying the requirement that
$$K|_{\partial\Sigma}\quad\text{vanishes along} \quad L\,,$$
while $I$ is a family of compatible almost complex structures on $M$ parametrized by $\Sigma$. A perturbation datum $(K,I)$ is said to be \tb{compatible} with the Floer data $\{(H^\theta,J^\theta)\}_\theta$ associated to the punctures of $\Sigma$ if
$$\epsilon_\theta^*K = H^\theta_t\,dt\quad\text{ and }\quad I(\epsilon_\theta(s,t)) = J^\theta_t \quad \text{for all }(s,t)\text{ and }\theta\,.$$

We can now define the Floer PDE. Assume that $\Sigma$ is a punctured Riemann surface, and we have fixed a choice of ends for it, as well as Floer data associated to every puncture, and a compatible perturbation datum $(K,I)$. We let $X_K$ be the $1$-form on $\Sigma$ with values in Hamiltonian vector fields on $M$, defined via
$$\omega(X_K(\xi),\cdot) = - dK(\xi)\quad\text{ for }\xi \in T\Sigma\,.$$
The \tb{Floer PDE} is the equation
$$\ol\partial_{K,I}u:=(du - X_K)^{0,1} = 0$$
for $u \in C^\infty(\Sigma,\partial\Sigma; M,L)$.

\subsection{Moduli spaces of solutions of the Floer PDE}\label{ss:moduli_spaces_sols_Floer_PDE}

Here we define solutions spaces of the Floer PDE, the corresponding moduli spaces, and describe their compactness properties.

\subsubsection{Cappings and action functionals}\label{sss:cappings_action_fcnls}

Let $D^2 \subset \C$ be the closed unit disk and let $\dot D^2 = D^2 -\{1\}$ be the punctured disk where the puncture is considered to be positive. The standard end is defined by
\begin{equation}\label{eqn:std_end_cappings_Lagr_HF}
\epsilon_{\std}\fc S^+ \to \dot D^2\,, \quad \epsilon_{\std}(z) = \frac{e^{\pi z} - i}{e^{\pi z} + i}\,.
\end{equation}

Let $y \fc ([0,1],\{0,1\}) \to (M,L)$ be a smooth curve. By definition, a \tb{capping} of $y$ is a b-smooth map
$$\wh y \in C^\infty_b(\dot D^2,\partial \dot D^2 ; M,L;y)$$
with respect to the standard end. A capping has a canonical extension to the compactified disk where we glue in an interval $[0,1]$ at the puncture along the end $\epsilon_{\std}$. Two cappings $\wh y,\wh y'$ are \tb{equivalent} if the concatenation of the canonical extension of $\wh y$ and of the canonical extension of $-\wh y'$ defines a disk with boundary on $L$ representing the trivial class in $\pi_2(M,L)$, where
$$-\wh y'(a+ib) = \wh y'(-a+ib)\,.$$
We use the notation
$$\wt y = [y,\wh y]$$
to denote the equivalence class of cappings of $y$ containing $\wh y$.

Similarly let $\dot S^2$ be the sphere punctured once, with the puncture being positive. The standard end for it is defined by
$$\epsilon_{\std} \fc C^+ \to \dot S^2\,,\quad \epsilon_{\std}(z) = \frac{e^{2\pi z} - i}{e^{2 \pi z} + i}\,,$$
where we view $S^2 = \C P^1 = \C \cup \{\infty\}$. Let $y \fc S^1 \to M$ be a smooth loop. A capping of $y$ is a b-smooth map
$$\wh y \in C^\infty_b(\dot S^2,M;y)$$
relative to the standard end. Such a capping has a canonical extension to the compactified sphere where we add a circle at infinity according to the standard end. We call two cappings $\wh y,\wh y'$ equivalent if the concatenation of the canonical extension of $\wh y$ and the canonical extension of $-\wh y'$ is a contractible sphere. Again, we denote $\wt y =[y,\wh y]$ the equivalence class of cappings containing $\wh y$.

We now define the action functionals. Let
$$\Omega_L = \{y \fc [0,1] \to M\,|\, [y] = 0 \in \pi_1(M,L)\}$$
and let
$$\wt\Omega_L = \{\wt y=[y,\wh y]\,|\, \wh y\text{ is a capping of }y\}\,.$$
We denote by $p \fc \wt \Omega_L \to \Omega_L$ the obvious projection. Let $H$ be a time-dependent Hamiltonian. The action functional associated to $H$ is
$$\cA_{H:L} \fc \wt\Omega_L \to \R\,,\quad \cA_{H:L}([y,\wh y]) = \int_0^1 H_t(y(t))\,dt-\int \wh y^*\omega\,.$$
A point $\wt y$ is critical for $\cA_{H:L}$ if and only if $y = p(\wt y)$ is a Hamiltonian orbit of $H$. If all the orbits of $H$ with boundary on $L$ are nondegenerate, there is a Conley-Zehnder index
$$m_{H:L} \fc \Crit \cA_{H:L} \to \Z\,.$$
We use the definition of this index, which satisfies the following shift property: for two cappings $\wh y,\wh y'$ of the same orbit $y$ we have
$$m_{H:L}([y,\wh y]) - m_{H:L}([y,\wh y']) = -\mu(\wh y\sharp -\wh y')\,.$$
In addition, we normalize it as follows. Assume $f$ is a $C^2$-small Morse function on $L$, and that we have identified a neighborhood of $L$ with a neighborhood of the zero section in $T^*L$; let $H$ be obtained by cutting off the pullback of $f$ to $T^*L$ outside the neighborhood. Then the Hamiltonian orbits of $H$ with boundary on $L$ are precisely the constant curves at the critical points of $f$. We let the Conley-Zehnder index of such an orbit together with the constant capping be equal the Morse index of the corresponding critical point. 

We use the abbreviated notations
$$|\wt y| := m_{H:L}(\wt y)\,,\quad |\wt y|':=n-m_{H:L}(\wt y)\,,$$
which should cause no confusion.

Similarly, let
$$\Omega = \{y \fc S^1 \to M\,|\,[y]=0\in\pi_1(M)\}$$
and
$$\wt \Omega = \{\wt y=[y,\wh y]\,|\, \wh y\text{ is a capping of }y\}\,.$$
We have the projection $p \fc \wt\Omega \to \Omega$. Let $H$ be a time-periodic Hamiltonian. The associated action functional is
$$\cA_H \fc \wt \Omega \to \R\,,\quad \cA_H([y,\wh y]) = \int_{S^1} H_t(y(t))\,dt - \int \wh y^*\omega\,.$$
Its critical points are $\wt y$ with $y=p(\wt y)$ being a periodic Hamiltonian orbit of $H$. If all the periodic orbits of $H$ are nondegenerate, there is a Conley-Zehnder index
$$m_{H} \fc \Crit \cA_{H} \to \Z\,.$$
This has the shift property
$$m_H([y,\wh y]) - m_H([y,\wh y']) = -2c_1(\wh y \sharp - \wh y')\,,$$
and is normalized to coincide with the Morse index of a critical point of $H$ provided $H$ is $C^2$-small, autonomous, and Morse, and the critical point is considered as a constant periodic orbit taken with the constant capping.

We use the abbreviated notations
$$|\wt y| := m_{H}(\wt y)\,,\quad |\wt y|':=2n-m_{H:L}(\wt y)\,.$$
\begin{rem}Below we have to use action functionals of both types, $\cA_{H:L}$ and $\cA_{H}$, on almost equal footing. In order to keep the notation less cumbersome, we will oftentimes denote both of them just by $\cA_{H}$. The context will make it clear on which space ($\wt\Omega$ or $\wt\Omega_L$) it is defined.
\end{rem}

\subsubsection{Solution spaces and moduli spaces}\label{sss:solution_spaces_moduli_spaces}

Until the end of \S\ref{s:HF}, we assume that $\wh\Sigma$ is either the sphere or the closed disk, and that there is a unique positive puncture $\theta$ and negative punctures $\{\theta_i\}_{i=1}^k$ (possibly $k=0$). Recall that $\Sigma = \wh \Sigma - \Theta$. Fix a choice of ends for $\Sigma$, nondegenerate Floer data $(H,J)$ associated to $\theta$ and $\{(H^i,J^i)\}_i$ associated to $\theta_i$, and a compatible perturbation datum $(K,I)$. Choose critical points $\wt y,\wt y_i$ of the action functionals corresponding to $H$ and $H^i$. We define the solution space
$$\cM_\Sigma(K,I;\{\wt y_i\}_i,\wt y) = \{u \in C^\infty_b(\Sigma,\partial\Sigma; M,L;\{y_i\}_i,y)\,|\, u\sharp \wh y_1\sharp\dots\sharp \wh y_k \in \wt y\,,\ol\partial_{K,I}u = 0\}\,,$$
that is the set of solutions of Floer's PDE with boundary conditions on $L$, asymptotics given by the orbits $y_i,y$ and homotopy class coming from the chosen equivalence classes of cappings for the orbits. Here $u\sharp \wh y_1 \sharp\dots \sharp \wh y_k$ denotes a map obtained by pregluing the b-smooth maps $u,\wh y_1,\dots,\wh y_k$ according to the obvious gluing tree. To every $u \in \cM_\Sigma(K,I;\{\wt y_i\},\wt y)$ there is associated the linearized operator (see \cite[Chapter 9]{Seidel_The_Book})
$$D_u \fc W^{1,p}(\Sigma,\partial \Sigma; E_u,F_u) \to L^p(\Sigma, \Omega^{01,}_\Sigma \otimes E_u)\,.$$
\begin{rem}\label{rem:relating_linearized_and_formal_linearized_ops}
In \cite[Chapter 8]{Seidel_The_Book} it is explained how to choose a connection $\nabla_u$ on $E_u$ for which this linearized operator coincides with the formal linearized operator $\nabla_u^{0,1}$ introduced \S\ref{ss:CROs_from_b_smooth_maps} for b-smooth maps. In the sequel we shall always reserve the notation $D_u$ to mean the linearized operator in case $u$ is a solution of the Floer PDE.
\end{rem}

The perturbation datum $(K,I)$ is called \tb{regular} if for every choice of critical points $\wt y_i,\wt y$ and every $u \in \cM_\Sigma(K,I;\{\wt y_i\}_i,\wt y)$ the operator $D_u$ is onto. In this case $\cM_\Sigma(K,I;\{\wt y_i\}_i,\wt y)$ is naturally a smooth manifold of dimension
$$\dim \cM_\Sigma(K,I;\{\wt y_i\}_i,\wt y) = |\wt y|' - \sum_{i=1}^k|\wt y_i|'\,.$$
The set of regular perturbation data $(K,I)$ compatible with the given Floer data is dense in the set of all compatible perturbation data \cite{Seidel_The_Book}.

We single out the special case of a translation-invariant perturbation datum: assume $\Sigma = S$ or $C$ and assume the perturbation datum has the form:
$$K(s,t) = H_t\,dt\,,I(s,t) = J_t$$
for all $(s,t)$. In this case the Floer PDE is the original equation for the negative gradient flow of the action functional, to wit
$$\langle \ol\partial_{K,I}u, \partial_s\rangle = \partial_s u + J(u) \big( \partial_t u - X_H(u) \big) = 0\,.$$
If $\wt y_\pm \in \Crit \cA_H$, we let
$$\wt\cM(H,J;\wt y_-,\wt y_+) := \cM_\Sigma(K,I;\wt y_-,\wt y_+)\,.$$
This space admits a natural $\R$-action by translation in the $s$ variable. We let
$$\cM(H,J;\wt y_-,\wt y_+)$$
denote the quotient if the action is free, otherwise we declare it to be empty. We call the Floer datum $(H,J)$ \tb{regular} if the corresponding translation-invariant perturbation datum is regular. In this case $\wt\cM(H,J;\wt y_-,\wt y_+)$ is a smooth manifold of dimension
$$|\wt y_+|' - |\wt y_-|' = |\wt y_-| - |\wt y_+|\,,$$
while $\cM(H,J;\wt y_-,\wt y_+)$ has dimension $|\wt y_-| - |\wt y_+| - 1$. We note that the set of regular data is dense in the set of all nondegenerate Floer data \cite{Floer_Hofer_Salamon_Transversality}.
\begin{rem}
Note that by definition a regular Floer datum is nondegenerate.
\end{rem}

\paragraph{Parametrized solution spaces}\label{par:parametrized_solution_spaces}

Next we treat families of surfaces and the Floer PDE on them. Assume for the moment that we have a smooth compact connected oriented surface with boundary $\wh\Sigma$ equipped with a finite set of punctures $\Theta$ and let $\Sigma= \wh\Sigma - \Theta$. Let $\cS \to \cR$ be a fiber bundle with fiber $\Sigma$, whose structure group is the group of orientation-preserving diffeomorphisms of $\wh\Sigma$ which are the identity on $\Theta$. In this case $\cS \to \cR$ extends to a fiber bundle $\wh\cS \to \cR$ with fiber $\wh\Sigma$ and the punctures of $\Sigma$ give rise to canonical smooth sections $\cR \to \wh\cS$. We identify a puncture with the corresponding section. We denote the fiber of $\cS \to \cR$ over $r \in \cR$ by $\Sigma_r$.

Assume there is a smooth family of conformal structures on the fibers of $\cS \to \cR$. A choice of ends for $\cS \to \cR$ is a family of fiberwise maps
$$\epsilon_\theta \fc \cR \times S^\pm \to \cS \quad \text{ or }\quad \epsilon_\theta \fc \cR \times C^\pm \to \cS$$
whose restrictions to the fibers over $r$ constitute a choise of ends for the Riemann surface $\Sigma_r$, such that $\epsilon_\theta$ is asymptotic to $\theta$.

Assume that we have fixed a choice of ends for $\cS$ and nondegenerate Floer data $\{(H^\theta,J^\theta)\}_\theta$ associated to the punctures of $\Sigma$. A perturbation datum on $\cS$ is a pair $(K,I)$ where $I$ is a family of compatible almost complex structures parametrized by $\cS$, while $K$ is a smooth family of $1$-forms on the vertical tangent bundle of $\cS \to \cR$ such that the restriction $(K_r,I_r)$ is a perturbation datum on $\Sigma_r$. The notion of compatibility of the perturbation datum with the Floer data extends to families in an obvious manner.

Let now $\wh\Sigma$ be either the sphere or the closed disk, and assume that it carries a unique positive puncture $\theta$ and negative punctures $\{\theta_i\}_{i=1}^k$. Fix a choice of ends for $\cS$, Floer data $(H,J),\{(H^i,J^i)\}_i$ associated to punctures of $\Sigma$ and a compatible perturbation datum $(K,I)$. Fix also critical points $\wt y,\wt y_i$ of the action functionals of $H$ and $H_i$. We define
$$\cM_{\cS}(K,I;\{\wt y_i\}_i,\wt y) = \{(r,u)\,|\,r\in\cR\,, u \in \cM_{\Sigma_r}(K_r,I_r;\{\wt y_i\}_i, \wt y)\}\,.$$
For $(r,u) \in \cM_{\cS}(K,I;\{\wt y_i\}_i,\wt y)$ we have the \tb{extended linearized operator}
$$D_{r,u} \fc T_r\cR \times W^{1,p}(\Sigma_r,\partial\Sigma_r; E_u,F_u) \to L^p(\Sigma_r,\Omega_{\Sigma_r}^{0,1}\otimes E_u)\,.$$
See \cite{Seidel_The_Book} for the precise definition. Note for future use that the restriction
$$D_{r,u}|_{0 \times W^{1,p}(\Sigma_r,\partial\Sigma_r; E_u,F_u)}$$
coincides with the linearized operator $D_u$ of $u \in \cM_{\Sigma_r}(K_r,I_r;\{\wt y_i\}_i,\wt y)$. We call the perturbation datum $(K,I)$ for $\cS$ \tb{regular} if for every choice of critical points $\wt y_i,\wt y$ and every $(r,u) \in \cM_{\cS}(K,I;\{\wt y_i\}_i,\wt y)$ the extended operator $D_{r,u}$ is onto. In this case $\cM_{\cS}(K,I;\{\wt y_i\}_i,\wt y)$ is a smooth manifold of dimension
$$\dim \cM_{\cS}(K,I;\{\wt y_i\}_i,\wt y) = \dim \cR + |\wt y|' - \sum_{i=1}^k|\wt y_i|'\,.$$
The set of regular perturbation data compatible with the given Floer data is dense in the set of all compatible perturbation data.

\begin{rem}
We note here that a regular solution of the Floer PDE is necessarily a b-smooth map \cite{Schwarz_PhD_thesis}. Therefore all the constructions regarding b-smooth maps in \S\ref{ss:b_smooth_maps_pregluing}, \S\ref{ss:CROs_from_b_smooth_maps} apply to them and their linearized operators.
\end{rem}

\subsubsection{Compactness and gluing}\label{sss:compactness_gluing}

Here we discuss the relevant compactness and gluing results for the above solution spaces. The general statement is that whenever a moduli space is zero-dimensional, it is compact, therefore a finite number of points, whereas when it is one-dimensional, it can be compactified into a compact $1$-dimensional manifold with boundary, where the boundary consists either of boundary points already present in the moduli space, or else of pairs of elements of $0$-dimensional moduli spaces. Moreover, the converse to compactness, called gluing, states that all suitable pairs are obtained in this way.

We start with the description of the relevant $0$-dimensional moduli spaces. There are three basic types of moduli spaces used in Floer homology:
\begin{enumerate}
 \item when the surface is either a strip or a cylinder and the perturbation datum is translation-invariant --- this leads to the definition of boundary operators;
 \item there is a single surface --- this is used to define various operations on Floer homology;
 \item there is a family of surfaces $\cS \to \cR$ with $\cR$ being $1$-dimensional --- this leads to relations between the operations and the boundary operators, such as chain homotopies and various algebraic identities.
\end{enumerate}

\paragraph{The case of translation-invariant perturbation datum on $S$ or $C$}

We treat the case $\Sigma = S$, the case of the cylinder $C$ being entirely similar. We fix a regular Floer datum $(H,J)$. The set $\wt\cM(H,J;\wt y_-,\wt y_+)$ is a smooth manifold of dimension
$$|\wt y_-| - |\wt y_+|\,.$$
When this difference is $1$, the quotient manifold $\cM(H,J;\wt y_-,\wt y_+)$ is $0$-dimensional and compact, therefore a finite set of points. When the difference is $2$, there are two cases:

Case I: $N_L \geq 3$ or $y_- \neq y_+$, and

Case II: $N_L = 2$ and $y_- = y_+$.

\noindent In case I, the manifold $\cM(H,J;\wt y_-,\wt y_+)$ admits a compactification $\ol\cM(H,J;\wt y_-,\wt y_+)$ whose boundary is
$$\partial\ol \cM(H,J;\wt y_-,\wt y_+) = \bigcup_{\wt y \in \Crit\cA_{H:L}} \cM(H,J;\wt y_-,\wt y) \times \cM(H,J;\wt y,\wt y_+)\,,$$
that is the only way to noncompactness is through Floer breaking.

Consider now case II. For $q \in L$, $A \in \pi_2(M,L,q)$, and an almost complex structure $J$, define
$$\wt \cM_1(J;q,A) = \{u \in C^\infty(D^2,S^1,1;M,L,q)\,|\, \ol \partial_J u = 0, [u] = A\}\,,$$
where $\ol\partial_J = \ol \partial_{0,J}$ is the Floer operator of the perturbation datum $(0,J)$ on $D^2$, that is with vanishing Hamiltonian term. We let $\cM_1(J;q,A)$ be the quotient of $\wt\cM_1(J;q,A)$ by the conformal automorphism group of $D^2$ preserving $1 \in S^1$. We call $J$ \tb{regular} if for every $q,A$, and $u \in \wt \cM_1(J;q,A)$ the linearized operator $D_u$ is onto. Since $L$ is monotone and $N_L \geq 2$, it follows that all disks in $\wt\cM_1(J;q,A)$ are simple \cite{Biran_Cornea_Quantum_structures_Lagr_submfds, Biran_Cornea_Rigidity_uniruling}, therefore the set of regular $J$ is dense. If $J$ is regular, $\cM_1(J;q,A)$ is a zero-dimensional manifold for generic $q$. The set of regular Floer data $(H,J)$ for which $J_0,J_1$ are regular is dense. If $(H,J)$ is such, we have in case II:
\begin{align*}
\partial\ol \cM(H,J;\wt y_-,\wt y_+) &= \bigcup_{\wt z \in\Crit\cA_{H:L}} \cM(H,J;\wt y_-, \wt z) \times \cM(H,J;\wt z,\wt y_+)\\
&\cup \cM_1(J_0; y(0), [\wh y_+ \sharp - \wh y_-]) \cup \cM_1(J_1; y(1), [\wh y_+ \sharp - \wh y_-])\,,
\end{align*}
where $y = y_- = y_+$. This means that in this particular case another possibility for noncompactness opens up, that of bubbling off of Maslov $2$ disks attached to endpoints of the Hamiltonian orbit $y$.

Bubbling off of holomorphic spheres of Chern number $1$ is also possible, however the set of points through which they pass has high codimension and therefore generically they do not appear in the boundary of the moduli space \cite{Hofer_Salamon_HF_Nov_rings}. Also see \cite{Hu_Lalonde_Relative_Seidel_morphism_Albers_map} .

The treatment in the case $\Sigma = C$ is entirely analogous, with the difference that the Floer datum $(H,J)$ is $1$-periodic in $t$ and for a generic datum there is no bubbling, the noncompactness being only due to Floer breaking.

\paragraph{The case of a single surface}

Assume $\wh\Sigma$ is the sphere or the closed disk and endow it with punctures, where exactly one puncture $\theta$ is positive, the other punctures $\{\theta_i\}_{i=1}^k$ being negative; let $\Sigma$ be the resulting punctured Riemann surface. Endow $\Sigma$ with a choice of ends, and fix regular Floer data $(H^i,J^i)$ and $(H,J)$ associated to the punctures $\theta_i$ and $\theta$. Fix also a regular perturbation datum $(K,I)$ compatible with the Floer data, and critical points $\wt y_i, \wt y$ of the action functionals of $H^i$, $H$, and consider
$$\cM_\Sigma(K,I;\{\wt y_i\}_i,\wt y)\,.$$
This is a smooth manifold of dimension $|\wt y|' - \sum_i |\wt y_i|'$. If this dimension is $0$, $\cM_\Sigma(K,I;\{\wt y_i\}_i,\wt y)$ is compact and therefore a finite number of points.

When $\cM_\Sigma(K,I;\{\wt y_i\}_i,\wt y)$ is $1$-dimensional, the only possible noncompactness is due to Floer breaking, namely we have
{
\begin{multline*}
\partial \ol \cM_\Sigma(K,I;\{\wt y_i\}_i,\wt y) =  \bigcup_{j=1}^{k} \bigcup_{\wt y_j' \in \Crit\cA_{H^j}}\cM(H^j,J^j;\wt y_j,\wt y_j') \times \cM_\Sigma(K,I;\{\wt y_i\}_{i\neq j},\wt y_j',\wt y) \\
\cup \bigcup_{\wt y' \in \Crit \cA_{H}}\cM_\Sigma(K,I;\{\wt y_i\}_i,\wt y') \times \cM(H,J;\wt y',\wt y)\,.
\end{multline*}
}

\paragraph{The case of a $1$-dimensional family}\label{par:compactness_families}

Here we have a family of punctured surfaces $\cS \to \cR$ as in \S\ref{par:parametrized_solution_spaces}. We will only need the cases when $\cR = [0,1]$ or $\cR=[0,\infty)$.

Assume first that $\cR = [0,1]$, that $\cS = \Sigma \times [0,1]$, where $\Sigma$ is either a punctured sphere or a punctured disk, that the punctures are $\Theta=\{\theta,\theta_1\dots,\theta_k\}$, where $\theta$ is the unique positive puncture, that we have fixed a choice of ends for $\cS$ which are constant near the boundary of $\cR$, and regular Floer data $(H,J)$ and $(H^i,J^i)$ associated to $\theta$ and $\theta_i$. Let $(K,I)$ be a regular compatible perturbation datum, which is constant near the boundary of $\cR$. Fix critical points $\wt y,\wt y_i$ of the action functionals of $H,H^i$. The set
$$\cM_\cS(K,I;\{\wt y_i\}_i,\wt y)$$
is a smooth manifold of dimension
$$\dim \cM_\cS(K,I;\{\wt y_i\}_i,\wt y) = 1 + |\wt y|' - \sum_i|\wt y_i|'\,.$$
When this dimension is zero, $\cM_\cS(K,I;\{\wt y_i\}_i,\wt y)$ is compact and therefore a finite number of points.

When it equals $1$, $\cM_\cS(K,I;\{\wt y_i\}_i,\wt y)$ can be compactified by adding Floer breaking, that is we have
\begin{multline*}
\partial\ol\cM_\cS(K,I;\{\wt y_i\}_i,\wt y) = \{0\}\times \cM_{\Sigma_0}(K_0,I_0;\{\wt y_i\}_i,\wt y) \cup \{1\} \times \cM_{\Sigma_1}(K_1,I_1;\{\wt y_i\}_i,\wt y) \\
\cup \bigcup_{j=1}^k \bigcup_{\wt y_j' \in \Crit\cA_{H^j}} \cM(H^j,J^j;\wt y_j,\wt y_j') \times \cM_{\cS}(K,I;\{\wt y_i\}_{i\neq j},\wt y_j',\wt y) \\
\cup \bigcup_{\wt y' \in \Crit\cA_{H}} \cM_{\cS}(K,I;\{\wt y_i\}_i,\wt y') \times \cM(H,J;\wt y',\wt y)\,.
\end{multline*}

When $\cR = [0,\infty)$, we require that the family $\cS$ and the choice of ends on it have a specific form. Namely, let $\Sigma^1,\Sigma^2$ be two punctured Riemann surfaces with puncture sets $\Theta_i = \{\theta^i,\theta^i_1,\dots,\theta^i_{k_i}\}$, $\theta^i$ being positive and the rest being negative. Fix a choice of ends for the $\Sigma^i$. Let $R_0 > 0$ and for $r \geq R_0$ let $\Sigma_r$ be obtained from gluing $\Sigma^1,\Sigma^2$, where the tree has two vertices corresponding to the surfaces $\Sigma^1,\Sigma^2$ and the unique edge between them is labeled by $(\theta^1,\theta^2_j)$ and the gluing length is $r$. Note that $\Sigma_r$ has $\theta^2,\{\theta^1_i\}_i,\{\theta^2_i\}_{i \neq j}$ as punctures. We require the family $\cS \to \cR = [0,\infty)$ to have fiber $\Sigma_r$ for $r \geq R_0$ and to have conformal structures and choices of ends to come from gluing as described.

Similarly, the choice of perturbation datum on $\cS$ comes from gluing. More precisely, assume we have chosen regular Floer data associated to the punctures of the $\Sigma^i$: $(H^i,J^i)$ for $\theta^i$ and $(H^{i,l},J^{i,l})_l$ for $\theta^i_l$, , $i=1,2$, such that $H^1 = H^{2,j}$, $J^1 = J^{2,j}$. Let $(K^i,I^i)$ be regular compatible perturbation data on $\Sigma^i$. There is an obvious perturbation datum $(K_r,I_r)$ on the glued surface $\Sigma_r$, since the perturbation data on the $\Sigma^i$ agree on the overlap, and we require the perturbation datum on $\cS$ to equal $(K_r,I_r)$ on the fibers $\Sigma_r$ over $r \geq R_0$. Note that the Floer data associated to punctures of $\Sigma_r$ (for all $r$) are $(H^2,J^2)$ for $\theta^2$ and $\{(H^{1,i},J^{1,i})\}_i$, $\{(H^{2,i},J^{2,i})\}_{i\neq j}$ for $\{\theta^1_i\}_i,\{\theta^2_i\}_{i \neq j}$.

Let therefore $\cS \to \cR = [0,\infty)$ be such a family. Endow $\cS$ with a choice of ends as above, with additional condition that they are locally constant near $0 \in \cR$. We already have a choice of Floer data for the punctures of $\cS$, which we now assume to be regular, and we let $(K,I)$ be a regular compatible perturbation datum, locally constant near $0\in\cR$, and which has the aforementioned form for $r \geq R_0$. The set of such $(K,I)$ is dense. Fix critical points $\wt y_2, \{\wt y_{1,i}\}_i,\{\wt y_{2,i}\}_{i \neq j}$ of the action functionals of $H^2,\{H^{1,i}\}_i,\{H^{2,i}\}_{i\neq j}$. Then the set
$$\cM_\cS(K,I;\{\wt y_{1,i}\}_i,\{\wt y_{2,i}\}_{i\neq j}, \wt y_2)$$
is a smooth manifold of dimension
$$\dim \cM_\cS(K,I;\{\wt y_{1,i}\}_i,\{\wt y_{2,i}\}_{i\neq j}, \wt y_2) = 1 + |\wt y_2|' - \sum_i|\wt y_{1,i}|' - \sum_{i\neq j}|\wt y_{2,i}|'\,.$$
When this dimension is zero, $\cM_\cS(K,I;\{\wt y_{1,i}\}_i,\{\wt y_{2,i}\}_{i\neq j}, \wt y_2)$ is compact and therefore a finite number of points. When it equals $1$, $\cM_\cS(K,I;\{\wt y_{1,i}\}_i,\{\wt y_{2,i}\}_{i\neq j}, \wt y_2)$ can be compactified by adding Floer breaking and breaking at the noncompact end of $\cR$, that is we have
\begin{multline*}
\partial\ol \cM_\cS(K,I;\{\wt y_{1,i}\}_i,\{\wt y_{2,i}\}_{i\neq j}, \wt y_2) = 
\{0\} \times \cM_{\Sigma_0}(K_0,I_0;\{\wt y_{1,i}\}_i,\{\wt y_{2,i}\}_{i\neq j}, \wt y_2) \\
\cup \bigcup_{\wt y_{2,j} \in \Crit\cA_{H^1}} \cM_{\Sigma^1}(K^1,J^1;\{\wt y_{1,i}\}_i,\wt y_{2,j}) \times \cM_{\Sigma^2}(K^2,I^2;\{\wt y_{2,i}\}_i,\wt y_2) \\
\cup \bigcup_{\wt y_2' \in \Crit \cA_{H^2}} \cM_{\cS}(K,I;\{\wt y_{1,i}\}_i,\{\wt y_{2,i}\}_{i\neq j}, \wt y_2') \times \cM(H^2,J^2;\wt y_2',\wt y_2)\\
\cup \bigcup_l \bigcup_{\wt y_{1,l}' \in \Crit \cA_{H^{1,l}}} \cM(H^{1,l},J^{1,l};\wt y_{1,l},\wt y_{1,l}') \times \cM_{\cS}(K,I;\{\wt y_{1,i}\}_{i\neq l},\wt y_{1,l}',\{\wt y_{2,i}\}_{i\neq j}, \wt y_2)\\
\cup \bigcup_{l\neq j} \bigcup_{\wt y_{2,l}' \in \Crit \cA_{H^{2,l}}} \cM(H^{2,l},J^{2,l};\wt y_{2,l},\wt y_{2,l}') \times \cM_{\cS}(K,I;\{\wt y_{1,i}\}_i,\{\wt y_{2,i}\}_{i\neq j,l}, \wt y_{2,l}', \wt y_2)\,.
\end{multline*}

\subsection{Orientations}\label{ss:orientatiosn}

\subsubsection{Canonical $\Z$-modules associated to critical points of action functionals}\label{sss:can_Z_modules_crit_pts_action_fcnls}

Recall the definition of a capping for a smooth curve $y \fc ([0,1],\{0,1\}) \to (M,L)$, \S\ref{sss:cappings_action_fcnls}. Fix an equivalence class of cappings $\wt y$ of $y$ and regard it as a topological space. Since a capping is a b-smooth map $\dot D^2 \to M$, it has a canonical extension to a continuous map defined on the compactification of $\dot D^2$ obtained by gluing a copy of $[0,1]$ along the standard end, \S\ref{ss:b_smooth_maps_pregluing}. This compactification is diffeomorphic to $D^2 \cap \{\RE z \leq 0\}$, and we view the canonical extension of $\wh y$ as a map $D^2 \cap \{\RE z \leq 0\} \to M$.

Let $C_{\wt y}$ be the space of continuous maps $D^2 \cap \{\RE z \leq 0\} \to M$ mapping the semicircle to $L$, the diameter to $y$, and belonging to the homotopy class dictated by $\wt y$. It is easy to see that the map taking a capping to the corresponding continuous extension is a homotopy equivalence between $\wt y$ and $C_{\wt y}$. We have the following obvious lemma.
\begin{lemma}\label{lem:fund_gp_space_of_cappings_isomorphic_to_rel_pi_3}
Fix $\wh y_0 \in C_{\wt y}$ and let $-\wh y_0 \fc D^2 \cap \{\RE z \geq 0\} \to M$ be defined by $-\wh y_0(s,t) = \wh y_0(-s,t)$. Then the concatenation map
$$C_{\wt y} \to \{w \fc (D^2,S^1,-i) \to (M,L,y(0)) \,|\, [w] = 0 \in \pi_2(M,L,y(0))\}\,, \quad  \wh y \mapsto \wh y \sharp - \wh y_0$$
is well-defined and is a homotopy equivalence. In particular the fundamental group $\pi_1(C_{\wt y}, \wh y_0)$ is isomorphic to $\pi_3(M,L,y(0))$, and since it is abelian, the isomorphism is independent of $\wh y_0$. \qed
\end{lemma}
\noindent This means that we have canonically identified the fundamental group of the space of cappings $\wt y$ with $\pi_3(M,L,y(0))$.

Let now $H \fc [0,1] \times M \to \R$ be a nondegenerate Hamiltonian, that is all its orbits with boundary on $L$ are nondegenerate, and fix an orbit $y$ of $H$. Pick an equivalence class of cappings $\wt y = [y,\wh y]$, that is a critical point of $\cA_{H:L}$. For any $\wh y$ we have an associated linearized operator $D_{\wh y}$, see \S\ref{ss:CROs_from_b_smooth_maps}. Its construction depends on various choices, such as an almost complex structure and a connection on the pullback bundle $\wh y^*TM$. Let $D_{\wt y}$ denote the collection of all the operators obtained in this way for all the cappings $\wh y$ in class $\wt y$ and all the auxiliary choices. Since the spaces of almost complex structures and connections are contractible, we see that the parameter space of $D_{\wt y}$ is homotopy equivalent to the space $\wt y$ of cappings in class $\wt y$. Therefore its fundamental group is canonically isomorphic to $\pi_3(M,L,y(0))$ by Lemma \ref{lem:fund_gp_space_of_cappings_isomorphic_to_rel_pi_3}. The determinant lines of the operators in the family $D_{\wt y}$ glue into the line bundle $\ddd (D_{\wt y})$. We have
\begin{lemma}\label{lem:first_Stiefel_Whitney_dD_wt_y_HF}
Relative to the canonical isomorphism of $\pi_3(M,L,y(0))$ with the fundamental group of the space of parameters over which the line bundle $\ddd(D_{\wt y})$ is defined, its first Stiefel--Whitney class equals
$$w_1(\ddd(D_{\wt y}))=w_2(TL) \circ \partial \fc \pi_3(M,L) \to \Z_2\,.$$
\end{lemma}
\begin{prf}
Let $(\wh y_\tau)_{\tau \in S^1}$ be a loop of cappings and lift it to a loop of operators $D_\tau = D_{\wh y_\tau}$. Recall that for a fixed capping $\wh y_0 \in \wt y$ we defined the reverse capping $- \wh y_0$ via $-\wh y_0(s,t) = \wh y_0(-s,t)$. To it there corresponds a Cauchy--Riemann operator $D_{-\wh y_0}$ on $D^2 - \{-1\}$. The latter surface has a negative end, and therefore we can form the glued operator $D_\tau \sharp D_{-\wh y_0}$ over $D^2$ for some gluing length. Combining the direct sum and the gluing isomorphisms, we obtain an isomorphism
$$\ddd(D_\tau) \otimes \ddd(D_{-\wh y_0}) \simeq \ddd(D_\tau \sharp D_{-\wh y_0})\,,$$
which is continuous in $\tau$, and which implies that the loop $(D_\tau)_\tau$ is orientable if and only if the loop $D_\tau \sharp D_{-\wh y_0}$ is. Let us therefore compute the first Stiefel--Whitney class of $\ddd(D_\tau \sharp D_{-\wh y_0})$ over $S^1$. The $w_1$ of a loop of Cauchy--Riemann operators on a compact surface has been computed, see for instance \cite{Seidel_The_Book, Georgieva_Orientability_problem_open_GW}. Since in our case the boundary condition of the loop of operators $\ddd(D_\tau \sharp D_{-\wh y_0})$ is stationary on the right half of the disk, we obtain
$$w_1\big((\ddd(D_\tau \sharp D_{-\wh y_0}))_\tau\big) = \langle w_2(TL), [U] \rangle\,,$$
where $U$ is the image in $L$ of the evaluation map
$$\partial D^2 \times S^1 \to L\,,\quad (\sigma,\tau) \mapsto (\wh y_\tau \sharp -\wh y_0)(\sigma)\,.$$
Unraveling the definitions, we see that the number $\langle w_2(TL), [U] \rangle$ equals $w_2(TL)\circ \partial$ evaluated on the loop in the space of contractible disks at $y(0)$ given by $(\wh y_\tau \sharp -\wh y_0)_{\tau \in S^1}$. \qed

\end{prf}

Due to assumption \tb{(O)}, we obtain that the determinant line $\ddd(D_{\wt y})$ is orientable. We let
$$C(\wt y)$$
be the free $\Z$-module of rank $1$ whose two generators are the two possible orientations of this determinant line bundle. Note that since $\wt y$ is connected, this definition makes sense.

Similarly, recall the definition of a capping for a smooth loop $y \fc S^1 \to M$. Pick a nondegenerate Hamiltonian $H \fc S^1 \times M \to \R$, that is all of its Hamiltonian orbits which are loops are nondegenerate. For a capping $\wh y$ of a periodic orbit $y$ of $H$ we have defined a linearized operator $D_{\wh y}$, see \S\ref{ss:CROs_from_b_smooth_maps}. Similarly to the above, we have the family $D_{\wt y}$ of all the linearized operators associated to cappings in a given class $\wt y$. We have
\begin{lemma}The determinant line $\ddd(D_{\wt y})$ is orientable.
\end{lemma}
\begin{prf}
It suffices to prove that the determinant bundle of a loop $D_\tau = D_{\wh y_\tau}$ of operators above a loop of cappings is orientable. Fix a capping $\wh y_0 \in \wt y$. We can form the glued operator $D_\tau \sharp D_{-\wh y_0}$ over $S^2$. The direct sum and the gluing isomorphisms combine to the isomorphism
$$\ddd(D_\tau) \otimes \ddd(D_{-\wh y_0}) \simeq \ddd(D_\tau \sharp D_{-\wh y_0})\,,$$
which is continuous in $\tau$, and which implies that we only have to prove the orientability of the determinant bundle $(\ddd(D_\tau \sharp D_{-\wh y_0}))_\tau$. However this latter bundle is a bundle of real Cauchy--Riemann operators on $S^2$, which is a closed Riemann surface. It is well-known that the set of real Cauchy--Riemann operators on a Hermitian vector bundle over a closed Riemann surface deformation retracts onto the subspace of complex-linear Cauchy--Riemann operators. The determinant lines of the complex linear operators are canonically oriented since their kernels and cokernels are complex vector spaces. It follows that the determinant bundle of the whole space of real Cauchy--Riemann operators over a closed Riemann surface is canonically oriented. This shows that our bundle over $S^1$ is orientable. \qed

\end{prf}

We therefore let
$$C(\wt y)$$
be the free $\Z$-module of rank $1$ whose two generators are the two possible orientations of this determinant line bundle. Again, since $\wt y$ is connected, this definition makes sense.

\subsubsection{Orientations and isomorphisms}\label{sss:orientations_isomorphisms}

Let $\Sigma$ be a punctured sphere or closed disk with punctures $\theta,\{\theta_i\}_i$ with $\theta$ being the only positive puncture. Assume we have chosen nondegenerate Floer data $(H,J)$, $(H^i,J^i)$ associated to $\theta,\theta_i$. Choose critical points $\wt y_i \in \Crit \cA_{H_i}$ and let $y$ be a Hamiltonian orbit of $H$. Let
$$w \in C^\infty_b(\Sigma,\partial\Sigma; M,L;\{y_i\}_i,y)\,.$$
Fix representative cappings $\wh y_i \in \wt y_i$. We can preglue the maps $\wh y_i$ and $w$, according to the obvious gluing tree, to form a new map, which is b-smooth and has $y$ as the unique asymptotic orbit. It therefore can be viewed as a capping for $y$, and so we denote it $\wh y$, and let $\wt y$ be its equivalence class. Let $D_{\wh y_i} \in D_{\wt y_i}$ be linearized operators for the cappings. We can glue these operators with a linearized operator $D_w$ cooresponding to $w$ so that the result can be deformed into a linearized operator $D_{\wh y}$ for the capping $\wh y$. We therefore have an isomorphism
$$\ddd(D_{\wh y}) \simeq \ddd\Big(D_w \oplus \bigoplus_i D_{\wh y_i}\Big)$$
which is the composition of the deformation and gluing isomorphisms.

The direct sum isomorphism (\S\ref{par:direct_sum_isos}) yields
$$\ddd\Big(D_w \oplus \bigoplus_i D_{\wh y_i}\Big) \simeq \ddd(D_w) \otimes \bigotimes_i\ddd(D_{\wh y_i})\,,$$
and it depends on the ordering of the punctures of $\Sigma$. Composing the two isomorphisms, we obtain the isomorphism
$$\ddd(D_{\wh y}) \simeq \ddd(D_w) \otimes \bigotimes_i\ddd(D_{\wh y_i})\,.$$
Passing to the families, we finally get
$$\ddd(D_{\wt y}) \simeq \ddd(D_w) \otimes \bigotimes_i\ddd(D_{\wt y_i})\,.$$
This means that there is a canonical bijection between orientations of $D_w$ and isomorphisms
$$\ddd(D_{\wt y}) \simeq \bigotimes_i\ddd(D_{\wt y_i})\,,$$
or equivalently isomorphisms
$$C(\wt y) \simeq \bigotimes_i C(\wt y_i)\,.$$
We emphasize that this bijection depends on the chosen ordering of the orbits $y_i$. Note that this bijection is continuous with respect to $w$.

Now we'll show how such isomorphisms correspond to orientations of solution spaces. There are two cases: the case of a single surface and the case of a family.

\paragraph{A single surface} Assume the Floer data associated to the punctures of $\Sigma$ are regular and choose a compatible regular perturbation datum $(K,I)$ on $\Sigma$, so that for every $u\in\cM_\Sigma(K,I;\{\wt y_i\}_i,\wt y)$ the linearized operator $D_u$ is onto and we have canonically
$$\ker D_u = T_u\cM_\Sigma(K,I;\{\wt y_i\}_i,\wt y)\,.$$
As we have just seen, isomorphisms $\bigotimes_i C(\wt y_i) \simeq C(\wt y)$ are in bijection with orientations of $D_u$ for any $u \in \cM_\Sigma(K,I;\{\wt y_i\}_i,\wt y)$, therefore we obtain: for every connected component of $\cM_\Sigma(K,I;\{\wt y_i\}_i,\wt y)$, there is a bijection between such isomorphisms and orientations of that connected component.

\paragraph{A family} Let $\cS \to \cR$ be a family of punctured Riemann surfaces with a single positive puncture $\theta$ and negative punctures $\{\theta_i\}_i$, and assume we have chosen a set of ends for it, a set of regular Floer data $(H,J)$, $(H^i,J^i)$ associated to $\theta,\theta^i$, and a compatible regular perturbation datum $(K,I)$. Fix $\wt y_i \in \Crit \cA_{H^i}$, $\wt y\in\Crit \cA_H$. For every $(r,u) \in \cM_\cS(K,I;\{\wt y_i\}_i,\wt y)$ we have canonically
$$\ker D_{r,u} = T_{(r,u)}\cM_\cS(K,I;\{\wt y_i\}_i,\wt y)\,.$$
The exact Fredholm triple
\begin{equation}\label{eqn:exact_triple_extended_linearized_op}
0 \to D_u \to D_{r,u} \to 0_{T_r\cR} \to 0
\end{equation}
leads to the canonical isomorphism
\begin{equation}\label{eqn:iso_exact_triple_extended_linearized_op}
\ddd(D_{r,u}) \simeq \ddd(D_u) \otimes \ddd(T_r\cR)\,.\end{equation}
As we have just seen, there is a bijection between isomorphisms $\bigotimes_i C(\wt y_i) \simeq C(\wt y)$ and orientations of $D_u$. Therefore there is a bijection between such isomorphisms and orientations of $\cM_\cS(K,I;\{\wt y_i\}_i,\wt y)$ \emph{relative to} $\cR$. In our applications below $\cR$ is always an interval in $\R$ and so it carries the positive orientation inherited from $\R$, which therefore implies that we have a canonical bijection between isomorphisms $\bigotimes_i C(\wt y_i) \simeq C(\wt y)$ and orientations of components of $\cM_\cS(K,I;\{\wt y_i\}_i,\wt y)$.

\subsubsection{Canonical orientations}\label{sss:canonical_ors}

Having defined the relevant moduli spaces of solutions of the Floer PDE, we now pass to the canonical orientations of the corresponding linearized operators. We distinguish the cases of a translation-invariant perturbation datum, a single surface, and a family.

\paragraph{Translation-invariant perturbation datum} Let $(H,J)$ be a regular Floer datum, which is time-periodic in case we consider periodic orbits, and fix $\wt y_\pm \in \Crit \cA_H$ of index difference $1$. For any $u \in \wt\cM(H,J;\wt y_-,\wt y_+)$ the linearized operator $D_u$ is onto and has index $1$, therefore its kernel is $1$-dimensional, and it is spanned by the infinitesimal translation $\partial_s$. We call the orientation
$$\partial_u := \partial_s \otimes 1^\vee$$
of $D_u$ \tb{canonical}.

\paragraph{A single surface} Let $\Sigma$ be a Riemann surface with $\theta$ being the only positive puncture and $\{\theta_i\}_i$ being the negative punctures, and assume we have chosen a set of ends for it, regular Floer data $(H,J)$, $(H^i,J^i)$ associated to $\theta,\theta_i$ and a compatible regular perturbation datum $(K,I)$; fix critical points $\wt y_i \in \Crit \cA_{H^i},\wt y \in \Crit \cA_H$. Assume the moduli space $\cM_\Sigma(K,I;\{\wt y_i\}_i,\wt y)$ is zero-dimensional and let $u$ be an element therein. The linearized operator $D_u$ is surjective and has index zero. Therefore it is an isomorphism and we let
$$\mfo_u = 1 \otimes 1^\vee \in \ddd(D_u)$$
be the positive orientation. We call this orientation $\mfo_u$ of $D_u$ {canonical}.

\paragraph{A family} \label{par:canonical_orientations_families} Assume we have a family of Riemann surfaces $\cS \to \cR$ where $\cR \subset \R$ is an interval. Let $\theta,\theta_i$ be the punctures of $\cS$ with $\theta$ being the only positive puncture, and assume we have chosen a set of ends for $\cS$, regular Floer data $(H,J),(H^i,J^i)$ associated to $\theta,\theta_i$, a regular compatible perturbation datum $(K,I)$, and critical points $\wt y_i \in \Crit \cA_{H^i},\wt y \in \Crit \cA_H$, such that $\cM_\cS(K,I;\{\wt y_i\}_i,\wt y)$ is zero-dimensional. Assume $(r,u) \in \cM_\cS(K,I;\{\wt y_i\}_i,\wt y)$. The operator $D_{r,u}$ is onto and has index zero, therefore it is an isomorphism. Recall the isomorphism \eqref{eqn:iso_exact_triple_extended_linearized_op}
$$\ddd(D_{r,u}) \simeq \ddd(D_u) \otimes \ddd(T_r\cR)\,.$$
Let $\partial_r$ be the positive orientation of $\cR$. We let $\mfo_u \in \ddd(D_u)$ be the orientation such that this isomorphism maps
$$1 \otimes 1^\vee \mapsto \mfo_u \otimes \partial_r\,.$$
We call this orientation $\mfo_u$ \tb{canonical}. Note that $D_u$ has index $-1$.

\subsubsection{Induced orientations}\label{sss:induced_ors}

Whenever we have a $1$-dimensional moduli space $\cM$, it can be compactified to a $1$-dimensional compact manifold $\ol\cM$ with boundary consisting of elements of $0$-dimensional moduli spaces. Here we show how the canonical orientations of the Fredholm operators corresponding to these $0$-dimensional spaces induce orientations on $\cM$. These computations will be used in \S\ref{sss:identities} to prove that the various operations in Floer homology satisfy suitable identities.

The notations and the treatment here parallel those of \S\ref{sss:compactness_gluing}, where the types of boundary points arising in compactification are described. Note that all Floer and perturbation data in sight are assumed to be regular and sufficiently generic so that compactness results of \S\ref{sss:compactness_gluing} apply.

\paragraph{The case of translation-invariant perturbation datum}\label{par:induced_orientation_translation_invt_pert_datum} Consider the $1$-dimensional moduli space $\cM(H,J;\wt y_-,\wt y_+)$ and a boundary point
$$\delta=([u],[v]) \in \cM(H,J;\wt y_-,\wt y) \times \cM(H,J;\wt y, \wt y_+)\,.$$
Let $\Delta$ be the connected component of $\ol\cM(H,J;\wt y_-,\wt y_+)$ such that $\delta \in \Delta$. Let $w \in \wt\cM(H,J;\wt y_-,\wt y_+)$ be obtained by gluing $u,v$ for some large gluing length. The canonical orientations of $D_u,D_v$ correspond, by \S\ref{sss:orientations_isomorphisms}, to isomorphisms
$$C(u)\fc C(\wt y_-) \simeq C(\wt y) \quad \text{and} \quad C(v) \fc C(\wt y) \simeq C(\wt y_+)\,.$$ The isomorphism
$$C(v) \circ C(u) \fc C(\wt y_-) \simeq C(\wt y_+)$$
corresponds to an orientation of $D_w$, which we will now compute. First, consider the isomorphism
\begin{equation}\label{eqn:iso_dDu_otimes_dDv_dDw_induced_or_translation_invt_pert_datum}
\ddd(D_v) \otimes \ddd(D_u) \simeq \ddd(D_w)
\end{equation}
which is the composition of the direct sum, gluing, and deformation isomorphisms. Let $\mfo_w \in \ddd(D_w)$ be the image of $\partial_v \otimes \partial_u$ under this isomorphism. We claim that $\mfo_w$ corresponds to the isomorphism $C(v) \circ C(u)$. Indeed, we have the following commutative diagram:
\begin{equation}\label{dia:computation_induced_ori_bdry_op_squared_HF}
\xymatrix{\ddd(D_{\wt y_-}) \ar[r] \ar[dr] & \ddd(D_v) \otimes \ddd(D_u) \otimes \ddd(D_{\wt y_-}) \ar[r] \ar[d] & \ddd(D_v) \otimes \ddd(D_{\wt y}) \ar[d]\\ & \ddd(D_w) \otimes \ddd(D_{\wt y_-}) \ar[r] & \ddd(D_{\wt y_+}) }
\end{equation}
where the leftmost horizontal arrow maps $\mfo_- \mapsto \partial_v \otimes \partial_u \otimes \mfo_-$ while the slanted arrow maps $\mfo_- \mapsto \mfo_w \otimes \mfo_-$. The triangle then commutes by the definition of $\mfo_w$. The square consists of combinations of direct sum, gluing, and deformation isomorphisms, and therefore commutes. Fix $\mfo_- \in \ddd(D_{\wt y_-})$ and let $\mfo = C(u)(\mfo_-)$ and $\mfo_+ = C(v)(\mfo)$. Then by the definition of $C(u), C(v)$ we know that this diagram maps
$$\xymatrix{\mfo_- \ar@{|->}[r] \ar@{|->}[dr] & \partial_v \otimes \partial_u \otimes \mfo_- \ar@{|->}[d] \ar@{|->}[r] & \partial_v \otimes \mfo \ar@{|->}[d] \\ & \mfo_w \otimes \mfo_- \ar@{|->} [r] & \mfo_+}$$
We see that the composition of the slanted arrow and the bottom arrow maps $\mfo_- \mapsto \mfo_+$. On the other hand, by the definition of the correspondence between isomorphisms $C(\wt y_-) \simeq C(\wt y_+)$ and orientations of $\ddd(D_w)$ we have that the isomorphism $C(v) \circ C(u)$, which maps $\mfo_- \mapsto \mfo_+$, corresponds to the orientation $\mfo_w$. It remains to explicitly compute $\mfo_w$.

Recall that the gluing map is a local diffeomorphism (see for instance \cite{Schwarz_PhD_thesis})
$$\wt\cM(H,J;\wt y_-,\wt y) \times \wt\cM(H,J;\wt y,\wt y_+) \to \wt\cM(H,J;\wt y_-,\wt y_+)$$
defined near $(u,v)$ and mapping a neighborhood of this point diffeomorphically to a neighborhood of $w$. The main property of this gluing map is that its differential, which is an isomorphism
$$dg \fc T_u\wt\cM(H,J;\wt y_-,\wt y) \oplus T_v\wt\cM(H,J;\wt y,\wt y_+) \simeq T_w \wt \cM(H,J;\wt y_-,\wt y_+)$$
is such that if we let $\ddd(dg) \fc \ddd(D_u \oplus D_v) \to \ddd(D_w)$ be the induced map on determinant lines, then the composition with the direct sum isomorphism
$$\ddd(D_v) \otimes \ddd(D_u) \xrightarrow{\oplus}\ddd(D_u \oplus D_v) \xrightarrow{\ddd(dg)} \ddd(D_w)$$
yields the isomorphism \eqref{eqn:iso_dDu_otimes_dDv_dDw_induced_or_translation_invt_pert_datum}, which also enters in the left vertical arrow in the diagram \eqref{dia:computation_induced_ori_bdry_op_squared_HF}. Here we have of course identified
$$\ddd(D_u) = \ddd(\ker D_u) = \ddd(T_u\wt\cM(H,J;\wt y_-,\wt y))$$
and similarly for $v,w$. Therefore we know that $\mfo_w$ is the image of $\partial_v \otimes \partial_u$ by the composition
$$\ddd(D_v) \otimes \ddd(D_u) \xrightarrow{\oplus} \ddd(D_v \oplus D_u) \xrightarrow{\ddd(dg)} \ddd(D_w)\,.$$
The first map maps $\partial_v \otimes \partial_u \mapsto \partial_v \wedge \partial_u$. Next, from the structure of the gluing map $g$ it follows that
$$dg(\partial_v + \partial_u) = \partial_w$$
while \label{footnote:differential_gluing_map} \footnote{This can be seen intuitively as follows. Fix a large gluing length $R$. Then the glued trajectory $w$ satisfies that $w(0)$ is close to $u(0)$ while $w(2R+1)$ is close to $v(0)$. If we now let $u',v'$ be trajectories which satisfy, say $u'(0) = u(-1)$, $v'(0) = v(1)$, so that the passage from $(u,v)$ to $(u',v')$ is in the direction of the vector $-\partial_u+\partial_v$, and we let $w'$ be the trajectory glued from $u',v'$ for the same gluing length, then $w'(0)$ is close to $u'(0)=u(-1)$, and $w'(2R+1)$ is close to $v'(0) = v(1)$, that is the points $w'(0),w'(2R+1)$ are futher apart than the points $w(0),w(2R+1)$, which means that $w'$ is faster than $w$, because it traverses a longer distance in the same amount of time. Therefore it spends less time near $\wt y$, meaning it is further away from $\delta$. Therefore $dg$ maps $-\partial_u + \partial_v$ to $\text{inward}_\delta$.}
$$dg(\partial_v - \partial_u) = \text{inward}_\delta\,,$$
where $\inward_\delta \in T_w\wt\cM(H,J;\wt y_-,\wt y_+)$ is a vector pointing away from $\delta$. Therefore
$$\mfo_w = \ddd(dg)(\partial_v \wedge \partial_u) = \ddd(dg)((\partial_v + \partial_u) \wedge (-\partial_v + \partial_u)) = \partial_w\wedge (-\text{inward}_\delta)\,.$$
This means that the isomorphism $C(v) \circ C(u)$ induces the orientation $-\partial_w\wedge \text{inward}_\delta$ on the connected component of $\wt\cM(H,J;\wt y_-,\wt y_+)$ containing $w$.

\paragraph{The case of a single surface}\label{par:induced_orientation_single_surf} Consider a $1$-dimensional moduli space $\cM_\Sigma (K,I;\{\wt y_i\}_i,\wt y)$. The boundary points correspond to Floer breaking, either at an incoming or at the outgoing end. Consider breaking at the outgoing end: 
$$\delta = (u,[v]) \in  \cM_\Sigma(K,I;\{\wt y_i\}_i,\wt y') \times \cM(H,J;\wt y',\wt y)\,.$$
Let $\Delta \subset \ol \cM_\Sigma (K,I;\{\wt y_i\}_i,\wt y)$ be the connected component with $\delta \in \Delta$. Let $w \in \Delta$ be obtained by gluing $u,v$ for some large gluing length. The canonical orientations of $D_u,D_v$ correspond to isomorphisms
$$C(u) \fc \bigotimes_i C(\wt y_i) \simeq C(\wt y')\quad\text{and}\quad C(v) \fc C(\wt y') \simeq C(\wt y)\,,$$
and the composition
$$C(v) \circ C(u) \fc \bigotimes_i C(\wt y_i) \simeq C(\wt y)$$
corresponds to an orientation of $D_w$, and therefore of $\Delta$. Let us compute this orientation. The differential of the gluing map is an isomorphism
$$dg \fc T_u\cM_\Sigma(K,I;\{\wt y_i\}_i,\wt y') \times T_v \wt\cM(H,J;\wt y',\wt y) \simeq T_w\cM_\Sigma (K,I;\{\wt y_i\}_i,\wt y)$$
mapping the spanning vector \footnote{This can be seen using the argument of the footnote on page \pageref{footnote:differential_gluing_map}.} $\partial_v$ to $\text{inward}_\delta$. Thus the induced isomorphism on determinant lines
$$\ddd(D_v) \otimes \ddd(D_u) \xrightarrow{\oplus} \ddd(D_v \oplus D_u) \xrightarrow{\ddd(dg)} \ddd(D_w)$$
maps $\partial_v \otimes \mfo_u \mapsto \text{inward}_\delta$.

The operator $D_v \oplus D_u \oplus \bigoplus_i D_{\wh y_i}$ glues into: $D_v \oplus D_{\wh y'}$, $D_w \oplus \bigoplus_i D_{\wh y_i}$, and $D_{\wh y}$. Using a combination of direct sum, linear gluing, and deformation isomorphisms, we obtain the commutative diagram
$$\xymatrix{\ddd(D_v) \otimes \ddd(D_u) \otimes \bigotimes_i\ddd(D_{\wt y_i}) \ar[dr] \ar[dd]^{dg \otimes \id}& &  \\
& \ddd(D_v) \otimes \ddd(D_{\wt y'}) \ar[r] \ar[dl] & \ddd(D_{\wt y}) \\
\ddd(D_w) \otimes \bigotimes_i\ddd(D_{\wt y_i}) \ar[rru]}$$
Pick $\mfo_i \in C(\wt y_i)$ and denote $\mfo' = C(u)\big(\bigotimes_i \mfo_i\big) \in C(\wt y')$ and $\mfo = C(v)(\mfo') \in C(\wt y)$. This diagram maps
$$\xymatrix{ \partial_v  \otimes \mfo_u \otimes \bigotimes_i \mfo_i \ar@{|->}[dr] \ar@{|->}[dd]^{dg \otimes \id}& &  \\
& \partial_v \otimes \mfo' \ar@{|->}[r] \ar@{|->}[dl] & \mfo \\
\text{inward}_\delta \otimes \bigotimes_i \mfo_i \ar@{|->}[urr]}$$
This means that the orientation of $D_w$ corresponding to the isomorphism $C(v)\circ C(u)$, which maps $\bigotimes_i\mfo_i \mapsto \mfo$, is given by $\text{inward}_\delta$. This is also the induced orientation on $\Delta$.

When the breaking happens at the $j$-th incoming end, we get a boundary point
$$\delta = ([u],v) \in \cM(H^j,J^j;\wt y_j,\wt y_j') \times \cM_\Sigma(K,I;\{\wt y_i\}_{i\neq j},\wt y_j',\wt y)\,.$$
Let again $\Delta$ be the connected component of $\ol \cM_\Sigma(K,I;\{\wt y_i\}_i,\wt y)$ with $\delta \in \Delta$. Let $w \in \cM_\Sigma(K,I;\{\wt y_i\}_i,\wt y)$ be obtained by gluing $u,v$ for some large gluing length. The canonical orientations of the operators $D_u,D_v$ correspond to isomorphisms
$$C(u) \fc C(\wt y_j) \simeq C(\wt y_j') \quad \text{and}\quad C(v) \fc \bigotimes_{i<j}C(\wt y_i) \otimes C(\wt y_j') \otimes \bigotimes_{i>j} C(\wt y_i) \simeq C(\wt y)\,.$$
The isomorphism
$$C(v) \circ (\id \otimes \dots \otimes C(u) \otimes \dots \otimes \id) \fc \bigotimes_i C(\wt y_i) \simeq C(\wt y)$$
corresponds to an orientation of $D_w$, and therefore of $\Delta$. Let us compute it. The differential of the gluing map is an isomorphism
$$dg \fc T_u \wt \cM(H^j,J^j;\wt y_j,\wt y_j') \times T_v \cM_\Sigma(K,I;\{\wt y_i\}_{i\neq j},\wt y_j',\wt y) \to T_w \cM_\Sigma(K,I;\{\wt y_i\}_i,\wt y)$$
mapping the spanning vector $\partial_u$ to $-\text{inward}_\delta$.

The operator $D_v \oplus D_u \oplus \bigoplus_i D_{\wh y_i}$ glues into $D_v \oplus \bigoplus_{i\neq j}D_{\wh y_i} \oplus D_{\wh y_j'}$, $D_w \oplus \bigoplus_i D_{\wh y_i}$, and $D_{\wh y}$. Noting that direct sum isomorphisms obey the Koszul rule with respect to the grading of the determinant lines, we have the following commutative diagram
$$\xymatrix{ & \ddd(D_v) \otimes \bigotimes_{i < j}\ddd(D_{\wt y_i}) \otimes \ddd (D_u) \otimes \ddd(D_{\wt y_j}) \otimes \bigotimes_{i>j}\ddd(D_{\wt y_i}) \ar [d]^R \ar[dl]\\
\ddd(D_{\wt y})    & \ddd(D_v) \otimes \ddd(D_u) \otimes \bigotimes_i \ddd(D_{\wt y_i}) \ar[d]^{\ddd(dg) \otimes \id} \ar[l]\\ 
& \ddd(D_w) \otimes \bigotimes_i \ddd(D_{\wt y_i}) \ar[lu]}$$
where $R$ differs from the mere exchange of factors by the Koszul sign $(-1)^{\sum_{i < j}|\wt y_i|'}$. Fix generators $\mfo_i \in C(\wt y_i)$ and let $\mfo = C(v)\circ (\id \otimes \dots \otimes C(u) \otimes \dots \otimes \id)(\bigotimes_i\mfo_i)\in C(\wt y)$. We have
\begin{equation}\label{dia:induced_orientation_single_surf_Floer_breaking_incoming_end}
\xymatrix{ & \mfo_v \otimes \bigotimes_{i < j} \mfo_i  \otimes \partial_u \otimes \mfo_j \otimes \bigotimes_{i>j} \mfo_i \ar@{|->} [d] \ar @{|->}[dl]\\
\mfo    &(-1)^{\sum_{i<j}|\wt y_i|'} \mfo_v \otimes \partial_u \otimes \bigotimes_i \mfo_i \ar @{|->}[d] \ar@{|->} [l]\\ 
& -(-1)^{\sum_{i<j}|\wt y_i|'} \text{inward}_\delta \otimes \bigotimes_i \mfo_i \ar @{|->}[lu]}
\end{equation}
This means that the isomorphism $C(v) \circ (\id \otimes \dots \otimes C(u) \otimes \dots \otimes \id)$ which maps $\bigotimes_i\mfo_i \mapsto \mfo$, corresponds to the orientation $-(-1)^{\sum_{i<j}|\wt y_i|'}\text{inward}_\delta$ of $D_w$. This is therefore the induced orientation on $\Delta$.

\paragraph{The case of a family}\label{par:induced_orientations_families} Let $\cS \to \cR$ be a family of punctured Riemann surfaces, and fix regular Floer data associated to its punctures and a regular compatible perturbation datum. We only deal with $\cR$ being a connected interval in $\R$ of the form $\cR = [0,1]$ or $\cR = [0,\infty)$. We orient $\cR$ by the orientation $\partial_r = 1 \in \ddd(T_r\cR) = \ddd(\R)$. We consider a $1$-dimensional moduli space $\cM_\cS(K,I;\{\wt y_i\}_i,\wt y)$ and the boundary of its compactification. There are three types of boundary points: the boundary points of the original moduli space before compactification, internal Floer breaking, and breaking at a noncompact end of $\cR$, see \S\ref{sss:compactness_gluing}.

We start with boundary points belonging to the moduli space itself: let
$$\delta = (r,u) \in \partial \cR \times \cM_{\Sigma_r}(K_r,I_r;\{\wt y_i\}_i,\wt y)\,.$$
The operator $D_u$ is canonically oriented by $1\otimes 1^\vee$, since it is an isomorphism. We have the isomorphism \eqref{eqn:iso_exact_triple_extended_linearized_op}
$$\ddd(D_{r,u}) \simeq \ddd(D_u) \otimes \ddd(T_r\cR)\,.$$
Since the operators $D_u,D_{r,u}$ are onto, by the normalization property (\S\ref{par:normalization_pty}) this isomorphism in fact comes from the short exact sequence
$$0 \to \ker D_u \to \ker D_{r,u} \xrightarrow{\pr} T_r\cR \to 0\,,$$
where $\pr \fc \ker D_{r,u} \to T_r\cR$ is the restriction of the projection $T_r\cR \oplus W^{1,p}(u) \to T_r\cR$. Since $\ker D_u = 0$, this $\pr$ is an isomorphism; by abuse of notation we denote $\partial_r$ the preimage of $\partial_r$ by $\pr$. It follows that
$$\partial_r \mapsto (1 \otimes 1^\vee) \otimes \partial_r\,.$$
Therefore the induced orientation on $D_{r,u}$ is given by $\partial_r$. This is therefore also the induced orientation on the connected component of $\cM_\cS(K,I;\{\wt y_i\}_i,\wt y)$ containing $(r,u)$. Note that if $r = 0$, this orientation coincides with $\inward_\delta$, while if $r = 1$, it coincides with $-\inward_\delta$.

Next we consider internal Floer breaking. This can happen at the outgoing end or at an incoming end. Consider first breaking at the outgoing end: let
$$\delta = ((r,u),[v])\in\cM_\cS(K,I;\{\wt y_i\}_i,\wt y') \times \cM(H,J;\wt y',\wt y)\,.$$
Let $\Delta$ be the connected component of $\ol\cM_\cS(K,I;\{\wt y_i\}_i,\wt y)$ with $\delta \in \Delta$. Let $(s,w) \in \Delta$ be obtained by gluing $(r,u)$ and $v$ for some large gluing length. The canonical orientations of $D_u,D_v$ correspond to isomorphisms
$$C(u) \fc C(\wt y') \simeq C(\wt y)\quad\text{and}\quad C(v)\fc\bigotimes_i C(\wt y_i) \simeq C(\wt y')\,.$$
The isomorphism
$$C(v) \circ C(u) \fc \bigotimes_i C(\wt y_i) \simeq C(\wt y)$$
corresponds to an orientation of $D_w$, and to an orientation of $\cM_\cS(K,I;\{\wt y_i\}_i,\wt y)$ via the isomorphism \eqref{eqn:iso_exact_triple_extended_linearized_op}. Let us compute these orientations. The differential of the gluing map is an isomorphism
$$dg \fc T_{(r,u)}\cM_\cS(K,I;\{\wt y_i\}_i,\wt y') \oplus T_v \wt \cM(H,J;\wt y',\wt y) \to T_{(s,w)}\cM_\cS(K,I;\{\wt y_i\}_i,\wt y)$$
mapping the spanning vector $\partial_v$ to $\text{inward}_\delta$. We have the exact square of Fredholm operators
$$\xymatrix{D_v \ar@{=}[r] \ar[d] & D_v \ar[r] \ar[d] & 0 \ar[d] \\ D_v \oplus D_u \ar[r] \ar[d] & D_v \oplus D_{r,u} \ar[r] \ar[d] & 0_{T_r\cR} \ar@{=}[d] \\ D_u \ar[r] & D_{r,u} \ar[r]  & 0_{T_r\cR}}$$
which induces the following commutative diagram:
$$\xymatrix{\ddd(D_v) \otimes \ddd(D_u) \otimes \ddd(T_r\cR) \ar[r] \ar[d]& \ddd(D_v) \otimes \ddd(D_{r,u}) \ar[d]\\
\ddd(D_v \oplus D_u) \otimes \ddd(T_r\cR) \ar[r]& \ddd(D_v\oplus D_{r,u})}$$
where the vertical arrows are direct sum isomorphisms while the horizontal arrows come from exact triples.

There is an isomorphism
\begin{equation}\label{eqn:iso_dDu_otimes_dDv_dDw_internal_Floer_breaking_families_outgoing_end}
\ddd(D_v \oplus D_u) \simeq \ddd(D_w)
\end{equation}
which is the combination of linear gluing and deformation isomorphisms. Since $s$ is close to $r$, we have the isomorphism
$$\ddd(D_v \oplus D_u) \otimes \ddd(T_r\cR) \simeq \ddd(D_w) \otimes \ddd(T_s\cR)\,.$$
The main property of the gluing map is the commutativity of the following diagram:
$$\xymatrix{\ddd(D_v \oplus D_u) \otimes \ddd(T_r\cR) \ar[r] \ar[d]& \ddd(D_v\oplus D_{r,u})\ar[d]^{\ddd(dg)}\\
\ddd(D_w) \otimes \ddd(T_s\cR) \ar[r]& \ddd(D_{s,w})}$$
where the horizontal arrows come from exact triples. Combining the two diagrams, we obtain the diagram on the left:

\begin{minipage}{8cm}
$$\xymatrix{\ddd(D_v) \otimes \ddd(D_u) \otimes \ddd(T_r\cR) \ar[r] \ar[d]& \ddd(D_v) \otimes \ddd(D_{r,u}) \ar[d]\\
\ddd(D_v \oplus D_u) \otimes \ddd(T_r\cR) \ar[r] \ar[d]& \ddd(D_v\oplus D_{r,u})\ar[d]^{\ddd(dg)}\\
\ddd(D_w) \otimes \ddd(T_s\cR) \ar[r]& \ddd(D_{s,w})}$$
\end{minipage}
\begin{minipage}{6cm}
$$\xymatrix{\partial_v \otimes \mfo_u  \otimes \partial_r \ar@{|->}[r] \ar@{|->}[d]& \partial_v \otimes (1\otimes 1^\vee)\ar@{|->}[d]\\
(\partial_v \wedge \mfo_u) \otimes \partial_r \ar@{|->}[r] \ar@{|->}[d]& \partial_v\ar@{|->}[d]\\
\mfo_w \otimes \partial_r \ar@{|->}[r]& \text{inward}_\delta}$$
\end{minipage}

\noindent Recall that we have the canonical orientations $\mfo_u \in \ddd(D_u)$ and $\partial_v \in \ddd(D_v)$. Let us denote by $\mfo_w$ the image of $\partial_v \otimes \mfo_u$ by the isomorphism \eqref{eqn:iso_dDu_otimes_dDv_dDw_internal_Floer_breaking_families_outgoing_end}. The diagram on the left maps orientations as shown in the diagram on the right, where the top horizontal arrow comes from the definition of $\mfo_u$, see \S \ref{par:canonical_orientations_families}, and the vertical arrows all come from the definitions. The goal of this computation is the bottom arrow, which as we can see maps $\mfo_w \otimes \partial_r \mapsto \text{inward}_\delta$.

Now the operator $D_v \oplus D_u \oplus \bigoplus_i D_{\wh y_i}$ glues into $D_w \oplus \bigoplus_i D_{\wh y_i}$, $D_v \oplus D_{\wh y'}$, and $D_{\wh y}$. Therefore, using a combination of direct sum, gluing, and deformation isomorphisms, we have the commutative diagram
$$\xymatrix{\ddd(D_v) \otimes \ddd(D_u) \otimes \bigotimes_i\ddd(D_{\wt y_i}) \ar[d]^{dg \otimes \id} \ar[r] & \ddd(D_{\wt y}) \\
\ddd(D_w) \otimes \bigotimes_i\ddd(D_{\wt y_i})  \ar[ru] }$$
Pick generators $\mfo_i \in C(\wt y_i)$ and let $\mfo = C(v)(C(u)(\bigotimes_i\mfo_i)) \in C(\wt y)$. This diagram maps
$$\xymatrix{\partial_v \otimes \mfo_u \otimes \bigotimes_i \mfo_i \ar@{|->}[d]^{dg \otimes \id} \ar@{|->}[r] & \mfo \\
\mfo_w \otimes \bigotimes_i \mfo_i \ar@{|->}[ru]}$$
Thus we see that the isomorphism $C(v) \circ C(u)$, which maps $\bigotimes_i \mfo_i \mapsto \mfo$, corresponds to the orientation $\mfo_w$ of $D_w$, which in turn corresponds to the orientation $\text{inward}_\delta$ induced on $\Delta$.

Assume now that the Floer breaking occurs at the $j$-th incoming end, and we have
$$\delta=([u],(r,v)) \in \cM(H^j,J^j;\wt y_j,\wt y_j') \times \cM_\cS(K,I;\{\wt y_i\}_{i\neq j},\wt y_j',\wt y)\,.$$
Let $\Delta$ be the connected component of $\ol\cM_\cS(K,I;\{\wt y_i\}_i,\wt y)$ with $\delta \in \Delta$. Let $(s,w)$ be obtained by gluing $u,(r,v)$ for some large gluing length. The canonical orientations $\partial_u \in \ddd(D_u)$ and $\mfo_v \in \ddd(D_v)$ correspond to isomorphisms
$$C(u) \fc C(\wt y_j) \simeq C(\wt y_j')\quad \text{and} \quad C(v) \fc \bigotimes_{i<j} C(\wt y_i) \otimes C(\wt y_j') \otimes \bigotimes_{i > j}C(\wt y_i) \simeq C(\wt y)\,.$$
The composition
$$C(v)\circ(\id \otimes \dots\otimes C(u) \otimes \dots \otimes \id) \fc \bigotimes_i C(\wt y_i) \simeq C(\wt y)$$
corresponds to an orientation of $D_w$, and via \eqref{eqn:iso_exact_triple_extended_linearized_op}, to an orientation of $\cM_\cS(K,I;\{\wt y_i\}_i,\wt y)$. Let us compute these orientations. The differential of the gluing map is an isomorphism
$$dg \fc T_u \wt \cM(H^j,J^j;\wt y_j,\wt y_j') \oplus T_{(r,v)}\cM_\cS(K,I;\{\wt y_i\}_{i\neq j},\wt y_j',\wt y) \to T_{(s,w)}\cM_\cS(K,I;\{\wt y_i\}_i,\wt y)$$
mapping the spanning vector $\partial_u$ to $-\text{inward}_\delta$. Similarly to the breaking at the outgoing end described above, we use the exact square of Fredholm operators
$$\xymatrix{D_v \ar[r] \ar[d] & D_v \oplus D_u \ar[r] \ar [d]  & D_u \ar@{=}[d] \\ D_{r,v} \ar[r] \ar[d]& D_{r,v} \oplus D_u \ar[r] \ar[d] & D_u \ar[d]\\ 0_{T_r\cR} \ar@{=}[r] & 0_{T_r\cR} \ar[r]& 0}$$
and the isomorphism
\begin{equation}\label{eqn:iso_dDu_otimes_dDv_dDw_internal_Floer_breaking_families_incoming_end}
\ddd(D_v \oplus D_u) \simeq \ddd(D_w)
\end{equation}
obtained as a combination of linear gluing and deformation isomorphisms, to obtain the commutative diagram on the left:

\begin{minipage}{8cm}
$$\xymatrix{\ddd(D_v) \otimes \ddd(D_u) \otimes \ddd(T_r\cR) \ar[r] \ar[d]& \ddd(D_{r,v}) \otimes \ddd(D_u) \ar[d]\\
\ddd(D_v \oplus D_u) \otimes \ddd(T_r\cR) \ar[r] \ar[d]& \ddd(D_{r,v}\oplus D_u)\ar[d]\\
\ddd(D_w) \otimes \ddd(T_s\cR) \ar[r]& \ddd(D_{s,w})}$$
\end{minipage}
\begin{minipage}{6cm}
$$\xymatrix{\mfo_v \otimes \partial_u \otimes \partial_r \ar@{|->}[r] \ar@{|->}[d]& -(1\otimes 1^\vee) \otimes \partial_u \ar@{|->}[d]\\
(\mfo_v \wedge\partial_u) \otimes \partial_r \ar@{|->}[r] \ar@{|->}[d]& -\partial_u\ar@{|->}[d]\\
\mfo_w \otimes \partial_r \ar@{|->}[r]& \text{inward}_\delta}$$
\end{minipage}

\noindent where the top square is induced by the exact Fredholm square while the bottom square comes from the compatibility of the gluing map with linear gluing. The top arrow is the composition of the interchange of factors, which includes the Koszul sign $(-1)^{\ind D_u \cdot\dim \cR} = -1$, and the isomorphism coming from the exact triple. Recall that we have canonical orientations $\mfo_v,\partial_u$ and let $\mfo_w$ be the image of $\mfo_v \otimes \partial_u$ by the isomorphism \eqref{eqn:iso_dDu_otimes_dDv_dDw_internal_Floer_breaking_families_incoming_end}. Then the diagram on the left maps as shown on the right.

The rest of the computation is almost identical to the case of Floer breaking at an incoming end for the case of a single surface: the isomorphism $C(v) \circ (\id \otimes \dots \otimes C(u) \otimes \dots \otimes \id)$ corresponds to the orientation of $D_w$ given by $(-1)^{\sum_{i<j}|\wt y_i|'}\mfo_w$ (that is, the Koszul sign is the same as in \eqref{dia:induced_orientation_single_surf_Floer_breaking_incoming_end}). As we have just computed, the orientation $\mfo_w$ corresponds to the orientation $\text{inward}_\delta$ of $\Delta$, therefore the total induced orientation on $\Delta$ is $(-1)^{\sum_{i<j}|\wt y_i|'} \text{inward}_\delta$.

Lastly, when $\cR = [0,\infty)$, we also have breaking at the noncompact end of $\cR$. We retain the notation of \S\ref{par:compactness_families} treating compactness and gluing in this case. Consider the $1$-dimensional moduli space $\cM_\cS(K,I;\{\wt y_{1,i}\}_i,\{\wt y_{2,i}\}_{i \neq j},\wt y_2)$ and let
$$\delta = (u,v) \in \cM_{\Sigma^1}(K^1,I^1;\{\wt y_{1,i}\}_i, \wt y_{2,j}) \times \cM_{\Sigma^2}(K^2,I^2;\{\wt y_{2,i}\}_i,\wt y_2)\,.$$
Let $\Delta$ be the connected component of $\ol\cM_\cS(K,I;\{\wt y_{1,i}\}_i,\{\wt y_{2,i}\}_{i \neq j},\wt y_2)$ such that $\delta \in \Delta$. Let $(r,w) \in \Delta$ be obtained by gluing $u,v$ for some large gluing length. The canonical orientations of $D_u,D_v$ correspond to isomorphisms
$$C(u) \fc \bigotimes_i C(\wt y_{1,i}) \simeq C(\wt y_{2,j})\quad\text{and} \quad C(v) \fc \bigotimes_i C(\wt y_{2,i}) \simeq C(\wt y_2)\,.$$
The composition
$$C(v) \circ (\id \otimes \dots \otimes C(u) \otimes \dots \otimes \id) \fc \bigotimes_{i < j} C(\wt y_{2,i}) \otimes \bigotimes_i C(\wt y_{1,i}) \otimes \bigotimes_{i > j} C(\wt y_{2,i}) \simeq C(\wt y_2)$$
corresponds to an orientation of $D_w$, and therefore to an orientation of $\Delta$, which we will now compute.

The gluing map is the trivial (zero) isomorphism
$$dg \fc T_u\cM_{\Sigma^1}(K^1,I^1;\{\wt y_{1,i}\}_i, \wt y_{2,j}) \oplus T_v \cM_{\Sigma^2}(K^2,I^2;\{\wt y_{2,i}\}_i,\wt y_2) \to T_w \cM_{\Sigma_r}(K_r,I_r;\{\wt y_{1,i}\}_i,\{\wt y_{2,i}\}_{i \neq j},\wt y_2)\,,$$
therefore the induced isomorphism
$$\ddd(D_v) \otimes \ddd(D_u) \simeq \ddd(D_w)\,,$$
which is the composition of the direct sum, linear gluing, and deformation isomorphisms, maps $\mfo_v \otimes \mfo_u = (1 \otimes 1^\vee) \otimes (1 \otimes 1^\vee)$ to $\mfo_w:= 1\otimes 1^\vee$, that is the positive orientation of $D_w$.

On the other hand, the operator
$$D_v \oplus D_u \oplus \bigoplus_i D_{\wh y_{1,i}} \oplus \bigoplus_{i \neq j} D_{\wh y_{2,i}}$$
glues into $D_w \oplus \bigoplus_i D_{\wh y_{1,i}} \oplus \bigoplus_{i \neq j} D_{\wh y_{2,i}}$, $D_v \oplus \bigoplus_i D_{\wh y_{2,i}}$, and $D_{\wh y_2}$. Using a combination of direct sum, linear gluing, and deformation isomorphisms, we get the following commutative diagram:
$$\xymatrix{\ddd(D_v) \otimes \bigotimes_{i < j} \ddd(D_{\wt y_{2,i}}) \otimes \ddd(D_u) \otimes \bigotimes_i \ddd(D_{\wt y_{1,i}}) \otimes \bigotimes_{i > j} \ddd(D_{\wt y_{2,i}}) \ar [d]^R \ar [dr]\\
\ddd(D_v) \otimes \ddd(D_u) \otimes \bigotimes_{i < j} \ddd(D_{\wt y_{2,i}}) \otimes \bigotimes_i \ddd(D_{\wt y_{1,i}}) \otimes \bigotimes_{i > j} \ddd(D_{\wt y_{2,i}}) \ar[d]^{\ddd(dg) \otimes \id} \ar[r] & \ddd(D_{\wt y_2})\\
\ddd(D_w) \otimes \bigotimes_{i < j} \ddd(D_{\wt y_{2,i}}) \otimes \bigotimes_i \ddd(D_{\wt y_{1,i}}) \otimes \bigotimes_{i > j} \ddd(D_{\wt y_{2,i}}) \ar[ur]}$$
where $R$ is the interchange of factors, including the Koszul sign, which is trivial since $\ind D_u = 0$. Pick generators $\mfo_{1,i} \in C(\wt y_{1,i})$ for all $i$, $\mfo_{2,i} \in C(\wt y_{2,i})$ for $i \neq j$, and let
$$\textstyle \mfo_2 = C(v)\big(\bigotimes_{i < j} \mfo_{2,i} \otimes C(u)\big(\bigotimes_i \mfo_{1,i}\big) \otimes \bigotimes_{i > j} \mfo_{2,i}) \in C(\wt y_2 \big)\,.$$
The diagram maps
$$\xymatrix{\mfo_v  \otimes \bigotimes_{i < j} \mfo_{2,i} \otimes \mfo_u \otimes \bigotimes_i \mfo_{1,i} \otimes \bigotimes_{i > j} \mfo_{2,j} \ar@{|->} [d]^R \ar@{|->}  [dr]\\
\mfo_v  \otimes \mfo_u \otimes \bigotimes_{i < j} \mfo_{2,i} \otimes \bigotimes_i \mfo_{1,i} \otimes \bigotimes_{i > j} \mfo_{2,j} \ar @{|->}  [d]^{\ddd(dg) \otimes \id} \ar @{|->} [r] & \mfo_2 \\
\mfo_w \otimes \bigotimes_{i < j} \mfo_{2,i} \otimes \bigotimes_i \mfo_{1,i} \otimes \bigotimes_{i > j} \mfo_{2,j} \ar @{|->}[ru]}$$
Therefore the isomorphism $C(v) \circ (\id \otimes \dots \otimes C(u) \otimes \dots \otimes \id)$ which maps
$$\textstyle\bigotimes_{i < j} \mfo_{2,i} \otimes \bigotimes_i \mfo_{1,i} \otimes \bigotimes_{i > j} \mfo_{2,i} \mapsto \mfo_2\,,$$
corresponds to the orientation $\mfo_w = 1 \otimes 1^\vee$ of $D_w$. Using the argument on page \pageref{par:induced_orientations_families} we see that the induced orientation on $\Delta$ is $\partial_r$, which evidently equals $-\inward_\delta$.

\subsection{Operations}\label{ss:operations}

Here we define the matrix elements of operations on Floer homology, which itself will be defined in the next subsection. The pattern is the same for all operations: given a collection of critical points $\wt y_i\in\Crit \cA_{H^i}$, $\wt y \in \Crit \cA_H$ the matrix element of an operation is a finite sum of the form
$$\sum_{u \in \cM(\{\wt y_i\}_i,\wt y)} C(u) \fc \bigotimes_i C(\wt y_i) \to C(\wt y)$$
where $\cM(\{\wt y_i\}_i,\wt y)$ is a $0$-dimensional moduli space of solutions of the Floer PDE and $C(u)$ is the isomorphism corresponding to the canonical orientation of the linearized operator $D_u$. Identities are proved by considering the compactified $1$-dimensional moduli spaces, whose boundary points are in bijection with summands of a desired identity. The main technical points here are the compactness and gluing results of \S\ref{sss:compactness_gluing}, the canonical orientations defined in \S\ref{sss:canonical_ors}, and the computation of induced orientations in \S\ref{sss:induced_ors}.

\subsubsection{Definition of operations}\label{sss:definition_ops}

There are three types of operations: boundary operators, multiplicative operators, and homotopy operators.

\paragraph{Boundary operators}\label{par:matrix_elts_boundary_op}

First we deal with the boundary operator in Lagrangian Floer homology. Let $(H,J)$ be a regular Floer datum for the strip $S = \R \times [0,1]$. Fix $\wt y_\pm \in \Crit \cA_H$ with $|\wt y_-| - |\wt y_+| = 1$. For every $u \in \wt \cM(H,J;\wt y_-,\wt y_+)$ the canonical orientation $\partial_u \in \ddd(D_u)$ corresponds to an isomorphism
$$C(u) \fc C(\wt y_-) \simeq C(\wt y_+)\,.$$
Clearly this only depends on $[u] \in \cM(H,J;\wt y_-,\wt y_+)$. The matrix element of the boundary operator $\partial_{H,J}$ is the homomorphism
$$\sum_{[u] \in \cM(H,J;\wt y_-,\wt y_+)} C(u) \fc C(\wt y_-) \to C(\wt y_+)\,.$$
This sum is finite since so is $\cM(H,J;\wt y_-,\wt y_+)$.

The matrix element of the boundary operator in Hamiltonian Floer homology is defined in an entirely analogous fashion, but the Floer datum is required to be $1$-periodic in time.

\paragraph{Multiplicative operators}\label{par:matrix_elts_operations}

There is a multiplicative operator corresponding to every punctured Riemann surface $\Sigma$ with exactly one positive puncture $\theta$ and negative punctures $\theta_i$, where $\wh \Sigma$ is the sphere or the closed disk. Fix a choice of ends for $\Sigma$, regular Floer data $(H,J),(H^i,J^i)$ associated to $\theta,\theta_i$, and a regular compatible perturbation datum $(K,I)$. For $\wt y \in \Crit \cA_{H^i},\wt y \in \Crit \cA_H$ with $|\wt y|' = \sum_i |\wt y_i|'$ pick $u \in \cM_\Sigma(K,I;\{\wt y_i\}_i,\wt y)$. The canonical orientation $\mfo_u \in \ddd(D_u)$ corresponds to an isomorphism
$$C(u) \fc \bigotimes_i C(\wt y_i) \simeq C(\wt y)\,.$$
The matrix element of the operation $\Phi_{\Sigma;K,I}$ is the homomorphism
$$\sum_{u \in \cM_\Sigma(K,I;\{\wt y_i\}_i,\wt y)} C(u) \fc \bigotimes_i C(\wt y_i) \to C(\wt y)\,.$$

\paragraph{Homotopy operators}\label{par:matrix_elts_htpy_ops}

Let $\cR = [0,1]$ or $[0,\infty)$ and let $\cS \to \cR$ be a family of punctured Riemann surfaces with one positive puncture $\theta$ and negative punctures $\theta_i$. Pick a choice of ends for $\cS$, regular Floer data $(H,J),(H^i,J^i)$ associated to $\theta,\theta_i$, and a regular compatible perturbation datum $(K,I)$ on $\cS$. Recall that in case $\cR = [0,\infty)$, the choice of ends and perturbation data has the special form described in \S\ref{par:compactness_families}. We keep a simpler numbering of Floer data and punctures, since the more specialized numbering described there is not needed for the purposes of the current definition.

Fix $\wt y_i \in \Crit \cA_{H^i}$, $\wt y \in \Crit \cA_H$ such that $1+|\wt y|' - \sum_i|\wt y_i|' = 0$. For $(r,u) \in \cM_\cS(K,I;\{\wt y_i\}_i,\wt y)$ the canonical orientation $\mfo_u \in \ddd(D_u)$ determines an isomorphism
$$C(u) \fc \bigotimes_i C(\wt y_i) \simeq C(\wt y)\,.$$
The matrix element of the homotopy operator $\Psi_{\cS;K,I}$ is the homomorphism
$$\sum_{(r,u) \in \cM_\cS(K,I;\{\wt y_i\}_i,\wt y)} C(u) \fc \bigotimes_i C(\wt y_i) \to C(\wt y)\,.$$

\subsubsection{Identities}\label{sss:identities}

Here we prove that various operations satisfy identities. There are three types of identities, expressing the fact that boundary operators square to zero, that multiplicative operators are chain maps, and that they satisfy algebraic identities, which in turn is proved using homotopy operators. Since every operation is defined in terms of matrix elements, in order to prove an identity, it suffices to prove that the corresponding combination of matrix elements of the operations involved equals zero. This is what we do here.

\paragraph{Boundary operators square to zero}\label{par:matrix_elts_bd_op_squared_vanish} Here we only prove this in the absence of bubbling. The remaining case is treated in \S\ref{ss:boundary_op_squares_zero_bubbling_HF}. It is enough to show that the matrix element of $\partial_{H,J}^2$ relative to $\wt y_\pm$ is zero, that is that
$$\sum_{\substack{\wt y \in \Crit \cA_H \\ [u]\in\cM(H,J;\wt y_-,\wt y) \\ [v] \in \cM(H,J;\wt y,\wt y_+)}}C(v)\circ C(u) = 0\,.$$
Consider the compactified moduli space $\ol \cM(H,J; \wt y_-,\wt y_+)$; the summands in the above sum are in bijection with its boundary points. For a boundary point $\delta$ let $C(\delta) \fc C(\wt y_-) \simeq C(\wt y_+)$ be the corresponding summand. It is enough to show that for each connected component $\Delta \subset \cM(H,J;\wt y_-,\wt y_+)$ having $\delta,\delta'$ as its boundary points, we have
$$C(\delta) + C(\delta') = 0\,.$$
Let $\wt \Delta \subset \wt \cM(H,J;\wt y_-,\wt y_+)$ be the component covering $\Delta$ under the quotient map $\wt\cM (H,J;\wt y_-, \wt y_+)\to \cM (H,J;\wt y_-, \wt y_+)$. By \S\ref{sss:orientations_isomorphisms}, we know that isomorphisms $C(\wt y_-) \simeq C(\wt y_+)$ are in bijection with orientations of $\wt \Delta$. We computed the orientation induced on $\wt \Delta$ from a boundary point in \S\ref{par:induced_orientation_translation_invt_pert_datum}. For $w \in \wt\Delta$ it is given by $-\partial_w \wedge \text{inward}_\delta$. Since clearly
$$\text{inward}_\delta = -\text{inward}_{\delta'}\,,$$
we see that $C(\delta) = -C(\delta')$.

\paragraph{Operations $\Phi_\Sigma$ are chain maps}\label{par:mult_ops_are_chain_maps}

We want to prove the vanishing of the matrix element
$$\sum_{\substack{\wt y' \in \Crit \cA_{H'} \\ u \in\cM_\Sigma(K,I;\{\wt y_i\}_i,\wt y') \\ [v] \in \cM(H,J;\wt y',\wt y)}} C(v)\circ C(u) 
-\sum_j (-1)^{\sum_{i<j}|\wt y_i|'} \sum_{\substack{\wt y_j' \in \Crit \cA_{H^j} \\ [u] \in\cM(H^j,J^j;\wt y_j,\wt y_j') \\ v \in \cM_\Sigma(K,I;\{\wt y_i\}_{i\neq j},\wt y_j',\wt y)}}  C(v)\circ (\id \otimes \dots \otimes C(u) \otimes \dots \otimes \id)$$
as a homomorphism
$$\bigotimes_jC(\wt y_j) \to C(\wt y)\,.$$
There is a bijection between the boundary points of the compactified $1$-dimensional moduli space $\ol\cM_\Sigma(K,I;\{\wt y_i\}_i,\wt y)$ and the summands of the above matrix element. For a boundary point $\delta \in \partial \ol \cM_\Sigma(K,I;\{\wt y_i\}_i,\wt y)$ we let $C(\delta)$ be the summand corresponding under this bijection to $\delta$, including the sign in front of it. It is enough to show that for each connected component $\Delta \subset \ol \cM_\Sigma(K,I;\{\wt y_i\}_i,\wt y)$ with $\partial \Delta = \{\delta,\delta'\}$ we have
$$C(\delta) + C(\delta') = 0\,.$$
Isomorphisms $\bigotimes_jC(\wt y_j) \to C(\wt y)$ are in bijection with orientations of $\Delta$. In \S\ref{par:induced_orientation_single_surf} we computed the orientations induced on $\Delta$ by the isomorphisms $C(\delta),C(\delta')$. The computations show precisely that these orientations are always opposite, whence the vanishing of the matrix element.

\paragraph{Chain homotopies for a single surface}\label{par:chain_htpies_single_surf} Consider the trivial family \footnote{Note however that the fiberwise conformal structure on $\cS$ is allowed to vary along $[0,1]$.} $\cS = \cR \times \Sigma$, $\cR = [0,1]$, where the ends and the perturbation datum are constant near $\partial\cR = \{0,1\}$. We want to prove the vanishing of the matrix element
\begin{multline*}
\sum_{u \in \cM_{\Sigma_1}(K_1,I_1;\{\wt y_i\}_i,\wt y)} C(u) - \sum_{u\in\cM_{\Sigma_0}(K_0,I_0; \{\wt y_i\}_i,\wt y)} C(u) -\\
- \sum_{\substack{\wt y' \in \Crit \cA_H \\ (r,u) \in \cM_{\cS}(K,I;\{\wt y_i\}_i,\wt y') \\ [v] \in \cM(H,J; \wt y',\wt y)}} C(v) \circ C(u) 
- \sum_j(-1)^{\sum_{i<j}|\wt y_i|'} \sum_{\substack{\wt y_j' \in \Crit \cA_{H^j} \\ [u] \in \cM(H^j,J^j;\wt y_j,\wt y_j') \\ (r,v) \in \cM_\cS(K,I;\{\wt y_i\}_{i\neq j},\wt y_j',\wt y)}} C(v) \circ C(u)
\end{multline*}
as a homomorphism
$$\bigotimes_jC(\wt y_j) \to C(\wt y)\,.$$
Again, the boundary points of the compactified $1$-dimensional moduli space $\ol\cM_\cS(K,I;\{\wt y_i\}_i,\wt y)$ are in bijection with the summands of the matrix element. If $C(\delta)$ denotes the summand corresponding to the boundary point $\delta$, including the sign in front of it, then it is enough to show that for every connected component $\Delta \subset \ol\cM_\cS(K,I;\{\wt y_i\}_i,\wt y)$ with $\partial \Delta = \{\delta,\delta'\}$ we have $C(\delta) + C(\delta') = 0$. Isomorphisms $\bigotimes_jC(\wt y_j) \to C(\wt y)$ correspond to orientations of $\Delta$. In \S\ref{par:induced_orientations_families} we computed the orientations induced on $\Delta$ by the isomorphisms $C(\delta),C(\delta')$. Examining these orientations we dedude $C(\delta) + C(\delta') = 0$ and therefore the vanishing of the matrix element.

\paragraph{Chain homotopies for a family of surfaces over $[0,\infty)$}\label{par:chain_htpies_families}

Keep the notations of \S\ref{par:compactness_families}. We wish to prove the vanishing of the matrix element
\begin{multline*}
\sum_{\substack{\wt y_{2,j} \in \Crit\cA_{H^1} \\ u\in \cM_{\Sigma^1}(K^1,J^1;\{\wt y_{1,i}\}_i,\wt y_{2,j}) \\ v \in \cM_{\Sigma^2}(K^2,I^2;\{\wt y_{2,i}\}_i,\wt y_2)}} C(v) \circ C(u)- \\
- \sum_{u \in \cM_{\Sigma_0}(K_0,I_0;\{\wt y_{1,i}\}_i,\{\wt y_{2,i}\}_{i\neq j}, \wt y_2)} C(u) -
\sum_{\substack{\wt y_2' \in \Crit \cA_{H^2} \\ (r,u)\in \cM_{\cS}(K,I;\{\wt y_{1,i}\}_i,\{\wt y_{2,i}\}_{i\neq j}, \wt y_2') \\ [v] \in \cM(H^2,J^2;\wt y_2',\wt y_2)}} C(v) \circ C(u)+ \\
- \sum_l \sum_{\substack{\wt y_{1,l}' \in \Crit \cA_{H^{1,l}} \\ [u] \in \cM(H^{1,l},J^{1,l};\wt y_{1,l},\wt y_{1,l}') \\ (r,v) \in \cM_{\cS}(K,I;\{\wt y_{1,i}\}_{i\neq l},\wt y_{1,l}',\{\wt y_{2,i}\}_{i\neq j}, \wt y_2)}} (-1)^{\sum_{i < j}|\wt y_{2,i}|' + \sum_{i < l} |\wt y_{1,i}|'} C(v) \circ C(u) \\
- \sum_{l\neq j} \sum_{\substack{\wt y_{2,l}' \in \Crit \cA_{H^{2,l}} \\ [u] \in \cM(H^{2,l},J^{2,l};\wt y_{2,l},\wt y_{2,l}') \\ (r,v) \in \cM_{\cS}(K,I;\{\wt y_{1,i}\}_i,\{\wt y_{2,i}\}_{i\neq j,l}, \wt y_{2,l}', \wt y_2)}} (-1)^{\text{sign}(l)} C(v) \circ C(u)
\end{multline*}
as a homomorphism
$$\bigotimes_{i < j} C(\wt y_{2,i}) \otimes \bigotimes_i C(\wt y_{1,i}) \otimes \bigotimes_{i > j} C(\wt y_{2,i}) \simeq C(\wt y_2)\,.$$
Here
$$\text{sign}(l) = \left\{\begin{array}{ll} \sum_{i < l} |\wt y_{2,i}|'  & \text { if } l < j \\ \sum_{i < l, i \neq j} |\wt y_{2,i}|' + \sum_{i=1}^m |\wt y_{1,i}|'  & \text { if } l > j \end{array} \right .$$
where $m$ is the number of negative ends of $\Sigma_1$.

The argument is identical to the above: the summands of the matrix element are in bijection with boundary points of the compactification $\ol \cM_\cS(K,I;\{\wt y_{1,i}\}_i,\{\wt y_{2,i}\}_{i\neq j}, \wt y_2)$; let $C(\delta)$ be the summand corresponding to the boundary point $\delta$, including any signs in front of it. It suffices to show that for any connected component $\Delta$ of the compactified space with $\partial \Delta = \{\delta,\delta'\}$ we have $C(\delta) + C(\delta') = 0$. The bijection between isomorphisms $\bigotimes_{i < j} C(\wt y_{2,i}) \otimes \bigotimes_i C(\wt y_{1,i}) \otimes \bigotimes_{i > j} C(\wt y_{2,i}) \simeq C(\wt y_2)$ and orientations of $\Delta$, plus the computation of induced orientations done in \S\ref{par:induced_orientations_families} imply  $C(\delta) + C(\delta') = 0$ and hence the vanishing of the desired matrix element.

\subsection{Floer homology}\label{ss:HF}

Here we define Lagrangian and Hamiltonian Floer homology, and the various algebraic operations thereupon.

\subsubsection{Lagrangian Floer homology}\label{sss:Lagr_HF}

Choose a regular Floer datum $(H,J)$ for the strip $S$, which is sufficiently generic so that the compactness results of \S\ref{sss:compactness_gluing} hold. Define the $\Z$-module
$$CF_*(H:L) = \bigoplus_{\wt y \in \Crit \cA_{H:L}}C(\wt y)\,.$$
This is graded by $m_{H:L}$.

\paragraph{The boundary operator} The boundary operator
$$\partial_{H,J} \fc CF_j(H:L) \to CF_{j-1}(H:L)$$
is defined via its matrix elements, see \S\ref{par:matrix_elts_boundary_op}. The vanishing of the matrix elements of $\partial_{H,J}^2$, \S\ref{par:matrix_elts_bd_op_squared_vanish} together with the results of \S\ref{ss:boundary_op_squares_zero_bubbling_HF} show that $\partial_{H,J}^2 = 0$ and therefore we can define the \tb{Lagrangian Floer homology}
$$HF_*(H,J:L)$$
as the homology of the \tb{Lagrangian Floer complex} $(CF_*(H:L),\partial_{H,J})$.

\paragraph{Continuation maps} Let $(H^i,J^i)$ be regular Floer data, $i=0,1$, associated to the ends of $S$, $(H^0,J^0)$ to the negative end, and $(H^1,J^1)$ to the positive end, and choose a regular perturbation datum $(K,I)$ on $S$ which is of the form
\begin{equation}\label{eqn:pert_datum_continuation_maps}
K(s,t)=H^s_t\,dt\,,\quad I(s,t) = J^s_t
\end{equation}
where $(H^s,J^s)_{s\in \R}$ is a smooth homotopy of Floer data, which is independent of $s$ for $s \notin (0,1)$. The corresponding \tb{continuation map}
$$\Phi_{S;K,I} \fc CF_j(H^0:L) \to CF_j(H^1:L)$$
is determined by its matrix elements, see \S\ref{par:matrix_elts_operations}. According to \S\ref{par:mult_ops_are_chain_maps}, it is a chain map:
$$\Phi_{S;K,I}\circ \partial_{H^0,J^0} = \partial_{H^1,J^1}\circ \Phi_{S;K,I}$$
and therefore it induces a map on homology
$$\Phi_{S;K,I} \fc HF_j(H^0,J^0:L) \to HF_j(H^1,J^1:L)\,.$$
If $(K',I')$ is another regular perturbation datum on $S$ as above, that is it corresponds to a different homotopy of Floer data, then the maps $\Phi_{S;K,I}$ and $\Phi_{S;K',I'}$ are chain homotopic. Indeed, consider the trivial family $\cS = S \times [0,1]$, where we associate the data $(H^i,J^i)$ to the ends of $S$ as above, and choose a regular compatible perturbation datum $(\ol K,\ol I)$ on $\cS$ which near $r=0$ equals $(K,I)$ and near $r=1$ equals $(K',I')$. It follows from \S\ref{par:chain_htpies_single_surf} that the operator
$$\Psi_{\cS;\ol K,\ol I} \fc CF_j(H^0:L) \to CF_{j+1}(H^1:L)$$
determined by its matrix elements \S\ref{par:matrix_elts_htpy_ops}, is a chain homotopy between $\Phi_{S;K,I}$ and $\Phi_{S;K',I'}$. Therefore the induced map on homology is independent of the choice of homotopy of Floer data; we denote it
$$\Phi_{H^0,J^0}^{H^1,J^1}\fc HF_*(H^0,J^0:L) \to HF_*(H^1,J^1:L)\,.$$

Consider now the trivial family $\cS = S \times [0,\infty)$, where for some $R_0 > 0$ the ends and a regular perturbation datum $(K,I)$ on $\cS$ come from gluing $\Sigma^1=\Sigma^2=S$ as described in \S\ref{par:compactness_families}, where the positive end of $\Sigma^1$ is glued to the negative end of $\Sigma^2$, the Floer data associated to the ends of $\Sigma^1$ are $(H^0,J^0),(H^1,J^1)$, the data associated to the ends of $\Sigma^2$ are $(H^1,J^1),(H^2,J^2)$, and there are perturbation data $(K^i,J^i)$ on $\Sigma_i$ of the form \eqref{eqn:pert_datum_continuation_maps}. Using the matrix elements defined in \S\ref{par:matrix_elts_htpy_ops} we can define a homotopy operator
$$\Psi_{\cS;K,I} \fc HF_j(H^0,J^0:L) \to HF_{j+1}(H^2,J^2:L)$$
between
$$\Phi_{H^0,J^0}^{H^2,J^2}\quad \text{and} \quad \Phi_{H^1,J^1}^{H^2,J^2} \circ \Phi_{H^0,J^0}^{H^1,J^1}\,.$$
This follows from the vanishing of the matrix elements, \S\ref{par:chain_htpies_single_surf}. This means that on homology we have
\begin{equation}\label{eqn:cocycle_id_continuation_maps_Lagr_HF}
\Phi_{H^0,J^0}^{H^2,J^2} = \Phi_{H^1,J^1}^{H^2,J^2} \circ \Phi_{H^0,J^0}^{H^1,J^1}\,.
\end{equation}

Let now $(H,J)$ be a regular Floer datum on $S$ and let $(K,I)$ be the corresponding translation-invariant perturbation datum. It can be seen that the continuation map $\Phi_{S;K,I}$ is the identity on chain level. This means that
$$\Phi_{H,J}^{H,J} = \id_{HF_*(H,J:L)}\,.$$
Combining with \eqref{eqn:cocycle_id_continuation_maps_Lagr_HF}, we obtain
$$\Phi_{H^0,J^0}^{H^1,J^1} = (\Phi_{H^1,J^1}^{H^0,J^0})^{-1}\,,$$
and in particular continuation maps are isomorphisms.

We can now define the \tb{abstract Floer homology} $HF_*(L)$ as the limit of the system of Floer homologies connected by the continuation isomorphisms.

\paragraph{The product}

Let $\Sigma$ be the disk with three boundary punctures $\theta_i$, $i=0,1,2$, arranged in positive cyclic order on $\partial D^2$; here $\theta_2$ is positive and the other punctures are negative. Endow $\Sigma$ with a choice of ends. Let $(H^i,J^i)$ be regular Floer data associated to the $\theta_i$ and let $(K,I)$ be a regular compatible perturbation datum on $\Sigma$. The matrix elements \S\ref{par:matrix_elts_operations} determine an operation
$$\Phi_{\Sigma;K,I} \fc CF_j(H^0:L) \otimes CF_k(H^1:L) \to CF_{j+k-n}(H^2:L)\,.$$
The vanishing of the matrix elements \S\ref{par:chain_htpies_single_surf} shows that $\Phi_{\Sigma;K,I}$ is a chain map, therefore we have a bilinear operation on homology
$$\Phi_{\Sigma;K,I} \fc HF_j(H^0,J^0:L) \otimes HF_k(H^1,J^1:L) \to HF_{j+k-n}(H^2,J^2:L)\,,$$
which is the \tb{product}. A change of auxiliary data, that is the choice of a conformal structure on $\Sigma$, the position of the punctures, the ends or perturbation datum, can be encoded in a family over $[0,1]$, which then produces a homotopy operator between the corresponding chain maps. This means that we have a well-defined operation
$$\star \fc HF_j(H^0,J^0:L) \otimes HF_k(H^1,J^1:L) \to HF_{j+k-n}(H^2,J^2:L)\,,$$
which only depends on the Floer data.

Gluing $\Sigma$ to a strip with perturbation data corresponding to a change of the Floer datum leads to a family $\cS \to [0,\infty)$, which can be used to produce a homotopy operator, which implies that on homology we have
$$\star \circ (\Phi \otimes \Phi) = \Phi \circ \star\,,$$
where $\Phi$ is shorthand notation for a continuation map, that is the product is compatible with the continuation maps. This in particular means that we have a well-defined product on the abstract Floer homology $HF_*(L)$.

\paragraph{Associativity of the product}

Let $\Sigma$ be the disk with four boundary punctures $\theta_i$, $i=0,1,2,3$, arranged in positive cyclic order, where $\theta_3$ is the only positive puncture. Fix a choice of ends for $\Sigma$, and let $(H^i,J^i)$ be regular Floer data associated to the $\theta_i$, and $(K,I)$ a regular compatible perturbation datum on $\Sigma$. The corresponding operation
$$\Phi_{\Sigma;K,I} \fc CF_j(H^0:L) \otimes CF_k(H^1:L) \otimes CF_l(H^2:L) \to CF_{j+k+l-2n}(H^3:L)$$
can be proved to be a chain map, using the same arguments as above, which means it descends to homology. It is then shown to be independent of the auxiliary data. Next one shows that it is compatible with the continuation morphisms, therefore it defines a ternary operation on $HF_*(L)$.

Note now that $\Sigma$ can glued from two surfaces $\Sigma^1=\Sigma^2$ which are disks with three boundary punctures of which only one is positive. This gluing can be done in two different ways: one can feed the positive puncture of $\Sigma_1$ into either one of the negative punctures of $\Sigma_2$. The resulting families $\cS \to [0,\infty)$ can be used to show that both compositions
$$\star \circ (\id \otimes \star)\quad\text{and}\quad \star \circ (\star \otimes \id)$$
are equal, on the level of homology, to the ternary operation we have just defined. In particular it means that the product $\star$ on $HF_*(L)$ is associative.

\paragraph{The unit}

Take now $\Sigma$ to be the disk with one positive puncture, endow it with an end, and fix a regular nondegenerate Floer datum $(H,J)$ and a compatible perturbation datum $(K,I)$. There is the corresponding nullary operation $\Phi_{\Sigma;K,I}$, which plays the role of the unit in Floer homology. Let us give a more detailed description of it. Fix $\wt y \in \Crit \cA_{H:L}$. The moduli space $\cM_{\Sigma}(K,I;\wt y)$ has dimension $|\wt y|' = n - |\wt y|$, therefore it is zero-dimensional whenever $|\wt y| = n$. Assume this. For any $u \in \cM_\Sigma(K,I;\wt y)$ the operator $D_u$ is an isomorphism therefore it is canonically oriented by the positive orientation $\mfo_u = 1 \otimes 1^\vee$. Orientations of $D_u$ are in bijection with isomorphisms $\bigotimes_\varnothing \simeq C(\wt y)$. The empty tensor product \footnote{If the reader does not like abstract nonsense, here is another way to do this: the linearized operator $D_u$ by definition belongs to the family $D_{\wt y}$, and since it is oriented, there is the corresponding generator of $C(\wt y)$; in our notation this is just $C(u)(1)$.} by definition is just the $\Z$-module $\Z$, therefore the canonical orientation of $D_u$ gives rise to an isomorphism $C(u) \fc \Z \simeq C(\wt y)$. Taking the sum:
$$\sum_{u \in \cM_\Sigma(K,I;\wt y)} C(u) \fc \Z \to C(\wt y)$$
gives the matrix element of the operation $\Phi_{\Sigma;K,I}$. Passing to the whole complex, we see that this amounts to a graded linear map
$$\Z[n] \to CF_*(H:L)\,,$$
where $\Z[n]$ denotes the graded abelian group which has $\Z$ in degree $n$ and $0$ everywhere else. This is the unit. Using the same methods as above it is shown to be a chain map, to be independent of the auxiliary data, and to be compatible with continuation morphisms. Therefore we have the canonical nullary operation
$$1 \fc \Z[n] \to HF_*(L)\,.$$

Let $\Sigma_\star$ be the disk with three boundary punctures, two negative and one positive puncture. We can glue $\Sigma$ to either one of the negative punctures of $\Sigma_\star$, and the resulting surface is isomorphic to the strip $S$ with one positive and one negative puncture. Using the same arguments as above, we see that this leads to the following identities:
$$\star \circ (1 \otimes \id) = \star \circ (\id \otimes 1) = \id$$
on $HF_*(L)$. This means precisely that the operation $1$ we have just defined is the unit for the product $\star$.

\subsubsection{Hamiltonian Floer homology}\label{sss:Hamiltonian_HF}

The treatment here is entirely parallel to the Lagrangian case, so we just outline the main arguments and results and establish notation.

Fix a regular nondegenerate Floer datum $(H,J)$ for the cylinder $C$, which is sufficiently generic so that the compactness results of \S\ref{sss:compactness_gluing} hold. Define the $\Z$-module
$$CF_*(H) = \bigoplus_{\wt y\in\Crit\cA_H}C(\wt y)\,.$$
This is graded by $m_H$.

\paragraph{The boundary operator}

The matrix elements of \S\ref{par:matrix_elts_boundary_op} assemble into the boundary operator
$$\partial_{H,J} \fc CF_j(H) \to CF_{j-1}(H)\,.$$
That $\partial_{H,J}^2 = 0$ follows from the vanishing of the corresponding matrix elements, see \S\ref{par:matrix_elts_bd_op_squared_vanish}. We let the \tb{Hamiltonian Floer homology}
$$HF_*(H,J)$$
be the homology of $(CF_*(H),\partial_{H,J})$.

\paragraph{Continuation maps}

These are defined in the exact same manner as in the Lagrangian case, and we leave the details to the reader. The upshot is that for any regular nondegenerate Floer data $(H^i,J^i)$, $i=0,1$, there is a well-defined morphism
$$\Phi_{H^0,J^0}^{H^1,J^1} \fc HF_*(H^0,J^0) \to HF_*(H^1,J^1)$$
which satisfies the cocycle identity and such that
$$\Phi_{H,J}^{H,J} = \id_{HF_*(H,J)}\,.$$
It follows that the continuation maps are isomorphisms. Therefore the abstract Floer homology $HF_*(M)$ can be defined as the limit of the system of Floer homologies for different Floer data, connected by the continuation isomorphisms.

\paragraph{The product}

This is defined analogously to the Lagrangian case, but there are some minor differences. The product here is defined using the surface $\Sigma$ which is the sphere $S^2$ with three punctures, two negative and one positive. It therefore gives rise to an operation on homology
$$* \fc HF_j(H^0,J^0) \otimes HF_k(H^1,J^1) \to HF_{j+k-2n}(H^2,J^2)$$
for any regular Floer data $(H^i,J^i)$ associated to the punctures of $\Sigma$. One shows that this is independent of the auxiliary data such as the choice of ends (since the space of ends around an interior puncture is connected), and the perturbation datum, and therefore the operation $*$ is well-defined on homology. It can also be shown to respect continuation maps, which means that there is a well-defined product on the abstract Floer homology $HF_*(M)$. Associativity is proved in the same way as in the Lagrangian case, by gluing two copies of $\Sigma$.

The difference from the Lagrangian case is the fact that the product $*$ is supercommutative. This can be seen as follows. Let $(K,I)$ be a regular compatible perturbation datum on $\Sigma$. There is an orientation-preserving diffeomorphism of $\Sigma$ which preserves the positive puncture and exchanges the negative ones. Let $(K',I')$ be the perturbation datum obtained by pushing forward the datum $(K,I)$ by this diffeomorphism. For every critical points $\wt y_i \in \Crit \cA_{H^i}$ we have a canonical identification of moduli spaces
$$\cM_\Sigma(K,I;\wt y_0,\wt y_1,\wt y_2) \simeq \cM_\Sigma(K',I';\wt y_1,\wt y_0,\wt y_2)\,,$$
where on the left we have the original conformal structure and ends while on the right we use the conformal structure and ends pushed forward by the diffeomorphism. This means that the corresponding matrix elements
$$\sum_{u \in \cM_\Sigma(K,I;\wt y_0,\wt y_1,\wt y_2)} C(u) \fc C(\wt y_0) \otimes C(\wt y_1) \to C(\wt y_2) \text{ and}\sum_{u \in \cM_\Sigma(K',I';\wt y_1,\wt y_0,\wt y_2)} C(u) \fc C(\wt y_1) \otimes C(\wt y_0) \to C(\wt y_2)$$
only differ by the Koszul sign $(-1)^{|\wt y_0|'|\wt y_1|'}$ which arises when we compose the direct sum isomorphisms
$$\ddd(D_{\wt y_0}) \otimes \ddd(D_{\wt y_1}) \simeq \ddd(D_{\wt y_0} \oplus D_{\wt y_1}) \simeq \ddd(D_{\wt y_1}) \otimes \ddd(D_{\wt y_0})\,,$$
see \S\ref{par:direct_sum_isos}. Therefore we can see the supercommutativity already on the chain level.

\paragraph{The unit}

This again is defined analogously to the Lagrangian case, the difference being that the unit is now a graded linear map
$$1: \Z[2n] \to HF_*(M)\,.$$

\subsubsection{Quantum module action}\label{sss:quantum_module_action}

This is defined as follows. Let $\Sigma_\bullet$ be the disk with two boundary punctures, one of which is positive, and one interior puncture which is negative. Endow $\Sigma_\bullet$ with a choice of ends, and let $(H^i,J^i)$ be regular Floer data associated to the punctures, the one with $i=1$ to the negative boundary puncture, $i=0$ to the interior puncture, and $i=2$ to the positive boundary puncture. Assume $(H^0,J^0)$ is $1$-periodic in time. Pick a regular compatible perturbation datum $(K,I)$. The matrix elements \S\ref{par:matrix_elts_operations} define an operation
$$\Phi_{\Sigma_\bullet;K,I} \fc CF_j(H^0) \otimes CF_k(H^1:L) \to CF_{j+k-2n} (H^2:L)\,.$$
As usual, one shows that this is a chain map, that the resulting map on homology is independent of the auxiliary choices, that it then respects the continuation maps, and therefore defines a bilinear operation on abstact homologies:
$$\bullet \fc HF_*(M) \otimes HF_*(L) \to HF_{*-2n}(L)\,.$$

Gluing $\Sigma_\bullet$ to different surfaces, we can prove the following identities:
$$\bullet \circ (1 \otimes \id) = \id\,,$$
that is the unit $1 \in HF_{2n}(M)$ acts as a unit;
$$\bullet \circ (* \otimes \id) = \bullet \circ (\id \otimes \bullet)\,,$$
which means $\bullet$ defines a module action; finally, we have
$$\bullet \circ (\id \otimes \star) = \star \circ (\bullet \otimes \id) = \star \circ (\id \otimes \bullet)\circ (R \otimes \id)\,, $$
where
$$R \fc HF_*(M) \otimes HF_*(L) \to HF_*(L) \otimes HF_*(M)$$
is the interchange of factors multiplied with the corresponding Koszul signs. This means that $HF_*(L)$ becomes a superalgebra over $HF_*(M)$ by means of $\bullet$.

One can also substitute the Lagrangian unit into $\bullet$ and get the so-called closed-open map, which is a degree $-n$ operation
$$\bullet \circ (\id \otimes 1) \fc HF_*(M) \to HF_{n-*}(L)\,,$$
which can be shown to be an algebra morphism.

\subsection{Arbitrary rings and local coefficients}\label{ss:arbitrary_rings_loc_coeffs}

The above definitions are made over the ground ring $\Z$. Given an arbitrary commutative ring $R$, we can form Floer complexes over $R$ by tensoring, $- \otimes_\Z R$. Thus we obtain the Floer homology over $R$, and all the above algebraic operations become $R$-linear.

In case $2=0$ in $R$, we can form the Floer complex without making assumption \tb{(O)}, as follows. For a nondegenerate Floer datum $(H,J)$ we define
$$CF_*(H:L) = \bigoplus_{\wt\gamma \in \Crit \cA_{H:L}}R \cdot \wt\gamma\,.$$
The boundary operator is given by counting the moduli spaces $\cM(H,J;\wt\gamma_-,\wt\gamma_+)$ modulo $2$. The algebraic structures are defined similarly.

Recall that a \tb{local system} on a topological space $X$ is a functor from the fundamental groupoid $\Pi_1(X)$ to the category of groups. In more elementary terms, it is given by assigning a group $G_x$ to every point $x \in X$ and an isomorphism $G_x \simeq G_y$ for every homotopy class of paths from $x$, where the isomorphisms behave coherently with respect to concatenation of paths. Similarly we can define a local system of $R$-modules.

Given a ground ring $R$, a \tb{flat $R$-bundle} over $\wt\Omega_L$ is by definition just a local system of $R$-modules over $\wt\Omega_L$. Let $\cE$ be such a bundle and for a path $u$ in $\wt\Omega_L$ running from $\wt\gamma_-$ to $\wt\gamma_+$ let
$$\cP_u \fc \cE_{\wt\gamma_-} \to \cE_{\wt\gamma_+}$$
be the corresponding parallel transport isomorphism. We can then form the Floer complex of $(H,J)$ \tb{twisted by} $\cE$, as follows. As an $R$-module, we have
$$CF_*(H:L;\cE) = \bigoplus_{\wt\gamma \in \Crit \cA_{H:L}} C(\wt\gamma) \otimes_R \cE_{\wt\gamma}\,.$$
The boundary operator has matrix elements
$$\sum_{[u] \in \cM(H,J;\wt\gamma_-,\wt\gamma_+)} C(u) \otimes_R \cP_u\,,$$
where we lift $u$ to a path in $\wt\Omega_L$. The corresponding Floer homology is denoted
$$HF_*(H,J:L;\cE)\,.$$
The continuation maps are similarly defined and we obtain the abstract Floer homology
$$HF_*(L;\cE)\,.$$

We will employ the following useful piece of notation. For $V$ a real $1$-dimensional vector space we let $|V|$ be its \tb{normalization}, which is the rank $1$ free $\Z$-module generated by its two possible orientations. If $\cL \to B$ is a real line bundle, we let $|\cL|$ be the flat $\Z$-bundle with fibers $|\cL_b|$ for $b \in B$. For a ring $R$, the $R$-normalization of $\cL$ is the flat locally free $R$-bundle of rank $1$ with fibers $|\cL_b|\otimes_\Z R$.

\subsection{Duality}\label{ss:duality_HF}

\subsubsection{Dual complexes and dual Hamiltonians}\label{sss:dual_cxs_dual_Hams}

We only treat the Lagrangian case in detail, the Hamiltonian case being entirely analogous. Fix a regular Floer datum $(H,J)$ on the strip $S$. We define
$$CF^*(H:L) = \bigoplus_{\wt y \in \Crit \cA_{H:L}}C(\wt y)^\vee$$
where
$$C(\wt y)^\vee = \Hom_\Z(C(\wt y),\Z)\,.$$
We grade this by $m_{H:L}$. Note for future use that there is a canonical isomorphism
$$C(\wt y)^\vee \equiv C(\wt y)\,,\quad c^\vee \mapsto c\,,$$
where $c \in C(\wt y)$ is a generator and $\langle c^\vee, c\rangle = 1$. The differential $\partial_{H,J}^\vee$ on $CF^*(H:L)$ is defined to be the dual of $\partial_{H,J}$: its matrix element
$$C(\wt y_+)^\vee \to C(\wt y_-)^\vee$$
is the dual of the matrix element of $\partial_{H,J}$ as a map $C(\wt y_-) \to C(\wt y_+)$ for $\wt y_\pm \in \Crit \cA_{H:L}$ of index difference $1$. We define another differential
$$\delta_{H,J} \fc CF^k(H,J:L) \to CF^{k+1}(H,J:L)\quad \text{by} \quad \delta_{H,J} = (-1)^{k-1}\partial_{H,J}^\vee\,.$$
We refer to the cochain complex
$$(CF^*(H,J:L),\delta_{H,J})$$
as the \tb{dual complex} and we let
$$HF^*(H,J:L)$$
be its cohomology, called the \tb{Floer cohomology} of $(H,J)$.

In an analogous fashion we can define the twisted dual complex
$$CF^*(H:L;\cE)$$
for a flat $\Z$-bundle $\cE$ over $\wt\Omega_L$ and the corresponding cohomology
$$HF^*(H,J:L;\cE)\,.$$

Let us define the dual Hamiltonian of $H$ to be
$$\ol H(t,x) = -H(1-t,x)\,.$$
This generates the flow obtained from the flow of $H$ by retracing it backward, that is $\phi^t_{\ol H} = \phi^{1-t}_H\phi_H^{-1}$. There is a bijection between orbits of $H$ and $\ol H$ given by
$$y \mapsto \ol y\,,\quad \ol y(t) = y(1-t)\,.$$
If $\wh y$ is a capping of $y$, then
$$\ol{\wh y} \fc \dot D^2 \to M \quad \text{defined by} \quad \ol{\wh y} (\sigma,\tau) = \wh y(\sigma,-\tau)$$
is a capping of $\ol y$ (using the same end for both maps). This establishes a bijection 
$$\Crit \cA_{H:L} \simeq \Crit \cA_{\ol H:L}\,,\quad \wt y=[y,\wh y] \mapsto \ol{\wt y}:=[\ol y,\ol{\wh y}]\,.$$
We have
$$m_{\ol H:L}(\ol{\wt y}) = n - m_{H:L}(\wt y)\quad \text{and} \quad \cA_{\ol H:L}(\ol{\wt y}) = - \cA_{H:L}(\wt y)\,.$$
This first relation is shown below \eqref{eqn:grading_dual_cx_Lagr_HF}.

We have the regular Floer datum $(\ol H,\ol J)$ where $\ol J_t(x) = J_{1-t}(x)$. Therefore the Floer complex
$$(CF_*(\ol H:L),\partial_{\ol H,\ol J})$$
is well-defined. For $\wt y_\pm \in \Crit \cA_{H:L}$ of index difference $1$ we have a canonical diffeomorphism between moduli spaces
$$\wt \cM(H,J;\wt y_-,\wt y_+) \simeq \wt \cM(\ol H,\ol J; \ol{\wt y}_+,\ol{\wt y}_-)\,,\quad u \mapsto \ol u\,,\text{ where }\ol u(s,t) = u(-s,1-t)\,.$$

For a capping $\wh y$ of an orbit $y$ of $H$ let us define the map
$$-\wh y \fc \dot D^2_- \to M\,,\quad -\wh y(\sigma,\tau) = \wh y(-\sigma,\tau)\,.$$
Here $\dot D^2_- = D^2 - \{-1\}$ has $-1$ as a negative puncture; we endow it with the standard negative end given by precomposing the standard end \eqref{eqn:std_end_cappings_Lagr_HF} with $(s,t) \mapsto (-s,-t)$. The map $\phi \fc \dot D^2 \to \dot D^2_-$, $\phi(z) = -z$ is a conformal isomorphism which also intertwines the standard ends. Therefore the linearized operators $D_{\ol{\wh y}}$ and $D_{-\wh y}$ are isomorphic in the sense of \S\ref{s:determinants}. It follows that their determinant lines are isomorphic as well, and in fact extend to a canonical isomorphism of determinant bundles
\begin{equation}\label{eqn:iso_det_lines_cappings_duality_Lagr_HF}
\ddd(D_{\ol{\wt y}}) \simeq \ddd(D_{-\wt y})\,.
\end{equation}
Next, the maps $\wh y,-\wh y$ have matching asymptotics and therefore can be preglued. It is not hard to see that the resulting map is homotopic through maps $(D^2,\partial D^2) \to (M,L)$ to the constant map at $y(0)$. Using the direct sum, linear gluing, and deformation isomorphisms, we obtain
$$\ddd(D_{\wh y}) \otimes \ddd(D_{-\wh y}) \simeq \ddd(D_{\wh y} \oplus D_{-\wh y}) \simeq \ddd(T_{y(0)}L)\,.$$
In particular $n = \ind 0_{T_{y(0)}L} = \ind D_{\wh y} + \ind D_{-\wh y}$ and thus
$$\ind D_{\ol {\wt y}} = \ind D_{-\wh y} = n - \ind D_{\wh y}\,,$$
which implies
\begin{equation}\label{eqn:grading_dual_cx_Lagr_HF}
m_{\ol H:L}(\ol{\wt y}) = n - \ind D_{\ol{\wt y}} = \ind D_{\wh y} = n - m_{H:L}(\wt y)\,.
\end{equation}
We thus have a canonical isomorphism
$$\ddd(D_{\wh y}) \otimes \ddd(D_{\ol{\wh y}}) \simeq \ddd(T_{y(0)}L)\,.$$
Tensoring with $\ddd(T_{y(0)}L)$ and noting that the square of a real line bundle is canonically oriented, we obtain
$$C(\wt y) \otimes |\cL_{\wt y}| \otimes C(\ol{\wt y}) \simeq \Z\,,$$
where $\cL$ is the flat $\Z$-bundle bundle on $\wt \Omega_L$ obtained by pulling back the bundle $|\ddd(TL)|$ on $L$ via the evaluation map $\wt \Omega_L \to L$, $\wt y \mapsto y(0)$.

Therefore we have the isomorphism
$$CF^{n-*}(H:L;\cL) \simeq CF_*(\ol H:L)$$
of graded $\Z$-modules. In fact this can be extended to a chain isomorphism, as follows. We have the following commutative diagram for $\wt y_\pm \in\Crit \cA_{H:L}$ of index difference $1$ and $u \in \wt\cM(H,J;\wt y_-,\wt y_+)$, obtained by employing the direct sum, linear gluing, and deformation isomorphisms:
$$\xymatrix{\ddd(TL) \ar@{=}[r]& \ddd(D_{\wh y_-} \oplus D_u \oplus D_{-\wh y_+}) & \ddd(D_{\wh y_+}) \otimes \ddd(D_{-\wh y_+}) \ar [l]\\
&\ddd(D_{\wh y_-}) \otimes \ddd(D_{-\wh y_-}) \ar[u] & \ddd(D_{\wh y_-}) \otimes \ddd(D_u) \otimes \ddd(D_{-\wh y_+}) \ar[u]^{(C(u) \otimes \id) \circ (R \otimes \id)} \ar[l]}$$
Where $R$ is the interchange of factors together with the Koszul sign $(-1)^{\ind D_u\cdot |\wt y_-|'} = (-1)^{|\wt y|'}$. Fix $\mfo_{y_-} \in C(\wt y_-)$ and $\mfo \in \ddd(TL)$ and let $\mfo_{y_+} = C(u)(\mfo_{y_-})$, and $\mfo_{-y_\pm} \in \ddd(D_{-\wt y_\pm})$ be such that the isomorphisms
$$\ddd(TL) \simeq \ddd(D_{\wt y_\pm}) \otimes \ddd(D_{-\wt y_\pm})$$
map $\mfo \mapsto \mfo_{y_\pm}\otimes \mfo_{-y_\pm}$. Then the above diagram maps
$$\xymatrix{\mfo & \mfo_{y_+} \otimes \mfo_{-y_+}  \ar@{|->} [l]\\
\mfo_{y_-} \otimes \mfo_{-y_-} \ar @{|->} [u] & (-1)^{|\wt y_-|'} \mfo_{y_-} \otimes \partial_u \otimes \mfo_{-y_+} \ar @{|->} [u] \ar @{|->} [l]}$$
The dual diagram, obtained by replacing all the operators with their ``minus-bar'' version, reads:
$$\xymatrix{\ddd(TL) \ar@{=}[r]& \ddd(D_{-\ol{\wh y}_-} \oplus D_{\ol u} \oplus D_{\ol{\wh y}_+}) & \ddd(D_{-\ol{\wh y}_+}) \otimes \ddd(D_{\ol{\wh y}_+}) \ar [l]\\
&\ddd(D_{-\ol{\wh y}_-}) \otimes \ddd(D_{\ol{\wh y}_-}) \ar[u] & \ddd(D_{-\ol{\wh y}_-}) \otimes \ddd(D_{\ol u}) \otimes \ddd(D_{\ol{\wh y}_+}) \ar[u] \ar[l]_-{\id \otimes C(\ol u)}}$$
Let now $\mfo_{\ol y_\pm} \in \ddd(D_{\ol{\wt y}_\pm})$ correspond to $\mfo_{-y_\pm} \in \ddd(D_{-\wt y_\pm})$ under the isomorphism \eqref{eqn:iso_det_lines_cappings_duality_Lagr_HF} and let $\mfo'_{\ol y_-} = C(\ol u) (\mfo_{\ol y_+})$. The dual diagram maps
$$\xymatrix{\mfo & \mfo_{-\ol y_+} \otimes \mfo_{\ol y_+}  \ar@{|->} [l]\\
\mfo_{-\ol y_-} \otimes \mfo_{\ol y_-} \ar @{|->} [u] & (-1)^{|\wt y_-|'} \mfo_{-\ol y_-} \otimes \partial_{\ol u} \otimes \mfo_{\ol y_+} \ar @{|->} [u] \ar @{|->} [l]}$$
whence it follows that $C(\ol u)$ maps $\mfo_{\ol y_+}$ to $(-1)^{|\wt y_-|'}\mfo_{\ol y_-}$. Thus the following diagram commutes:
$$\xymatrix{C(\ol{\wt y}_+) \ar[d]^{C(\ol u)} \ar[r] & [C(\wt y_+) \otimes \cL_{\wt y_+}]^\vee \ar [d] ^{(-1)^{|\wt y_-|'} [C(u) \otimes \cP_u]^\vee}\\
C(\ol{\wt y}_-) \ar[r] & [C(\wt y_-) \otimes \cL_{\wt y_-}]^\vee}$$
This means that we have obtained a canonical isomorphism of chain complexes
\begin{equation}\label{eqn:duality_iso_chain_cxs_Lagr_HF}
(CF^{n-*}(H:L;\cL),\delta_{H,J}\otimes \cP) = (CF_*(\ol H:L),\partial_{\ol H,\ol J})\,,
\end{equation}
where $\cP$ is the parallel transport operator of $\cL$. In particular we obtain a canonical isomorphism of homologies
$$HF^{n-*}(H,J:L;\cL) \simeq HF_*(\ol H, \ol J:L)\,.$$
This is the \tb{duality isomorphism}.

\subsubsection{Augmentation}\label{sss:augmentation_HF}

As we saw above in \S\ref{sss:unit_Lagr_QH}, the unit in Floer homology can be viewed as a chain map
$$\Z[n] \to CF_*(\ol H:L)\,.$$
Therefore we have the dual map
$$CF^*(\ol H:L) \to \Z[n]\,,$$
which is also a chain map. Since the former cochain complex is canonically isomorphic to $CF_{n-*}(H:L;\cL)$ by the duality isomorphism \eqref{eqn:duality_iso_chain_cxs_Lagr_HF}, we get the chain map
$$CF_*(H:L;\cL) \to \Z\text{ (in degree zero)}\,.$$
This is the \tb{augmentation map}. Since it is a chain morphism, we get the induced map on Floer homology
$$HF_*(H,J:L;\cL) \to \Z\,.$$
Since the unit commutes with continuation maps, the same is true of the augmentation and therefore we obtain the canonical augmentation map on the abstact Floer homology
$$\epsilon \fc HF_*(L;\cL) \to \Z\,.$$

\section{Quantum homology}\label{s:QH}

In this section we define Lagrangian quantum homology as well as the quantum homology of our symplectic manifold $M$, and various algebraic operations on them. Lagrangian quantum homology was constructed by Biran--Cornea, see \cite{Biran_Cornea_Quantum_structures_Lagr_submfds, Biran_Cornea_Rigidity_uniruling} and references therein. The present section uses the analytical results of these papers. Our contribution here is the precise construction of relevant Cauchy--Riemann operators, their determinant lines and the relation between these, and the ensuing construction of the canonical quantum complex over an arbitrary ground ring, without using relative Pin-structures, including the minimal condition for which such a construction is possible (assumption \tb{(O)}), and also a proof of the fact that the quantum boundary operator squares to zero in the case $N_L = 2$, see \S\ref{s:boundary_op_squares_zero_bubbling}.

\S\ref{ss:boundary_gluing} is concerned with gluing of Riemann surfaces and Cauchy--Riemann operators on them at a boundary or an interior point; the material presented here also appears in condensed form in \cite{Seidel_The_Book}. In \S\ref{ss:Lagr_QH} we define the quantum complex of a quantum datum, the corresponding boundary operator, and prove it squares to zero. We also define the quantum product and the corresponding unit. \S\ref{ss:arbitrary_rings_loc_coeffs_Lagr_QH}, \S\ref{ss:duality_Lagr_QH} are concerned with arbitrary coefficients and duality in quantum homology. \S\ref{ss:QH_of_M} defines the quantum homology of $M$ as well as its module action on the Lagrangian quantum homology. In \S\ref{ss:spectral_seqs} we describe the natural spectral sequence which starts with the (twisted) Morse complex and converges to the quantum homology.

\subsection{Gluing at a boundary or an interior point}\label{ss:boundary_gluing}

There is one technical aspect of Cauchy--Riemann operators not covered in \S\ref{s:HF}, namely gluing of Riemann surfaces, Cauchy--Riemann operators on them, and pregluing smooth maps, at a boundary or at an interior point rather than at a puncture. This subsection collects the necessary definitions and facts regarding these operations.

\paragraph{Riemann surfaces and smooth maps}

Given two Riemann surfaces $\Sigma_i$, $i=1,2$ and points $z_i \in \Sigma_i$, both either boundary or interior, we can form the glued surface $\Sigma_1 \sharp \Sigma_2$ by choosing collars around the points $z_i$, a gluing length, and identifying punctured neighborhoods of the $z_i$ according to the collars and the gluing length. If $u_i \fc (\Sigma_i,\partial \Sigma_i) \to (M,L)$ are smooth maps, $z_i \in \Sigma$ are both either boundary or interior, and $u_1(z_1) = u_2(z_2)$, we can use the expression of the $u_i$ near $z_i$ via exponential maps, similarly to what we described in \S\ref{ss:b_smooth_maps_pregluing}, to preglue $u_1$ and $u_2$ to a smooth map $u_1 \sharp u_2 \fc \Sigma_1 \sharp \Sigma_2 \to M$. Of course, if the $\Sigma_i$ have punctures and the $u_i$ are b-smooth, the resulting surface $\Sigma_1 \sharp \Sigma_2$ will inherit the punctures, and the preglued map $u_1\sharp u_2$ will be b-smooth as well.

\paragraph{Cauchy--Riemann operators}

Gluing Cauchy--Riemann operators at a boundary or interior point consists of two stages. At the first stage we assume given two Riemann surfaces $\Sigma_i$, $i=1,2$, two points $z_i \in \Sigma_i$, both boundary or both interior, and Hermitian bundle pairs $(E_i,F_i) \to (\Sigma_i,\partial \Sigma_i)$ endowed with connections $\nabla_i$. These give rise to the associated Cauchy--Riemann operators
$$D_i := \ol\partial_{\nabla_i} = \nabla_i^{0,1} \fc W^{1,p}(\Sigma_i, \partial \Sigma_i; E_i, F_i) \to L^p(\Sigma_i, \Omega_{\Sigma_i}^{0,1} \otimes E_i)\,.$$
These are Fredholm.\footnote{In case the $\Sigma_i$ have punctures, we need to assume in addition that the $D_i$ are admissible and nondegenerate, however we suppress these details for the sake of clarity.} Assume we are given a unitary isomorphism $(E_1)_{z_1} \simeq (E_2)_{z_2}$, which in case the $z_i$ are boundary also identifies $(F_1)_{z_1}$ with $(F_2)_{z_2}$. Assume that the $z_i$ are boundary. Using the identification, we can produce the following exact Fredholm triple:
$$\xymatrix{0 \ar[r] & W^{1,p}_1 \sharp W^{1,p}_2 \ar[r] \ar[d]^{D_1 \sharp D_2} & W^{1,p}_1 \oplus W^{1,p}_2 \ar[rr]^-{\ev_{z_1} - \ev_{z_2}} \ar[d]^{D_1 \oplus D_2} & & (F_1)_{z_1} \ar[r] \ar[d] & 0\\
0 \ar[r] & L^p_1 \oplus L^p_2 \ar@{=}[r] & L^p_1 \oplus L^p_2 \ar[rr] & & 0}$$
where
$$W^{1,p}_i = W^{1,p}(\Sigma_i, \partial \Sigma_i; E_i, F_i)\,,\quad L^p_i = L^p(\Sigma_i, \Omega_{\Sigma_i}^{0,1} \otimes E_i)\,,$$
$$W^{1,p}_1 \sharp W^{1,p}_2 = \{(\xi_1,\xi_2) \in W^{1,p}_1 \oplus W^{1,p}_2\,|\,\xi_1(z_1) = \xi_2(z_2)\}\,,\quad \text{and}\quad D_1 \sharp D_2 = (D_1 \oplus D_2)|_{W^{1,p}_1 \sharp W^{1,p}_2}\,.$$
This exact triple gives rise to the isomorphism
\begin{equation}\label{eqn:boudnary_gluing_iso}
\ddd(D_1 \sharp D_2) \otimes \ddd((F_1)_{z_1}) \to \ddd(D_1 \oplus D_2)\,,
\end{equation}
which we refer to as the \tb{boundary gluing isomorphism}. Note that it depends on the ordering of the operators. It can be checked that if we exchange the operators $D_1, D_2$, which amounts to replacing the map $\ev_{z_1} - \ev_{z_2}$ with $\ev_{z_2} - \ev_{z_1}$, then the above isomorphism is multiplied by $(-1)^{\dim (F_1)_{z_1}}$.

If the $z_i$ are interior, a similar argument yields the isomorphism
$$\ddd(D_1 \sharp D_2) \otimes \ddd((E_1)_{z_1}) \to \ddd(D_1 \oplus D_2)\,,$$
and if we use the canonical orientation of $E_1$, then we simply get
$$\ddd(D_1 \sharp D_2) \to \ddd(D_1 \oplus D_2)\,.$$
This isomorphism is independent of the ordering of $D_1,D_2$.

At the second stage we produce an operator on the glued surface $\Sigma_1 \sharp \Sigma_2$. First we need to glue the Hermitian bundle pairs by identifying them over the collars; this can be done using the given identification $(E_1)_{z_1} \simeq (E_2)_{z_2}$. Then we can deform the operator $D_1 \sharp D_2$, for example by deforming the connections, so that the resulting operators over the collars coincide relative to the identification of the bundles. This then defines an operator on the glued bundle pair. The resulting deformation of Cauchy--Riemann operators produces an isomorphism between $\ddd(D_1 \sharp D_2)$ and the determinant line of the glued operator. By slightly abusing notation, we denote the glued operator by the same symbol $D_1 \sharp D_2$. We refer to either one of these operators as being \tb{boundary glued} from $D_1,D_2$.

\paragraph{Smooth maps and their linearizations} Finally, assume we have two smooth maps $u_i \fc (\Sigma_i,\partial \Sigma_i) \to (M,L)$ with $u_1(z_1) = u_2(z_2)$, and consider the preglued map $u_1 \sharp u_2$. The linearizations $D_{u_1}$ and $D_{u_2}$ can be glued according to the previous paragraph, since we of course have an identification $(E_{u_1})_{z_1} = (E_{u_2})_{z_2} = T_{u_1(z_1)}M$ and similarly for $(F_i)_{z_i} = T_{u_1(z_1)}L$. The glued operator $D_{u_1} \sharp D_{u_2}$ can then be deformed into the linearization $D_{u_1\sharp u_2}$. In particular we have the canonical deformation isomorphism
$$\ddd(D_{u_1} \sharp D_{u_2}) \simeq \ddd(D_{u_1 \sharp u_2})\,.$$

\subsection{Lagrangian quantum homology}\label{ss:Lagr_QH}

Given a Morse function $f$ on $L$ we let $\Crit f$ be the set of its critical points. For $q \in \Crit f$ we denote its Morse index by $|q|_f$, or usually just by $|q|$. If $\rho$ is a Riemannian metric on $L$, we let $\cS_f(q),\cU_f(q)$ be the stable and unstable manifolds of $f$ with respect to $\rho$ at a critical point $q$. Usually we will drop the subscript $_f$ if the function is clear from the context. We call $(f,\rho)$ a Morse-Smale pair if every stable manifold of $f$ intersects every unstable manifold of $f$ transversely. In this case the set
$$\wt\cM(q,q') = \cU(q) \cap \cS(q') \subset L$$
is naturally a smooth manifold of dimension $|q| - |q'|$ for any $q,q' \in \Crit f$.

A \tb{quantum homology datum} for $L$ is a triple $\cD = (f,\rho,J)$, where $f \in C^\infty(L)$ is a Morse function, $\rho$ a Riemannian metric on $L$, such that $(f,\rho)$ is a Morse-Smale pair, and $J$ is an $\omega$-compatible almost complex structure. We call $\cD$ regular if $J$ is chosen generically, in the sense that the various pearly moduli spaces defined below are transversely cut out, in the sense of Biran--Cornea \cite{Biran_Cornea_Quantum_structures_Lagr_submfds, Biran_Cornea_Rigidity_uniruling}.

\subsubsection{Generators and the complex as a $\Z$-module}\label{sss:generators_cx_Lagr_QH}

Fix a quantum homology datum $\cD = (f,\rho, J)$ for $L$. The $\Z$-module underlying the quantum complex of $L$ is defined as
$$QC_*(\cD:L) = \bigoplus_{\substack{q \in \Crit f\\ A \in \pi_2(M,L,q)}}C(q,A)$$
where $C(q,A)$ is a certain rank $1$ free $\Z$-module associated to the pair $(q,A)$, which will be defined below. The grading is determined by requiring the elements of $C(q,A)$ to be homogeneous of degree $|q| - \mu(A)$.

In order to define the module $C(q,A)$ for $q \in \Crit f$ and $A \in \pi_2(M,L,q)$, we need a preliminary construction. First let us fix once and for all an arbitrary background connection $\nabla$ on $M$. Set
$$C^\infty_A = \{u \in C^\infty(D^2,S^1,1;M,L,q)\,|\, [u] = A\}\,.$$
For $u \in C^\infty_A$ we have the bundle pair $(E_u,F_u) \to (D^2,S^1)$ given by $E_u = u^*TM$, $F_u = (u|_{S^1})^*TL$. This carries the Hermitian structure $\omega_u = u^*\omega$, $J_u = u^*J$. We let $\nabla_u = u^*\nabla$ be the induced connection on $E_u$. We then have the associated Cauchy--Riemann operator
$$D_u = \nabla_u^{0,1} \fc W^{1,p}(D^2,S^1;E_u,F_u) \to L^p(D^2; \Omega_{D^2}^{0,1} \otimes E_u)$$
where we have extended it to the Sobolev completions. This operator is Fredholm with
$$\ind D_u = n + \mu(A)\,.$$
Thus we have the family of Fredholm operators
$$D_A = (D_u)_{u \in C^\infty_A}\,,$$
and the corresponding determinant line bundle \footnote{In contrast to the situation in \S\ref{ss:CROs_from_b_smooth_maps} where we had to include the auxiliary choices of an almost complex structure and a connection into the parameter space, here there is no need to do that.} $\det D_A$ over $C^\infty_A$. If we let $C_A$ be the space of continuous disks in class $A$, the inclusion $C^\infty_A \to C_A$ is a homotopy equivalence. The latter space can be identified with the space of continuous maps $([0,1]^2,\partial[0,1]^2,\{0,1\} \times [0,1] \cup [0,1] \times \{0\}) \to (M,L,q)$ in class $A$. Such continuous maps can be concatenated using the first coordinate of $[0,1]^2$ (this corresponds to the product in $\pi_2(M,L,q)$). For a map $u \fc [0,1]^2 \to M$ we let $-u \fc [0,1]^2 \to M$ be defined via $-u(s,t) = u(1-s,t)$. We have the following foundational lemma.
\begin{lemma} Fix $u_0 \in C^\infty_A$ and view it as a map of the square to $M$ as just described. The concatenation map $C_A \to C_0$, $u \mapsto u\sharp - u_0$, is a homotopy equivalence, where $C_0$ is the space of contractible disks at $q$. In particular the fundamental group of $C^\infty_A$ is isomorphic to $\pi_1(C_0,q)$, which in turn is canonically isomorphic to $\pi_3(M,L,q) = \pi_1(C_0,q)$, since $\pi_3(M,L,q)$ is abelian. Moreover, relative to this identification, the first Stiefel--Whitney class of the line bundle $\det D_A$ satisfies
$$w_1(\det D_A) = w_2(TL) \circ \partial \fc \pi_3(M,L,q) \to \Z_2$$
where $\partial \fc \pi_3(M,L,q) \to \pi_2(L,q)$ is the boundary map. \qed
\end{lemma}
\noindent The proof of the homotopy part of the statement is left to the reader as an exercise. The second part is proved very similarly to the proof of Lemma \ref{lem:first_Stiefel_Whitney_dD_wt_y_HF} concerning the case of Floer homology.

This lemma shows that $\det D_A$ is orientable if and only if $w_2(TL) \circ \partial$ vanishes, that is if assumption \tb{(O)} holds. We assume this vanishing from this point on.

There is an evaluation map \footnote{This is well-defined because $p > 2$ implies the elements of $W^{1,p}$ are continuous.}
$$\ev_1 \fc W^{1,p}(D^2,S^1;E_u,F_u) \to T_qL\,, \quad \xi \mapsto \xi(1)\,,$$
and we let $D_u \sharp T\cS(q)$ be the restriction of $D_u$ to $\ev_1^{-1}(T_q\cS(q))$. Sometimes we will employ the full notation $D_u \sharp T_q\cS_f(q)$. We have
\begin{lemma}\label{lem:family_ops_orientable_def_cx_QH}
The family of operators $D_A \sharp T\cS(q) := (D_u \sharp T\cS(q))_{u \in C^\infty_A}$ is orientable.
\end{lemma}
\begin{prf}
We have the natural exact Fredholm triple
$$0 \to D_u \sharp T\cS(q) \to D_u \oplus 0_{T_q\cS(q)} \xrightarrow{\ev_1 - \,\text{inclusion}} 0_{T_qL} \to 0\,,$$
Together with the direct sum isomorphism, it yields the isomorphisms
$$\ddd(D_u) \otimes \ddd(T\cS(q)) \simeq \ddd(D_u \oplus T\cS(q)) \simeq \ddd(D_u\sharp T\cS(q)) \otimes \ddd(T_qL)$$
which are continuous with respect to $u$. We see that we have obtained an isomorphism of line bundles
$$\ddd(D_A) \simeq \ddd(D_A \sharp T\cS(q))\,,$$
whence the bundle $\ddd(D_A \sharp T\cS(q))$ is orientable. \qed
\end{prf}

We can now define $C(q,A)$: it is the rank $1$ free $\Z$-module generated by the two possible orientations of the line bundle $\ddd(D_A \sharp T\cS(q))$. Note that this makes sense because $C^\infty_A$ is connected. The definition of the graded $\Z$-module $QC_*(\cD:L)$ is thereby completed.

\subsubsection{Boundary operator}\label{sss:boundary_op_Lagr_QH}

To define the boundary operator, we need first to define the spaces of pearls. We follow \cite{Biran_Cornea_Quantum_structures_Lagr_submfds, Biran_Cornea_Rigidity_uniruling}. Fix $q,q' \in \Crit f$. We will define the space of parametrized pearls $\wt \cP_k(q,q')$ for $k \geq 0$. We define\
$$\wt \cP_0(q,q')$$
to be the space of parametrized negative gradient lines of $f$ from $q$ to $q'$. It is naturally identified with $\wt\cM(q,q')$  via $u \mapsto u(0)$. This is a smooth manifold of dimension $|q| - |q'|$. It admits a natural $\R$-action and we let
$$\cP_0(q,q') = \wt \cP_0(q,q')/\R$$
be the quotient space if this action is free. When it is not free, which is the case if and only if $q = q'$, we let $\cP_0(q,q') = \varnothing$. Note that $\dim \cP_0(q,q') = |q| - |q'| - 1$.

Assume now $k > 0$. We have the evaluation map
$$\ev \fc (C^\infty(D^2,S^1;M,L))^k \to L^{2k}\,,\quad u = (u_1,\dots,u_k) \mapsto (u_1(-1),u_1(1),\dots,u_k(-1),u_k(1))\,.$$
We denote by $\phi_f^t \fc L \to L$ the time-$t$ flow map of the negative gradient $-\nabla_\rho f$. We let
$$Q_{f,\rho} = \{(x,\phi_f^t(x))\,|\, x \notin \Crit f\,, t > 0\} \subset L^2$$
be the flow manifold. We usually abbreviate this to $Q$. Note for further use that this has dimension $n+1$. We will also need the extended flow manifold
$$\ol Q = \{(x,\phi_f^t(x))\,|\, x \notin \Crit f\,, t\geq 0\}\,,$$
which is a manifold with boundary. This $\ol Q$ carries a natural hyperplane distribution, which we denote $\Gamma \subset T\ol Q$; this is just the collection of graphs of the differentials of the flow maps:
$$\Gamma_{(x,\phi_f^t(x))} = \{(X,\phi_{f*}^t(X))\,|\,X \in T_xL\} \subset T_{(x,\phi_f^t(x))}\ol Q\,.$$
For future use note that if $(x,y) \in Q$, then we have a natural basis of $T_{(x,y)}Q/\Gamma_{(x,y)}$, given by the coset of the vector $(-\nabla_\rho f(x),0)$:
\begin{equation}\label{eqn:basis_vector_for_TQ_mod_Gamma}
e_{x,y}:=(-\nabla_\rho f(x),0) + \Gamma_{(x,y)}\in T_{(x,y)}Q/\Gamma_{(x,y)}\,.
\end{equation}

Let $\wt\cM(L;J)$ be the space of parametrized nonconstant $J$-holomorphic disks with boundary on $L$. The \tb{pearly spaces} are then defined to be
$$\wt\cP_k(q,q') = \ev^{-1}(\cU(q) \times Q^{k-1} \times \cS(q')) \cap (\wt\cM(L;J))^k \subset (\wt \cM(L;J))^k\,,\quad \text{ and }\quad \wt\cP(q,q') = \bigcup_{k \geq 0} \wt\cP_k(q,q')\,.$$

We identify the group of conformal automorphisms of $(D^2,\partial D^2,\pm 1)$ with $\R$, as follows. The formula \eqref{eqn:std_end_cappings_Lagr_HF} defines a biholomorphism between $\R \times [0,1]$ and $D^2 - \{\pm 1\}$. The group $\R$ acts on $\R \times [0,1]$ by translations to the right. The biholomorphism thus induces an isomorphism $\R \to \Aut(D^2,S^1,\pm 1)$, the latter being the group of conformal automorphisms of $D^2$ preserving $\pm 1$. We make $\R$ act on the set of smooth maps $(D^2,S^1) \to (M,L)$ via
$$\tau \cdot u = u(\cdot + \tau,\cdot)$$
relative to the coordinates $(s,t)$ on $D^2 - \{\pm 1\} \simeq \R \times [0,1]$.

This gives rise to an action of $\R^k$ on $\wt\cP_k(q,q')$ by reparametrizations on each disk. We let $\cP_k(q,q')$ and $\cP(q,q')$ be the corresponding quotient spaces.

For $u \in (C^\infty (D^2,S^1;M,L))^k$ let
$$\mu(u):=\sum_i \mu(u_i)\,.$$
Biran--Cornea proved \cite{Biran_Cornea_Quantum_structures_Lagr_submfds, Biran_Cornea_Rigidity_uniruling}:
\begin{prop}
There is a subset of $\cJ(M,\omega)$ of the second category such that for every $J$ in the subset, for each $k \geq 0$ and for each pair $q,q' \in \Crit f$ the space $\cP_k(q,q')$ is a smooth manifold of local dimension at $[u] \in \cP_k(q,q')$
$$\dim_{[u]}\cP_k(q,q') = |q| - |q'| + \mu(u) - 1$$
whenever this number does not exceed $1$. \qed
\end{prop}
For the rest of \S\ref{sss:boundary_op_Lagr_QH} we assume that $J$ is chosen so that the proposition holds. We define for future use the elements
\begin{equation}\label{eqn:std_basis_vector_pearly_space_QH}
\partial_{u_i} \in T_u\wt\cP_k(q,q')
\end{equation}
given by the infinitesimal action of $\R$ on $u_i$.

The boundary operator
$$\partial_\cD \fc QC_*(\cD:L) \to QC_{*-1}(\cD:L)$$
will be defined in terms of its matrix elements, which are homomorphisms
$$C(q,A) \to C(q',A')$$
for $|q| - \mu(A) = |q'| - \mu(A') + 1$. We now proceed to define these matrix elements.

Fix $q,q' \in \Crit f$ and $A \in \pi_2(M,L,q)$. For any $u \in \ev^{-1}(\cU(q) \times \ol Q{}^{k-1} \times \cS(q'))$ there is a natural way of constructing a homotopy class $A \sharp u \in \pi_2(M,L,q')$. Indeed, each disk $u_j$ can be viewed as a continuous map $([0,1]^2, \partial[0,1]^2) \to (M,L)$ which maps $\{0\} \times [0,1]$ to $u_j(-1)$ and $\{1\} \times [0,1]$ to $u_j(1)$. Every piece of gradient trajectory connecting either a critical point $q$ or $q'$ with $u_1(-1)$ or with $u_k(1)$, or $u_j(1)$ with $u_{j+1}(-1)$, can be viewed as a continuous map $([0,1]^2, \partial[0,1]^2) \to (M,L)$ which is independent of the second variable. Now take all these maps defined on $[0,1]^2$ and concatenate them using the first coordinate, with each other, and with a representative of $A$, also viewed as such a continuous map, where the concatenation order is dictated by the linear structure of the string $u$. This works in case $k=0$ as well.

For $u \in \wt\cP_k(q,q')$ which satisfies $|q| - |q'| + \mu(u) - 1 = 0$, we will construct an isomorphism
$$C(u) \fc C(q,A) \to C(q',A\sharp u)\,.$$
Fix now $A' \in \pi_2(M,L,q')$ such that $|q| - \mu(A) = |q'| - \mu(A') + 1$. Then the corresponding matrix element of the boundary operator is the sum
$$\sum_{\substack{[u] \in \cP(q,q'): \\ A\sharp u = A'}}C(u) \fc C(q,A) \to C(q',A')\,.$$
A compactness argument shows that this sum is finite.

It therefore remains to define the isomorphism $C(u)$. We describe this is detail. The same basic construction will be used again below to define the algebraic structures.

Before we describe the isomorphism $C(u)$ for $u \in \wt\cP_0(q,q')$, we establish a general correspondence between orientations of the connected component of $\wt\cM(q,q')$ containing a gradient trajectory $u$ and isomorphisms $\ddd(D_A\sharp T\cS(q)) \simeq \ddd(D_{A'} \sharp T\cS(q'))$, where $q,q' \in \Crit f$ are two arbitrary critical points, and $A'$ is obtained from $A$ by transferring it along $u$.

First note that we have the following natural exact sequence
\begin{equation}\label{eqn:exact_seq_moduli_space_gradient_lines_stable_unstable}
0 \to T_{u(0)}\wt\cM(q,q') \to T_{u(0)} \cS(q') \to T_{u(0)}L / T_{u(0)} \cU(q) \to 0\,,
\end{equation}
where the penultimate arrow is the composition of the inclusion followed by the quotient map. The exactness follows from the fact that its kernel is precisely
$$T_{u(0)}\wt\cM(q,q') = T_{u(0)} \cS(q') \cap T_{u(0)} \cU(q)\,,$$
because $\wt\cM(q,q') = \cS(q') \cap \cU(q)$. Note that the bundle $TL/T\cU(q)$ is trivial over $\cU(q)$ (the latter being contractible), and we have naturally $T_qL / T_q \cU(q) = T_q\cS(q)$. Combining this fact with the natural isomorphism of determinant lines \eqref{eqn:iso_det_lines_exact_seq_vector_spaces} we obtain the isomorphism
\begin{equation}\label{eqn:iso_oris_moduli_space_gradient_lines_stable_mfds}
\ddd(T\cS(q')) \simeq \ddd(T_{u(0)}\wt\cM(q,q')) \otimes \ddd(T\cS(q))
\end{equation}
This means that there is a natural bijection between orientations of the connected component of $\wt\cM(q,q')$ containing $u$ and isomorphisms $\ddd(T\cS(q)) \simeq \ddd(T\cS(q'))$.

Using exact triples as in the proof of Lemma \ref{lem:family_ops_orientable_def_cx_QH}, we obtain canonical isomorphisms
$$\ddd(D_A) \otimes \ddd(T\cS(q)) \simeq \ddd(D_A \sharp T\cS(q)) \otimes \ddd(T_qL)\,,$$
$$\ddd(D_{A'}) \otimes \ddd(T\cS(q')) \simeq \ddd(D_{A'} \sharp T\cS(q')) \otimes \ddd(T_{q'}L)\,.$$
Using the natural isomorphisms $\ddd(D_A) \simeq \ddd(D_{A'})$ and $\ddd(T_qL) \simeq \ddd(T_{q'}L)$ obtained by transferring orientations along $u$, we see that there is a natural bijection between the set of isomorphisms $\ddd(D_A \sharp T\cS(q)) \simeq \ddd(D_{A'} \sharp T\cS(q'))$ and the set of isomorphisms $\ddd(T\cS(q)) \simeq \ddd(T\cS(q'))$, which in turn is in bijection with orientations of the connected component of $\wt\cM(q,q')$ containing $u$. This is the general correspondence we wanted to construct.

Now if $u \in \wt\cP_0(q,q')$ and $|q| = |q'| + 1$, $T_{u(0)}\wt\cM(q,q') = T_u\wt\cP_0(q,q')$ is canonically oriented by the translation vector $\dot u(0) = \partial_u$. We have the corresponding isomorphism $\ddd(D_A \sharp T\cS(q)) \simeq \ddd(D_{A'} \sharp T\cS(q'))$, or equivalently the isomorphism $C(q,A) \simeq C(q',A')$. This isomorphism is $C(u)$.

We now pass to the definition of $C(u)$ for $u \in \wt \cP_k(q,q')$ with $k \geq 1$. We need to define some additional Banach spaces. For a smooth map $v \fc (D^2,S^1) \to (M,L)$ we use the abbreviations
$$W^{1,p}(v) = W^{1,p}(D^2,S^1;v^*TM,(v|_{S^1})^*TL)\quad \text{ and } \quad L^p(v) = L^p(D^2;\Omega^{0,1}_{D^2}\otimes v^*TM)\,.$$
For a string of maps $v = (v_1,\dots,v_k) \in (C^\infty(D^2,S^1;M,L))^k$ we let
$$W^{1,p}(v) = \bigoplus_i W^{1,p}(v_i) \quad \text{ and } \quad L^p(v) = \bigoplus_i L^p(v_i)\,.$$
We define for $v \in (C^\infty(D^2,S^1;M,L))^k$ with $\ev(v) \in \cU(q) \times \ol Q{}^{k-1} \times \cS(q')$ the spaces
$$Z_\Gamma^v = \{\xi = (\xi_1,\dots,\xi_k) \in W^{1,p}(v)\,|\, (\xi_i(1),\xi_{i+1}(-1)) \in \Gamma_{(v_i(1),v_{i+1}(-1))}\text{ for }i < k\}\,,$$
$$X_\Gamma^v = \{\xi \in Z_\Gamma^v \,|\, \xi_k(1) \in T_{v_k(1)}\cS(q')\}\,,$$
$$Y_\Gamma^v = \{\xi \in X_\Gamma^v\,|\, \xi_1(-1) \in T_{v_1(-1)}\cU(q)\}\,,$$
and if $\ev (v) \in \cU(q) \times Q^{k-1} \times \cS(q')$, we define in addition
\begin{multline*}
Y_Q^v = \{\xi \in W^{1,p}(v)\,|\, \xi_1(-1) \in T_{v_1(-1)}\cU(q)\,,(\xi_i(1),\xi_{i+1}(-1)) \in T_{(v_i(1),v_{i+1}(-1))}Q\text{ for }i < k\,,\\\xi_k(1) \in T_{v_k(1)}\cS(q')\}\,.
\end{multline*}
For such $v$ we will also employ the notation
$$T_vQ:= \bigoplus_{i < k} T_{(v_i(1),v_{i+1}(-1))}Q$$
and
\begin{equation}\label{eqn:definition_of_TQ_mod_Gamma}
T_vQ/\Gamma:= \bigoplus_{i < k} T_{(v_i(1),v_{i+1}(-1))}Q/\Gamma_{(v_i(1),v_{i+1}(-1))}\,.
\end{equation}
The latter space has a natural basis
\begin{equation}\label{eqn:canonical_basis_of_TQ_mod_Gamma}
e^v_i:=e_{v_i(1),v_{i+1}(-1)} \in T_{(v_i(1),v_{i+1}(-1))}Q/\Gamma_{(v_i(1),v_{i+1}(-1))}
\end{equation}
for $i=1,\dots,k-1$, where we use the notation \eqref{eqn:basis_vector_for_TQ_mod_Gamma}. The natural evaluation map $\ev \fc Y^v_Q \to T_vQ/\Gamma$ induces a short exact sequence of Banach spaces:
$$0 \to Y^v_\Gamma \to Y^v_Q \xrightarrow{\ev} T_vQ \to 0\,.$$

The isomorphism $C(u)$ is constructed in two stages. First, we establish a canonical bijection between the orientations of $\ddd(D_u|_{Y^u_\Gamma})$ and orientations of $\ddd(T_u \wt\cP_k(q,q')) \otimes \ddd(T_uQ/\Gamma)$. Then we establish a canonical bijection between the orientations of $\ddd(D_u|_{Y^u_\Gamma})$ and isomorphisms $\ddd(D_A\sharp T\cS(q)) \simeq \ddd(D_{A'} \sharp T\cS(q'))$ where $A' = A \sharp u$, or equivalently isomorphisms $C(q,A) \simeq C(q',A')$. Selecting the orientation
$$\textstyle (-1)^{k+1}\bigwedge_i\partial_{u_i} \otimes \bigwedge_ie^u_i \in \ddd(T_u \wt\cP_k(q,q')) \otimes \ddd(T_uQ/\Gamma)\,,$$
we get, using these bijections, the desired isomorphism $C(u) \fc C(q,A) \simeq C(q',A')$.

To establish the first bijection consider the exact Fredholm triple
$$\xymatrix{0 \ar[r] & Y_\Gamma^u \ar[r] \ar[d]^{D_u|_{Y^u_\Gamma}} & Y_Q^u \ar[r]^{\ev} \ar[d]^{D_u|_{Y^u_Q}} & T_uQ/\Gamma\ar[r] \ar[d] & 0 \\ 0 \ar[r] & L^p(u) \ar@{=}[r] & L^p(u) \ar[r] & 0}$$
where $D_u$ is the linearization of $u$. Note that $T_u\wt\cP_k(q,q') = \ker D_u|_{Y_Q^u}$, and that $D_u|_{Y_Q^u}$ is surjective. Then the exact triple yields the isomorphism
\begin{equation}\label{eqn:correspondence_ors_pearly_space_TQ_mod_Gamma_and_D_u}
\ddd(T_u\wt\cP_k(q,q')) = \ddd(D_u|_{Y_Q^u}) \simeq \ddd(D_u|_{Y_\Gamma^u}) \otimes \ddd(T_uQ/\Gamma)
\end{equation}
and by tensoring with $\ddd(T_uQ/\Gamma)$ we see that there is indeed a canonical bijection between orientations of $\ddd(D_u|_{Y^u_\Gamma})$ and orientations of $\ddd(T_u\wt\cP_k(q,q')) \otimes \ddd(T_uQ/\Gamma)$.

We now construct the second bijection. We need some preliminary constructions and lemmata. Let $X_i,Y_i$, $i=1,2$, be Banach spaces, $D_i \in \cF(X_i,Y_i)$, and let $V$ be a finite-dimensional vector space. Assume $\theta_i \fc X_i \to V$ are surjective linear continuous maps. Let $W \subset V \oplus V$ be a subspace of dimension $\dim W = \dim V$. Then we can define the space $X_1 \sharp_W X_2 = (\theta_1 \oplus \theta_2)^{-1}(W) \subset X_1 \oplus X_2$, and the operator
$$D_1 \sharp_W D_2 = (D_1 \oplus D_2)|_{X_1 \sharp_W X_2}\,.$$
We have the following exact triple
$$\xymatrix{0 \ar[r] & X_1\sharp_W X_2 \ar[r] \ar[d]^{D_1\sharp_W D_2} & X_1 \oplus X_2 \oplus W \ar[rr]^-{\theta_1 \oplus \theta_2 - \, \text{inclusion}} \ar[d]^{D_1 \oplus D_2 \oplus 0_W} & & V \oplus V \ar[r] \ar[d] & 0\\
0 \ar[r] & Y_1 \oplus Y_2 \ar[r] & Y_1 \oplus Y_2 \ar[rr] & &  0}$$
which together with the direct sum isomorphism yields the isomorphism
$$\ddd(D_1 \oplus D_2) \otimes \ddd(W) \simeq \ddd(D_1 \sharp_W D_2) \otimes \ddd(V \oplus V)\,.$$
This shows that $(\ddd(W))_W$ and $(\ddd(D_1 \sharp_W D_2))_W$ are isomorphic as line bundles over the Grassmannian $G_{\dim V}(V \oplus V)$. In particular a path in this Grassmannian between two subspaces $W,W'$ induces an isomorphism $\ddd(D_1 \sharp_W D_2) \simeq \ddd(D_1 \sharp_{W'} D_2)$.

Assume now that $V$ splits as the direct sum $W_1 \oplus W_2$ and let $W = (W_1\oplus 0) \oplus (0\oplus W_2) \subset V \oplus V$. There is a path of subspaces between $W$ and the diagonal $\Delta_V$, given by the images of the embeddings
$$V = W_1 \oplus W_2 \to V \oplus V\,,\quad (w_1,w_2) \mapsto (w_1 + tw_2, tw_1 + w_2)$$
for $t \in [0,1]$. We have the induced isomorphism $\ddd(W) \simeq \ddd(\Delta_V)$ and therefore the isomorphism
\begin{equation}\label{eqn:iso_abstract_deformation_incidence_condition_W}
\ddd(D_1 \sharp_W D_2) \simeq \ddd(D_1 \sharp _{\Delta_V} D_2)\,,
\end{equation}
which will be frequently used in the sequel. We say that this isomorphism is obtained by \tb{deforming the incidence condition}.

Now we turn to a special case. Assume $w \in C^\infty_A$ with $A \in \pi_2(M,L,q)$, in particular $w(1) = q$. Assume further that $u' \in (C^\infty(D^2,S^1;M,L))^k$, $k \geq 1$, satisfies $\ev(u') \in \{q\} \times \Delta_L^{k-1} \times \{q'\}$. We have the operators
$$D_w \fc W^{1,p}(w) \to L^p(w)\quad \text{and} \quad D_{u'}|_{X_\Gamma^{u'}} \fc X_\Gamma^{u'} \to L^p(u')\,,$$
and the surjective continuous homomorphisms
$$\theta_1 \fc W^{1,p}(w) \to T_qL\,,\; \xi \mapsto \xi(1)\,;\quad  \theta_2 \fc X_\Gamma^{u'} \to T_qL\,,\; (\xi_1,\dots,\xi_k) \mapsto \xi_1(-1)\,.$$
Note that
$$D_w|_{\theta_1^{-1}(T_q\cS(q))} = D_w \sharp T\cS(q)\,,$$
that
$$\theta_2^{-1}(T_q\cU(q)) = Y_\Gamma^{u'}\,,$$
and that
$$(D_{u'}|_{X_\Gamma^{u'}})|_{\theta_2^{-1}(T_q\cU(q))} = D_{u'}|_{Y_\Gamma^{u'}}\,.$$
Therefore if we let $W = (T_q\cS(q) \oplus 0) \oplus (0 \oplus T_q\cU(q)) \subset T_qL \oplus T_qL$, then
$$D_w \sharp_W D_{u'}|_{X_\Gamma^{u'}} = D_w\sharp T\cS(q) \oplus D_{u'}|_{Y_\Gamma^{u'}}\,.$$
On the other hand, the operator $D_w \sharp_{\Delta_{T_qL}} D_{u'}|_{X_\Gamma^{u'}}$ is precisely what we denoted by $D_w \sharp D_{u'}|_{X_\Gamma^{u'}}$ in \ref{ss:boundary_gluing}, and which is the precursor to boundary gluing of the operators $D_w, D_{u'}|_{X_\Gamma^{u'}}$. Combining these considerations with the above isomorphism \eqref{eqn:iso_abstract_deformation_incidence_condition_W}, and with the direct sum isomorphism, we obtain the isomorphism
\begin{equation}\label{eqn:iso_correspondence_ors_D_u_Y_Gamma_isos}
\ddd(D_{u'}|_{Y_\Gamma^{u'}}) \otimes \ddd(D_w \sharp T\cS(q)) \simeq \ddd(D_w \sharp D_{u'}|_{X_\Gamma^{u'}})\,.
\end{equation}
We have the following lemma.
\begin{lemma}\label{lem:bijection_isos_C_q_A_ors_D_u_prime}
The isomorphism \eqref{eqn:iso_correspondence_ors_D_u_Y_Gamma_isos} induces a bijection between orientations of $\ddd(D_{u'}|_{Y_\Gamma^{u'}})$ and isomorphisms $C(q,A) \to C(q',A')$ where $A' = A \sharp u'$. It is continuous in $u'$.
\end{lemma}
\begin{prf}
The continuity follows from the corresponding property of direct sum and deformation isomorphisms.

We know that concatenating $w$ with the constituent disks of $u'$ gives us a representative of $A'$. Moreover, the operator $D_w \sharp D_{u'}|_{X_\Gamma^{u'}}$, after performing boundary gluing and deformation as described in \S\ref{ss:boundary_gluing}, yields a representative of the family $D_{A'}\sharp T\cS(q')$. This can be seen as follows. The operator $D_w \sharp D_{u'}|_{Z^{u'}_\Gamma}$ satisfies the necessary incidence conditions and therefore can be glued to an operator in the family $D_{A'}$ since the surface on which it is defined is precisely the concatenation of the disks $w,u'_1,\dots,u'_k$, that is $w\sharp u_1 \sharp \dots \sharp u_k$, which is a disk in class $A'$. The operator $D_w \sharp D_{u'}|_{X_\Gamma^{u'}}$ is the restriction of $D_w \sharp D_{u'}|_{Z_\Gamma^{u'}}$ to the subspace where $\xi_k \in W^{1,p}(u_k)$ satisfies $\xi_k(1) \in T_{q'}\cS(q')$. Therefore the glued operator obtained from $D_w \sharp D_{u'}|_{Z_\Gamma^{u'}}$, when restricted to the same subspace, yields an operator in the family $D_{A'}\sharp T\cS(q')$. This restricted operator is also be obtained by gluing the operator $D_w \sharp D_{u'}|_{X_\Gamma^{u'}}$.

Passing to the families, we see that the isomorphism \eqref{eqn:iso_correspondence_ors_D_u_Y_Gamma_isos} yields an isomorphism of line bundles
$$\ddd(D_{u'}|_{Y_\Gamma^{u'}}) \otimes \ddd(D_A \sharp T\cS(q)) \simeq \ddd(D_{A'} \sharp T\cS(q'))\,.$$
This shows that there is a canonical bijection between orientations of $\ddd(D_{u'}|_{Y_\Gamma^{u'}})$ and isomorphisms $\ddd(D_A \sharp T\cS(q)) \simeq \ddd(D_{A'} \sharp T\cS(q'))$, or equivalently, isomorphisms $C(q,A) \simeq C(q',A')$, as claimed. \qed
\end{prf}

In order to obtain a bijection between orientations of $D_u|_{Y_\Gamma^u}$ and isomorphisms $C(q,A) \simeq C(q',A')$, we deform $u$ into $u' \in (C^\infty(D^2,S^1;M,L))^k$ with $\ev(u') \in \{q\} \times \Delta_L^{k-1} \times \{q'\}$. The isomorphism induced by the deformation
\begin{equation}\label{eqn:deformation_iso_from_u_to_u_prime}
\ddd(D_u|_{Y^u_\Gamma}) \simeq \ddd(D_{u'}|_{Y^{u'}_\Gamma})
\end{equation}
and the bijection between the orientations of $\ddd(D_{u'}|_{Y^{u'}_\Gamma})$ and isomorphisms $C(q,A) \simeq C(q',A')$ from Lemma \ref{lem:bijection_isos_C_q_A_ors_D_u_prime} yields the desired bijection. Therefore we have to find a way to deform $u$ into $u'$ with the desired properties. This can be done as follows. Consider the disk $u_1$ and deform it slightly so that a neighborhood of $-1$ maps to $u_1(-1)$. Then deform the constant map on this neighborhood to a map covering the piece of gradient trajectory going from $q$ to $u_1(-1)$. Perform a similar deformation of every $u_j$, $j > 1$ so that the resulting deformation maps a neighborhood of $-1$ to the piece of gradient trajectory from $u_{j-1}(1)$ to $u_j(-1)$. Finally, deform $u_k$ additionally so that a neighborhood of $1$ covers the piece of gradient trajectory from $u_k(1)$ to $q'$. All this can be done so that the deformation and the maps involved are smooth.

We have therefore completed the definition of the boundary operator.

\begin{thm}\label{thm:boundary_op_squares_zero_QH}
$\partial_\cD^2 = 0$.
\end{thm}

\noindent The proof of this theorem is quite involved and will occupy the rest of \S\ref{sss:boundary_op_Lagr_QH}.

\begin{prf}
Fix two critical points $q,q'' \in \Crit f$ and classes $A \in \pi_2(M,L,q), A'' \in \pi_2(M,L,q'')$ such that $|q| - \mu(A) = |q''| - \mu(A'') + 2$. We need to prove the vanishing of the corresponding matrix element of $\partial_\cD^2$, which equals
\begin{equation}\label{eqn:vanishing_matrix_elts_boundary_op_QH_squared}
\sum_{q' \in \Crit f} \sum_{\substack{A' \in \pi_2(M,L,q'): \\ |q| - \mu(A) = |q'| - \mu(A') + 1}} \sum_{\substack{([u],[v]) \in \cP(q,q') \times \cP(q',q''):\\A\sharp u = A', A'\sharp v = A''}} C(v) \circ C(u) \fc C(q,A) \to C(q'',A'')\,.
\end{equation}

Let us denote by $\cP^1(q,q'')$ the $1$-dimensional part of $\cP(q,q'')$. Biran--Cornea \cite{Biran_Cornea_Quantum_structures_Lagr_submfds, Biran_Cornea_Rigidity_uniruling} gave a description of the compactification $\ol\cP{}^1(q,q'')$ for $\cP^1(q,q'')$. The compactness of $\cP^1(q,q'')$ fails in one of the following ways: one of the gradient trajectories undergoes Morse breaking; two holomorphic disks collide, that is one of the gradient lines connecting the disks shrinks to a point; a holomorphic disk breaks into two; or a Maslov $2$ holomorphic disk bubbles off at a point. The case involving bubbling, which only happens when $N_L = 2$, will be handled in \S\ref{ss:boundary_op_squares_zero_bubbling_QH}; here we only consider the case in which bubbling does not occur.

Let us call $u \in \wt\cM(L,J) ^ m$, $m > 0$ a \tb{degenerate pearly trajectory} between $q,q''$ if there is $j < m$ such that $\ev(u) \in \cU(q) \times Q^{j-1} \times \Delta_L \times Q^{m-j-1} \times \cS(q'')$, that is it is an ordinary pearly trajectory, except that two of the holomorphic disks touch.

A degenerate pearly trajectory appears as a boundary point in the compactification of exactly two components of $\cP^1(q,q'')$: one in which the two disks are separated by a positive-time gradient trajectory, and the other one where one of the holomorphic disks breaks into the two touching disks in $u$. We redefine $\ol\cP{}^1(q,q'')$ to be the disjoint union of all the compactified components of $\cP^1(q,q'')$, and where we identify two boundary points if they correspond to the same degenerate pearly trajectory. This endows $\ol\cP{}^1(q,q'')$ with the structure of a compact $1$-dimensional topological manifold with boundary whose points are pairs of pearly trajectories $([u],[v]) \in \cP(q,q') \times \cP(q',q'')$ with $\dim_{[u]}\cP(q,q') = \dim_{[v]}\cP(q',q'') = 0$, therefore we see that the summands of \eqref{eqn:vanishing_matrix_elts_boundary_op_QH_squared} are in bijection with the boundary $\partial\ol\cP{}^1(q,q'')$.

For $\delta = ([u],[v]) \in \partial \ol\cP{}^1(q,q'')$ let $C(\delta) = C(v) \circ C(u)$ be the summand corresponding to $\delta$. It is enough to show that for every connected component $\Delta \subset \ol\cP{}^1(q,q'')$ with boundary $\partial\Delta = \{\delta,\delta'\}$ we have $C(\delta) + C(\delta') = 0$. The component $\Delta$ either has holomorphic disks or it doesn't. If it has no holomorphic disks, by the discussion on page \pageref{eqn:exact_seq_moduli_space_gradient_lines_stable_unstable} we know that for $w \in \wt\cP_0(q,q'')$ with $[w] \in \Delta$ there is a bijection between isomorphisms $C(q,A) \simeq C(q'',A'')$ and orientations of $T_w\wt\cP(q,q'')$.

If $\Delta$ does contain holomorphic disks, we have the following situation. Let $\{w^t\}_{t \in [0,1]}$ be a continuous parametrization of $\Delta$ with $w^0 = \delta$, $w^1 = \delta'$. There are finitely many instances
$$t_0 = 0 < t_1 < \dots < t_m = 1$$
such that $w^t \notin \cP^1(q,q'')$ if and only if $t = t_i$ for some $i$, that is $w^{t_i}$ is a denerate trajectory. By abuse of notation, for $t \neq t_i$ let $w^t \in \wt\cP(q,q'')$ be a continuous family of representatives of $w^t \in \cP^1(q,q'')$. Lemma \ref{lem:bijection_isos_C_q_A_ors_D_u_prime} together with the isomorphism \eqref{eqn:correspondence_ors_pearly_space_TQ_mod_Gamma_and_D_u} show that there is a bijection between isomorphisms $C(q,A) \simeq C(q'',A'')$, orientations of $D_{w^t}|_{Y_\Gamma^{w^t}}$, and orientations of $\ddd(T_{w^t}\wt\cP(q,q'')) \otimes \ddd(T_{w^t}Q/\Gamma)$, which is continuous in $t$ within each interval $(t_i,t_{i+1})$. The theorem then follows from the next lemma.
\begin{lemma}\label{lem:induced_or_pearly_spaces_boundary_op_squared}
Assume $\delta = ([u],[v]) \in \cP_k(q,q') \times \cP_l(q',q'')$. Then under the above bijection, the isomorphism $C(\delta)$ corresponds to:
\begin{itemize}
\item If $k=l=0$, the orientation $- \partial_w \wedge \text{\rm inward}_\delta \in \ddd(T_w\wt\cP(q,q''))$;
\item Otherwise, the orientation
\begin{equation}\label{eqn:orientation_induced_on_pearly_space_TQ_mod_Gamma_by_C_delta}
\textstyle - \bigwedge_i \partial_{w^t_i} \wedge \text{\rm inward}_\delta \otimes \bigwedge_i e^{w^t}_i \in \ddd(T_{w^t}\wt\cP(q,q'')) \otimes \ddd(T_{w^t}Q/\Gamma)
\end{equation}
whenever $t$ is not one of the degenerate values $t_j$. Here $\text{\rm inward}_\delta \in T_{w^t}\wt\cP(q,q'')$ denotes a vector which is transverse to the infinitesimal action of the automorphism group, and which points away from the boundary point $\delta$.
\end{itemize}
\end{lemma}
\noindent Indeed, assume the lemma. Let $\delta' = ([u'],[v']) \in \cP_{k'}(q,q_1) \times \cP_{l'}(q_1,q'')$. According to the lemma, if $k = l = 0$, then, of course $k' = l' = 0$, and the isomorphism $C(\delta')$ corresponds to the orientation
$$-\partial_w \wedge\text{inward}_{\delta'} \in T_{w}\wt \cP(q,q'')\,.$$
Since $\text{inward}_{\delta} = -\text{inward}_{\delta'}$, we see that in this case $C(\delta) + C(\delta') = 0$.

Otherwise the lemma says that $C(\delta')$ corresponds to
$$\textstyle - \bigwedge_i \partial_{w^t_i} \wedge \text{\rm inward}_{\delta'} \otimes \bigwedge_i e^{w^t}_i \in \ddd(T_{w^t}\wt\cP(q,q'')) \otimes \ddd(T_{w^t}Q/\Gamma)$$
for $t \in (t_{m-1},t_m)$. On the other hand, $C(\delta)$ corresponds to the orientation 
$$\textstyle - \bigwedge_i \partial_{w^t_i} \wedge \text{\rm inward}_\delta \otimes \bigwedge_i e^{w^t}_i \in \ddd(T_{w^t}\wt\cP(q,q'')) \otimes \ddd(T_{w^t}Q/\Gamma)$$
for the same range of $t$. Noting that $\text{inward}_\delta = -\text{inward}_{\delta'}$, we see that these two orientations are opposite, which proves that $C(\delta) + C(\delta') = 0$ and therefore that $\partial_\cD^2 = 0$. \qed
\end{prf}
\noindent This proves Theorem \ref{thm:boundary_op_squares_zero_QH}, modulo Lemma \ref{lem:induced_or_pearly_spaces_boundary_op_squared}.

\begin{prf}[of Lemma \ref{lem:induced_or_pearly_spaces_boundary_op_squared}] In case $k = l = 0$, this is a standard computation in Morse theory, see for instance Schwarz's book \cite{Schwarz_Morse_H_book}. In order to put that computation in the context of the present paper, we note that we have the following commutative diagram:
$$\xymatrix{\ddd(T_v\wt\cP(q',q'')) \otimes \ddd(T_u\wt\cP(q,q')) \otimes \ddd(T\cS(q)) \ar[r] \ar[d] & \ddd(T_v\wt\cP(q',q'')) \otimes \ddd(T\cS(q'))\ar[d] \\ \ddd(T_w\wt\cP(q,q'')) \otimes \ddd(T\cS(q)) \ar[r] & \ddd(T\cS(q''))} $$
where all the arrows except the left one come from the isomorphism \eqref{eqn:iso_oris_moduli_space_gradient_lines_stable_mfds}, whereas the left arrow comes from the direct sum isomorphism and the differential of the gluing map, which is an isomorphism
$$T_u\wt\cP(q,q') \times T_v\wt\cP(q',q'') \simeq T_w\wt\cP(q,q'')\quad \text{mapping}\quad \partial_v + \partial_u \mapsto \partial_w \text{ and }\partial_u - \partial_v \mapsto -\inward_\delta\,.$$
We can now tensor this diagram with the identity on $\ddd(D_A)$, then replace $\ddd(D_A)$ with $\ddd(D_{A'})$ and $\ddd(D_{A''})$ where needed, and replace the identity isomorphisms with the corresponding deformation isomorphisms. The exact triple isomorphisms as in the proof of Lemma \ref{lem:family_ops_orientable_def_cx_QH} then will yield the following diagram:
$$\resizebox{\textwidth}{!}{\xymatrix{\ddd(T_v\wt\cP(q',q'')) \otimes \ddd(T_u\wt\cP(q,q')) \otimes \ddd(D_A\sharp T\cS(q)) \otimes \ddd(T_qL) \ar[r] \ar[d] & \ddd(T_v\wt\cP(q',q'')) \otimes \ddd(D_{A'} \sharp T\cS(q')) \otimes \ddd(T_{q'}L) \ar[d] \\ \ddd(T_w\wt\cP(q,q'')) \otimes \ddd(D_A \sharp T\cS(q)) \otimes \ddd(T_qL) \ar[r] & \ddd(D_{A''} \sharp T\cS(q'')) \otimes \ddd(T_{q''}L)}} $$
Here we can now erase the terms $\ddd(T_\cdot L)$ and obtain the diagram
$$\xymatrix{\ddd(T_v\wt\cP(q',q'')) \otimes \ddd(T_u\wt\cP(q,q')) \otimes \ddd(D_A\sharp T\cS(q)) \ar[r] \ar[d] & \ddd(T_v\wt\cP(q',q'')) \otimes \ddd(D_{A'} \sharp T\cS(q')) \ar[d] \\ \ddd(T_w\wt\cP(q,q'')) \otimes \ddd(D_A \sharp T\cS(q)) \ar[r] & \ddd(D_{A''} \sharp T\cS(q''))} $$
which tells us that the isomorphism $C(\delta) = C(v) \circ C(u)$, which is obtained by chasing the diagram along through the top right corner, corresponds to the orientation of $T_w\wt\cP(q,q'')$ obtained from the differential of the gluing map, that is the orientation $-\partial_w \wedge \inward_\delta$.

Assume now that at least one of the numbers $k,l$ is nonzero. The claim then will be proved by induction on $j$. The base of the induction consists of showing that the isomorphism $C(\delta)$ corresponds to the orientation
$$\textstyle - \bigwedge_i \partial_{w^t_i} \wedge \text{\rm inward}_\delta \otimes \bigwedge_i e^{w^t}_i \in \ddd(T_{w^t}\wt\cP(q,q'')) \otimes \ddd(T_{w^t}Q/\Gamma)$$
for $t \in (0,t_1)$. It suffices to show this for some $t$ small enough due to continuity. Let therefore $w = w^t$ for some small $t > 0$. There is a canonical bijection between isomorphisms $C(q,A) \simeq C(q'',A'')$, orientations of $\ddd(D_w|_{Y_\Gamma^w})$, and orientations of $T_w\wt\cP(q,q'') \otimes \ddd(T_wQ/\Gamma)$, see Lemma \ref{lem:bijection_isos_C_q_A_ors_D_u_prime} and equation \eqref{eqn:correspondence_ors_pearly_space_TQ_mod_Gamma_and_D_u}. Therefore the isomorphism $C(\delta) = C(v) \circ C(u)$ corresponds to a certain orientation of $\ddd(D_w|_{Y_\Gamma^w})$. Let us find this orientation.

We will separate the proof into two cases. The first case is when both $k,l > 0$. The second case is treated below.

Consider the family of operators
$$D_{w^t}|_{Y_\Gamma^{w^t}} \fc Y_\Gamma^{w^t} \to L^p(w^t)\,;$$
these form a Fredholm morphism between the Banach bundles $(Y_\Gamma^{w^t})_{t \in (0,t_1)}$ and $(L^p(w^t))_{t\in(0,t_1)}$. Note that since $w^t$ is obtained by Morse gluing at $q'$ from $u,v$, it follows that the Banach bundle $Y_\Gamma^{w^t}$ over $(0,t_1)$ can be extended to a bundle over $[0,t_1)$ by adding the fiber over $0$ which is precisely $Y_\Gamma^u \oplus Y_\Gamma^v$. This follows from the following lemma in Morse theory, left to the reader as an exercise:
\begin{lemma}\label{lem:convergence_in_Grassmannian_Morse_theory}
Consider the exterior sum $TL \boxplus TL$ as a rank $2n$ vector bundle over $L \times L$ and consider the associated Grassmann bundle of $n$-dimensional subspaces $G_n(TL \boxplus TL)$. Then in this space $\Gamma_{(w^t_k(1),w^t_{k+1}(-1))}$, that is the graph of the differential of the flow map of $-\nabla_\rho f$ connecting $w^t_k(1)$ with $w^t_{k+1}(-1)$, converges to $T_{u_k(1)}\cS(q') \oplus T_{v_1(-1)}\cU(q')$ as $t \to 0$. \qed
\end{lemma}
\noindent The bundle $L^p(w^t)$ extends to $t=0$ by $L^p(u) \oplus L^p(v)$. The Fredholm morphism $D_{w^t}|_{Y_\Gamma^{w^t}}$ also extends over $0$, where it coincides with $D_u|_{Y_\Gamma^u} \oplus D_v|_{Y_\Gamma^v}$.

We can deform $w^t$ into $(w^t)' \in (C^\infty(D^2,S^1;M,L))^{k+l}$ with $\ev((w^t)') \in \{q\} \times \Delta_L^{k+l-1} \times \{q''\}$. Similarly, we can deform $u$ and $v$ into $u',v'$ with $\ev(u') \in \{q\} \times \Delta_L^{k-1} \times \{q'\}$ and $\ev(v') \in \{q'\} \times \Delta_L^{l-1} \times \{q''\}$, according to the prescription after equation \eqref{eqn:deformation_iso_from_u_to_u_prime}. Note that these deformations can be done so that $(w^t)'$ tends to $w':=(u_1',\dots,u_k',v_1',\dots,v_l')$ as $t \to 0$. These considerations show that the operator $D_{w^t}|_{Y_\Gamma^{w^t}}$ deforms into the operators $D_u|_{Y_\Gamma^u} \oplus D_v|_{Y_\Gamma^v}$, $D_{(w^t)'}|_{Y_\Gamma^{(w^t)'}}$, and $D_{w'}|_{Y_\Gamma^{w'}}$. Note as well that there is a deformation between $D_u|_{Y_\Gamma^u} \oplus D_v|_{Y_\Gamma^v}$ and $D_{u'}|_{Y_\Gamma^{u'}} \oplus D_{v'}|_{Y_\Gamma^{v'}}$, while the latter operator can be deformed as described above into $D_{w'}|_{Y_\Gamma^{w'}}$ by deforming the incidence condition at $q'$.

This implies that the families $D_A \sharp T\cS(q) \oplus D_{w^t}|_{Y_\Gamma^{w^t}}$, $D_A \sharp T\cS(q) \oplus D_u|_{Y_\Gamma^u} \oplus D_v|_{Y_\Gamma^v}$, $D_A \sharp T\cS(q) \oplus D_{(w^t)'}|_{Y_\Gamma^{(w^t)'}}$, $D_A \sharp T\cS(q) \oplus D_{w'}|_{Y_\Gamma^{w'}}$, and $D_A \sharp T\cS(q) \oplus D_{u'}|_{Y_\Gamma^{u'}} \oplus D_{v'}|_{Y_\Gamma^{v'}}$ are all deformations of one another. Moreover, the space of parameters of these deformations, which consists of the interval $[0,t_1)$ multiplied by the parameter used to deform $u$ into $u'$ and so on, is contractible, being a product of intervals. Therefore all these operators have mutually canonically isomorphic determinant lines, the isomorphisms being deformation isomorphisms. Furthermore, by deforming the incidence conditions at $q$ and $q'$, we see that $D_A \sharp T\cS(q) \oplus D_{u'}|_{Y_\Gamma^{u'}} \oplus D_{v'}|_{Y_\Gamma^{v'}}$ deforms into $D_A \sharp D_{u'}|_{Z_\Gamma^{u'}} \sharp D_{v'}|_{X_\Gamma^{v'}}$, while $D_A \sharp T\cS(q) \oplus D_{w'}|_{Y_\Gamma^{w'}}$ deforms into $D_A \sharp D_{w'}|_{X_\Gamma^{w'}}$. Both the operators $D_A \sharp D_{u'}|_{Z_\Gamma^{u'}} \sharp D_{v'}|_{X_\Gamma^{v'}}$ and $D_A \sharp D_{w'}|_{X_\Gamma^{w'}}$ can be glued to yield representatives of the family $D_{A''} \sharp T\cS(q'')$. These considerations, together with the direct sum isomorphisms, yield the commutative diagram
$$\resizebox{\textwidth}{!}{\xymatrix{\ddd(D_v|_{Y_\Gamma^v}) \otimes \ddd(D_u|_{Y_\Gamma^u}) \otimes \ddd(D_A\sharp T\cS(q)) \ar[r] \ar[d] &\ddd(D_{v'}|_{Y_\Gamma^{v'}}) \otimes \ddd(D_{u'}|_{Y_\Gamma^{u'}}) \otimes \ddd(D_A\sharp T\cS(q)) \ar[r] \ar[d] & \ddd(D_{A''} \sharp T\cS(q''))\\
\ddd(D_w|_{Y_\Gamma^w}) \otimes \ddd(D_A\sharp T\cS(q)) \ar[r] &\ddd(D_{w'}|_{Y_\Gamma^{w'}}) \otimes \ddd(D_A\sharp T\cS(q)) \ar[ru]}}$$
Let $\mfo \in \ddd(D_A\sharp T\cS(q))$, $\mfo'' = C(\delta)(\mfo)$. The isomorphisms $C(u),C(v)$ correspond to orientations $\mfo_u \in \ddd(D_u|_{Y^u_\Gamma})$, $\mfo_v \in \ddd(D_v|_{Y^v_\Gamma})$. Consider the isomorphism
$$\ddd(D_v|_{Y^v_\Gamma}) \otimes \ddd(D_u|_{Y^u_\Gamma}) \to \ddd(D_w|_{Y^w_\Gamma})$$
obtained by composing the direct sum isomorphism with the deformation isomorphism described above, and let $\mfo_w \in \ddd(D_w|_{Y^w_\Gamma})$ be the image of $\mfo_v \otimes \mfo_u$ under this isomorphism. Letting $\mfo_{u'} \in \ddd(D_{u'}|_{Y^{u'}_\Gamma})$, $\mfo_{v'} \in \ddd(D_{v'}|_{Y^{v'}_\Gamma})$, $\mfo_{w'} \in \ddd(D_{w'}|_{Y^{w'}_\Gamma})$ be the orientations corresponding to $\mfo_u,\mfo_v, \mfo_w$ under the deformation isomorphisms \eqref{eqn:deformation_iso_from_u_to_u_prime}, we see that the diagram maps
$$\xymatrix{\mfo_v \otimes \mfo_u \otimes \mfo \ar@{|->}[r] \ar@{|->}[d] & \mfo_{v'} \otimes \mfo_{u'} \otimes \mfo \ar@{|->}[r] \ar@{|->}[d] & \mfo''\\
\mfo_w \otimes \mfo \ar@{|->}[r] & \mfo_{w'} \otimes \mfo \ar@{|->}[ru]}$$
Since the composition of the two bottom arrows maps $\mfo_w \otimes \mfo \mapsto \mfo''$, we see that the isomorphism $C(\delta)$, which maps $\mfo \mapsto \mfo''$, corresponds to the orientation $\mfo_w$.

The next step is showing that this orientation $\mfo_w$ corresponds to the orientation \eqref{eqn:orientation_induced_on_pearly_space_TQ_mod_Gamma_by_C_delta}.

We have the following commutative diagram
\begin{equation}\label{dia:computation_induced_ori_bdry_op_squared_QH}
\resizebox{\textwidth}{!}{\xymatrix{\ddd(D_v|_{Y_\Gamma^v}) \otimes \ddd(T_vQ/\Gamma) \otimes \ddd(D_u|_{Y_\Gamma^u}) \otimes \ddd(T_uQ/\Gamma) \otimes \ddd(T^k_wQ/\Gamma) \ar[d] \ar[r] & \ddd(D_v|_{Y_Q^v}) \otimes \ddd(D_u|_{Y_Q^u}) \otimes \ddd(T^k_wQ/\Gamma) \ar[d] \\
\ddd(D_v|_{Y_\Gamma^v} \oplus D_u|_{Y_\Gamma^u}) \otimes \ddd(T_vQ/\Gamma \oplus T_uQ/\Gamma) \otimes \ddd(T^k_wQ/\Gamma) \ar[d] \ar[r] & \ddd(D_v|_{Y_Q^v} \oplus D_u|_{Y_Q^u}) \otimes \ddd(T^k_wQ/\Gamma) \ar[d] \\
\ddd(D_w|_{Y_\Gamma^w}) \otimes \ddd(T_w^{\neg k}Q/\Gamma) \otimes \ddd(T^k_wQ/\Gamma) \ar[d] \ar[r] & \ddd(D_w|_{Y_Q^{w,\neg k}}) \otimes \ddd(T^k_wQ/\Gamma) \ar[d]\\
\ddd(D_w|_{Y_\Gamma^w}) \otimes \ddd(T_wQ/\Gamma) \ar[r]& \ddd(D_w|_{Y_Q^w})} }
\end{equation}
where we use the following notations:
$$T_w^kQ/\Gamma = T_{(w_k(1),w_{k+1}(-1))}Q/\Gamma_{(w_k(1),w_{k+1}(-1))}\,,$$
$$T_w^{\neg k}Q/\Gamma = \bigoplus_{j\neq k}T_{(w_j(1),w_{j+1}(-1))}Q/\Gamma_{(w_j(1),w_{j+1}(-1))}\,,$$
$$Y_Q^{w,\neg k} = \{\xi \in Y_Q^w\,|\, (\xi_k(1),\xi_{k+1}(-1)) \in \Gamma_{w_k(1),w_{k+1}(-1)}\}\,.$$
The top square is obtained as follows. We have the exact Fredholm square
$$\xymatrix{D_v|_{Y_\Gamma^v} \ar[r] \ar[d] & D_v|_{Y_Q^v} \ar[r] \ar[d] & 0_{T_vQ/\Gamma} \ar[d] \\
D_v|_{Y_\Gamma^v} \oplus D_u|_{Y_\Gamma^u} \ar[r] \ar[d] & D_v|_{Y_Q^v} \oplus D_u|_{Y_Q^u} \ar[r] \ar[d] & 0_{T_vQ/\Gamma} \oplus 0_{T_uQ/\Gamma} \ar[d] \\
D_u|_{Y_\Gamma^u} \ar[r] & D_u|_{Y_Q^u} \ar[r] & 0_{T_uQ/\Gamma}}$$
The top square is then the commutative square corresponding to this exact Fredholm square, see \S\ref{par:exact_squares_pty}, tensored with the identity on $\ddd(T_w^kQ/\Gamma)$.

The bottom square is the commutative square corresponding to the exact Fredholm square
$$\xymatrix{D_w|_{Y_\Gamma^w} \ar@{=}[r] \ar[d] & D_w|_{Y_\Gamma^w} \ar[r] \ar [d] & 0 \ar[d] \\
D_w|_{Y_Q^{w,\neg k}} \ar[r] \ar[d] & D_w|_{Y_Q^w} \ar[r] \ar[d] & 0_{T^k_wQ/\Gamma} \ar@{=}[d]\\
0_{T^{\neg k}_wQ/\Gamma} \ar[r] & 0_{T_wQ/\Gamma} \ar[r] & 0_{T^k_wQ/\Gamma}}$$
It remains to describe the middle square. For every $t>0$ we have the exact Fredholm triple
$$0 \to D_{w^t}|_{Y_\Gamma^{w^t}} \to D_{w^t}|_{Y_Q^{w^t,\neg k}} \to 0_{T_{w^t}^{\neg k}Q/\Gamma} \to 0$$
yielding the isomorphism
\begin{equation}\label{eqn:iso_det_lines_w_t_close_to_boundary_point_pearly_space}
\ddd(D_{w^t}|_{Y_Q^{w^t,\neg k}}) \simeq \ddd(D_{w^t}|_{Y_\Gamma^{w^t}}) \otimes \ddd(T_{w^t}^{\neg k}Q/\Gamma)
\end{equation}
which is continuous in $t$. As $t \to 0$, the space $T_{w^t}^{\neg k}Q/\Gamma$ naturally converges to $T_uQ/\Gamma \oplus T_vQ/\Gamma$. Moreover, a Banach bundle argument similar to the one just before Lemma \ref{lem:convergence_in_Grassmannian_Morse_theory} shows that $D_{w^t}|_{Y_Q^{w^t,\neg k}}$ deforms into $D_u|_{Y_Q^u} \oplus D_v|_{Y_Q^v}$ as $t \to 0$. Therefore we have natural isomorphisms
$$\ddd(D_{w^t}|_{Y_\Gamma^{w^t}}) \simeq \ddd(D_u|_{Y_\Gamma^u} \oplus D_v|_{Y_\Gamma^v}) \quad \text{ and }\quad \ddd(D_{w^t}|_{Y_Q^{w^t,\neg k}}) \simeq \ddd(D_u|_{Y_Q^u} \oplus D_v|_{Y_Q^v})\,.$$
The middle square is obtained by substituting these isomorphisms into \eqref{eqn:iso_det_lines_w_t_close_to_boundary_point_pearly_space}, and tensoring with the identity on $\ddd(T_w^kQ/\Gamma)$.

We will see shortly that the diagram \eqref{dia:computation_induced_ori_bdry_op_squared_QH} maps
\begin{equation}\label{dia:calculation_induced_or_boundary_op_squared_how_diagram_maps}
\xymatrix{\mfo_v \otimes \bigwedge_i e^v_i \otimes \mfo_u \otimes \bigwedge_i e^u_i \otimes e^w_k \ar@{|->}[r] \ar@{|->}[d] & (-1)^{k+l} \bigwedge_i \partial_{v_i} \otimes \bigwedge_i \partial_{u_i} \otimes e^w_k \ar@{|->}[d] \\
(-1)^{l-1} (\mfo_v \wedge \mfo_u) \otimes (\bigwedge_i e^v_i \wedge \bigwedge_i e^u_i) \otimes e^w_k \ar@{|->}[r] \ar@{|->}[d] & (-1)^{k+l}(\bigwedge_i \partial_{v_i} \wedge \bigwedge_i \partial_{u_i}) \otimes e^w_k \ar@{|->}[d] \\
(-1)^{k(l-1)} \mfo_w \otimes \bigwedge_{i \neq k} e^w_i \otimes e^w_k \ar@{|->}[r] \ar@{|->}[d] & (-1)^{kl+k+l} \bigwedge_i \partial_{w_i} \otimes e^w_k \ar@{|->}[d] \\
(-1)^{(k+1)(l-1)} \mfo_w \otimes \bigwedge_i e^w_i \ar@{|->}[r] &  (-1)^{kl + k + l} \bigwedge_i \partial_{w_i} \wedge \text{inward}_{\delta}}
\end{equation}
but first let us deduce the desired result for $k,l \neq 0$. We see that the bottom arrow maps
$$\textstyle \mfo_w \otimes \bigwedge_i e^w_i \mapsto -\bigwedge_i \partial_{w_i} \wedge \text{inward}_{\delta}\,,$$
which means that $\mfo_w$ corresponds to the orientation
$$\textstyle - \bigwedge_i \partial_{w_i} \wedge \text{inward}_{\delta} \otimes \bigwedge_i e^w_i \in \ddd(T_w\wt\cP(q,q'')) \otimes \ddd(T_wQ/\Gamma)\,.$$
On the other hand, we saw above that $C(\delta)$ corresponds to $\mfo_w$. Therefore $C(\delta)$ corresponds to the orientation \eqref{eqn:orientation_induced_on_pearly_space_TQ_mod_Gamma_by_C_delta}, as claimed.

Let us explain the diagram \eqref{dia:calculation_induced_or_boundary_op_squared_how_diagram_maps}. Here we can compute the top horizontal arrow, as well as all the vertical arrows using definitions and the normalization property, see \S\ref{par:normalization_pty}. Our goal is the bottom horizontal arrow. Since the diagram commutes, we can compute the remaining three horizontal arrows, including the bottom one.

In the top square the left arrow consists of direct sum isomorphisms, together with the interchange isomorphism, which is responsible for the sign $(-1)^{l-1}$, since $\ind D_u|_{Y_\Gamma^u} = 1$, $\dim T_vQ/\Gamma = l-1$. The right arrow can be computed using the normalization property, since $D_u|_{Y_Q^u},D_v|_{Y_Q^v}$ are surjective. Its top horizontal arrow is obtained as the tensor product of isomorphisms \eqref{eqn:correspondence_ors_pearly_space_TQ_mod_Gamma_and_D_u} for $u,v$.

In the middle square, the right arrow comes from the fact that under deformation, $\partial_{u_i} \mapsto \partial_{w_i}$ while $\partial_{v_i} \mapsto \partial_{w_{k+i}}$, and also because $\bigwedge_i\partial_{v_i} \wedge \bigwedge_i\partial_{u_i} = (-1)^{kl} \bigwedge_i\partial_{u_i} \wedge \bigwedge_i\partial_{v_i}$. The left arrow is a combination of the deformation isomorphism sending $\mfo_v \wedge \mfo_u \mapsto \mfo_w$ (this is how we defined $\mfo_w$), and the fact that the deformation isomorphism sends $e_i^u \mapsto e_i^w$, $e_i^v \mapsto e_{k+i}^w$, and that $\bigwedge_i e_i^v \wedge \bigwedge_i e_i^u = (-1)^{(k-1)(l-1)} \bigwedge_i e_i^u \wedge \bigwedge_i e_i^v$.

In the bottom square, the left arrow comes from the normalization property together with the equality $\bigwedge_{i \neq k} e_i^w \wedge e_k^w = (-1)^{l-1}\bigwedge_i e_i^w$. The right arrow comes from the normalization property together with the fact that the map $\ker D_w|_{Y_Q^w} \to T^k_wQ/\Gamma$ maps $\text{inward}_\delta \mapsto e_k^w$, because the vector $e_k^w$ corresponds to shrinking the segment of gradient line between the points $w_k(1),w_{k+1}(-1)$, which evidently corresponds to moving the pearly trajectory away from the boundary point $\delta$.

Now we treat the case when precisely one of the numbers $k,l$ vanishes. Assume first $k = 0$. Using a combination of the techniques for the Morse case and the above treatment, we can obtain the following commutative diagram:
$$\xymatrix{\ddd(D_v|_{Y_\Gamma^v}) \otimes \ddd(T_u\wt\cP(q,q')) \otimes \ddd(D_A\sharp T\cS(q)) \ar[r] \ar[d] &\ddd(D_{v}|_{Y_\Gamma^{v}}) \otimes \ddd(D_{A'}\sharp T\cS(q'))  \ar[d] \\
\ddd(D_w|_{Y_\Gamma^w}) \otimes \ddd(D_A\sharp T\cS(q)) \ar[r] & \ddd(D_{A''} \sharp T\cS(q''))}$$
in which the top arrow is obtained from the correspondence described on page \pageref{eqn:exact_seq_moduli_space_gradient_lines_stable_unstable}, the right and the bottom arrows come from direct sum, deformation, and gluing isomorphisms, and the left arrow comes from the differential of the gluing map on pearly spaces. Recall that $C(v)$ corresponds to the orientation $\mfo_v$, $C(u)$ to the orientation $\partial_u$, and let $\mfo_w \in \ddd(D_w|_{Y_\Gamma^w})$ be the orientation which is the image of $\mfo_v \otimes \partial_u$ by the isomorphism defining the left arrow. It follows that the isomorphism $C(v) \circ C(u) \fc C(q,A) \simeq C(q'',A'')$ corresponds to $\mfo_w$, and our task now is to compute this orientation. We have the following commutative diagram
\begin{equation}\label{dia:calculation_induced_or_boundary_op_squared_k_or_l_zero} 
\xymatrix{\ddd(D_v|_{Y_\Gamma^v}) \otimes \ddd(T_vQ/\Gamma) \otimes \ddd(T_u\wt\cP(q,q')) \ar[r] \ar[d] & \ddd(D_v|_{Y_Q^v}) \otimes \ddd(T_u\wt\cP(q,q')) \ar[d] \\ 
\ddd(D_w|_{Y_\Gamma^w}) \otimes \ddd(T_wQ/\Gamma) \ar[r] & \ddd(D_w|_{Y_Q^w})}
\end{equation}
which maps
$$\xymatrix{\mfo_v \otimes \bigwedge_ie^v_i \otimes \partial_u \ar@{|->}[r] \ar@{|->}[d]& (-1)^{l+1} \bigwedge_i \partial_{v_i} \otimes \partial_u \ar@{|->}[d] \\ (-1)^{l-1}\mfo_w \otimes \bigwedge_i e^w_i \ar@{|->}[r] & (-1)^l \bigwedge_i \partial_{w_i} \wedge \text{inward}_\delta}$$
We see that $\mfo_w$ corresponds to the orientation $-\bigwedge_i \partial_{w_i} \wedge \text{inward}_\delta \otimes \bigwedge_ie^w_i$, as claimed. To explain the diagram \eqref{dia:calculation_induced_or_boundary_op_squared_k_or_l_zero}, note that the top arrow comes from the definition of the boundary operator, the left arrow acquires the Koszul sign, while the right arrow comes from the fact that the deformation isomorphism maps $\partial_{v_i} \mapsto \partial_{w_i}$ and $\partial_u \mapsto - \inward_\delta$. The proof for the case $k > 0$, $l = 0$, follows the same scheme.

Now we prove the induction step. This consists of showing the following: assume that the isomorphism $C(\delta)$ corresponds to the orientation
$$\textstyle -\bigwedge_i \partial_{w^t_i} \wedge \text{\rm inward}_\delta \otimes \bigwedge_i e^{w^t}_i \in \ddd(T_{w^t}\wt\cP(q,q'')) \otimes \ddd(T_{w^t}Q/\Gamma)$$
for $t \in (t_j,t_{j+1})$; then for $t\in(t_{j+1},t_{j+2})$ it corresponds to the orientation
$$\textstyle - \bigwedge_i \partial_{w^t_i} \wedge \text{\rm inward}_\delta \otimes \bigwedge_i e^{w^t}_i \in \ddd(T_{w^t}\wt\cP(q,q'')) \otimes \ddd(T_{w^t}Q/\Gamma)\,.$$
Without loss of generality we may assume that the degeneration of the pearly trajectory as $t \nearrow t_{j+1}$ consists of the shrinking of one of the gradient trajectories to a point. We will use $w^t$ for $t$ close to $t_{j+1}$ but smaller than or equal to it, and denote $w' = w^t$ for $t$ close to $t_{j+1}$ but strictly larger than it. Assume that the disks in $w^t$ which collide as $t \nearrow t_{j+1}$ carry numbers $r,r+1$. In this case the $r$-th disk of $w'$ is obtained by gluing the two colliding disks $w^{t_{j+1}}_r, w^{t_{j+1}}_{r+1}$.

We have the following commutative diagram for $t \leq t_{j+1}$ close to $t_{j+1}$.
\begin{equation}\label{dia:sign_computation_crossing_disk_collision_boundary_op_squares_zero}
\xymatrix{\ddd(D_{w^t}|_{Y_\Gamma^{w^t}}) \otimes \ddd(T_{w^t}Q/\Gamma) \ar[r] & \ddd(D_{w^t}|_{Y_Q^{w^t}})\\
\ddd(D_{w^t}|_{Y_\Gamma^{w^t}}) \otimes \ddd(T_{w^t}^{\neg r}Q/\Gamma) \otimes \ddd(T_{w^t}^rQ/\Gamma) \ar[r] \ar[u] \ar[d] & \ddd(D_{w^t}|_{Y_Q^{{w^t},\neg r}}) \otimes \ddd(T_{w^t}^rQ/\Gamma) \ar[u] \ar[d] \\
\ddd(D_{w'}|_{Y_\Gamma^{w'}}) \otimes \ddd(T_{w'}Q/\Gamma) \otimes \ddd(T_{w^t}^rQ/\Gamma) \ar[r] & \ddd(D_{w'}|_{Y_Q^{w'}}) \otimes \ddd(T_{w^t}^rQ/\Gamma)}
\end{equation}
The top square is the commutative square corresponding to the following exact Fredholm square:
$$\xymatrix{D_{w^t}|_{Y_\Gamma^{w^t}} \ar@{=}[r] \ar[d] & D_{w^t}|_{Y_\Gamma^{w^t}} \ar[r] \ar [d] & 0 \ar[d] \\
D_{w^t}|_{Y_Q^{{w^t},\neg r}} \ar[r] \ar[d] & D_{w^t}|_{Y_Q^{w^t}} \ar[r] \ar[d] & 0_{T^r_{w^t}Q/\Gamma} \ar@{=}[d]\\
0_{T^{\neg r}_{w^t}Q/\Gamma} \ar[r] & 0_{T_{w^t}Q/\Gamma} \ar[r] & 0_{T^r_{w^t}Q/\Gamma}}$$
The bottom square is obtained as follows. Since the pearly trajectory $w'$ is obtained from $w:=w^{t_{j+1}}$ by gluing, see \cite{Biran_Cornea_Quantum_structures_Lagr_submfds, Biran_Cornea_Rigidity_uniruling}, we have a natural isomorphism
$$\ddd(D_w|_{Y_Q^{w,\neg r}}) \simeq \ddd(D_{w'}|_{Y_Q^{w'}})$$
expressing the fact that the tangent space to the space of degenerate pearls at $w$, which is precisely the kernel of $D_w|_{Y_Q^{w,\neg r}}$, maps isomorphically onto $T_{w'}\wt\cP(q,q'') = \ker D_{w'}|_{Y_Q^{w'}}$ by the differential of the gluing map. The operators $D_w|_{Y_Q^{w,\neg r}}, D_{w'}|_{Y_Q^{w'}}$ also appear in the following exact Fredholm triples:
\begin{equation}\label{eqn:exact_triple_w_crossing_disk_collision_boundary_op_squares_zero}
0 \to D_w|_{Y_\Gamma^w} \to D_w|_{Y_Q^{w,\neg r}} \to 0_{T_w^{\neg r}Q/\Gamma} \to 0\,.
\end{equation}
\begin{equation}\label{eqn:exact_triple_w_prime_crossing_disk_collision_boundary_op_squares_zero}
0 \to D_{w'}|_{Y_\Gamma^{w'}} \to D_{w'}|_{Y_Q^{w'}} \to 0_{T_{w'}Q/\Gamma} \to 0\,.
\end{equation}
We note that the operator $D_w|_{Y_\Gamma^w}$ can be boundary glued at the point where the disks $w_r,w_{r+1}$ touch, and the result can be deformed into $D_{w'}|_{Y_\Gamma^{w'}}$. Also, as $t \searrow t_{j+1}$, the space $T_{(w')^t}Q/\Gamma$ tends to $T_w^{\neg r}Q/\Gamma$. It is then the feature of the gluing map that the following diagram commutes:
$$\xymatrix{\ddd(D_w|_{Y_\Gamma^w}) \otimes \ddd(T_w^{\neg r}Q/\Gamma) \ar[r] \ar[d] & \ddd(D_w|_{Y_Q^{w,\neg r}}) \ar[d] \\
\ddd(D_{w'}|_{Y_\Gamma^{w'}}) \otimes \ddd(T_{w'}Q/\Gamma) \ar[r] & \ddd(D_{w'}|_{Y_Q^{w'}})}$$
where the horizontal arrows come from the exact triples \eqref{eqn:exact_triple_w_crossing_disk_collision_boundary_op_squares_zero}, \eqref{eqn:exact_triple_w_prime_crossing_disk_collision_boundary_op_squares_zero}, while the vertical arrows come from gluing and deformation. The bottom square of the diagram \eqref{dia:sign_computation_crossing_disk_collision_boundary_op_squares_zero} is obtained by taking the middle horizontal arrow, letting $t \nearrow t_{j+1}$, except in the term $\ddd(T_{w^t}^rQ/\Gamma)$, where $t$ is kept fixed, applying the square to the resulting isomorphism, and finally tensoring with $\ddd(T_{w^t}^rQ/\Gamma)$.

Now the isomorphism $C(\delta)$ corresponds to an orientation $\mfo_{w^t} \in \ddd(D_{w^t}|_{Y_\Gamma^{w^t}})$, an orientation $\mfo_w \in \ddd(D_w|_{Y_\Gamma^w})$, and an orientation $\mfo_{w'} \in \ddd(D_{w'}|_{Y_\Gamma^{w'}})$. It is not difficult to see that the deformation isomorphism $\ddd(D_{w^t}|_{Y_\Gamma^{w^t}}) \simeq \ddd(D_w|_{Y^w_\Gamma})$ maps $\mfo_{w^t} \mapsto \mfo_w$, and the above isomorphism $\ddd(D_w|_{Y_\Gamma^w}) \simeq \ddd(D_{w'}|_{Y_\Gamma^{w'}})$ maps $\mfo_w\mapsto \mfo_{w'}$. Assume $w$ has $s$ disks. Then the diagram \eqref{dia:sign_computation_crossing_disk_collision_boundary_op_squares_zero} maps
$$\xymatrix{(-1)^{s-r-1} \mfo_{w^t} \otimes \bigwedge_i e^{w^t}_i \ar@{|->}[r] & (-1)^{s-r} \bigwedge_i\partial_{w^t_i} \wedge \text{inward}_\delta\\
\mfo_{w^t} \otimes \bigwedge_{i\neq r} e^{w^t}_i \otimes e^{w^t}_r \ar@{|->}[d] \ar@{|->}[r] \ar@{|->}[u] & (-1)^{s-r} \bigwedge_i\partial_{w^t_i} \otimes e^{w^t}_r \ar@{|->}[u] \ar@{|->}[d]\\
\mfo_{w'} \otimes \bigwedge_i e^{w'}_i \otimes e^{w^t}_r \ar@{|->}[r] & -\bigwedge_i \partial_{w'_i} \wedge \text{inward}_\delta \otimes e^{w^t}_r}$$
We will derive this shortly, but now let us see how it implies the claim. We see that the assumption that the orientation $\mfo_{w^t}$ corresponds to the orientation
$$\textstyle - \bigwedge_i \partial_{w^t_i} \wedge \text{inward}_\delta \otimes \bigwedge_i e^{w^t}_i$$
implies that the orientation $\mfo_{w'}$ corresponds to the orientation 
$$\textstyle - \bigwedge_i \partial_{w'_i} \wedge \text{inward}_\delta \otimes \bigwedge_i e^{w'}_i\,,$$
which is what the bottom arrow tells us. This means that the sign stays the same, which was precisely what we wanted to prove. 

Let us now explain the diagram. Here we can compute the top arrow and all the vertical arrows, while the remaining two horizontal arrows are obtained by commutativity. The goal of the computation is the bottom horizontal arrow. In the top square, the left arrow comes from the normalization property \S\ref{par:normalization_pty}, together with the fact that we have to interchange some factors in the wedge product, which gives the sign. The top arrow comes from the assumption on the orientation $\mfo_w$. The right arrow comes from the exact triple
$$0 \to D_{w^t}|_{Y_Q^{{w^t},\neg r}}  \to D_{w^t}|_{Y_Q^{w^t}} \to 0_{T^r_{w^t}Q/\Gamma} \to 0\,,$$
where we note that $\text{inward}_\delta \in \ker D_{w^t}|_{Y_Q^{w^t}}$ maps to $e^{w^t}_r \in T^r_{w^t}Q/\Gamma$, because both correspond to shrinking the piece of gradient trajectory separating the disks $w^t_r$, $w^t_{r+1}$.

In the bottom square the left arrow comes from the deformation isomorphisms, which map $\mfo_{w^t} \mapsto \mfo_{w'}$ and $e_i^{w^t} \mapsto e_i^{w'}$ for $i < r$, $e_i^{w^t} \mapsto e_{i-1}^{w'}$ for $i > r$. The right arrow comes from the deformation isomorphism, together with the following computation. It is straightforward to check that the differential of the gluing map, which is an isomorphism $\ker(D_{w^t}|_{Y_Q^{{w^t},\neg r}}) \simeq \ker(D_{w'}|_{Y_Q^{w'}})$, maps $\partial_{w^t_i} \mapsto \partial_{w'_i}$ for $i < r$, $\partial_{w^t_r} + \partial_{w^t_{r+1}} \mapsto \partial_{w'_r}$, $-\partial_{w^t_r} + \partial_{w^t_{r+1}} \mapsto \text{inward}_\delta$, and $\partial_{w^t_i} \mapsto \partial_{w'_{i-1}}$ for $i > r+1$. Therefore the induced isomorphism $\ddd(D_{w^t}|_{Y_Q^{{w^t},\neg r}}) \simeq \ddd(D_{w'}|_{Y_Q^{w'}})$ maps
$$\textstyle \bigwedge_i \partial_{w^t_i} \mapsto \bigwedge_{i \leq r}\partial_{w'_i} \wedge \text{inward}_\delta \wedge \bigwedge_{i > r+1}\partial_{w'_{i-1}} = (-1)^{s-r-1} \bigwedge_i \partial_{w'_i} \wedge \text{inward}_\delta\,.$$

The proof of the lemma, and therefore of Theorem \ref{thm:boundary_op_squares_zero_QH}, is now complete. \qed
\end{prf}

This theorem allows us to define the quantum homology
$$QH_*(\cD:L)$$
as the homology of the complex $(QC_*(\cD:L),\partial_\cD)$.

\subsubsection{Product}\label{sss:product_Lagr_QH}

Fix quantum data $\cD_i = (f_i,\rho,J)$ for $L$ for $i=0,1,2$, where we assume that the pairs $(f_i,\rho)$ are Morse--Smale and $J$ is chosen so that the $\cD_i$ are regular. The \tb{quantum product} is a bilinear map
\begin{equation}\label{eqn:product_QH_chain_level}
\star \fc QC_k(\cD_0:L) \otimes QC_l(\cD_1:L) \to QC_{k+l-n}(\cD_2:L)\,.
\end{equation}
As with the boundary operator, this is defined by its matrix elements which are homomorphisms
$$C(q_0,A_0) \otimes C(q_1,A_1) \to C(q_2,A_2)$$
for $q_i \in \Crit f_i$, $A_i \in \pi_2(M,L,q_i)$ such that $|q_0| + |q_1| - |q_2| - \mu(A_0) - \mu (A_1) + \mu(A_2) - n = 0$.
In order to define these matrix elements, we need first to describe the spaces of pearly triangles. Therefore fix $q_i \in \Crit f_i$. Recall that $\wt\cM(L,J)$ denotes the space of parametrized nonconstant $J$-holomorphic disks with boundary on $L$. We denote by $\wt\cM^\circ(L,J)$ the space of all parametrized $J$-holomorphic disks with boundary on $L$, including constant ones. Fix $k_i \geq 0$ for $i=0,1,2$. We have the evaluation map
$$\ev \fc (C^\infty(D^2,S^1;M,L))^{k_0 + k_1 + k_2 + 1}\to L^{2(k_0+k_1+k_2)+3}$$
given by
\begin{multline*}
\ev(U=(u^0,u^1,u^2,u)) = (u^0_1(-1),u^0_1(1),\dots,u^0_{k_0}(-1),u^0_{k_0}(1),u(1);\\
 u^1_1(-1),\dots,u^1_{k_1}(1),u(e^{2\pi i/3}); u(e^{4\pi i/3}),u^2_1(-1),\dots,u^2_{k_2}(1))\,,
\end{multline*}
where $u^i \in (C^\infty(D^2,S^1;M,L))^{k_i}$ and $u \in C^\infty(D^2,S^1;M,L)$. We let
\begin{multline*}
\wt\cP_{k_0,k_1,k_2}(q_0,q_1;q_2) = \ev^{-1}(\cU_{f_0}(q_0) \times Q_{f_0,\rho}^{k_0} \times \cU_{f_1}(q_1) \times Q_{f_1,\rho}^{k_1} \times Q_{f_2,\rho}^{k_2} \times \cS_{f_2}(q_2))\cap\\ (\wt\cM(L,J))^{k_0} \times (\wt\cM(L,J))^{k_1} \times (\wt\cM(L,J))^{k_2} \times \wt\cM^\circ(L,J)\,. 
\end{multline*}
This is the space of parametrized \tb{pearly triangles}. We have a natural action of $\R^{k_0+k_1+k_1}$ on this space and we let $\cP_{k_0,k_1,k_2}(q_0,q_1;q_2)$ be the quotient. We also define
$$\wt\cP(q_0,q_1;q_2) = \bigcup_{k_0,k_1,k_2 \geq 0} \wt\cP_{k_0,k_1,k_2}(q_0,q_1;q_2) \quad \text{ and }\quad \cP(q_0,q_1;q_2) = \bigcup_{k_0,k_1,k_2 \geq 0} \cP_{k_0,k_1,k_2}(q_0,q_1;q_2)\,.$$
For $U \in (C^\infty(D^2,S^1;M,L))^{k_0+k_1+k_2 +1}$ we let $\mu(U)$ be the sum of the Maslov numbers of the constituent disks of $U$. We have the following result by Biran--Cornea \cite{Biran_Cornea_Quantum_structures_Lagr_submfds, Biran_Cornea_Rigidity_uniruling}:
\begin{prop}
For Morse-Smale pairs $(f_i,\rho)_{i=0,1,2}$ there is a subset of $\cJ(M,\omega)$ of the second category such that for each $J$ in the subset, each triple of critical points $q_i \in \Crit f_i$ and each triple of nonnegative integers $k_i$ the space $\cP_{k_0,k_1,k_2}(q_0,q_1;q_2)$ is a smooth manifold of local dimension at $[U]$
$$|q_0| + |q_1| - |q_2| + \mu(U) - n$$
whenever this number is at most $1$. \qed
\end{prop}

We proceed to the definition of the matrix element. Fix $A_i \in \pi_2(M,L,q_i)$ for $i = 0,1$. For $U \in (C^\infty(D^2,S^1;M,L))^{k_0+k_1+k_2+1}$ with
$$\ev(U) \in \cU_{f_0}(q_0) \times Q_{f_0,\rho_0}^{k_0} \times \cU_{f_1}(q_1) \times Q_{f_1,\rho_1}^{k_1} \times Q_{f_2,\rho_2}^{k_2} \times \cS_{f_2}(q_2)$$
there is an obvious way to construct a class $A_0 \sharp A_1 \sharp U \in \pi_2(M,L,q_2)$ by concatenating representatives of $A_0,A_1$ with the constituent disks of $U$ according to the gradient trajectories connecting the evaluation points of the disks. For $U \in \wt\cP(q_0,q_1;q_2)$ satisfying $|q_0| + |q_1| - |q_2| + \mu(U) - n = 0$ we will construct an isomorphism
$$C(U) \fc C(q_0,A_0) \otimes C(q_1,A_1) \simeq C(q_2,A_0\sharp A_1\sharp U)\,.$$
For a class $A_2 \in \pi_2(M,L,q_2)$ with $|q_0|+|q_1| - |q_2| - \mu(A_0) - \mu(A_1) + \mu(A_2) - n = 0$, the matrix element is defined to be
\begin{equation}\label{eqn:matrix_elts_product_QH}
\sum_{\substack{[U] \in \cP(q_0,q_1;q_2):\\A_0 \sharp A_1 \sharp U = A_2}} C(U)\fc C(q_0,A_0) \otimes C(q_1,A_1) \simeq C(q_2,A_2)\,.
\end{equation}
It remains to define the isomorphism $C(U)$. To this end we define additional Banach spaces. For $V =(v^0,v^1,v^2,v) \in (C^\infty(D^2,S^1;M,L))^{k_0+k_1+k_2+1}$ with
$$\ev(V) \in \cU_{f_0}(q_0) \times \ol Q{}_{f_0,\rho}^{k_0} \times \cU_{f_1}(q_1) \times \ol Q{}_{f_1,\rho}^{k_1} \times \ol Q{}_{f_2,\rho}^{k_2} \times \cS_{f_2}(q_2)$$
we define
\begin{multline*}
X_\Gamma^V = \{\Xi = (\xi^0,\xi^1,\xi^2,\xi) \in W^{1,p}(v^0) \oplus W^{1,p}(v^1) \oplus W^{1,p}(v^2) \oplus W^{1,p}(v)\,|\\
\xi^2_{k_2}(1) \in T_{v^2_{k_2}(1)}\cS_{f_2}(q_2)\,; (\xi^m_j(1),\xi^m_{j+1}(-1)) \in \Gamma_{(v^m_j(1),v^m_{j+1}(-1))},\\
(\xi^m_{k_m}(1),\xi(e^{2\pi i m/3})) \in \Gamma_{(v^m_{k_m}(1),v(e^{2\pi i m/3}))}\text{ for }m=0,1 \text{ and }j < k_m;\\
(\xi^2_j(1),\xi^2_{j+1}(-1)) \in \Gamma_{(v^2_j(1),v^2_{j+1}(-1))}\text{ for }j=1,\dots,k_2; (\xi(e^{4\pi i/3}),\xi^2_1(-1)) \in \Gamma_{(v(e^{4\pi i/3}),v^2_1(-1))}\}
\end{multline*}
$$Y_\Gamma^V = \{\Xi \in X_\Gamma^V\,|\, \xi^m_1(-1) \in T_{v^m_1(-1)}\cU_{f_m}(q_m)\text{ for }m=0,1\}$$
If $V$ satisfies
$$\ev(V) \in \cU_{f_0}(q_0) \times Q_{f_0,\rho}^{k_0} \times \cU_{f_1}(q_1) \times Q_{f_1,\rho}^{k_1} \times Q_{f_2,\rho}^{k_2} \times \cS_{f_2}(q_2)\,,$$
we define in addition
\begin{multline*}
Y_Q^V = \{\Xi = (\xi^0,\xi^1,\xi^2,\xi) \in W^{1,p}(v^0) \oplus W^{1,p}(v^1) \oplus W^{1,p}(v^2) \oplus W^{1,p}(v)\,|\\
\xi^2_{k_2}(1) \in T_{v^2_{k_2}(1)}\cS_{f_2}(q_2)\,; (\xi^m_j(1),\xi^m_{j+1}(-1)) \in T_{(v^m_j(1),v^m_{j+1}(-1))}Q_{f_m,\rho},\\
(\xi^m_{k_m}(1),\xi(e^{2\pi i m/3}) \in T_{(v^m_{k_m}(1),v(e^{2\pi i m/3}))}Q_{f_m,\rho}\text{ for }m=0,1 \text{ and }j < k_m;\\
(\xi^2_j(1),\xi^2_{j+1}(-1)) \in T_{(v^2_j(1),v^2_{j+1}(-1))}Q\text{ for }j=1,\dots,k_2; (\xi(e^{4\pi i/3}),\xi^2_1(-1)) \in T_{(v(e^{4\pi i/3}),v^2_1(-1))}Q_{f_2,\rho}\}
\end{multline*}
and
\begin{multline*}
T_VQ = \bigoplus_{j=1}^{k_0-1}T_{(v^0_j(1),v^0_{j+1}(-1))}Q_{f_0,\rho} \oplus T_{(v^0_{k_0}(1),v(1))}Q_{f_0,\rho} \oplus \\
\bigoplus_{j=1}^{k_1-1}T_{(v^1_j(1),v^1_{j+1}(-1))}Q_{f_1,\rho} \oplus T_{(v^1_{k_1}(1),v(e^{2\pi i/3}))}Q_{f_1,\rho} \oplus \\
T_{(v(e^{4\pi i/3}),v^2_1(-1))}Q_{f_2,\rho} \oplus \bigoplus_{j=1}^{k_2-1}T_{(v^2_j(1),v^2_{j+1}(-1))}Q_{f_2,\rho}\,.
\end{multline*}
Also we define the space $T_VQ/\Gamma$ in a manner similar to the definition \eqref{eqn:definition_of_TQ_mod_Gamma}; its dimension equals $k_0+k_1+k_2$. Note that this space has a basis defined similarly to \eqref{eqn:canonical_basis_of_TQ_mod_Gamma}, and whose elements we denote by $e_i^{v^j}$ for $i=0,1,2$ and $j=1,\dots,k_i$.

Similarly to the case of the boundary operator above, the isomorphism $C(U)$ is defined in two stages. At the first stage we construct a bijection between orientations of $D_U|_{Y_\Gamma^U}$ and orientations of the line $\ddd(T_U\wt\cP(q_0,q_1;q_2)) \otimes \ddd(T_UQ/\Gamma)$. Then we construct a bijection between orientations of $D_U|_{Y_\Gamma^U}$ and isomorphisms $C(q_0,A_0) \otimes C(q_1,A_1) \simeq C(q_2,A_2)$. Once we have these bijections, the desired isomorphism $C(U)$ is the one corresponding to the following orientation:
$$\textstyle (-1)^{k_2}\bigwedge_i\partial_{u^0_i} \wedge \bigwedge_i \partial_{u^1_i} \wedge \bigwedge_i \partial_{u^2_i} \otimes \bigwedge_i e^{u^0}_i \wedge \bigwedge_i e^{u^1}_i \wedge \bigwedge_i e^{u^2}_i \in \ddd(T_U\wt\cP_{k_0+k_1+k_2}(q_0,q_1;q_2)) \otimes \ddd(T_UQ/\Gamma)\,.$$

The first bijection is constructed as follows. The exact Fredholm triple
$$0 \to D_U|_{Y_\Gamma^U} \to D_U|_{Y_Q^U} \to 0_{T_UQ/\Gamma} \to 0$$
induces an isomorphism
$$\ddd(T_U\wt\cP(q_0,q_1;q_2)) = \ddd(D_U|_{Y_Q^U}) \simeq \ddd(D_U|_{Y_\Gamma^U}) \otimes \ddd(T_UQ/\Gamma)\,,$$
where the first equality is due to the fact that $T_U \wt\cP(q_0,q_1;q_2) = \ker D_U|_{Y_Q^U}$, which follows from the definition of the pearly triangles. Tensoring with $\ddd(T_UQ/\Gamma)$, we see that indeed there is a bijection between orientations of $\ddd(D_U|_{Y_\Gamma^U})$ and orientations of $\ddd(T_U\wt\cP(q_0,q_1;q_2)) \otimes \ddd(T_UQ/\Gamma)$.

The second bijection is constructed as follows. Consider the family of operators
$$D_{A_0} \sharp T\cS_{f_0}(q_0) \oplus D_{A_1} \sharp T\cS_{f_1}(q_1) \oplus D_U|_{Y_\Gamma^U}\,.$$
We first deform $U$ into $U' \in (C^\infty(D^2,S^1;M,L))^{k_0+k_1+k_2+1}$ with
$$\ev(U') \in \{q_0\} \times \Delta_L^{k_0} \times \{q_1\} \times \Delta_L^{k_1} \times \Delta_L^{k_2} \times \{q_2\}$$
just like we did in \S\ref{sss:boundary_op_Lagr_QH} when defining the boundary operator. Our operator therefore deforms into
$$D_{A_0} \sharp T\cS_{f_0}(q_0) \oplus D_{A_1} \sharp T\cS_{f_1}(q_1) \oplus D_{U'}|_{Y_\Gamma^{U'}}\,.$$
Next, we can deform the incidence conditions at $q_0,q_1$ in the sense of the isomorphism \eqref{eqn:iso_abstract_deformation_incidence_condition_W}, to arrive at the boundary glued operator $D_{A_0} \sharp D_{A_1} \sharp D_{U'}|_{X_\Gamma^{U'}}$, which after deformation yields a representative of the family $D_{A_2} \sharp T\cS_{f_2}(q_2)$. These deformations, together with the direct sum isomorphism, yield a string of isomorphisms
\begin{align*}
\ddd(D_U|_{Y_\Gamma^U}) \otimes \ddd(D_{A_0} \sharp T\cS_{f_0}(q_0)) \otimes \ddd(D_{A_1} \sharp T\cS_{f_1}(q_1)) &\simeq \ddd(D_{A_0} \sharp T\cS_{f_0}(q_0) \oplus D_{A_1} \sharp T\cS_{f_1}(q_1) \oplus D_U|_{Y_\Gamma^U})\\
&\simeq \ddd(D_{A_0} \sharp T\cS_{f_0}(q_0) \oplus D_{A_1} \sharp T\cS_{f_1}(q_1) \oplus D_{U'}|_{Y_\Gamma^{U'}})\\
&\simeq \ddd(D_{A_0} \sharp D_{A_1} \sharp D_{U'}|_{X_\Gamma^{U'}})\\
&\simeq \ddd(D_{A_2} \sharp T\cS_{f_2}(q_2))
\end{align*}
whose composition indeed shows that there is a bijection between isomorphisms $C(q_0,A_0) \otimes C(q_1,A_1) \simeq C(q_2,A_2)$ and orientations of $D_U|_{Y_\Gamma^U}$.

We have therefore completed the definition of the matrix elements of the product \eqref{eqn:matrix_elts_product_QH} and hence we have defined the product as a bilinear operation \eqref{eqn:product_QH_chain_level}.

We now prove
\begin{thm}\label{thm:product_QH_is_chain_map}
The operation $\star$ is a chain map. More precisely, we have 
$$\partial_{\cD_2} \circ \star =  \star \circ (\partial_{\cD_0} \otimes \id + (-1)^{n - k} \id \otimes \partial_{\cD_1}) \fc QC_k(\cD_0:L) \otimes QC_l(\cD_1:L) \to QC_{k+l-n-1}(\cD_2:L)\,.$$
\end{thm}

\begin{prf}
It suffices to prove the vanishing of the matrix element
\begin{multline}\label{eqn:vanishing_matrix_elts_to_prove_product_is_chain_map}
\sum_{q_2' \in \Crit f_2} \sum_{\substack{A_2' \in \pi_2(M,L,q_2'):\\ |q_2'| - \mu(A_2') = |q_2| - \mu(A_2) + 1}} \sum_{\substack{([U],[w]) \in \cP(q_0,q_1;q_2') \times \cP(q_2',q_2): \\ A_0 \sharp A_1 \sharp U = A_2', A_2' \sharp w = A_2}} C(w) \circ C(U) - \\
- \sum_{q_0' \in \Crit f_0} \sum_{\substack{A_0' \in \pi_2(M,L,q_0'): \\ |q_0| - \mu(A_0) = |q_0'| - \mu(A_0') + 1}} \sum_{\substack{([w],[U]) \in \cP(q_0,q_0') \times \cP(q_0',q_1;q_2): \\ A_0 \sharp w = A_0',A_0' \sharp A_1 \sharp U = A_2}} C(U) \circ (C(w) \otimes \id)\\
- (-1)^{n-k} \sum_{q_1' \in \Crit f_1} \sum_{\substack{A_1' \in \pi_2(M,L,q_1'): \\ |q_1| - \mu(A_1) = |q_1'| - \mu(A_1') + 1}} \sum_{\substack{([w],[U]) \in \cP(q_1,q_1') \times \cP(q_0,q_1';q_2): \\ A_1 \sharp w = A_1', A_0 \sharp A_1' \sharp U = A_2}} C(U) \circ (\id \otimes C(w))
\end{multline}
as a homomorphism
$$C(q_0,A_0) \otimes C(q_1,A_1) \to C(q_2,A_2)\,.$$
Let us denote by $\cP^1(q_0,q_1;q_2)$ the $1$-dimensional part of the space of pearly triangles. Biran--Cornea \cite{Biran_Cornea_Quantum_structures_Lagr_submfds, Biran_Cornea_Rigidity_uniruling} described the structure of the compactification $\ol\cP{}^1(q_0,q_1;q_2)$. The space $\cP^1(q_0,q_1;q_2)$ fails to be compact in one of three ways: either one of the gradient trajectories undergoes Morse breaking, or one of the gradient trajectories shrinks to a point, or a holomorphic disk breaks into two. Note that the breaking can happen at the core and in the resulting degenerate triangle one of the disks (the one carrying three marked points) may be constant. We redefine $\ol\cP{}^1(q_0,q_1;q_2)$ to be the union of the compactified connected components of $\cP^1(q_0,q_1;q_2)$ where two boundary points of different components are identified if they represent the same degenerate pearly triangle, where by a degenerate pearly triangle we mean a pearly triangle where precisely one of the gradient trajectories connecting two holomorphic disks has length zero. Thus $\ol\cP{}^1(q_0,q_1;q_2)$ has the structure of a compact $1$-dimensional topological manifold with boundary. Its boundary points represent Morse breaking and are in an obvious bijection with the summands of the matrix element \eqref{eqn:vanishing_matrix_elts_to_prove_product_is_chain_map}. For $\delta \in \partial \ol\cP{}^1(q_0,q_1;q_2)$ let $C(\delta)$ denote the summand of \eqref{eqn:vanishing_matrix_elts_to_prove_product_is_chain_map} corresponding to $\delta$ (together with the sign). It therefore suffices to prove the following: for any connected component $\Delta \subset \ol\cP{}^1(q_0,q_1;q_2)$ with $\partial \Delta = \{\delta,\delta'\}$ we have $C(\delta) + C(\delta') = 0$.

We now proceed to the computation of the various summands. Fix a point $\delta \in \ol\cP{}^1(q_0,q_1;q_2)$ and let $[X] \in \cP^1(q_0,q_1;q_2)$ lie close to $\delta$. As we saw, the isomorphism $C(\delta) \fc C(q_0,A_0) \otimes C(q_1,A_1) \simeq C(q_2,A_2)$ determines an orientation of $D_X|_{Y_\Gamma^X}$, which in turn corresponds to an orientation of $\ddd(T_X\wt\cP(q_0,q_1;q_2)) \otimes \ddd(T_XQ/\Gamma)$. For $V = (v^0,v^1,v^2;v) \in \wt\cP(q_0,q_1;q_2)$ we abbreviate
$$\textstyle \bigwedge_i \partial_{v^0_i} \wedge \bigwedge_i \partial_{v^1_i} \wedge \bigwedge_i \partial_{v^2_i} \wedge \text{inward}_\delta \otimes \bigwedge_i e^{v^0}_i \wedge \bigwedge_i e^{v^1}_i \wedge \bigwedge_i e^{v^2}_i \in \ddd(T_V \wt \cP(q_0,q_1;q_2)) \otimes \ddd(T_VQ/\Gamma)$$
to
$$\textstyle \bigwedge_i\partial_{V_i} \wedge \text{inward}_\delta\otimes \bigwedge_ie^V_i\,,$$
where $\text{inward}_\delta \in T_V \wt \cP(q_0,q_1;q_2)$ is a tangent vector directed away from the boundary point $\delta$. This should cause no confusion. We also refer to this particular orientation as the \tb{standard orientation} of the line $\ddd(T_V \wt \cP(q_0,q_1;q_2)) \otimes \ddd(T_VQ/\Gamma)$.

We have the following lemma.
\begin{lemma}\label{lem:computation_induced_oris_product_is_chain_map}
We have the following cases:
\begin{itemize}
\item If $\delta = ([U],[w]) \in \cP_{k_0,k_1,k_2}(q_0,q_1;q_2') \times \cP_r(q_2',q_2)$, the isomorphism $C(\delta) = C(w) \circ C(U)$ corresponds to $(-1)^{k_0+k_1}$ times the standard orientation of $\ddd(T_X\wt\cP(q_0,q_1;q_2)) \otimes \ddd(T_XQ/\Gamma)$.
\item If $\delta = ([w],[U]) \in \cP_r(q_0,q_0') \times \cP_{k_0,k_1,k_2}(q_0',q_1;q_2)$, the isomorphism $C(\delta) = -C(U) \circ (C(w) \otimes \id)$ corresponds to $(-1)^{k_0+k_1+r}$ times the standard orientation of $\ddd(T_X\wt\cP(q_0,q_1;q_2)) \otimes \ddd(T_XQ/\Gamma)$.
\item If $\delta = ([w],[U]) \in \cP_r(q_1,q_1') \times \cP_{k_0,k_1,k_2}(q_0,q_1';q_2)$, then the isomorphism $C(\delta) = -(-1)^{n-k} C(U) \circ (\id \otimes C(w))$ corresponds to $(-1)^{k_0+k_1+r}$ times the standard orientation of $\ddd(T_X\wt\cP(q_0,q_1;q_2)) \otimes \ddd(T_XQ/\Gamma)$.
\end{itemize}
\end{lemma}
\noindent Lemma \ref{lem:computation_induced_oris_product_is_chain_map} is proved below. In order to complete the proof of the vanishing of the matrix element \eqref{eqn:vanishing_matrix_elts_to_prove_product_is_chain_map}, we also need to keep track of the change of orientations when crossing disk collision/breaking points in $\ol\cP{}^1(q_0,q_1;q_2)$. This is described in the following lemma:
\begin{lemma}\label{lem:sign_change_product_is_chain_map}
Let an isomorphism $C \fc C(q_0,A_0) \otimes C(q_1,A_1) \simeq C(q_2,A_2)$ be given and let $[V],[W] \in \cP^1(q_0,q_1;q_2)$ be two points lying on two different sides of a degenerate pearly triangle in the space $\ol\cP{}^1(q_0,q_1;q_2)$, and close to it. Suppose that the isomorphism $C$ corresponds to the orientation
$$\textstyle\bigwedge_i \partial_{V_i} \wedge \eta_V \otimes \bigwedge_i e^V_i \in \ddd(T_V\wt\cP(q_0,q_1;q_2)) \otimes \ddd(T_VQ/\Gamma)$$
where $\eta_V \in T_V\wt\cP(q_0,q_1;q_2)$ is an arbitrary vector transverse to the infinitesimal action of the automorphism group at $V$. Then $C$ corresponds to the orientation
$$\textstyle\epsilon \bigwedge_i \partial_{W_i} \wedge \eta_W \otimes \bigwedge_i e^W_i \in \ddd(T_W\wt\cP(q_0,q_1;q_2)) \otimes \ddd(T_WQ/\Gamma)\,,$$
where $\eta_W \in T_W\wt\cP(q_0,q_1;q_2)$ points in the same direction \footnote{Note that this is well-defined!} as $\eta_V$ and $\epsilon \in \{\pm 1\}$ is a sign which equals $-1$ if the passage from $V$ to $W$ happens through disk breaking/collision in one of legs $0,1$ of the triange, and it equals $1$ if the breaking/collision happens in leg $2$ of the triangle.
\end{lemma}
\noindent This lemma is proved below. Let us see how Lemmas \ref{lem:computation_induced_oris_product_is_chain_map}, \ref{lem:sign_change_product_is_chain_map} allow us to complete the proof of Theorem \ref{thm:product_QH_is_chain_map}. There are various combinatorial types of possible components $\Delta \subset \ol\cP{}^1(q_0,q_1;q_2)$. If $\partial \Delta = \{\delta,\delta'\}$, the proof that $C(\delta) + C(\delta') = 0$ follows an identical argument for all of these types, therefore we will only give a full proof for one of them: suppose $\delta = ([U],[w]) \in \cP_{k_0,k_1,k_2}(q_0,q_1;q_2') \times \cP_r(q_2',q_2)$ and $\delta' = ([w'],[U']) \in \cP_{r'}(q_1,q_1') \times \cP_{k_0',k_1',k_2'}(q_0,q_1';q_2)$. Let $X^t \in \wt\cP(q_0,q_1;q_2)$ be a continuous family of pearly triangles such that $[X^t]$ gives a continuous parametrization of $\Delta \cap \cP^1(q_0,q_1;q_2)$, that is $X^t$ is defined for all but a finite number of values of $t \in [0,1]$. According to Lemma \ref{lem:computation_induced_oris_product_is_chain_map}, the isomorphism $C(\delta)$ corresponds to
$$\textstyle (-1)^{k_0+k_1} \bigwedge_i\partial_{X^t_i} \wedge \text{inward}_\delta \otimes \bigwedge_i e^{X^t}_i \in \ddd(T_{X^t}\wt\cP(q_0,q_1;q_2)) \otimes \ddd(T_{X^t}Q/\Gamma)$$
for small positive $t$. As $t$ grows from $0$ to $1$, $X^t$ undergoes a number of jumps which correspond to instances of disk collision/breaking. Let us denote by $n_i$ the number of disk collision/breaking instances which take place in leg $i$ of the triangle. Then clearly we have
$$n_0 \equiv k_0 + k_0' \mod 2\,,\quad n_1 \equiv k_1+r' + k_1' \mod 2\,,\quad n_2 \equiv k_2 + r + k_2' \mod 2\,.$$
Therefore for $t$ close to $1$ the isomorphism $C(\delta)$ corresponds to the orientation
$$\textstyle (-1)^{k_0+k_1+n_0+n_1} \bigwedge_i\partial_{X^t_i} \wedge \text{inward}_\delta \otimes \bigwedge_i e^{X^t}_i = (-1)^{k_0' + k_1' + r'} \bigwedge_i\partial_{X^t_i} \wedge \text{inward}_\delta \otimes \bigwedge_i e^{X^t}_i\,,$$
as follows from Lemma \ref{lem:sign_change_product_is_chain_map}. On the other hand, Lemma \ref{lem:computation_induced_oris_product_is_chain_map} implies that the isomorphism $C(\delta')$ corresponds to the orientation
$$\textstyle (-1)^{k_0'+k_1'+r'} \bigwedge_i\partial_{X^t_i} \wedge \text{inward}_{\delta'} \otimes \bigwedge_i e^{X^t}_i \in \ddd(T_{X^t}\wt\cP(q_0,q_1;q_2)) \otimes \ddd(T_{X^t}Q/\Gamma)$$
for $t$ close to $1$. Since clearly $\text{inward}_\delta = -\text{inward}_{\delta'}$, we see that the orientations are opposite and therefore $C(\delta) + C(\delta') = 0$. Similar arguments prove the vanishing of this sum for other combinatorial types of the components $\Delta$. This finishes the proof of the theorem. \qed
\end{prf}

This means that, modulo Lemmas \ref{lem:computation_induced_oris_product_is_chain_map}, \ref{lem:sign_change_product_is_chain_map}, we have defined a bilinear operation on homology:
$$\star \fc QH_k(\cD_0:L) \otimes QH_l(\cD_1:L) \to QH_{k+l-n}(\cD_2:L)\,.$$

\noindent We now prove Lemma \ref{lem:computation_induced_oris_product_is_chain_map}.
\begin{prf}[of Lemma \ref{lem:computation_induced_oris_product_is_chain_map}] We only prove the lemma assuming $r > 0$. The remaining case can be handled using arguments similar to those of the proof of Lemma \ref{lem:induced_or_pearly_spaces_boundary_op_squared}.

Consider the first case. Using arguments similar to those appearing in the proof of Lemma \ref{lem:induced_or_pearly_spaces_boundary_op_squared}, we obtain the following commutative diagram:
\begin{equation}\label{dia:induced_or_product_chain_map_leg_2}
\resizebox{\textwidth}{!}{\xymatrix{\ddd(D_w|_{Y_\Gamma^w}) \otimes \ddd(T_wQ/\Gamma) \otimes \ddd(D_U|_{Y_\Gamma^U}) \otimes \ddd(T_UQ/\Gamma) \otimes \ddd(T_{x^2}^{k_2+1}Q/\Gamma) \ar[r] \ar[d] & \ddd(D_w|_{Y_Q^w}) \otimes \ddd(D_U|_{Y_Q^U}) \otimes \ddd(T_{x^2}^{k_2+1}Q/\Gamma) \ar[d] \\
\ddd(D_w|_{Y_\Gamma^w} \oplus D_U|_{Y_\Gamma^U}) \otimes \ddd(T_wQ/\Gamma \oplus T_UQ/\Gamma) \otimes \ddd(T_{x^2}^{k_2+1}Q/\Gamma) \ar[r] \ar[d] & \ddd(D_w|_{Y_Q^w} \oplus D_U|_{Y_Q^U}) \otimes \ddd(T_{x^2}^{k_2+1}Q/\Gamma) \ar[d] \\
\ddd(D_X|_{Y_\Gamma^X}) \otimes \ddd(T_X^{x^2:\neg k_2+1}Q/\Gamma) \otimes \ddd(T_{x^2}^{k_2+1}Q/\Gamma) \ar[r] \ar[d] & \ddd(D_X|_{Y_Q^{X,x^2:\neg k_2+1}}) \otimes \ddd(T_{x^2}^{k_2+1}Q/\Gamma) \ar[d] \\
\ddd(D_X|_{Y_\Gamma^X}) \otimes \ddd(T_XQ/\Gamma) \ar[r] & \ddd(D_X|_{Y_Q^X})}}
\end{equation}
where
$$T_{x^2}^{k_2+1}Q/\Gamma = T_{(x^2_{k_2}(1),x^2_{k_2+1}(-1))}Q/\Gamma_{(x^2_{k_2}(1),x^2_{k_2+1}(-1))}\,,$$
$$T_X^{x^2:\neg k_2+1}Q/\Gamma = \ker \big(T_XQ/\Gamma \to T_{x^2}^{k_2+1}Q/\Gamma\big)\,,$$
$$Y_Q^{X,x^2:\neg k_2+1} = \{\Xi = (\xi^0,\xi^1,\xi^2,\xi) \in Y_Q^X\,|\, (\xi^2_{k_2}(1),\xi^2_{k_2+1}(-1)) \in \Gamma_{(x^2_{k_2}(1),x^2_{k_2+1}(-1))}\}\,.$$
The top and bottom squares of the diagram are obtained from suitable exact Fredholm squares while the middle square is obtained from gluing. Now, by definition, the isomorphism $C(w)$ corresponds to the orientation $\mfo_w \in \ddd(D_w|_{Y_\Gamma^w})$ which in turns corresponds to the orientation $(-1)^{r+1}\bigwedge_i \partial_{w_i} \otimes \bigwedge_i e^w_i \in \ddd(D_w|_{Y_Q^w}) \otimes \ddd(T_wQ/\Gamma)$; the isomorphism $C(U)$ corresponds to the orientation $\mfo_U \in \ddd(D_U|_{Y_\Gamma^U})$ which in turn corresponds to the orientation $(-1)^{k_2}\bigwedge_i \partial_{U_i} \otimes \bigwedge_i e^U_i \in \ddd(D_U|_{Y_Q^U}) \otimes \ddd(T_UQ/\Gamma)$. The direct sum isomorphism, composed with the gluing isomorphism, gives the isomorphism
\begin{equation}\label{eqn:iso_det_lines_composition_product_boundary_op_leg_2}
\ddd(D_w|_{Y_\Gamma^w}) \otimes \ddd(D_U|_{Y_\Gamma^U}) \simeq \ddd(D_X|_{Y_\Gamma^X})\,;
\end{equation}
let us denote by $\mfo_X$ the image of $\mfo_w \otimes \mfo_U$ under this isomorphism. Lemma \ref{lem:bijection_isos_C_q_A_ors_D_u_prime} together with the deformation isomorphism says that the isomorphism $C(w) \circ C(U)$ corresponds to an orientation of $\ddd(D_X|_{Y_\Gamma^X})$; it is not hard to see that this orientation is precisely $\mfo_X$ --- indeed, it follows from associativity of gluing and arguments involving deformation of incidence consditions that the following diagram commutes:
$$\xymatrix{\ddd(D_w|_{Y_\Gamma^w}) \otimes \ddd(D_U|_{Y_\Gamma^U}) \otimes \ddd(D_{A_0} \sharp T\cS(q_0)) \otimes \ddd(D_{A_1} \sharp T\cS(q_1)) \ar[r] \ar[d] & \ddd(D_{A_2}\sharp T\cS(q_2)) \\
\ddd(D_X|_{Y_\Gamma^X}) \otimes \ddd(D_{A_0} \sharp T\cS(q_0)) \otimes \ddd(D_{A_1} \sharp T\cS(q_1)) \ar[ur] }$$
where the arrows pointing to the right come from deforming the trajectories so that their disks touch, and then applying Lemma \ref{lem:bijection_isos_C_q_A_ors_D_u_prime}, therefore they correspond to the isomorphism $C(w) \circ C(U)$, while the vertical arrow is the isomorphism \eqref{eqn:iso_det_lines_composition_product_boundary_op_leg_2} tensored with identity.

Now diagram \eqref{dia:induced_or_product_chain_map_leg_2} maps
$$\resizebox{\textwidth}{!}{\xymatrix{\mfo_w \otimes \bigwedge_i e^w_i \otimes \mfo_U \otimes \bigwedge_i e^U_i \otimes e^{x^2}_{k_2+1} \ar@{|->}[r] \ar@{|->}[d] &  (-1)^{r+1}\bigwedge_i\partial_{w_i} \otimes (-1)^{k_2} \bigwedge_i \partial_{U_i} \otimes e^{x^2}_{k_2+1} \ar@{|->}[d]\\ 
(\mfo_w \wedge \mfo_U) \otimes (\bigwedge_i e^w_i \wedge \bigwedge_i e^U_i) \otimes e^{x^2}_{k_2+1} \ar@{|->}[r] \ar@{|->}[d] &  (-1)^{k_2+r+1} (\bigwedge_i\partial_{w_i} \wedge \bigwedge_i \partial_{U_i}) \otimes e^{x^2}_{k_2+1} \ar@{|->}[d] \\
(-1)^{(r-1)(k_0+k_1+k_2)}\mfo_X \otimes (\bigwedge_i e^{x^0}_i \wedge \bigwedge_i e^{x^1}_i \wedge \bigwedge_{i \neq k_2+1} e^{x^2}_i) \otimes e^{x^2}_{k_2+1}  \ar@{|->}[r] \ar@{|->}[d] &  (-1)^{k_2+r + 1 + r(k_0+k_1+k_2)} \bigwedge_i\partial_{X_i} \otimes e^{x^2}_{k_2+1} \ar@{|->}[d] \\
(-1)^{(r-1)(k_0+k_1+k_2 + 1)}\mfo_X \otimes \bigwedge_i e^X_i \ar@{|->}[r] &  (-1)^{k_2 + r + 1 + r(k_0+k_1+k_2)} \bigwedge_i\partial_{X_i} \wedge \text{inward}_\delta}}$$
We will explain in a moment why it is so. Assuming this, we see from the bottom arrow, which is the goal of this computation, that the orientation $\mfo_X$, and therefore the isomorphism $C(\delta) = C(w) \circ C(U)$, corresponds to the orientation
$$\textstyle (-1)^{k_0+k_1}\bigwedge_i\partial_{X_i} \wedge \text{inward}_\delta \otimes \bigwedge_i e^X_i\,,$$
as claimed.

Let us briefly explain the diagram. In the top square the top arrow comes from the definition of the orientations $\mfo_w,\mfo_U$, the left arrow is just the direct sum isomorphism composed with the interchange isomorphism, which does not produce a sign since $\ind D_U|_{Y_\Gamma^U} = 0$. The right arrow comes from the normalization property \S\ref{par:normalization_pty}. In the middle square we use the fact that $\partial_{u^i_j} \mapsto \partial_{x^i_j}$ for all relevant $i,j$, while $\partial_{w_i}\mapsto \partial_{x^2_{k_2+i}}$, and the additional sign comes from the interchange of factors in the wedge product. In the left arrow we similarly have $e_j^{u^i} \mapsto e_j^{x^i}$ and $e_i^w \mapsto e_{k_2+i+1}^{x^2}$, and the additional sign comes from interchange of factors in the wedge product. In the bottom square, in the right arrow we use the fact that $e_{k_2+1}^{x^2}$ lifts to $\text{inward}_\delta \in \ker D_X|_{Y_Q^X}$, while in the left arrow the additional sign again comes from interchange of factors in the wedge product.

Consider now the second case. We have the following commutative diagram, obtained by methods similar to those above.
\begin{equation}\label{dia:induced_or_product_chain_map_leg_0}
\resizebox{\textwidth}{!}{\xymatrix{\ddd(D_U|_{Y_\Gamma^U}) \otimes \ddd(T_UQ/\Gamma) \otimes \ddd(D_w|_{Y_\Gamma^w}) \otimes \ddd(T_wQ/\Gamma) \otimes \ddd(T_{x^0}^rQ/\Gamma) \ar[r] \ar[d] &  \ddd(D_U|_{Y_Q^U}) \otimes \ddd(D_w|_{Y_Q^w}) \otimes \ddd(T_{x^0}^rQ/\Gamma) \ar[d] \\
\ddd(D_U|_{Y_\Gamma^U} \oplus D_w|_{Y_\Gamma^w}) \otimes \ddd(T_UQ/\Gamma \oplus T_wQ/\Gamma ) \otimes \ddd(T_{x^0}^rQ/\Gamma) \ar[r] \ar[d] & \ddd(D_U|_{Y_Q^U} \oplus D_w|_{Y_Q^w}) \otimes \ddd(T_{x^0}^rQ/\Gamma) \ar[d] \\
\ddd(D_X|_{Y_\Gamma^X}) \otimes \ddd(T_X^{x^0:\neg r}Q/\Gamma) \otimes \ddd(T_{x^0}^rQ/\Gamma) \ar[r] \ar[d] & \ddd(D_X|_{Y_Q^{X,x^0:\neg r}}) \otimes \ddd(T_{x^0}^rQ/\Gamma) \ar[d] \\
\ddd(D_X|_{Y_\Gamma^X}) \otimes \ddd(T_XQ/\Gamma) \ar[r] & \ddd(D_X|_{Y_Q^X})}}
\end{equation}
By definition, the isomorphism $C(w)$ corresponds to orientations
$$\textstyle(-1)^{r+1}\bigwedge_i \partial_{w_i} \otimes \bigwedge_i e^w_i \in \ddd(D_w|_{Y_Q^w}) \otimes \ddd(T_wQ/\Gamma)$$
and $\mfo_w \in \ddd(D_w|_{Y_\Gamma^w})$, while $C(U)$ corresponds to orientations
$$\textstyle(-1)^{k_2}\bigwedge_i \partial_{U_i} \otimes \bigwedge_i e^U_i \in \ddd(D_U|_{Y_Q^U}) \otimes \ddd(T_UQ/\Gamma)$$
and $\mfo_U \in \ddd(D_U|_{Y_\Gamma^U})$. We have the isomorphism
$$\ddd(D_U|_{Y_\Gamma^U}) \otimes \ddd(D_w|_{Y_\Gamma^w}) \simeq \ddd(D_X|_{Y_\Gamma^X})$$
obtained by composing the direct sum isomorphism with deformation and gluing; let $\mfo_X$ be the image of $\mfo_U \otimes \mfo_w$ under this isomorphism. This isomorphism enters in the vertical arrow of the following commutative diagram:
$$\xymatrix{\ddd(D_U|_{Y_\Gamma^U}) \otimes \ddd(D_w|_{Y_\Gamma^w}) \otimes  \ddd(D_{A_0} \sharp T\cS(q_0)) \otimes \ddd(D_{A_1} \sharp T\cS(q_1)) \ar[r] \ar[d] & \ddd(D_{A_2}\sharp T\cS(q_2)) \\
\ddd(D_X|_{Y_\Gamma^X}) \otimes \ddd(D_{A_0} \sharp T\cS(q_0)) \otimes \ddd(D_{A_1} \sharp T\cS(q_1)) \ar[ur] }$$
where the horizontal arrows correspond to the isomorphism $C(U) \circ (C(w) \otimes \id)$. It follows that $\mfo_X$ is precisely the orientation corresponding to this isomorphism by an obvious modification of Lemma \ref{lem:bijection_isos_C_q_A_ors_D_u_prime} applied to a deformed triangle where the disks touch. Therefore the diagram \eqref{dia:induced_or_product_chain_map_leg_0} maps
$$\resizebox{\textwidth}{!}{\xymatrix{\mfo_U \otimes \bigwedge_i e^U_i \otimes \mfo_w \otimes \bigwedge_i e^w_i \otimes e^{x^0}_r \ar@{|->}[r] \ar@{|->}[d] & (-1)^{k_2} \bigwedge_i\partial_{U_i} \otimes (-1)^{r+1}\bigwedge_i \partial_{w_i} \otimes e^{x^0}_r \ar@{|->}[d]  \\
(-1)^{k_0+k_1+k_2}(\mfo_U \wedge \mfo_w) \otimes (\bigwedge_i e^U_i \wedge \bigwedge_i e^w_i) \otimes e^{x^0}_r \ar@{|->}[r] \ar@{|->}[d] & (-1)^{k_2+r+1} (\bigwedge_i\partial_{U_i} \wedge \bigwedge_i \partial_{w_i}) \otimes e^{x^0}_r \ar@{|->}[d] \\
(-1)^{r(k_0+k_1+k_2)} \mfo_X \otimes (\bigwedge_{i \neq r} e^{x^0}_i \wedge \bigwedge_i e^{x^1}_i \wedge \bigwedge_i e^{x^2}_i) \otimes e^{x^0}_r \ar@{|->}[r] \ar@{|->}[d] & (-1)^{k_2+r + 1 + r(k_0+k_1+k_2)} \bigwedge_i\partial_{X_i} \otimes e^{x^0}_r \ar@{|->}[d] \\
(-1)^{(r+1)(k_0+k_1+k_2)} \mfo_X \otimes \bigwedge_i e^X_i \ar@{|->}[r] & (-1)^{k_2 + r + 1 + r(k_0+k_1+k_2)} \bigwedge_i\partial_{X_i} \wedge \text{inward}_{\delta}}}$$
whence it follows that $\mfo_X$ corresponds to the isomorphism $C(U) \circ (C(w) \otimes \id)$ and to the orientation
$$\textstyle (-1)^{k_0 + k_1 + r + 1}\bigwedge_i\partial_{X_i} \wedge \text{inward}_\delta \otimes \bigwedge_i e^X_i\,,$$
as claimed. Note that $C(\delta) = -C(U) \circ (C(w) \otimes \id)$ by definition.

We now turn the the third case. We similary have the following commutative diagram:
\begin{equation}\label{dia:induced_or_product_chain_map_leg_1}
\resizebox{\textwidth}{!}{\xymatrix{\ddd(D_U|_{Y_\Gamma^U}) \otimes \ddd(T_UQ/\Gamma) \otimes \ddd(D_w|_{Y_\Gamma^w}) \otimes \ddd(T_wQ/\Gamma) \otimes \ddd(T_{x^1}^rQ/\Gamma) \ar[d] \ar[r] &  \ddd(D_U|_{Y_Q^U}) \otimes \ddd(D_w|_{Y_Q^w}) \otimes \ddd(T_{x^1}^rQ/\Gamma) \ar[d] \\
\ddd(D_U|_{Y_\Gamma^U} \oplus D_w|_{Y_\Gamma^w}) \otimes \ddd(T_UQ/\Gamma \oplus T_wQ/\Gamma ) \otimes \ddd(T_{x^1}^rQ/\Gamma) \ar[d] \ar[r] & \ddd(D_U|_{Y_Q^U} \oplus D_w|_{Y_Q^w}) \otimes \ddd(T_{x^1}^rQ/\Gamma)  \ar[d]\\
\ddd(D_X|_{Y_\Gamma^X}) \otimes \ddd(T_X^{x^1:\neg r}Q/\Gamma) \otimes \ddd(T_{x^1}^rQ/\Gamma) \ar[d] \ar[r] & \ddd(D_X|_{Y_Q^{X,x^1:\neg r}}) \otimes \ddd(T_{x^1}^rQ/\Gamma)  \ar[d] \\
\ddd(D_X|_{Y_\Gamma^X}) \otimes \ddd(T_XQ/\Gamma) \ar[r] & \ddd(D_X|_{Y_Q^X})}}
\end{equation}
We have: the isomorphism $C(w)$ corresponds to orientations
$$\textstyle(-1)^{r+1}\bigwedge_i \partial_{w_i} \otimes \bigwedge_i e^w_i \in \ddd(D_w|_{Y_Q^w}) \otimes \ddd(T_wQ/\Gamma)$$
and $\mfo_w \in \ddd(D_w|_{Y_\Gamma^w})$, while $C(U)$ corresponds to orientations 
$$\textstyle(-1)^{k_2}\bigwedge_i \partial_{U_i} \otimes \bigwedge_i e^U_i \in \ddd(D_U|_{Y_Q^U}) \otimes \ddd(T_UQ/\Gamma)$$
and $\mfo_U \in \ddd(D_U|_{Y_\Gamma^U})$. We have the isomorphism
$$\ddd(D_U|_{Y_\Gamma^U}) \otimes \ddd(D_w|_{Y_\Gamma^w}) \simeq \ddd(D_X|_{Y_\Gamma^X})$$
obtained by composing the direct sum isomorphism with deformation and gluing; let $\mfo_X$ be the image of $\mfo_U \otimes \mfo_w$ under this isomorphism. This isomorphism enters in the bottom vertical arrow of the following commutative diagram:
$$\xymatrix{\ddd(D_U|_{Y_\Gamma^U}) \otimes \ddd(D_{A_0} \sharp T\cS(q_0)) \otimes \ddd(D_w|_{Y_\Gamma^w}) \otimes  \ddd(D_{A_1} \sharp T\cS(q_1)) \ar[rd] \ar[d]^{R} &  \\
\ddd(D_U|_{Y_\Gamma^U}) \otimes \ddd(D_w|_{Y_\Gamma^w}) \otimes  \ddd(D_{A_0} \sharp T\cS(q_0)) \otimes \ddd(D_{A_1} \sharp T\cS(q_1)) \ar[r] \ar[d] & \ddd(D_{A_2}\sharp T\cS(q_2)) \\
\ddd(D_X|_{Y_\Gamma^X}) \otimes \ddd(D_{A_0} \sharp T\cS(q_0)) \otimes \ddd(D_{A_1} \sharp T\cS(q_1)) \ar[ur] & }$$
where the top and bottom right arrows correspond to the isomorphism $C(U) \circ (\id \otimes C(w))$, the top vertical arrow is the interchange of factors times the Koszul sign
$$(-1)^{\ind D_{A_0} \sharp T\cS(q_0)\cdot \ind D_w|_{Y_\Gamma^w}} = (-1)^{n-k}\,.$$
It follows that the isomorphism $C(U) \circ (\id \otimes C(w))$ corresponds to the orientation $(-1)^{n-k}\mfo_X$. The diagram \eqref{dia:induced_or_product_chain_map_leg_1} maps
$$\xymatrix{\mfo_U \otimes \bigwedge_i e^U_i \otimes \mfo_w \otimes \bigwedge_i e^w_i \otimes e^{x^1}_r  \ar@{|->}[r] \ar@{|->}[d] & (-1)^{k_2} \bigwedge_i\partial_{U_i} \otimes (-1)^{r+1}\bigwedge_i \partial_{w_i} \otimes e^{x^1}_r \ar@{|->}[d] \\
(-1)^{k_0+k_1+k_2}(\mfo_U \wedge \mfo_w) \otimes (\bigwedge_i e^U_i \wedge \bigwedge_i e^w_i) \otimes e^{x^1}_r  \ar@{|->}[r] \ar@{|->}[d] & (-1)^{k_2 + r + 1} (\bigwedge_i\partial_{U_i} \wedge \bigwedge_i \partial_{w_i}) \otimes e^{x^1}_r \ar@{|->}[d] \\
(-1)^{r(k_1+k_2) + k_0} \mfo_X \otimes (\bigwedge_i e^{x^0}_i \wedge \bigwedge_{i \neq r} e^{x^1}_i \wedge \bigwedge_i e^{x^2}_i) \otimes e^{x^1}_r  \ar@{|->}[r] \ar@{|->}[d] & (-1)^{k_2 + r + 1 + r(k_1+k_2)} \bigwedge_i\partial_{X_i} \otimes e^{x^1}_r \ar@{|->}[d] \\
(-1)^{(r+1)(k_1+k_2) + k_0} \mfo_X \otimes \bigwedge_i e^X_i  \ar@{|->}[r] & (-1)^{k_2 + r + 1 + r(k_1+k_2)} \bigwedge_i\partial_{X_i} \wedge \text{inward}_{\delta}}$$
whence it follows that $\mfo_X$, which corresponds to the isomorphism $(-1)^{n-k}C(U) \circ (\id \otimes C(w))$, corresponds to the orientation
$$\textstyle (-1)^{k_0 + k_1 + r + 1}\bigwedge_i\partial_{X_i} \wedge \text{inward}_\delta \otimes \bigwedge_i e^X_i\,.$$
It follows that the isomorphism $C(U) \circ (\id \otimes C(w))$ corresponds to the orientation
$$\textstyle (-1)^{n-k+ k_0 + k_1 + r + 1}\bigwedge_i\partial_{X_i} \wedge \text{inward}_\delta \otimes \bigwedge_i e^X_i\,,$$
as claimed. Note that by definition $C(\delta) = -(-1)^{n-k}C(U) \circ (\id \otimes C(w))$. This finishes the proof of Lemma \ref{lem:computation_induced_oris_product_is_chain_map}. \qed
\end{prf}

It remains to prove Lemma \ref{lem:sign_change_product_is_chain_map}.
\begin{prf}
We can assume without loss of generality that the two disks which touch in the degenerate triangle are separated by a positive length gradient trajectory in $V$, and therefore in $W$ the two disks are glued together into a single one. Assume first that at least one of the disks lies in leg $0$ of the triangle, while the other disk may either belong to the same leg or be the core. Therefore in $V$ the two disks are separated by a piece of gradient trajectory, and let us assume its number is $j$. Then we have the following commutative diagram:
\begin{equation}\label{dia:sign_change_product_is_chain_map_leg_0}
\xymatrix{\ddd(D_V|_{Y_\Gamma^V}) \otimes \ddd(T_VQ/\Gamma) \ar[r] & \ddd(D_V|_{Y_Q^V})\\
\ddd(D_V|_{Y_\Gamma^V}) \otimes \ddd(T_V^{v^0:\neg j}Q/\Gamma) \otimes \ddd(T_{v^0}^jQ/\Gamma) \ar[r] \ar[u] \ar[d] & \ddd(D_V|_{Y_Q^{V,v^0:\neg j}}) \otimes \ddd(T_{v^0}^jQ/\Gamma) \ar[u] \ar[d] \\
\ddd(D_W|_{Y_\Gamma^W}) \otimes \ddd(T_WQ/\Gamma) \otimes \ddd(T_{v^0}^jQ/\Gamma) \ar[r] & \ddd(D_W|_{Y_Q^W}) \otimes \ddd(T_{v^0}^jQ/\Gamma)}
\end{equation}
where the top square corresponds to the exact Fredholm square
$$\xymatrix{D_V|_{Y_\Gamma^V} \ar@{=}[r] \ar[d] & D_V|_{Y_\Gamma^V} \ar[r] \ar[d] & \ar[d] 0\\
D_V|_{Y_Q^{V,v^0:\neg j}} \ar[r] \ar[d] & D_V|_{Y_Q^V} \ar[r] \ar[d] & T_{v^0}^jQ/\Gamma \ar@{=}[d] \\
T_V^{v^0:\neg j}Q/\Gamma \ar[r] & T_VQ/\Gamma \ar[r] & T_{v^0}^jQ/\Gamma}$$
To obtain the bottom square, we note that, similarly to what we had during the proof of Lemma \ref{lem:induced_or_pearly_spaces_boundary_op_squared}, we have the exact triples
$$0 \to D_V|_{Y_\Gamma^V} \to D_V|_{Y_Q^{V,v^0:\neg j}} \to T_V^{v^0:\neg j}Q/\Gamma \to 0\,.$$
$$0 \to D_W|_{Y_\Gamma^W} \to D_W|_{Y_Q^W} \to T_WQ/\Gamma \to 0\,.$$
We have canonical isomorphisms
$$\ddd(D_V|_{Y_\Gamma^V}) \simeq \ddd(D_W|_{Y_\Gamma^W})\,,\; \ddd(D_V|_{Y_Q^{V,v^0:\neg j}}) \simeq \ddd(D_W|_{Y_Q^W})\,,\;\text{and}\;\ddd(T_V^{v^0:\neg j}Q/\Gamma) \simeq \ddd(T_WQ/\Gamma)\,,$$
obtained as follows. The pearly triangle $V$ can be deformed into the degenerate triangle, in which two disks are then glued to obtain $W$. During this process the operator $D_V|_{Y_\Gamma^V}$ undergoes deformation and linear gluing, with $D_W|_{Y_\Gamma^W}$ as the result; similarly, the space $T_V^{v^0:\neg j}Q/\Gamma$ deforms into $T_WQ/\Gamma$. If we let $V'$ be the degenerate pearly triangle, we have the corresponding linearized operator $D_{V'}|_{Y_Q^{V'}}$, which is surjective, and moreover the deformation of $V$ into $V'$ yields an isomorphism
$$\ker(D_V|_{Y_Q^{V,v^0:\neg j}}) \simeq \ker(D_{V'}|_{Y_Q^{V'}})\,.$$
In addition, the differential of the gluing map yields an isomorphism
$$\ker(D_{V'}|_{Y_Q^{V'}}) \simeq \ker(D_W|_{Y_Q^W})\,.$$
In total we obtain an isomorphism
$$\ddd(D_V|_{Y_Q^{V,v^0:\neg j}}) \simeq \ddd(D_W|_{Y_Q^W})\,.$$
It is then the feature of the gluing map that the following diagram commutes:
$$\xymatrix{\ddd(D_V|_{Y_\Gamma^V}) \otimes \ddd(T_V^{v^0:\neg j}Q/\Gamma) \ar[r] \ar[d] & \ddd(D_V|_{Y_Q^{V,v^0:\neg j}}) \ar[d] \\
\ddd(D_W|_{Y_\Gamma^W}) \otimes \ddd(T_WQ/\Gamma) \ar[r] & \ddd(D_W|_{Y_Q^W})}$$
with the horizontal arrows coming from the above exact triples and the vertical arrows being the above canonical isomorphisms.

Assume now that the isomorphism $C$ corresponds to orientations
$$\textstyle \mfo_V \in \ddd(D_V|_{Y_\Gamma^V})\,,\; \bigwedge_i \partial_{V_i} \wedge \eta_V \otimes \bigwedge_ie^V_i \in \ddd(T_V\wt\cP(q_0,q_1;q_2)) \otimes T_VQ/\Gamma\,,$$
and to orientations
$$\textstyle \mfo_W \in \ddd(D_W|_{Y_\Gamma^W})\,,\; \epsilon \bigwedge_i \partial_{W_i} \wedge \eta_W \otimes \bigwedge_ie^W_i \in \ddd(T_W\wt\cP(q_0,q_1;q_2)) \otimes T_WQ/\Gamma\,,$$
and let us compute the sign $\epsilon$. For the convenience of the computation, and without loss of generality, we assume that $\eta_V$ is directed from $V$ toward the degenerate triangle $V'$. Note that the isomorphism $\ddd(D_V|_{Y_\Gamma^V}) \simeq \ddd(D_W|_{Y_\Gamma^W})$ maps $\mfo_V \mapsto \mfo_W$. Then the diagram \eqref{dia:sign_change_product_is_chain_map_leg_0} maps:
$$\xymatrix{(-1)^{k_0+k_1+k_2-j}\mfo_V \otimes \bigwedge_i e^V_i \ar@{|->}[r] & (-1)^{k_0+k_1+k_2 - j} \bigwedge_i\partial_{V_i}\wedge \eta_V\\
\mfo_V \otimes \bigwedge_{i \neq j} e^{v^0}_i \wedge \bigwedge_i e^{v^1}_i \wedge \bigwedge_i e^{v_2}_i \otimes e^{v^0}_j \ar@{|->}[r] \ar@{|->}[u] \ar@{|->}[d] & (-1)^{k_0+k_1+k_2 - j} \bigwedge_i\partial_{V_i} \otimes e^{v^0}_j \ar@{|->}[u] \ar@{|->}[d]\\
\mfo_W \otimes \bigwedge_i e^W_i \otimes e^{v^0}_j \ar@{|->}[r] & -\bigwedge_i\partial_{W_i} \wedge \eta_W \otimes e^{v^0}_j}$$
The bottom arrow tells us that the orientation $\mfo_W$ corresponding to the isomorphism $C$ also corresponds to the orientation $-\bigwedge_i\partial_{W_i} \wedge \eta_W \otimes \bigwedge_ie^W_i$, therefore $\epsilon = -1$ as claimed.

The only computation that needs explanation is the right bottom arrow. The differential of the gluing map sends: $\partial_{v^{1,2}_i} \mapsto \partial_{w^{1,2}_i}$ for all $i$, $\partial_{v^0_i} \mapsto \partial_{w^0_i}$ for $i < j$. If $j < k_0$, it sends $\partial_{v^0_j} + \partial_{v^0_{j+1}} \mapsto \partial_{w^0_j}$ and $-\partial_{v^0_j} + \partial_{v^0_{j+1}} \mapsto \eta_W$, and $\partial_{v^0_i} \mapsto \partial_{w^0_{i-1}}$ for $i > j+1$, therefore
$$\textstyle \bigwedge_i\partial_{V_i} \mapsto -(-1)^{k_0+k_1+k_2-j}\bigwedge_i\partial_{W_i} \wedge \eta_W\,.$$
In case $j=k_0$, the gluing map sends $\partial_{v^0_{k_0}} \mapsto -\eta_W$, and therefore we have
$$\textstyle \bigwedge_i\partial_{V_i} \mapsto - (-1)^{k_0+k_1+k_2-j}\bigwedge_i\partial_{W_i} \wedge \eta_W\,.$$

The computation of $\epsilon$ in the case when the breaking/collision happens in leg $1$ of the triangle is entirely analogous. Let us therefore compute $\epsilon$ in case the breaking happens in leg $2$. Using identical arguments, we obtain the following commutative diagram, where we assume that the piece of gradient trajectory of $V$ which shrinks to $0$ at the collision point bears number $j$:
$$\xymatrix{\ddd(D_V|_{Y_\Gamma^V}) \otimes \ddd(T_VQ/\Gamma) \ar[r] & \ddd(D_V|_{Y_Q^V})\\
\ddd(D_V|_{Y_\Gamma^V}) \otimes \ddd(T_V^{v^2:\neg j}Q/\Gamma) \otimes \ddd(T_{v^2}^jQ/\Gamma) \ar[r] \ar[u] \ar[d] & \ddd(D_V|_{Y_Q^{V,v^2:\neg j}}) \otimes \ddd(T_{v^2}^jQ/\Gamma) \ar[u] \ar[d] \\
\ddd(D_W|_{Y_\Gamma^W}) \otimes \ddd(T_WQ/\Gamma) \otimes \ddd(T_{v^2}^jQ/\Gamma) \ar[r] & \ddd(D_W|_{Y_Q^W}) \otimes \ddd(T_{v^2}^jQ/\Gamma)}$$
It maps
$$\xymatrix{(-1)^{k_2-j}\mfo_V \otimes \bigwedge_i e^V_i \ar@{|->}[r] & (-1)^{k_2 - j} \bigwedge_i\partial_{V_i}\wedge \eta_V\\
\mfo_V \otimes \bigwedge_i e^{v^0}_i \wedge \bigwedge_i e^{v^1}_i \wedge \bigwedge_{i \neq j} e^{v_2}_i \otimes e^{v^2}_j \ar@{|->}[r] \ar@{|->}[d] \ar@{|->}[u] & (-1)^{k_2 - j} \bigwedge_i\partial_{V_i} \otimes e^{v^2}_j \ar@{|->}[d] \ar@{|->}[u] \\
\mfo_W \otimes \bigwedge_i e^W_i \otimes e^{v^2}_j \ar@{|->}[r] & \bigwedge_i\partial_{W_i} \wedge \eta_W \otimes e^{v^0}_j}$$
We see from the bottom arrow that the orientation $\mfo_W$ corresponding to the isomorphism $C$ also corresponds to the orientation $\bigwedge_i\partial_{W_i} \wedge \eta_W \otimes \bigwedge_ie^W_i$, therefore $\epsilon = 1$ as claimed.

Again, the bottom right arrow is obtained as follows. The differential of the gluing map sends: $\partial_{v^{0,1}_i} \mapsto \partial_{w^{0,1}_i}$ for all $i$, $\partial_{v^2_i} \mapsto \partial_{w^2_{i-1}}$ for $i > j$. If $j > 1$, it sends $\partial_{v^2_{j-1}} + \partial_{v^2_{j}} \mapsto \partial_{w^2_{j-1}}$ and $-\partial_{v^2_{j-1}} + \partial_{v^2_{j}} \mapsto \eta_W$, and $\partial_{v^2_i} \mapsto \partial_{w^2_{i}}$ for $i < j-1$, therefore
$$\textstyle \bigwedge_i\partial_{V_i} \mapsto (-1)^{k_2 - j}\bigwedge_i\partial_{W_i} \wedge \eta_W\,.$$
In case $j=1$, the gluing map sends $\partial_{v^2_1} \mapsto \eta_W$, and therefore we have
$$\textstyle \bigwedge_i\partial_{V_i} \mapsto (-1)^{k_2-j}\bigwedge_i\partial_{W_i} \wedge \eta_W\,.$$
The proof of the lemma is now complete.
\qed
\end{prf}

We have thus completed the definition of the product on quantum homology.

\subsubsection{Unit}\label{sss:unit_Lagr_QH}

Assume we have a regular quantum datum $\cD = (f,\rho,J)$. We will now construct an element in $QH_n(\cD:L)$ which serves as a unit for the quantum product. This is defined as follows. Let $q$ be a maximum of $f$. Then the stable manifold $\cS(q)$ is just the singleton $\{q\}$. Let $0 \in \pi_2(M,L,q)$ be the zero class and consider the operator $D_0 \sharp T\cS(q)$. We can canonically orient this operator, as follows. Let $w$ be the constant disk at $q$; it clearly represents the class $0$. The corresponding Cauchy--Riemann operator $D_w$ is just the standard Doulbeault operator on the trivial bundle pair $(T_qM,T_qL) \to (D^2,S^1)$ with the Hermitian structure $(\omega_q,J_q)$. This operator is surjective and its kernel consists of constant sections with values in $T_qL$. Now the operator $D_w \sharp T\cS(q)$ is its restriction to the subspace of sections of this bundle pair vanishing at $1 \in D^2$. Clearly this restricted operator has trivial kernel, and since it has index $0$, it must be an isomorphism. Therefore we have the canonical positive orientation of this operator by $1 \otimes 1^\vee$, which induces an orientation on the family $D_0 \sharp T\cS(q)$, and therefore an element of $C(q,0)$, which we denote $1_q$. The \tb{unit} is now the element
$$1_\cD = \sum_{\substack{q \in \Crit f:\\|q| = n}}1_q \in QC_n(\cD:L)\,.$$
We claim that this element is a cycle. Indeed, the Morse boundary operator vanishes on it, and from index considerations in $\partial_\cD(1_\cD)$ there are no terms involving nonconstant holomorphic disks.

To see that the class $1_\cD \in QH_n(\cD:L)$ is the unit with respect to the quantum product, we choose $f_1 = f_2 = f'$ (this is possible \cite{Biran_Cornea_Quantum_structures_Lagr_submfds, Biran_Cornea_Rigidity_uniruling}), that is the same function for the first and second slots, and show that the map
$$1_\cD \star - \fc QC_k (f',\rho,J:L) \to QC_k(f',\rho,J:L)$$
is the identity, so that in fact we have identity already on chain level. First we claim that any pearly triangle belonging to the zero-dimensional part of $\cP(q,q';q'')$ where $q \in \Crit f$, $q',q'' \in \Crit f'$ must be constant, meaning that the only holomorphic disk is the core, that it is the constant disk, and that the lengths of all the gradient trajectories are zero. Indeed, let $U$ be such a triangle. First we will show that the zeroth leg contains no disks. If there were a disk, it would be nonconstant, and since lying in the unstable manifold of $q$ is an open condition, after quotienting out the action of the automorphism group of the triange, we would obtain a positive-dimensional space, which is a contradiction. The same argument shows that the core must be constant, which means that we can obtain a pearly trajectory from $q'$ to $q''$, and that it has index zero. This is only possible if there are no nonconstant disks and that the resulting gradient trajectory from $q'$ to $q''$ is constant, which is what we claimed.

Next, for a generically chosen $f'$, all of its critical points lie in the union of the unstable manifolds of the maxima of $f$, which means that all the constant pearly triangles indeed appear and that we obtain the identity map on $QC_k(f',\rho,J:L)$ as a result of multiplying by $1_\cD$, proving that it indeed acts as the unit.

\subsection{Arbitrary rings and twisted coefficients}\label{ss:arbitrary_rings_loc_coeffs_Lagr_QH}

Analogously to the case of Floer homology, we can use an arbitrary ground ring $R$ if $L$ satisfies assumption \tb{(O)}, or otherwise use a ring in which $2=0$. Also, given a flat $R$-bundle $\cE$ over $\wt\Omega_L$ we can define quantum homology of $L$ twisted by $\cE$, which we denote $QH_*(\cD:L;\cE)$. We only need to note that the pairs $(q,A)$ with $q \in L$ and $A \in \pi_2(M,L,q)$ naturally give rise to points of the space $\wt\Omega_L$. Details are left to the reader.

\subsection{Duality}\label{ss:duality_Lagr_QH}

The treatment of duality in quantum homology is very similar to the case of Floer homology, and we follow it closely. The goal here is to establish a canonical chain isomorphism
$$QC_*(\ol\cD:L) \equiv QC^{n-*}(\cD:L;\cL)\,,$$
where $\ol\cD$ is the dual quantum datum, $\cL \to \wt\Omega_L$ is the flat $\Z$-bundle obtained as the normalization of the pullback of $\ddd(TL)$ to $\wt\Omega_L$ via the evaluation map $\wt\gamma \mapsto \gamma(0)$.
First we define quantum cohomology. Fix a regular quantum homology datum $\cD = (f,\rho,J)$ and define
$$C(q,A)^\vee = \Hom_\Z(C(q,A),\Z)\,.$$
Put
$$QC^*(\cD:L) = \bigoplus_{\substack{q \in \Crit f \\ A \in \pi_2(M,L,q)}}C(q,A)^\vee\,.$$
This is graded by assigning the elements of $C(q,A)^\vee$ the degree $|q| - \mu(A)$. The matrix coefficients of the dual differential $\partial_\cD^\vee$ are dual to the matrix elements of $\partial_\cD$ as maps $C(q',A')^\vee \to C(q,A)^\vee$. We define another differential
$$\delta_\cD \fc QH^k(\cD:L) \to QH^{k+1}(\cD:L) \quad\text{via} \quad \delta_\cD = (-1)^{k-1}\partial_\cD^\vee\,.$$
The cochain complex
$$(QC^*(\cD:L),\delta_\cD)$$
is the \tb{quantum cochain complex} and its cohomology is the \tb{quantum cohomology}
$$QH^*(\cD:L)\,.$$

The dual quantum datum is defined by $\ol\cD=(-f,\rho,J)$. The functions $f, -f$ have the same critical points. For $w \in C^\infty(D^2,S^1,1;M,L,q)$ we let $\ol w$ be defined as $\ol w(\sigma,\tau) = w(\sigma,-\tau)$. Then $\ol w$ represents the class $[w]^{-1} \in \pi_2(M,L,q)$. We also define $-w \in C^\infty(D^2,S^1,-1;M,L,q)$ via $-w(\sigma,\tau) = w(-\sigma,\tau)$. Clearly $-w = \ol w \circ \phi$ where $\phi \fc D^2 \to D^2$ is defined by $\phi(z) = -z$. This map is a conformal isomorphism therefore it induces an isomorphism of determinant lines $\ddd(D_{\ol w}) = \ddd(D_{-w})$. We also have the isomorphisms
\begin{equation}\label{eqn:iso_det_lines_duality_Lagr_QH}
\ddd(D_w \sharp T\cS_f(q)) \otimes \ddd(D_{-w}\sharp T\cU_f(q)) \simeq \ddd(D_w \sharp T\cS_f(q) \oplus D_{-w} \sharp T\cU_f(q)) \simeq \ddd(D_w \sharp D_{-w}) \simeq \ddd(T_qL)\,,
\end{equation}
where the second isomorphism comes from deforming the incidence condition at $q$ from
$$(T\cS_f(q)\oplus 0) \oplus (0 \oplus T\cU_f(q)) \subset T_qL \oplus T_qL \quad \text{to} \quad \Delta_{T_qL}\,,$$
and the third isomorphism comes from deforming the operator $D_w \sharp D_{-w}$ into the operator $D_0$ and using the canonical isomorphism $\ddd(D_0) = \ddd(T_qL)$.

Since we have $\cU_f(q) = \cS_{-f}(q)$, this yields, combined with the isomorphism \eqref{eqn:iso_det_lines_duality_Lagr_QH}, the following:
$$\ddd(D_w \sharp T\cS_f(q)) \otimes \ddd(D_{\ol w} \sharp T\cS_{-f}(q)) \simeq \ddd(T_qL)\,.$$
This means that we have a canonical isomorphism
$$C_f(q,A) \otimes C_{-f}(q,A^{-1}) \otimes \cL_q = \Z\,.$$
This implies that we have canonically
$$QC_*(\ol\cD:L) = QC^{n-*}(\cD:L;\cL)$$
as modules, where we observe that the elements of $C_{-f}(q,A^{-1})$ have degree
$$|q|_{-f} - \mu(A^{-1}) = n-|q|_f + \mu(A)\,,$$
which is $n$ minus the degree of $C_f(q,A)$. We will now obtain an identification of the differentials. Fix pairs $(q_\pm,A_\pm)$ so that $|q_-| - \mu(A_-) = |q_+| - \mu(A_+) + 1$, and let $u \in \wt\cP(q_-,q_+)$ be such that $A_- \sharp u = A_+$. Represent the classes $A_\pm$ by maps $w_\pm$. We have the following commutative diagram, obtained by employing the direct sum, gluing, and deformation isomorphisms:
$$\xymatrix{\ddd(D_{w_-}\sharp T\cS_f(q_-)) \otimes \ddd(D_u|_{Y_\Gamma^u}) \otimes \ddd(D_{-w_+} \sharp T\cS_{-f}(q_+)) \ar[r] \ar[d]^{(C(u) \otimes \id)\circ (R\otimes \id)}& \ddd(D_{w_-}\sharp T\cS_f(q_-)) \otimes \ddd(D_{-w_-} \sharp T\cS_{-f}(q_-)) \ar[d] \\ \ddd(D_{w_+}\sharp T\cS_f(q_+)) \otimes \ddd(D_{-w_+} \sharp T\cS_{-f}(q_+) ) \ar[r] & \ddd(TL)}$$
where $R$ is the interchange of factors including the Koszul sign
$$(-1)^{\ind D_u|_{Y_\Gamma^u} \cdot \ind D_{w_-}\sharp T\cS_f(q_-)} = (-1)^{n-|q_-|_{f}+\mu(A_-)}\,.$$
Fix $\mfo \in \ddd(TL)$ where we trivialize $\ddd(TL)$ along the lower boundary of $u$, viewed as a degenerate strip with boundary on $L$. Fix $\mfo_{q_-} \in C_f(q_-,A_-)$, let $\mfo_{q_+} = C(u)(\mfo_{q_-}) \in C_f(q_+,A_+)$ and let $\mfo_{-q_\pm} \in \ddd(D_{-w_\pm} \sharp T\cS_{-f}(q_\pm))$ be such that $\mfo_{q_\pm} \otimes \mfo_{-q_\pm} \mapsto \mfo$ via the isomorphism \eqref{eqn:iso_det_lines_duality_Lagr_QH}. The diagram then maps
$$\xymatrix{(-1)^{n - |q_-| + \mu(A_-)}\mfo_{q_-} \otimes \mfo_u \otimes \mfo_{-q_+} \ar@{|->}[r] \ar@{|->}[d]  & \mfo_{q_-} \otimes \mfo_{-q_-} \ar@{|->}[d] \\ \mfo_{q_+} \otimes \mfo_{-q_+} \ar@{|->}[r] & \mfo}$$
Dualizing in a manner similar to the treatment of duality in Floer homology, see \S\ref{ss:duality_HF}, we obtain the diagram
$$\xymatrix{\ddd(D_{-\ol w_-}\sharp T\cS_f(q_-)) \otimes \ddd(D_{\ol u}|_{Y_\Gamma^{\ol u}}) \otimes \ddd(D_{\ol w_+} \sharp T\cS_{-f}(q_+)) \ar^-{\id \otimes C(\ol u)}[r] \ar[d]& \ddd(D_{- \ol w_-}\sharp T\cS_f(q_-)) \otimes \ddd(D_{\ol w_-} \sharp T\cS_{-f}(q_-)) \ar[d] \\ \ddd(D_{-\ol w_+}\sharp T\cS_f(q_+)) \otimes \ddd(D_{\ol w_+} \sharp T\cS_{-f}(q_+)) \ar[r] & \ddd(TL)}$$
Let $\ol \mfo_{q_\pm} \in \ddd(D_{\ol w_\pm} \sharp T\cS_{-f}(q_\pm))$, $\ol \mfo_{-q_\pm} \in \ddd(D_{-\ol w_\pm} \sharp T\cS_f(q_\pm))$ be obtained from $\mfo_{-q_\pm}$, $\mfo_{q_\pm}$, respectively, by dualization. The latter diagram then maps
$$\xymatrix{(-1)^{n - |q_-| + \mu(A_-)} \ol \mfo_{-q_-} \otimes \mfo_{\ol u} \otimes \ol \mfo_{q_+} \ar@{|->}[r] \ar@{|->}[d]  & \ol\mfo_{-q_-} \otimes \ol\mfo_{q_-} \ar@{|->}[d]\\
\ol\mfo_{-q_+} \otimes \ol\mfo_{q_+} \ar@{|->}[r]  & \mfo   }$$
Thus we see that $C(\ol u) (\ol \mfo_{q_+}) = (-1)^{n - |q_-| + \mu(A_-)} \ol \mfo_{q_-}$. This means that the following diagram commutes:
$$\xymatrix{C_{-f}(q_+,A_+^{-1})  \ar[r]  \ar[d]^{C(\ol u)} & C_f(q_+,A_+)^\vee \otimes \cL_{q_+} \ar[d]^{(-1)^{n - |q_-| + \mu(A_-)}C(u)^\vee \otimes \cP}\\ C_{-f}(q_-,A_-^{-1}) \ar[r] & C_f(q_-,A_-)^\vee \otimes \cL_{q_-}}$$
where $\cP$ is the parallel transport isomorphism on the line bundle $\cL$, see \S\ref{ss:duality_HF}. This means that we have established a canonical isomorphism of chain complexes:
$$(QC_*(\ol\cD:L),\partial_{\ol \cD}) = (QC^{n-*}(\cD:L;\cL),\delta_{\cD} \otimes \cP)\,.$$
It therefore induces the \tb{duality isomorphism} on homology:
\begin{equation}\label{eqn:duality_Lagr_QH}
QH_*(\ol\cD:L) = QH^{n-*}(\cD:L;\cL)\,.
\end{equation}

\subsubsection{Augmentation}

Similarly to the case of Floer homology, we can view the unit as a graded map
$$1 \fc \Z[n] \to QH_*(\ol\cD:L)\,.$$
The duality isomorphism \eqref{eqn:duality_Lagr_QH} means that we obtain a graded map
$$\Z[n] \to QH^{n-*}(\cD:L;\cL)\,,$$
and by dualizing we obtain
$$QH_*(\cD:L;\cL) \to \Z\,,$$
which is the \tb{augmentation map}.

\subsection{Quantum homology of $M$}\label{ss:QH_of_M}

We call a triple $\cD = (f,\rho,J)$ a quantum datum for $M$ if $(f,\rho)$ is a Morse--Smale pair on $M$ and $J$ is an $\omega$-compatible almost complex structure on $M$. We call it regular if the various moduli spaces below are transversely cut out. This is the case for a generic $J$.

\subsubsection{Generators, the complex as a module, and the boundary operator}\label{sss:generators_cx_bd_op_QH_of_M}

Fix a regular quantum datum $\cD = (f,\rho,J)$ for $M$. For a critical point $q \in \Crit f$ and a homotopy class $A \in \pi_2(M,q)$ we can construct the family of operators $D_A$, just like in the Lagrangian case. Members of $D_A$ are formal linearized operators $D_u$ of smooth maps $u \fc (S^2,1) \to (M,q)$ in class $A$ with respect to some auxiliary connection $\nabla$ on $M$. We have the following foundational lemma.
\begin{lemma}\label{lem:D_A_has_canonical_orientation}
The family $D_A$ possesses a canonical orientation.
\end{lemma}
\begin{prf}
This is a family of real Cauchy--Riemann operators on a closed Riemann surface. The set of real Cauchy--Riemann operators retracts onto the subset of complex linear operators, which have canonically oriented determinant lines. It follows that the determinant line bundle of the set of real Cauchy--Riemann operators, and therefore the determinant line of $D_A$, is canonically oriented. \qed
\end{prf}

We also have the family $D_A \sharp T\cS(q)$, defined analogously to the Lagrangian case. Using an exact Fredholm triple analogous to the ones appearing in the proof of Lemma \ref{lem:family_ops_orientable_def_cx_QH}, we see that this family is orientable as well, and thus we can define
$$C(q,A)$$
to be the rank $1$ free abelian group whose two generators are its two possible orientations. We let
$$QC_*(\cD) = \bigoplus_{\substack{q \in \Crit f \\ A \in \pi_2(A,q)}} C(q,A)\,.$$
This is graded by assigning the elements of $C(q,A)$ the degree $|q| - 2c_1(A)$. The boundary operator $\partial_\cD \fc QC_*(\cD) \to QC_{*-1}(\cD)$ is just the ordinary Morse boundary operator, enhanced by the homotopy classes. More precisely, it is defined as follows. We first define the Morse boundary operator. When $A = 0 \in \pi_2(M,q)$, we have canonically $\ddd(D_A) = \ddd(T_qM)$, and therefore canonically $\ddd(D_A \sharp T\cS(q)) = \ddd(T\cS(q))$. For a Morse trajectory $u \in \wt\cM(q,q')$ of index $1$ we have a natural exact sequence, see \eqref{eqn:exact_seq_moduli_space_gradient_lines_stable_unstable}:
$$0 \to \R \partial_u \to T\cS(q') \to T\cS(q) \to 0$$
whence
$$\ddd(T\cS(q')) \simeq \ddd(\R \partial_u) \otimes \ddd(T\cS(q))\,.$$
Substituting the positive orientation of $\R$ we get the isomorphism
$$C(u) \fc C(q,0) \simeq C(q',0)\,.$$
For general $A$, this isomorphism induces an isomorphism
$$C(u) \fc C(q,A) \simeq C(q',A')$$
where $A' \in \pi_2(M,q')$ is obtained by transferring the class $A$ to $q'$ along $u$. Indeed, from the relations
$$\ddd(D_A) \otimes \ddd(T\cS(q)) \simeq \ddd(D_A \sharp T\cS(q)) \otimes \ddd(T_qM)$$
$$\ddd(D_{A'}) \otimes \ddd(T\cS(q')) \simeq \ddd(D_{A'} \sharp T\cS(q')) \otimes \ddd(T_{q'}M)$$
we see that to induce such an isomorphism, it suffices to produce isomorphisms $\ddd(D_A) \simeq \ddd(D_{A'})$ and $\ddd(T_qM) \simeq \ddd(T_{q'}M)$. The former is obtained by deformation induced by moving the base point along $u$, while the latter comes from the fact that $M$ is oriented. The boundary operator is then given by its matrix elements. The matrix element between $(q,A)$ and $(q',A')$ where $|q| = |q'| + 1$ and $c_1(A') = c_1(A)$ is given by the sum
$$\sum_{[u] \in \cM(q,q'): A\sharp u = A'}C(u) \fc C(q,A) \to C(q',A')\,.$$
\begin{thm}\label{thm:boundary_op_squares_zero_QH_of_M}
The boundary operator satisfies $\partial_\cD^2 = 0$.
\end{thm}
\noindent We can therefore define the quantum homology of $M$, $QH_*(\cD)$, as the homology of $(QC_*(\cD),\partial_\cD)$.

\begin{prf}[of Theorem \ref{thm:boundary_op_squares_zero_QH_of_M}]This immediately follows from the parallel proof in Morse theory, coupled with the observation that transferring a class $A \in \pi_2(M,q)$ along paths ending at $q''$, which are homotopic with fixed endpoints, yields the same class $A'' \in \pi_2(M,q'')$. In the Morse-theoretic proof one uses the compactified $1$-dimensional moduli space of gradient trajectories, and transferring a class along the trajectories comprising the two boundary points of a connected component of it yields the same class at the other critical point, therefore the proof goes through, even though the quantum boundary operator distinguishes homotopy classes of spheres. \qed
\end{prf}

\subsubsection{Product}\label{sss:product_QH_of_M}

To define the product we need to define moduli spaces of spiked spheres. Fix data quantum data $\cD_i = (f_i,\rho,J)$, $i=0,1,2$. Let $\wt\cM^\circ(J)$ be the space of parametrized $J$-holomorphic spheres in $M$, including constant ones. We have the evaluation map
$$\ev \fc C^\infty(S^2,M) \to M^3\,, \quad u \mapsto (u(0),u(1),u(\infty))\,,$$
where we view $S^2 = \C P^1$. The space of spiked spheres is
$$\cP(q_0,q_1;q_2) = \ev^{-1}(\cU(q_0) \times \cU(q_1) \times \cS(q_2)) \cap \wt\cM^\circ(J)\,.$$
For generic $J$ it is a smooth manifold of local dimension at $u$
$$|q_0| + |q_1| - |q_2| +2c_1(u) - 2n$$
provided this number is $\leq 1$. For $u$ with $\ev(u) \in \cU(q_0) \times \cU(q_1) \times \cS(q_2)$ define the space
$$Y^u = \{\xi \in W^{1,p}(u) \,|\, \xi(0) \in T_{u(0)}\cU(q_0),\xi(1) \in T_{u(1)}\cU(q_1),\xi(\infty) \in T_{u(\infty)}\cS(q_2)\}\,.$$
Assume $\dim_u \cP(q_0,q_1;q_2) = 0$. Then the linearized operator $D_u|_{Y^u}$ has index zero and is surjective, therefore it possesses the canonical positive orientation $\mfo_u = 1\otimes 1^\vee$. Deform $u$ into $u' \in C^\infty(S^2,M)$ satisfying $\ev(u') = (q_0,q_1,q_2)$. Then we have a canonical isomorphism
$$\ddd(D_u|_{Y^u}) \simeq \ddd(D_{u'}|_{Y^{u'}})\,.$$
Using direct sum, deformation, and linear gluing isomorphisms, combined with arguments involving deformation of incidence conditions at $q_0,q_1$ (see the proof of Lemma \ref{lem:computation_induced_oris_product_is_chain_map} of the Lagrangian case), we get an isomorphism
$$\ddd(D_{u'}|_{Y^{u'}}) \otimes \ddd(D_{A_0}\sharp T\cS(q_0)) \otimes \ddd(D_{A_1}\sharp T\cS(q_1))\simeq \ddd(D_{A_2}\sharp T\cS(q_2))\,,$$
where $A_2 = A_0 \sharp A_1 \sharp u$, and composing this with the deformation isomorphism $\ddd(D_u|_{Y^u}) \simeq \ddd(D_{u'}|_{Y^{u'}})$, we get a bijection between orientations of $D_u|_{Y^u}$ and isomorphisms $C(q_0,A_0) \otimes C(q_1,A_1) \simeq C(q_2,A_2)$. Thus the standard orientation $\mfo_u \in \ddd(D_u|_{Y^u})$ gives rise to the isomorphism
$$C(u) \fc C(q_0,A_0) \otimes C(q_1,A_1) \simeq C(q_2,A_2)\,.$$
The matrix element of the product is then
$$\sum_{\substack{u \in \cP(q_0,q_1;q_2): \\ A_0\sharp A_1\sharp u = A_2}} C(u) \fc C(q_0,A_0) \otimes C(q_1,A_1) \simeq C(q_2,A_2)\,.$$
We have thus defined a bilinear operation
$$* \fc QC_k(\cD_0) \otimes QC_l(\cD_1) \to QC_{k+l-2n}(\cD_2)\,.$$

\begin{thm}
The operation $*$ is a chain map. More precisely
$$\partial_{\cD_2} \circ * = * \circ (\partial_{\cD_0} \otimes \id + (-1)^{2n-k} \id \otimes \partial_{\cD_1}) \fc QC_k(\cD_0) \otimes QC_l(\cD_1) \to QC_{k+l-2n-1}(\cD_2)\,.$$
\end{thm}
\begin{prf}
This is proved in complete analogy with the Lagrangian case, except that there is no holomorphic curve breaking or collision to take into account. The $1$-dimensional moduli space $\cP(q_0,q_1;q_2)$ can be compactified by adding Morse breaking. One then computes the induced orientations on $\cP$ coming from isomorphisms corresponding to the constituent trajectories of boundary points to conclude. \qed
\end{prf}

\subsubsection{Unit}\label{sss:unit_QH_of_M}

This is defined analogously to the Lagrangian case. Let $\cD = (f,\rho,J)$ be a regular quantum datum. Like in the Lagrangian case, we have the elements $1_q \in C(q,0)$ for every maximum $q \in \Crit f$ where $0 \in \pi_2(M,q)$ is the zero class. Their sum
$$1_\cD = \sum_{\substack{q \in \Crit f: \\ |q| = 2n}}1_q \in QC_{2n}(\cD)$$
is the unit. Again, like in the Lagrangian case, one checks that $1_\cD$ actually is a unit on chain level, and therefore in homology.

\subsubsection{Quantum module action}\label{sss:quantum_module_action_QH}

Fix a compatible almost complex structure $J$ and quantum data $\cD_i = (f_i,\rho,J)$, $i=0,1$ for $L$, and a quantum datum $\cD = (f,\rho',J)$ for $M$, such that all of them are regular. The \tb{quantum module action} is a bilinear operation
$$\bullet \fc QC_k(\cD) \otimes QC_l(\cD_0:L) \to QC_{k+l-2n}(\cD_1:L)\,.$$
This is defined via its matrix elements which are homomorphisms
$$C(q,A) \otimes C(q_0,A_0) \to C(q_1,A_1)$$
for $q \in \Crit f,q_i \in \Crit f_i$, $A \in \pi_2(M,q)$, $A_i \in \pi_2(M,L,q_i)$, such that $|q| + |q_0| - |q_1| - 2c_1(A) - \mu(A_0) + \mu(A_1) - 2n = 0$. In order to define these, we need to describe additional pearly moduli spaces. Fix $q_i \in \Crit f_i$ and $q \in \Crit f$. Recall the spaces $\wt\cM^\circ(L,J)$, $\wt\cM(L,J)$ of $J$-holomorphic disks with boundary on $L$, respectively nonconstant $J$-holomorphic disks. For $k_0,k_1 \geq 0$ we have the evaluation map
$$\ev \fc (C^\infty(D^2,S^1;M,L))^{k_0+k_1+1} \to L^{2(k_0+k_1)+3}$$
defined by
\begin{multline*}
\ev(U = (u^0,u^1,u)) = (u(0);u^0_1(-1),u^0_1(1),\dots,u^0_{k_0}(-1),u^0_{k_0}(1),u(-1);\\ u(1),u^1_1(-1),u^1_1(1),\dots,u^1_{k_1}(-1),u^1_{k_1}(1))\,,
\end{multline*}
where $u^i \in (C^\infty(D^2,S^1;M,L))^{k_i}$. We let
\begin{multline*}
\wt\cP_{k_0,k_1}(q,q_0;q_1) = \ev^{-1}(\cU_f(q) \times \cU_{f_0}(q_0) \times Q_{f_0,\rho}^{k_0} \times Q_{f_1,\rho}^{k_1} \times \cS_{f_1}(q_1)) \cap \\ (\wt\cM(L,J))^{k_0} \times (\wt\cM(L,J))^{k_1} \times \wt\cM^\circ(L,J)\,. 
\end{multline*}
For the sake of convention we call this space the space of spiked pearls. There is a natural $\R^{k_0+k_1}$-action on this space and we let $\cP_{k_0+k_1}(q,q_0;q_1)$ be the quotient. We also define
$$\wt\cP(q,q_0;q_1) = \bigcup_{k_0,k_1 \geq 1}\wt\cP_{k_0,k_1}(q,q_0;q_1) \quad \text{ and }\quad \cP(q,q_0;q_1) = \bigcup_{k_0,k_1 \geq 1} \cP_{k_0,k_1}(q,q_0;q_1)\,.$$
We have \cite{Biran_Cornea_Quantum_structures_Lagr_submfds, Biran_Cornea_Rigidity_uniruling}:
\begin{prop}
Fix Morse-Smale pairs $(f_i,\rho)_{i=0,1}$ and $(f,\rho')$. Then there is a subset of $\cJ(M,\omega)$ of the second category so that for each $J$ in this subset, for all $q_i \in \Crit f_i, q \in \Crit f$ the space $\cP(q,q_0;q_1)$ is a smooth manifold of local dimension at $[U]$ equal to
$$|q| + |q_0| - |q_1| + \mu(U) - 2n$$
provided this number is at most zero. \qed
\end{prop}

We proceed with the definition of the matrix element. Fix $A \in \pi_2(M,q)$ and $A_0 \in \pi_2(M,L,q_0)$. There is an obvious way to produce a class $A\sharp A_0 \sharp U \in \pi_2(M,L,q_1)$ for $U \in (C^\infty(D^2,S^1;M,L))^{k_0+k_1+1}$ with $\ev(U) \in \cU_f(q) \times \cU_{f_0}(q_0) \times Q_{f_0,\rho_0}^{k_0} \times Q_{f_1,\rho_1}^{k_1} \times \cS_{f_1}(q_1)$ by concatenating representatives of $A,A_0$ with the constituent disks of $U$ along the pieces of gradient trajectories connecting the evaluation points of the disks. For $U \in \wt\cP(q,q_0;q_1)$ satisfying $|q| + |q_0| - |q_1| + \mu(U) - 2n = 0$ we will construct an isomorphism
$$C(U) \fc C(q,A) \otimes C(q_0,A_0) \simeq C(q_1,A\sharp A_0 \sharp U)\,.$$
For $A_1 \in \pi_2(M,L,q_1)$ with $|q| + |q_0| - |q_1| - 2c_1(A) - \mu(A_0) + \mu(A_1) - 2n = 0$, the matrix element of $\bullet$ is then the sum
$$\sum_{\substack{[U] \in \cP(q,q_0;q_1): \\ A\sharp A_0 \sharp U = A_1}} C(U) \fc C(q,A) \otimes C(q_0,A_0) \simeq C(q_1,A_1)\,.$$
We therefore have to define the isomorphism $C(U)$. We need some additional spaces. For $U$ with $\ev(U) \in \cU_f(q) \times \cU_{f_0}(q_0) \times \ol Q{}_{f_0,\rho}^{k_0} \times \ol Q{}_{f_1,\rho}^{k_1} \times \cS_{f_1}(q_1)$ let us define
\begin{multline*}
Y_\Gamma^U = \{\Xi =(\xi^0,\xi^1,\xi) \in W^{1,p}(u^0) \oplus W^{1,p}(u^1) \oplus W^{1,p}(u)\,|\,\\
\xi(0) \in T_{u(0)}\cU_{f}(q); (\xi^0_j(1),\xi^0_{j+1}(-1)) \in \Gamma_{(u^0_j(1),u^0_{j+1}(-1))}\text{ for }j < k_0;\\
(\xi^1_j(1),\xi^1_{j+1}(-1)) \in \Gamma_{(u^1_j(1),u^1_{j+1}(-1))}\text{ for }j \geq 1;\\
\xi^0(-1) \in T_{u^0(-1)}\cU_{f_0}(q_0), (\xi^0_{k_0},\xi(-1)) \in \Gamma_{(u^0_{k_0}(1),u(-1))};\\
\xi^1_{k_1}(1) \in T_{u^1_{k_1}(1)}\cS_{f_1}(q_1),(\xi(1),\xi^1_1(-1)) \in \Gamma_{(u(1),u^1_{k_1}(-1))}\}
\end{multline*}
and for $U$ with $\ev(U) \in \cU_f(q) \times \cU_{f_0}(q_0) \times Q_{f_0,\rho}^{k_0} \times Q_{f_1,\rho}^{k_1} \times \cS_{f_1}(q_1)$ let us define
\begin{multline*}
Y_Q^U = \{\Xi =(\xi^0,\xi^1,\xi) \in W^{1,p}(u^0) \oplus W^{1,p}(u^1) \oplus W^{1,p}(u)\,|\,\\
\xi(0) \in T_{u(0)}\cU_{f}(q); (\xi^0_j(1),\xi^0_{j+1}(-1)) \in T_{(u^0_j(1),u^0_{j+1}(-1))}Q_{f_0,\rho} \text{ for }j < k_0;\\
(\xi^1_j(1),\xi^1_{j+1}(-1)) \in T_{(u^1_j(1),u^1_{j+1}(-1))}Q_{f_1,\rho} \text{ for }j \geq 1;\\
\xi^0(-1) \in T_{u^0(-1)}\cU_{f_0}(q_0), (\xi^0_{k_0},\xi(-1)) \in T_{(u^0_{k_0}(1),u(-1))}Q_{f_0,\rho};\\
\xi^1_{k_1}(1) \in T_{u^1_{k_1}(1)}\cS_{f_1}(q_1),(\xi(1),\xi^1_1(-1)) \in T_{(u(1),u^1_{k_1}(-1))}Q_{f_1,\rho}\}
\end{multline*}
The isomorphism $C(U)$ is defined in an entirely similar fashion to the isomorphism entering the definition of the Lagrangian quantum product. Namely, we construct two canonical bijections: the first one is between orientations of $D_U|_{Y_\Gamma^U}$ and orientations of $\ddd(T_U\wt\cP(q,q_0;q_1)) \otimes \ddd(T_UQ/\Gamma)$, while the second one is between orientations of $D_U|_{Y_\Gamma^U}$ and isomorphisms $C(q,A) \otimes C(q_0,A_0) \simeq C(q_1,A_1)$.

The first bijection comes from the exact triple
$$0 \to D_U|_{Y_\Gamma^U} \to D_U|_{Y_Q^U} \to 0_{T_UQ/\Gamma} \to 0$$
which yields the isomorphism
$$\ddd(D_U|_{Y_Q^U}) \simeq \ddd(D_U|_{Y_\Gamma^U}) \otimes \ddd(T_UQ/\Gamma)\,,$$
whence the desired bijection.

The construction of the second bijection follows the same lines as for the Lagrangian product. Namely, we can deform $U$ into $U' \in (C^\infty(D^2,S^1;M,L))^{k_0+k_1+1}$ with $\ev(U') \in \{q\} \times \{q_0\} \times \Delta_L^{k_0+k_1} \times \{q_1\}$, which induces an isomorphism $\ddd(D_U|_{Y_\Gamma^U}) \simeq \ddd(D_{U'}|_{Y_\Gamma^{U'}})$. Then we can deform the operator family $D_A \sharp T\cS_f(q) \oplus D_{A_0}\sharp T\cS_{f_0} \oplus D_{U'}|_{Y_\Gamma^{U'}}$ by changing the incidence conditions at $q,q_0$ into an operator \footnote{The space $X_\Gamma^{U}$ is defined just like $Y_\Gamma^U$ but with no condition on $\xi(0)$ and $\xi^1_1(-1)$.} $D_A \sharp D_{A_0} \sharp D_{U'}|_{X_\Gamma^{U'}}$, which upon gluing and deformation yields a representative of the family $D_{A_1} \sharp T\cS_{f_1}(q_1)$. Together with the direct sum isomorphisms, we obtain the isomorphism
$$\ddd(D_{U}|_{Y_\Gamma^{U}}) \otimes \ddd(D_A \sharp T\cS_f(q)) \otimes \ddd(D_{A_0} \sharp T\cS_{f_0})  \simeq \ddd(D_{A_1} \sharp T\cS_{f_1}(q_1))\,,$$
which shows that indeed there is a bijection between orientations of $D_{U}|_{Y_\Gamma^{U}}$ and isomorphisms $C(q,A) \otimes C(q_0,A_0) \simeq C(q_1,A_1)$.

The isomorphism $C(U)$ now corresponds to the orientation
$$\textstyle (-1)^{k_1}\bigwedge_i\partial_{U_i} \otimes \bigwedge_i e^U_i \in \ddd(T_U\wt\cP(q,q_0;q_1)) \otimes \ddd(T_UQ/\Gamma)$$
where we abbreviated $\bigwedge_i\partial_{U_i} = \bigwedge_i\partial_{u^0_i}\wedge \bigwedge_i\partial_{u^1_i}$ and $\bigwedge_i e^U_i = \bigwedge_i e^{u^0}_i \wedge \bigwedge_i e^{u^1}_i$.

We have
\begin{thm}
The operation $\bullet$ is a chain map. More precisely:
$$\partial_{\cD_1} \circ \bullet = \bullet \circ (\partial_{\cD} \otimes \id + (-1)^{2n-k} \id \otimes \partial_{\cD_0}) \fc QC_k(\cD) \otimes QC_l(\cD_0:L) \to QC_{k+l-2n-1}(\cD_1:L)$$
\end{thm}
\begin{prf}
Unlike the proof in the case of Lagrangian quantum product, due to transversality issues in case $N_L = 2$ (see \cite{Biran_Cornea_Quantum_structures_Lagr_submfds, Biran_Cornea_Rigidity_uniruling} and reference therein), one needs to use more general objects in order to prove the asserted relation, namely one needs to use moduli spaces of spiked pearls where the center is allowed to carry a Hamiltonian perturbation. The proof in \cite{Biran_Cornea_Quantum_structures_Lagr_submfds, Biran_Cornea_Rigidity_uniruling} then proceeds as follows: first one defines a perturbed variant of the operation $\bullet$ using the perturbed moduli spaces; it is then straightforward to show that this perturbed operation is indeed a chain map. Next, one proves that for all sufficiently small perturbations, there is a canonical bijection between the $0$-dimensional components of perturbed and unperturbed moduli spaces of spiked pearls, and therefore that the two operations are identical.

Of course, in the original proof of Biran--Cornea orientations are not taken into account, as all the counts are modulo 2. In \cite{Biran_Cornea_Lagr_top_enumerative_invts} the authors bring orientations into play, but we will not follow that argument, since our methods do not presuppose additional choices such as an orientation of $L$ and (relative S)Pin structures, and since the proof can be carried through without such choices.

The method used in the present paper to show, for example, that the Lagrangian quantum product defines a chain map, can be used with obvious minimal changes to accomodate the situation where the core carries a Hamiltonian perturbation. The same exact strategy can be used to show that if we define $\bullet$ using moduli spaces of spiked pearls in which the center carries a Hamiltonian perturbation, then it is a chain map, and this is based on transversality and gluing results of Biran--Cornea \cite{Biran_Cornea_Quantum_structures_Lagr_submfds, Biran_Cornea_Rigidity_uniruling}. What remains to be shown is that the perturbed operation and the unperturbed one coincide in case the perturbation is small enough. Let $\cP^0(q,q_0;q_1)$ and $\cP^0(q,q_0;q_1;H)$ denote the $0$-dimensional components of the moduli spaces of spiked pearls and perturbed spiked pearls, respectively. Biran--Cornea prove that for all $H$ sufficiently small there is a canonical bijection $\cP^0(q,q_0;q_1;H) \simeq \cP^0(q,q_0;q_1)$. We have to show that if $U_H \in \cP^0(q,q_0;q_1;H)$ and $U \in \cP^0(q,q_0;q_1)$ correspond under this bijection, then the isomorphisms $C(U_H)$, $C(U)$ are equal. This is however obvious, since the maps $U_H,U$ are close, therefore can be deformed into each other, which means there is a canonical isomorphism $\ddd(D_U|_{Y_\Gamma^U}) = \ddd(D_{U_H}|_{Y_\Gamma^{U_H}})$. We have a similar isomorphism $\ddd(T_U\wt\cP(q,q_0;q_1)) \otimes \ddd(T_UQ/\Gamma) \simeq \ddd(T_{U_H}\wt\cP(q,q_0;q_1;H)) \otimes \ddd(T_{U_H}Q/\Gamma)$, which means that the isomorphisms $C(U_H),C(U)$, which correspond to specific orientations of these spaces, coincide. \qed
\end{prf}

Thus we have a well-defined bilinear operation on homology
$$\bullet \fc QH_k(\cD) \otimes QH_l(\cD_0:L) \to QH_{k+l-2n}(\cD_1:L)\,.$$
This is called the quantum module action.

\subsection{Spectral sequences}\label{ss:spectral_seqs}

The quantum complexes admit natural filtrations by the Maslov or Chern numbers, and these give rise to spectral sequences.

We only consider the Lagrangian case. The Lagrangian quantum complex corresponding to a quantum datum $\cD = (f,\rho,J)$ is
$$QC_*(\cD:L) = \bigoplus_{ \substack{x \in \Crit f \\ A \in \pi_2(M,L,x)}} C(x,A)\,.$$
Let us define the increasing filtration
$$F_pQC_*(\cD:L) = \bigoplus_{x \in \Crit f} \bigoplus_{\substack{A \in \pi_2(M,L,x) \\ \mu(A) \geq - pN_L}} C(x,A)\,.$$
It follows from the definition of $\partial_\cD$ that it preserves this filtration. Therefore we have the associated spectral sequence whose zeroth page is
$$E^0_{p,q} = F_p QC_{p+q}(\cD:L) / F_{p-1} QC_{p+q}(\cD:L)$$
and so
$$E^0_{p,*} = F_p QC_*(\cD:L) / F_{p-1} QC_*(\cD:L) \simeq \bigoplus_{x \in \Crit f} \bigoplus_{\substack{A \in \pi_2(M,L,x) \\ \mu(A) = - pN_L}} C(x,A)\,.$$
The boundary operator $\partial^0$ on $E^0_{p,*}$ comes from the Morse boundary operator of $f$, twisted by the local system $\xi_p$, where
$$\xi_p(x) = \{A \in \pi_2(M,L,x)\,|\,\mu(A) = -pN_L\}\,,$$
that is $(E^0_{p,*},\partial^0) \simeq (CM_*(f;\xi_p), \partial^f_{\xi_p})$, where $CM_*$ means the Morse complex, and $\partial^f_{\xi_p}$ is the Morse boundary operator of $f$ twisted by $\xi_p$. Therefore the first page is the twisted homology:
$$E^1_{p,q} \simeq H_{p+q - pN_L}(L;\xi_p)\,.$$
Thus we have
\begin{thm}
The filtration by the Maslov index of disks, $F_*QC_*(\cD:L)$ induces a spectral sequence whose first page is isomorphic to the singular homology of $L$ twisted by the local systems $\xi_p$, and which converges to $QH_*(\cD:L)$. Moreover, using different quantum data, this spectral sequence can be seen to be multiplicative in an obvious sense. \qed
\end{thm}

\section{PSS isomorphisms}\label{s:PSS}

This section is dedicated to the definition and properties of the PSS isomorphisms. The idea of the construction, at least over $\Z_2$, is contained in \cite{Biran_Cornea_Quantum_structures_Lagr_submfds, Biran_Cornea_Rigidity_uniruling} and the references therein, however to the best of our knowledge, this is the first time the construction is carried out in detail over an arbitrary ground ring.

In \S\ref{ss:definition_PSS} we define the PSS morphisms on chain level via their matrix elements, prove that they are chain maps and that they are independent of the auxiliary data such as a perturbation datum. \S\ref{ss:PSS_maps_properties} covers the main properties of the PSS maps, namely the fact that they are indeed isomorphisms and that they respect the natural algebraic structures on Floer and quantum homology.

\subsection{Definition}\label{ss:definition_PSS}

We start with the Lagrangian Floer and quantum homologies. Let $(H,J)$ be a regular Floer datum for $L$ and $\cD = (f,\rho,I_0)$ a regular quantum datum for $L$. The \tb{PSS map}
$$\PSS_{H,J}^\cD \fc CF_*(H:L) \to QC_*(\cD:L)$$
will be defined through its matrix elements, to define which we need new moduli spaces which combine solutions of Floer's PDE and pearly trajectories. Denote $D^2_- = D^2 - \{-1\}$ and consider $-1$ as a negative puncture. Endow it with the standard negative end, and associate the Floer datum $(H,J)$ to it. Choose a regular perturbation datum $(K,I)$ on $D^2_-$ which is compatible with $(H,J)$ and which satisfies $K = 0$, $I = I_0$ near $1$. Let
$$\cM_-(\wt\gamma) = \{u \in C^\infty_b(D^2_-,\partial D^2_-; M,L; \gamma)\,|\,\ol\partial_{K,I}u = 0\,,[\wh\gamma\sharp u] = 0\in \pi_2(M,L)\}\,.$$
Since the perturbation datum is regular, this is a smooth manifold of dimension $|\wt\gamma|$. We have the evaluation map
$$\ev \fc \cM_-(\wt\gamma) \times (\wt\cM(L,I_0))^k \to L^{2k+1}$$
given by
$$\ev(u = (u_0,\dots,u_k)) = (u_0(1),u_1(-1),u_1(1),\dots,u_k(-1),u_k(1))\,.$$
Fix $q \in \Crit f$ and define
$$\wt\cP_k(\wt\gamma,q) = \ev^{-1}(Q^k \times \cS(q))\,.$$
There is a natural $\R^k$-action on this space and we let $\cP_k(\wt\gamma,q)$ be the quotient. We also denote
$$\wt\cP(\wt\gamma,q) = \bigcup_{k\geq 0} \wt\cP_k(\wt\gamma,q)\quad \text{ and }\quad \cP(\wt\gamma,q) = \bigcup_{k \geq 0} \cP_k(\wt\gamma,q)\,.$$
As Biran--Cornea indicate \cite{Biran_Cornea_Quantum_structures_Lagr_submfds, Biran_Cornea_Rigidity_uniruling}, it can be shown that for fixed $(H,J)$ and $\cD$ there is a subset of $\cJ(M,\omega)$ of second category such that for each $I_0$ in it $\cP(\wt\gamma,q)$ is a smooth manifold of local dimension at $[u]$ equal to $|\wt\gamma| - |q| + \mu(u)$ whenever this number does not exceed $1$.

We now proceed with the definition of the matrix elements. Given
$$u \in C^\infty_b(D^2_-,\partial D^2_-;M,L,\gamma) \times (C^\infty(D^2,S^1;M,L))^k$$
with $\ev(u) \in Q^k \times \cS(q)$ there is an obvious way of producing a class $\wt\gamma \sharp u \in \pi_2(M,L,q)$ by concatenating $\wh\gamma$ with $u_0$ and the disks in $u$ along pieces of gradient trajectories connecting the evaluation points. For $u \in \wt\cP(\wt\gamma,q)$ with $|\wt\gamma| - |q| + \mu(u) = 0$ and any class of cappings $\wt\gamma'$ for $\gamma$ we will construct an isomorphism
$$C(u) \fc C(\wt\gamma') \to C(q,\wt\gamma' \sharp u)$$
For $A \in \pi_2(M,L,q)$ the matrix element of the PSS map is then
$$\sum_{\substack{[u] \in \cP(\wt\gamma,q):\\ \wt\gamma' \sharp u = A}} C(u) \fc C(\wt\gamma') \to C(q,A)\,.$$
It remains to define the isomorphism $C(u)$. This is entirely analogous to all the cases described above: there are natural bijections between orientations of \footnote{The spaces $Y_\Gamma^u$, $Y_Q^u$, $X_\Gamma^u$, $T_uQ/\Gamma$ are defined in an obvious manner.} $D_u|_{Y_\Gamma^u}$, isomorphisms $C(\wt\gamma') \simeq C(q,A)$ for $A = \wt\gamma' \sharp u$, and orientations of $\ddd(T_u\wt\cP(\wt\gamma,q)) \otimes \ddd(T_uQ/\Gamma)$. The isomorphism $C(u)$ is then the one corresponding to the orientation
$$\textstyle(-1)^k\bigwedge_i \partial_{u_i} \otimes \bigwedge_i e^u_i \in \ddd(T_u\wt\cP(\wt\gamma,q)) \otimes \ddd(T_uQ/\Gamma)\,,$$
where $k$ is the number of holomorphic disks in $u$. The bijections are constructed as follows: the first one is obtained using deformation and gluing arguments akin to those used in the definition of quantum homology and various algebraic structures on it. The second one comes from the exact Fredholm triple
$$0 \to D_u|_{Y_\Gamma^u} \to D_u|_{Y_Q^u} \to 0_{T_uQ/\Gamma} \to 0 \,.$$

We have
\begin{thm}\label{thm:PSS_is_chain_map_Floer_to_quantum}
The operator $\PSS_{H,J}^\cD$ thus defined is a chain map:
$$\partial_\cD \circ \PSS_{H,J}^\cD = \PSS_{H,J}^\cD \circ \partial_{H,J}\,.$$
\end{thm}
\begin{prf}
The argument is almost identical to the above proofs, so we just briefly outline the main steps. Let $\cP^1(\wt\gamma,q)$ denote the union of $1$-dimensional components of $\cP(\wt\gamma,q)$. This space can be compactified by adding points of the following types: Floer breaking at the negative end of $u_0$; Morse breaking of one of the gradient trajectories; splitting of a holomorphic disk $u_i$, $i > 0$, into two; splitting of $u_0$ into an element of $\cM_-(\wt \gamma)$ and a holomorphic disk; and collision of $u_i,u_{i+1}$ for some $i \geq 0$. Just as in case of pearly spaces, we can define $\ol\cP{}^1(\wt\gamma,q)$ to be the union of compactified components of $\cP^1(\wt\gamma,q)$ where points corresponding to collision/breaking belonging to two different components are identified if they represent the same degenerate trajectory. Thus $\ol\cP{}^1(\wt\gamma,q)$ is endowed with the structure of a $1$-dimensional compact topological manifold with boundary. The boundary points correspond to Floer and Morse breaking and are in an obvious bijection with the summands of the matrix element of
$$\partial_\cD \circ \PSS_{H,J}^\cD - \PSS_{H,J}^\cD \circ \partial_{H,J}$$
going from $C(\wt\gamma')$ to $C(q,A)$. One needs to show that for every connected component of $\ol\cP{}^1(\wt\gamma,q)$ the summands corresponding to its two boundary points cancel out. 

This is achieved as follows. Let $[w] \in \cP^1(\wt\gamma,q)$. Then isomorphisms $C(\wt\gamma') \to C(q,\wt\gamma'\sharp w)$ are in bijection with orientations of $D_w|_{Y_\Gamma^w}$ and of $\ddd(T_w\wt\cP(\wt\gamma,q)) \otimes \ddd(T_wQ/\Gamma)$. The following can be shown:
\begin{itemize}
\item Suppose $\delta = ([u],[v]) \in \cM(H,J;\wt\gamma,\wt\gamma') \times \cP_k(\wt\gamma',q)$ is a boundary point of $\ol\cP{}^1(\wt\gamma,q)$ and $[w]$ lies close to $\delta$. Then the orientation corresponding to $C(v) \circ C(u)$ is $- \bigwedge_i \partial_{w_i} \wedge \text{inward}_\delta \otimes \bigwedge_i e^w_i \in \ddd(T_w\wt\cP(\wt\gamma,q)) \otimes \ddd(T_wQ/\Gamma)$.
\item Suppose $\delta = ([u],[v]) \in \cP_k(\wt\gamma,q') \times \cP_l(q',q)$ and $[w]$ is close to $\delta$. Then the isomorphism $C(v) \circ C(u)$ corresponds to the orientation $\bigwedge_i \partial_{w_i} \wedge \text{inward}_\delta \otimes \bigwedge_i e^w_i \in \ddd(T_w\wt\cP(\wt\gamma,q)) \otimes \ddd(T_wQ/\Gamma)$.
\item If $C \fc C(\wt\gamma') \simeq C(q,A)$ is a fixed isomorphism and $u,v \in \wt\cP^1(\wt \gamma,q)$ are two points lying on different sides of a degenerate trajectory in $\ol\cP{}^1(\wt\gamma,q)$ close to it, and $C$ corresponds to the orientation $\bigwedge_i \partial_{u_i}\wedge \eta_u \otimes \bigwedge_i e^u_i$, where $\eta_u \in T_u\wt\cP(\wt\gamma,q)$ is a vector transverse to the infinitesimal action of the automorphism group at $u$, then $C$ corresponds to the orientation $\bigwedge_i \partial_{v_i}\wedge \eta_v \otimes \bigwedge_i e^v_i$, where $\eta_v \in T_v\wt\cP(\wt\gamma,q)$ points in the same direction as $\eta_u$. In other words, the sign in front of the standard orientation is unchanged when crossing a breaking/collision point.
\end{itemize}
Let us prove, for example, that if there a connected component $\Delta$ with boundary points $\delta = ([u],[v]) \in \cM(H,J;\wt\gamma,\wt\gamma') \times \cP_k(\wt\gamma',q)$ and $\delta' = ([u'],[v']) \in \cP_{k'}(\wt\gamma,q') \times \cP_{l'}(q',q)$, then $C(v) \circ C(u) = C(v') \circ C(u')$. Let $\{[w^t]\}_{t\in (0,1)}$ be a parametrization of $\Delta$ at nondegenerate points, such that
$$[w^t] \xrightarrow[t \to 0]{} \delta\,.$$
The isomorphism $C(v) \circ C(u)$ corresponds to the orientation
$$\textstyle - \bigwedge_i \partial_{w^t_i} \wedge \text{inward}_\delta \otimes \bigwedge_i e^{w^t}_i$$
for $t$ close to $0$. We see that $C(v) \circ C(u)$ also corresponds to the orientation
$$\textstyle - \bigwedge_i \partial_{w^t_i} \wedge \text{inward}_\delta \otimes \bigwedge_i e^{w^t}_i$$
for $t$ close to $1$. On the other hand the isomorphism $C(v') \circ C(u')$ corresponds to the orientation
$$\textstyle \bigwedge_i \partial_{w^t_i} \wedge \text{inward}_{\delta'} \otimes \bigwedge_i e^{w^t}_i$$
for $t$ close to $1$. Since $\text{inward}_\delta = - \text{inward}_{\delta'}$, we see that the two orientations are equal, whence the equality of the isomorphisms $C(v) \circ C(u) = C(v') \circ C(u')$.

The computation of the above orientations is done in analogy with the proofs appearing in \S\ref{s:QH}: one uses commutative diagrams of determinant lines coming from exact Fredholm squares and from gluing, and uses the structure of the respective gluing maps to conclude. \qed
\end{prf}

The opposite PSS map
$$\PSS^{H,J}_\cD \fc QC_*(\cD:L) \to CF_*(H:L)$$
is constructed analogously. We will describe its matrix elements. To define them, we need opposite moduli spaces. Recall that $\dot D^2 = D^2 - \{1\}$ and consider $1$ as a positive puncture. Endow it with the standard end, and associate the Floer datum $(H,J)$ to it. Choose a regular perturbation datum $(K,I)$ on $\dot D^2$ which is compatible with $(H,J)$ and which satisfies $K = 0$, $I = I_0$ near $-1$. Let
$$\cM_+(\wt\gamma) = \{u \in C^\infty_b(\dot D^2,\partial \dot D^2; M,L; \gamma)\,|\,\ol\partial_{K,I}u = 0\,,\wh\gamma \text{ is equivalent to }u \text{ as a capping of }\gamma\}\,.$$
Since the perturbation datum is regular, this is a smooth manifold of dimension $|\wt\gamma|' = n - |\wt\gamma|$. We have the evaluation map
$$\ev \fc (\wt\cM(L,I_0))^k \times \cM_+(\wt\gamma) \to L^{2k+1}$$
given by
$$\ev(u = (u_1,\dots,u_k;u_0)) = (u_1(-1),u_1(1),\dots,u_k(-1),u_k(1);u_0(-1))\,.$$
Fix $q \in \Crit f$ and define
$$\wt\cP_k(q,\wt\gamma) = \ev^{-1}(\cU(q) \times Q^k)\,.$$
There is a natural $\R^k$-action on this space and we let $\cP_k(q,\wt\gamma)$ be the quotient. We also denote
$$\wt\cP(q,\wt\gamma) = \bigcup_{k\geq 0} \wt\cP_k(q,\wt\gamma)\quad \text{ and }\quad \cP(q,\wt\gamma) = \bigcup_{k \geq 0} \cP_k(q,\wt\gamma)\,.$$
Again, it can be shown that for fixed $(H,J)$ and $\cD$ there is a subset of $\cJ(M,\omega)$ of the second category such that for each $I_0$ in it $\cP(q,\wt\gamma)$ is a smooth manifold of local dimension at $[u]$ equal to $|q| - |\wt\gamma| + \mu(u)$ whenever this number does not exceed $1$.

We now proceed with the definition of the matrix elements. Given
$$u \in (C^\infty(D^2,S^1;M,L))^k \times C^\infty_b(\dot D^2,\partial \dot D^2;M,L;\gamma)$$
with $\ev(u) \in \cU(q) \times Q^k$ there is an obvious way of producing a class of cappings $A \sharp u$ of $\gamma$ for $A \in \pi_2(M,L,q)$ by concatenating a representative of $A$ with the disks of $u$ and the capping $u_0$ of $\gamma$ along pieces of gradient trajectories connecting the evaluation points. For $u \in \wt\cP(q,\wt\gamma)$ with $|q| - |\wt\gamma| + \mu(u) = 0$ we will construct an isomorphism
$$C(u) \fc C(q,A) \to C(A\sharp u)\,.$$
The matrix element is then
$$\sum_{\substack{[u] \in \cP(q,\wt\gamma): \\ A\sharp u = \wt\gamma'}} C(u) \fc C(q,A) \to C(\wt\gamma')\,.$$
It remains to define the isomorphism $C(u)$. This is identical to the above construction: there are natural bijections between orientations of $D_u|_{Y_\Gamma^u}$, isomorphisms $C(q,A) \simeq C(\wt\gamma)$ for $\wt\gamma = A\sharp u$, and orientations of $\ddd(T_u\wt\cP(q,\wt\gamma)) \otimes \ddd(T_uQ/\Gamma)$. The isomorphism $C(u)$ is then the one corresponding to the orientation
$$\textstyle\bigwedge_i \partial_{u_i} \otimes \bigwedge_i e^u_i \in \ddd(T_u\wt\cP(\wt\gamma,q)) \otimes \ddd(T_uQ/\Gamma)\,.$$
The bijections are constructed as follows: the first one is obtained using deformation and gluing arguments akin to those used in the definition of quantum homology and various algebraic structures on it. The second one comes from the exact Fredholm triple $0 \to D_u|_{Y_\Gamma^u} \to D_u|_{Y_Q^u} \to 0_{T_uQ/\Gamma} \to 0$.

We have
\begin{thm}
The operator $\PSS^{H,J}_\cD$ is a chain map:
$$\partial_{H,J} \circ \PSS^{H,J}_\cD = \PSS^{H,J}_\cD \circ \partial_\cD\,.$$
\end{thm}
\begin{prf}
This is proved in the exact same way as Theorem \ref{thm:PSS_is_chain_map_Floer_to_quantum}, the only difference being in the orientations induced by the boundary points and the fact that there is a sign change when crossing a breaking/collision point. To summarize the argument, $\ol\cP{}^1(q,\wt\gamma)$ has the structure of a $1$-dimensional compact topological manifold with boundary consisting of points corresponding to Floer or Morse breaking. If $[w] \in \cP^1(q,\wt\gamma)$, then isomorphisms $C(q,A) \simeq C(A\sharp w)$ are in bijection with orientations of $D_w|_{Y_\Gamma^w}$ and of $\ddd(T_w\wt\cP(q,\wt\gamma)) \otimes \ddd(T_wQ/\Gamma)$. The following can be shown:
\begin{itemize}
\item Suppose $\delta = ([u],[v]) \in \cP_k(q,\wt\gamma') \times \cM(H,J;\wt\gamma',\wt\gamma)$ is a boundary point of $\ol\cP{}^1(q,\wt\gamma)$ and $[w]$ lies close to $\delta$. Then the orientation corresponding to $C(v) \circ C(u)$ is $(-1)^k\bigwedge_i \partial_{w_i} \wedge \text{inward}_\delta \otimes \bigwedge_i e^w_i \in \ddd(T_w\wt\cP(\wt\gamma,q)) \otimes \ddd(T_wQ/\Gamma)$.
\item Suppose $\delta = ([u],[v]) \in \cP_k(q,q') \times \cP_l(q',\wt\gamma)$ and $[w]$ is close to $\delta$. Then the isomorphism $C(v) \circ C(u)$ corresponds to the orientation $(-1)^{k+l+1}\bigwedge_i \partial_{w_i} \wedge \text{inward}_\delta \otimes \bigwedge_i e^w_i \in \ddd(T_w\wt\cP(\wt\gamma,q)) \otimes \ddd(T_wQ/\Gamma)$.
\item If $C \fc C(q,A) \simeq C(\wt\gamma')$ is a fixed isomorphism and $u,v \in \wt\cP^1(\gamma,q)$ are two points lying on different sides of a degenerate trajectory in $\ol\cP{}^1(\wt\gamma,q)$ close to it, and $C$ corresponds to the orientation $\bigwedge_i \partial_{u_i}\wedge \eta_u \otimes \bigwedge_i e^u_i$, where $\eta_u \in T_u\wt\cP(\wt\gamma,q)$ is a vector transverse to the infinitesimal action of the automorphism group at $u$, then $C$ corresponds to the orientation $-\bigwedge_i \partial_{v_i}\wedge \eta_v \otimes \bigwedge_i e^v_i$, where $\eta_v \in T_v\wt\cP(\wt\gamma,q)$ points in the same direction as $\eta_u$. In other words, the sign in front of the standard orientation flips when crossing a breaking/collision point.
\end{itemize}
From this we can deduce the vanishing of the matrix elements of
$$\partial_{H,J} \circ \PSS^{H,J}_\cD - \PSS^{H,J}_\cD \circ \partial_\cD\,. \quad\qed$$
\end{prf}

The definition of PSS maps between Hamiltonian Floer homology and the quantum homology of $M$ proceeds in a similar, though much simpler, manner. We describe it here briefly, mainly to establish notation for later use.

Fix a regular nondegenerate periodic Floer datum $(H,J)$ and a regular quantum datum $\cD = (f,\rho,I_0)$ for $M$. We start with the definition of
$$\PSS_{H,J}^\cD \fc CF_*(H) \to QC_*(\cD)\,.$$
Let $S^2_- = S^2 - \{-1\}$, where we view $S^2 = \C P^1$, with $-1$ being a negative puncture. Endow it with the standard negative end and associate $(H,J)$ to it. Pick a regular perturbation datum $(K,I)$ on $S^2_-$ which is compatible with $(H,J)$ and which satisfies $K = 0$, $I = I_0$ near $1$. For $\wt x \in \Crit \cA_H$ define
$$\cM_-(\wt x) = \{u \in C^\infty_b(S^2_-,M;x)\,|\,\ol\partial_{K,I}u=0\,,[\wh x \sharp u] = 0 \in \pi_2(M)\}\,.$$
This is a smooth manifold of dimension $|\wt x|$. We have the evaluation map
$$\ev \fc \cM_-(\wt x) \to M\,,\quad \ev(u) = u(1)\,.$$
For $q \in \Crit f$ let
$$\cP(\wt x, q) = \ev^{-1}(\cS(q))\,.$$
This is a smooth manifold of dimension $|\wt x| - |q|$.

We now define the matrix elements of the PSS map. Given $u \in C^\infty_b(S^2_-,M;x)$ with $\ev(u) \in \cS(q)$ we can produce a class $\wt x \sharp u \in \pi_2(M,q)$ by concatenating a representative of $\wt x$ with $u$ and transferring it to $q$ using the piece of gradient trajectory connecting $u(1)$ with $q$. For $u \in \cP(\wt x,q)$ and any class of cappings $\wt x'$ of $x$ we will define an isomorphism
$$C(u) \fc C(\wt x') \to C(q,\wt x' \sharp u)$$
The matrix element of the PSS map is then
$$\sum_{\substack{u \in \cP(\wt x,q): \\ \wt x' \sharp u = A}} C(u) \fc C(\wt x') \to C(q,A)\,.$$
To define the isomorphism $C(u)$ we note that, as usual, there is a bijection between orientations \footnote{The space $Y^u$ is defined in an obvious manner.} of $D_u|_{Y^u}$ and isomorphisms $C(\wt x') \simeq C(q,\wt x' \sharp u)$. On the other hand, $D_u|_{Y^u}$ is an index zero surjective operator, therefore an isomorphism. The isomorphism $C(u)$ is the one corresponding to the positive canonical orientation $1 \otimes 1^\vee \in \ddd(D_u|_{Y^u})$. The bijection is constructed using the gluing and deformation isomorphisms.

This defines the PSS map, and we have
\begin{thm}
The PSS map is a chain map:
$$\partial_\cD \circ \PSS_{H,J}^\cD = \PSS_{H,J}^\cD \circ \partial_{H,J}\,.$$
\end{thm}
\begin{prf}
The $1$-dimensional part $\cP^1(\wt x,q)$ can be compactified by adding Floer and Morse breaking. It suffices to note that a boundary point corresponding to Floer breaking induces the outward orientation on $\cP^1(\wt x, q)$ while if it corresponds to Morse breaking, the induced orientation is inward. Since there is no breaking/collision, the fact that for an element $w \in \cP^1(\wt x,q)$ there is a bijection between isomorphisms $C(\wt x') \simeq C(q,\wt x'\sharp w)$ and orientations of $\ddd(D_w|_{Y^w}) = \ddd(T_w\cP(\wt x,q))$, suffices to establish the vanishing of the matrix elements of the difference $\partial_\cD \circ \PSS_{H,J}^\cD = \PSS_{H,J}^\cD \circ \partial_{H,J}$. \qed
\end{prf}

Similarly, we have the opposite PSS map
$$\PSS^{H,J}_\cD \fc QC_*(\cD) \to CF_*(H)\,.$$
Let $S^2_+ = S^2 - \{1\}$ with the puncture $\{1\}$ being positive. Endow it with the standard end and associate $(H,J)$ to it. Pick a regular perturbation datum $(K,I)$ on $S^2_+$ which is compatible with $(H,J)$ and which satisfies $K = 0$, $I = I_0$ near $-1$. For $\wt x \in \Crit \cA_H$ define
$$\cM_+(\wt x) = \{u \in C^\infty_b(S^2_+,M;x)\,|\,\ol\partial_{K,I}u=0\,,u \in \wt x \text{ as a capping}\}\,.$$
This is a smooth manifold of dimension $|\wt x|' = 2n - |\wt x|$. We have the evaluation map
$$\ev \fc \cM_+(\wt x) \to M\,,\quad \ev(u) = u(-1)\,.$$
For $q \in \Crit f$ let
$$\cP(q,\wt x) = \ev^{-1}(\cU(q))\,.$$
This is a smooth manifold of dimension $|q| - |\wt x|$.

We now define the matrix elements of the PSS map. Given $u \in C^\infty_b(S^2_+,M;x)$ with $\ev(u) \in \cU(q)$, and $A \in \pi_2(M,q)$ we can produce a class of cappings $A \sharp u$ of $x$ by concatenating a representative of $A$ with $u$ along the piece of gradient trajectory connecting $q$ with $u(-1)$. For $u \in \cP(q,\wt x)$ and $A \in \pi_2(M,q)$ we will define an isomorphism
$$C(u) \fc C(q,A) \simeq C(A \sharp u)\,.$$
The matrix element of the PSS map is then
$$\sum_{\substack{u \in \cP(q, \wt x): \\ A \sharp u = \wt x'}}C(u) \fc C(q,A) \to C(\wt x')\,.$$
To define the isomorphism $C(u)$ note that there is a bijection between orientations of $D_u|_{Y^u}$ and isomorphisms $C(q,A) \simeq C(A \sharp u)$. On the other hand, $D_u|_{Y^u}$ is an index zero surjective operator, therefore an isomorphism. The isomorphism $C(u)$ is the one corresponding to the positive canonical orientation $1 \otimes 1^\vee \in \ddd(D_u|_{Y^u})$. The bijection is constructed using the gluing and deformation isomorphisms.

This defines the PSS map, and we can prove
\begin{thm}
The PSS map is a chain map:
$$\partial_{H,J} \circ \PSS^{H,J}_\cD = \PSS^{H,J}_\cD \circ \partial_\cD\,.\quad \qed$$
\end{thm}

\subsubsection{Independence of auxiliary data}\label{sss:indep_aux_data_PSS}

The idea is to use moduli spaces of mixed pearls parametrized over $[0,1]$ and to use a homotopy of auxiliary data, such as conformal structures and perturbation data over the interval. We then define a chain homotopy between the PSS maps corresponding to the perturbation data over $0,1$ using the $0$-dimensional moduli spaces, and then the $1$-dimensional moduli spaces are used in order to show that this is indeed a chain homotopy. This combines techniques from \S\ref{s:QH} and the present section.

We will give details only for the Lagrangian case, the other one being entirely similar, albeit much simpler. More precisely, let $(H,J)$ be a regular Floer datum, $\cD = (f,\rho, I_0)$ a regular quantum homology datum for $L$, and let $(K,I) = \{(K^r,I^r)\}_{r \in [0,1]}$ be a homotopy of perturbation data on $D^2_-$ compatible with $(H,J)$, which is stationary for $r$ close to $0,1$. Assume that both $(K^0,I^0)$, $(K^1,I^1)$ are regular and the homotopy is regular as well. Any two perturbation data can be connected by such a homotopy. For $\wt\gamma \in \Crit \cA_{H:L}$ define
$$\cM_-(K,I;\wt\gamma) = \{(r,u)\,|\, r\in[0,1]\,,u\in C^\infty_b(D^2_-,\partial D^2_-;M,L;\gamma)\,,\ol\partial_{K^r,I^r}u = 0\,,[\wh\gamma\sharp u] = 0 \in \pi_2(M,L)\}\,.$$
This is a smooth manifold of dimension $|\wt\gamma| + 1$. We have the evaluation map
$$\ev \fc \cM_-(K,I;\wt\gamma) \times (\wt\cM(L,I_0))^k \to L^{2k+1}$$
defined via
$$\ev(u = (u_0,u_1,\dots, u_k)) = (u_0(1),u_1(-1),\dots,u_k(1))$$
and we let
$$\wt\cP_k(K,I;\wt\gamma,q) = \ev^{-1}(Q^k \times \cS(q))\,.$$
We have a natural $\R^k$-action on this space and we let $\cP_k(K,I;\wt\gamma,q)$ be the quotient. We let $\wt\cP(K,I;\wt\gamma,q)$ and $\cP(K,I;\wt\gamma,q)$ be the respective unions of these spaces over $k \geq 0$. For generic $I_0$ the space $\cP(K,I;\wt\gamma,q)$ is a smooth manifold of local dimension at $u$ equal to $|\wt\gamma| - |q| + \mu(u) + 1$ if this number does not exceed $1$.

Using the $0$-dimensional part of this space one can define a homomorphism
$$\Psi\equiv \Psi_{K,I} \fc CF_*(H:L) \to QC_{*+1}(\cD:L)\,.$$
This is defined via its matrix elements. To define these it suffices, using the above methods, to orient the operator $D_u|_{Y_\Gamma^u}$ for $(r,u) \in \wt\cP(K,I;\wt\gamma,q)$ such that $|\wt\gamma| - |q| + \mu(u) + 1 = 0$. This orientation then induces an isomorphism $C(\wt\gamma') \simeq C(q,\wt\gamma'\sharp u)$, and summing up over the elements of the $0$-dimensional part of $\cP(K,I;\wt\gamma,q)$ produces the desired matrix elements. In order to orient $D_u|_{Y_\Gamma^u}$ we note that there is an exact Fredholm triple
$$0 \to D_u|_{Y_\Gamma^u} \to D_{r,u}|_{Y_Q^u} \to 0_{T_r[0,1] \oplus T_uQ/\Gamma} \to 0$$
which induces a bijection between orientations of $D_u|_{Y_\Gamma^u}$ and orientations of $\ddd(T_{(r,u)}\wt\cP(K,I;\wt\gamma,q)) \otimes \ddd(T_r[0,1] \oplus T_uQ/\Gamma)$. Now pick the orientation $\bigwedge_i\partial_{u_i} \otimes \bigwedge_i e^u_i \wedge \partial_r$.

Using the $1$-dimensional part of $\cP(K,I;\wt\gamma,q)$ and chasing suitable commutative diagrams as above, we can prove
\begin{thm}
The map $\Psi_{K,I}$ defines a chain homotopy between the PSS maps constructed using the perturbation data $(K^0,I^0)$ and $(K^1,I^1)$. \qed
\end{thm}

Analogously one proves that the opposite PSS maps are also independent of the perturbation datum used. This allows us to see that there are well-defined morphisms
$$\PSS_{H,J}^\cD \fc HF_*(H,J:L) \to QH_*(\cD:L)\,,\quad \text{ and }\quad \PSS_\cD^{H,J} \fc QH_*(\cD:L) \to HF_*(H,J:L)\,.$$

\subsection{Properties}\label{ss:PSS_maps_properties}

\subsubsection{PSS maps and the continuation morphisms}\label{sss:PSS_maps_cont_isomorphisms}

Let $(H^i,J^i)$, $i=0,1$ be regular Floer data associated to the puncture of $D^2_-$, and let $(H^s,J^s)_{s\in\R}$ be a regular homotopy between these which is stationary outside $s \in (0,1)$. Let $\cD = (f,\rho,I_0)$ be a regular quantum homology datum. We will prove here that the PSS maps respect Floer continuation morphisms, more precisely that
$$\PSS_{H^1,J^1}^\cD \circ \Phi_{H^0,J^0}^{H^1,J^1} = \PSS_{H^0,J^0}^\cD \fc HF_*(H^0,J^0:L) \to QH_*(\cD:L)\,.$$
In order to prove this we will construct a suitable chain homotopy. First we need some additional moduli spaces. Let $(K,I) = \{(K^r,I^r)\}_{r\geq 0}$ be a regular perturbation datum on the trivial family $[0,\infty) \times D^2_-$, which is stationary for $r$ near $0$, which is compatible with $(H^0,J^0)$, which satisfies $K=0,I=I_0$ near $1 \in D^2_-$, and in addition for $r$ large enough we have the following expression in the coordinates $(s,t) \in (-\infty,0] \times [0,1]$ on the end of $D^2_-$:
$$K^r(s,t) = H^{s+r+1}_t\,dt\,, \quad I^r(s,t) = J^{s+r+1}_t\,.$$

For $\wt\gamma \in \Crit\cA_{H^0:L}$ let
$$\cM_-(K,I;\wt\gamma) = \{(r,u)\,|\, r\in[0,\infty)\,,u\in C^\infty_b(D^2_-,\partial D^2_-;M,L;\gamma)\,,\ol\partial_{K^r,I^r}u = 0\,,[\wh\gamma\sharp u] = 0 \in \pi_2(M,L)\}\,.$$
This is a smooth manifold of dimension $|\wt\gamma| + 1$. We have the evaluation map
$$\ev \fc \cM_-(K,I;\wt\gamma) \times (\wt\cM(L,I_0))^k \to L^{2k+1}$$
defined via
$$\ev(u = (u_0,u_1,\dots, u_k)) = (u_0(1),u_1(-1),\dots,u_k(1))$$
and for $q \in \Crit f$ we let
$$\wt\cP_k(K,I;\wt\gamma,q) = \ev^{-1}(Q^k \times \cS(q))\,.$$
We have a natural $\R^k$-action on this space and we let $\cP_k(K,I;\wt\gamma,q)$ be the quotient. We let $\wt\cP(K,I;\wt\gamma,q)$ and $\cP(K,I;\wt\gamma,q)$ be the respective unions of these spaces over $k \geq 0$. For generic $I_0$ the space $\cP(K,I;\wt\gamma,q)$ is a smooth manifold of local dimension at $u$ equal to $|\wt\gamma| - |q| + \mu(u) + 1$ if this number does not exceed $1$.

Using the $0$-dimensional part of this space one can define a homomorphism
$$\Psi\equiv \Psi_{K,I} \fc CF_*(H^0:L) \to QC_{*+1}(\cD:L)\,.$$
This is defined via its matrix elements. To define these it suffices, using the above methods, to orient the operator $D_u|_{Y_\Gamma^u}$ for $(r,u) \in \wt\cP(K,I;\wt\gamma,q)$ such that $|\wt\gamma| - |q| + \mu(u) + 1 = 0$. This orientation then induces an isomorphism $C(\wt\gamma') \simeq C(q,\wt\gamma'\sharp u)$, and summing up over the elements of the $0$-dimensional part of $\cP(K,I;\wt\gamma,q)$ produces the desired matrix elements. In order to orient $D_u|_{Y_\Gamma^u}$ we note that there is an exact Fredholm triple
$$0 \to D_u|_{Y_\Gamma^u} \to D_{r,u}|_{Y_Q^u} \to 0_{T_r[0,\infty) \oplus T_uQ/\Gamma} \to 0$$
which induces a bijection between orientations of $D_u|_{Y_\Gamma^u}$ and orientations of $\ddd(T_{(r,u)}\wt\cP(K,I;\wt\gamma,q)) \otimes \ddd(T_r[0,\infty) \oplus T_uQ/\Gamma)$. Now pick the orientation $\bigwedge_i\partial_{u_i} \otimes \bigwedge_i e^u_i \wedge \partial_r$.

The $1$-dimensional part of $\cP(K,I;\wt\gamma,q)$ compactifies by adding Morse breaking, Floer breaking at $\wt\gamma$, which corresponds to the continuation morphism, as well as disk collision and breaking. It is then possible to show, using the same methods as above, that the following is true:
\begin{thm}
The map $\Psi_{K,I}$ defines a chain homotopy between $\PSS_{H^0,J^0}^\cD$ and $\PSS_{H^1,J^1}^\cD \circ \Phi_{H^0,J^0}^{H^1,J^1}$. \qed
\end{thm}

The proof that the opposite PSS morphism, as well as the PSS morphisms for quantum homology of $M$, respect continuation morphisms, proceeds in the same manner. Therefore the various PSS morphisms assemble into canonical morphisms
$$\PSS^\cD \fc HF_*(L) \to QH_*(\cD:L)\,,\quad \PSS_\cD \fc QH_*(\cD:L) \to HF_*(L)\,,$$
and
$$\PSS^\cD \fc HF_*(M) \to QH_*(\cD)\,,\quad \PSS_\cD \fc QH_*(\cD) \to HF_*(M)\,.$$

\subsubsection{PSS maps are isomorphisms}\label{sss:PSS_are_isomorphisms}

We show here that the compositions
$$\PSS_{H,J}^\cD \circ \PSS_\cD^{H,J} \quad \text{ and } \quad \PSS_\cD^{H,J} \circ \PSS_{H,J}^\cD$$
are chain homotopic to identity.

Let us start with the second composition. We will in fact produce a chain homotopy between $\PSS_\cD^{H,J} \circ \PSS_{H,J}^\cD$ and a continuation morphism from $(H,J)$ to itself. For this we need some new moduli spaces. For $k \geq 1$ and $\wt\gamma_\pm \in \Crit \cA_{H:L}$ we have the evaluation map
$$\ev \fc \cM_-(\wt\gamma_-) \times (\wt\cM(L,I_0))^{k-1} \times \cM_+(\wt\gamma_+) \to L^{2k}\,,$$
$$\ev(u = (u_0,u_1,\dots,u_{k-1},u_k)) = (u_0(1),u_1(-1),u_1(1),\dots,u_{k-1}(-1),u_{k-1}(1),u_k(-1))\,,$$
and we let
$$\wt\cP_k(\wt\gamma_-,\wt\gamma_+) = \ev^{-1}(Q^k)\,.$$
There is a natural $\R^{k-1}$-action on this space and we let $\cP_k(\wt\gamma_-,\wt\gamma_+)$ be the quotient.

We define $\wt\cP_0(\wt\gamma_-,\wt\gamma_+) \equiv \cP_0(\wt\gamma_-,\wt\gamma_+) = \cM_\cS(K,I;\wt\gamma_-,\wt\gamma_+)$ where $\cS = [r_0,\infty) \times S$ is a trivial family for some $r_0 > 0$, and where $(K,I)$ is a regular perturbation datum which equals $(K^r,I^r) = (0,I_0)$ for $s \in [-r,r]$, equals $(H_t\,dt, J_t)$ for $s$ outside $(-r-1,r+1)$, and interpolates between the two for the rest of values of $s$. We let
$$\wt\cP(\wt\gamma_-,\wt\gamma_+) = \bigcup_{k \geq 0}\wt\cP_k(\wt\gamma_-,\wt\gamma_+)\,,\quad \cP(\wt\gamma_-,\wt\gamma_+) = \bigcup_{k \geq 0}\cP_k(\wt\gamma_-,\wt\gamma_+)\,.$$
The set $\cP(\wt\gamma_-,\wt\gamma_+)$ is a smooth manifold of local dimension at $u$ equal to $|\wt\gamma_-| - |\wt\gamma_+| + \mu(u) + 1$ if this number does not exceed $1$. It is now clear how to proceed with the definition of the chain homotopy. For $u \in \wt\cP_k(\wt\gamma_-,\wt\gamma_+)$ with $|\wt\gamma_-| - |\wt\gamma_+| + \mu(u) + 1 = 0$, we need to orient $D_u|_{Y_\Gamma^u}$ in order to obtain an isomorphism $C(u) \fc C(\wt\gamma_-) \simeq C(\wt\gamma_-\sharp u)$. For $k = 0$ this proceeds as described in \S\ref{s:HF}. For $k \geq 1$ we have the exact triple
$$0 \to D_u|_{Y_\Gamma^u} \to D_u|_{Y_Q^u} \to T_uQ/\Gamma \to 0\,.$$
The operator $D_u|_{Y_Q^u}$ is onto with $(k-1)$-dimensional kernel which is nothing but $T_u\wt\cP_k(\wt\gamma_-,\wt\gamma_+)$. We orient $D_u|_{Y_\Gamma^u}$ by the orientation corresponding to
$$\textstyle\bigwedge_i\partial_{u_i} \otimes \bigwedge_ie^u_i \in \ddd(T_u\wt\cP(\wt\gamma_-,\wt\gamma_+)) \otimes \ddd(T_uQ/\Gamma)\,.$$
Similarly to the above, one proves
\begin{thm}
The map $\Psi$ thus defined is a chain homotopy between $\PSS_\cD^{H,J}\circ \PSS_{H,J}^\cD$ and the continuation map $\Phi_{H,J}^{H,J}$ defined with the help of the perturbation datum $(K^{r_0},I^{r_0})$. \qed
\end{thm}

We now will show that the composition $\PSS^\cD_{H,J} \circ \PSS_\cD^{H,J}$ is homotopic to the identity. For this we need new moduli spaces. Recall that $D^2 - \{\pm1\}$ is biholomorphic to $\R \times [0,1]$. Denote by $(s,t)$ the coordinates on $D^2 - \{\pm1\}$ induced by this biholomorphism. Consider a regular perturbation datum $(K,I)$ on the trivial family $\cS = [0,\infty) \times D^2$ which satisfies: for $r$ close to $0$, $K = 0, I = I_0$, for $r$ large $(K^r(s,t), I^r(s,t))$ equals $(H_t\,dt, J_t)$ for $s \in [-r,r]$, equals $(0,I_0)$ for $s \notin (-r-1,r+1)$, and interpolates between the two for the rest of values of $s$. Consider the space
$$\cM_\cS(K,I) = \{(r,u) \in [0,\infty) \times C^\infty(D^2,S^1;M,L)\,|\, \ol\partial_{K^r,I^r}u = 0\}\,.$$
This is a smooth manifold of local dimension at $u$ equal to $n + \mu(u) + 1$. We have the evaluation map
$$\ev \fc (\wt\cM(L,I_0))^{k_0} \times \cM_\cS(K,I) \times (\wt\cM(L,I_0))^{k_1} \to L^{2(k_0+k_1+1)}$$
defined by
$$\ev(U = (u^0_1,\dots,u^0_{k_0};u;u^1_1,\dots,u^1_{k_1})) = (u^0_1(-1),\dots,u^0_{k_0}(1);u(-1),u(1);u^1_1(-1),\dots,u^1_{k_1}(1))\,.$$
For $q,q'\in\Crit f$ define
$$\wt\cP_{k_0,k_1}(K,I;q,q') = \ev^{-1}(\cU(q) \times Q^{k_0+k_1} \times \cS(q'))\,.$$
There is a natural $\R^{k_0+k_1}$-action on this space and we let $\cP_{k_0,k_1}(K,I;q,q')$ be the quotient. We also define
$$\wt\cP(K,I;q,q') = \bigcup_{k_0,k_1\geq 0}\wt\cP_{k_0,k_1}(K,I;q,q')\quad \text{ and } \quad \cP(K,I;q,q') = \bigcup_{k_0,k_1\geq 0}\cP_{k_0,k_1}(K,I;q,q')\,.$$
The space $\cP(K,I;q,q')$ is a smooth manifold of local dimension at $U$ equal to $|q| - |q'| + \mu(U) + 1$ whenever this number does not exceed $1$.

In order to define the desired chain homotopy, we need only orient the operator $D_U|_{Y_\Gamma^U}$ for $U \in \wt\cP_{k_0,k_1}(K,I;q,q')$ such that $|q| - |q'| + \mu(U) + 1 = 0$. We have the exact Fredholm triple
$$0 \to D_U|_{Y_\Gamma^U} \to D_{r,U}|_{T_r[0,\infty) \times Y_Q^U} \to T_r[0,\infty) \times T_UQ/\Gamma \to 0$$
whence the isomorphism
$$\ddd(D_{r,U}|_{T_r[0,\infty) \times Y_Q^U}) \simeq \ddd(D_U|_{Y_\Gamma^U}) \otimes \ddd(T_UQ/\Gamma)\,.$$
We orient the operator $D_U|_{Y_\Gamma^U}$ by the orientation corresponding to the orientation
$$\textstyle (-1)^{k_0} \bigwedge_i\partial_{U_i} \otimes \bigwedge_ie^U_i \wedge \partial_r\,.$$
Using the $1$-dimensional moduli spaces, one proves
\begin{thm}
The map $\Psi$ thus defined is a chain homotopy between $\PSS^\cD_{H,J} \circ \PSS^{H,J}_\cD$ and the identity. \qed
\end{thm}
\begin{prf}
We need only show that for $r = 0$ one gets the identity. This follows from dimension considerations and the fact that the central disk is $I_0$-holomorphic in this case. Unless it is constant, it admits a $1$-parameter family of reparametrizations which contradicts the fact that the corresponding pearly space is $0$-dimensional. Therefore the disk is constant and $q=q'$. \qed
\end{prf}

Therefore we have shown that the compositions
$$\PSS_{H,J}^\cD \circ \PSS_\cD^{H,J} \quad \text{ and } \quad \PSS_\cD^{H,J} \circ \PSS_{H,J}^\cD$$
equal identity maps on homology, and therefore they are inverse isomorphisms.

The proof of analogous facts for the Hamiltonian Floer homology and the quantum homology of $M$ proceeds in a similar, though much simpler, manner.

\subsubsection{Continuation maps for quantum homology}\label{sss:continuation_maps_QH}

It is possible to define continuation morphisms on quantum homology directly using the approach of Biran--Cornea \cite{Biran_Cornea_Quantum_structures_Lagr_submfds, Biran_Cornea_Rigidity_uniruling}. But since we are ultimately interested in spectral invariants \cite{Leclercq_Zapolsky_Spectral_invts_monotone_Lags}, it is less important to have an independent good definition of quantum homology and therefore we choose another path.

Let $\cD,\cD'$ be two data for Lagrangian quantum homology. We define the {continuation morphism}
$$\Phi_\cD^{\cD'} \fc QH_*(\cD:L) \to QH_*(\cD':L) \quad \text{by} \quad \Phi_\cD^{\cD'} = \PSS_{H,J}^{\cD'} \circ \PSS_{\cD}^{H,J}$$
for a regular Floer datum $(H,J)$ on $L$. The previous subsection implies that this definition is independent of $(H,J)$ since the PSS maps respect Floer continuation maps. Since PSS maps are isomorphisms, so are the quantum homology continuation maps. In particular we can now define the abstract quantum homology $QH_*(L)$ as the limit of the system of homologies $QH_*(\cD:L)$ connected by the isomorphisms $\Phi_\cD^{\cD'}$.

Also we can define the abstract PSS isomorphisms
$$\PSS \fc HF_*(L) \to QH_*(L) \quad \text{ and }\quad \PSS \fc QH_*(L) \to HF_*(L)\,.$$

\subsubsection{PSS maps respect the algebraic structures}\label{sss:PSS_maps_respect_alg_structs}

We start with the product in Lagrangian Floer and quantum homologies. Here we show that for any regular Floer data $(H^i,J^i)$ and regular quantum data $\cD_i = (f_i,\rho,I_0)$ for $L$, $i=0,1,2$ we have
$$\PSS_{H^2,J^2}^{\cD_2} \circ \star \circ (\PSS_{\cD_0}^{H^0,J^0} \otimes \PSS_{\cD_1}^{H^1,J^1}) = \star \fc QH_*(\cD_0:L) \otimes QH_*(\cD_1:L) \to QH_*(\cD_2:L)\,,$$
where $\star$ denotes the product on both the Floer and quantum homology. This implies that the PSS maps respect the multiplicative structure on both sides. 

We will exhibit a chain homotopy between $\PSS_{H^2,J^2}^{\cD_2} \circ \star \circ (\PSS_{\cD_0}^{H^0,J^0} \otimes \PSS_{\cD_1}^{H^1,J^1})$ and $\star$. In order to do so we define new moduli spaces. We choose a regular perturbation datum $(K,I)$ on the trivial family $\cS = [0,\infty) \times D^2$, as follows. Fix a choice of ends on the punctured surface
$$D^2-\{e^{2\pi i j/3}\,|\,j=0,1,2\}\,,$$
such that the puncture $e^{4\pi i/3}$ is positive and the other two are negative. We require the perturbation datum to satisfy the following conditions: for $r$ close to $0$, $K^r = 0, I^r = I_0$, for $r$ large, $(K^r,I^r) = (H^j_t\,dt,J^j_t)$ for $\epsilon_js \in [0,r]$, $(K^r,I^r) = (0,I_0)$ for $\epsilon_js \in [r+1,\infty)$, and interpolating between the two for $\epsilon_js \in (r,r+1)$ on the $j$-th end; here $\epsilon_j$ is the sign of the $j$-th end, that is $\epsilon_j = -1$ for $j=0,1$ and $\epsilon_2 = 1$. We let
$$\cM_\cS(K,I) = \{(r,u) \in [0,\infty) \times C^\infty(D^2,S^1;M,L)\,|\,\ol\partial_{K^r,I^r}u = 0\}\,.$$
This is a smooth manifold of local dimension at $u$ equal to $n + \mu(u) + 1$.

We have the evaluation map
$$\ev \fc (\wt\cM(L,I_0))^{k_0+k_1+k_2} \times \cM_\cS(K,I) \to L^{2(k_0+k_1+k_2)+3}$$
defined by
\begin{multline*}
\ev(U=(u^0,u^1,u^2;u)) = (u^0_1(-1),\dots,u^0_{k_0}(1),u(1); u^1_1(-1),\dots,u^1_{k_1}(1),u(e^{2\pi i/3});\\ u(e^{4\pi i/3}),u^2_1(-1),\dots,u^2_{k_2}(1))\,.
\end{multline*}
For $q_i \in \Crit f_i$ we let
$$\wt\cP_{k_0,k_1,k_2}(K,I;q_0,q_1;q_2) = \ev^{-1}(\cU_{f_0}(q_0) \times Q^{k_0}_{f_0,\rho} \times \cU_{f_1}(q_1) \times Q^{k_1}_{f_1,\rho} \times Q^{k_2}_{f_2,\rho} \times \cS_{f_2}(q_2))\,.$$
We let $\cP_{k_0,k_1,k_2}(K,I;q_0,q_1;q_2)$ be the quotient of this space by the natural $\R^{k_0+k_1+k_2}$-action, and put
$$\wt \cP(K,I;q_0,q_1;q_2) = \bigcup_{k_0,k_1,k_2\geq 0}\wt \cP_{k_0,k_1,k_2}(K,I;q_0,q_1;q_2)\;,$$
$$\cP(K,I;q_0,q_1;q_2) = \bigcup_{k_0,k_1,k_2\geq 0} \cP_{k_0,k_1,k_2}(K,I;q_0,q_1;q_2)\,.$$
The space $\cP(K,I;q_0,q_1,q_2)$ is a smooth manifold of local dimension at $U$ equal to $|q_0| + |q_1| - |q_2| + \mu(U) + 1 - n$ whenever this number does not exceed $1$. In order to define the chain homotopy, we need to orient the operator $D_U|_{Y_\Gamma^U}$ if $U \in \wt\cP_{k_0,k_1,k_2}(K,I;q_0,q_1;q_2)$ is such that $|q_0| + |q_1| - |q_2| + \mu(U) + 1 - n= 0$. The exact Fredholm triple $0 \to D_U|_{Y_\Gamma^U} \to D_{r,U}|_{T_r[0,\infty)\times Y_Q^U} \to T_r[0,\infty) \times T_UQ/\Gamma \to 0$ gives rise to the isomorphism
$$\ddd(D_{r,U}|_{T_r[0,\infty) \times Y_Q^U}) \simeq \ddd(D_U|_{Y_\Gamma^U}) \otimes \ddd(T_UQ/\Gamma)\,.$$
We orient $D_U|_{Y_\Gamma^U}$ by the orientation corresponding to the orientation
$$\textstyle(-1)^{k_0+k_1}\bigwedge_i\partial_{U_i} \otimes \bigwedge_i e^U_i \wedge \partial_r\,.$$
Using the compactified $1$-dimensional moduli space, it is possible to show
\begin{thm}
The map $\Psi$ thus defined is a chain homotopy between $\PSS_{H^2,J^2}^{\cD_2} \circ \star \circ (\PSS_{\cD_0}^{H^0,J^0} \otimes \linebreak\PSS_{\cD_1}^{H^1,J^1})$ and $\star$. \qed
\end{thm}
\begin{prf}
Since for $r = 0$ the perturbation datum coincides with $(0,I_0)$, we are counting holomorphic disks and therefore indeed we obtain the Lagrangian quantum product at the corresponding end. \qed
\end{prf}

The proof that the PSS maps intertwine the two other algebraic structures, namely the product on $QH_*(M)$ and $HF_*(M)$ and the quantum module structure, is entirely analogous. The case of the quantum homology of $M$ is much simpler since there are no signs involved and there is no breaking/collision. In the case of the quantum module structure one needs to introduce the sign $(-1)^{k_0}$ when defining the chain homotopy, where $k_0$ is the number of disks in leg $0$ of the spiked pearly trajectory.

The upshot is that the PSS maps are isomorphisms between the quantum and Floer homologies, and they intertwine the algebraic structures on both sides, which include the product structures on the quantum and Floer homologies and the quantum module action. Note that the PSS maps automatically respect the unit elements, and that the inverse PSS maps respect the algebraic structures as well.

\section{Boundary operators in the presence of bubbling}\label{s:boundary_op_squares_zero_bubbling}

Here we show that the boundary operators on Lagrangian Floer and quantum homologies square to zero when there is bubbling present. By index considerations this only happens when $N_L = 2$.

\subsection{$\partial_{H,J}^2 = 0$}\label{ss:boundary_op_squares_zero_bubbling_HF}

Let $(H,J)$ be a regular Floer datum for $L$, where $N_L = 2$. We wish to show that $\partial_{H,J}^2 = 0$. Since the boundary operator is determined by its matrix elements, it suffices to show that all the matrix elements of its square vanish, that is whenever $\wt\gamma_\pm \in \Crit \cA_{H:L}$ are such that $|\wt\gamma_-| = |\wt\gamma_+| + 2$, then we have
$$\sum_{\substack{\wt\delta \in \Crit \cA_{H:L} \\ |\wt\delta| = |\wt\gamma_-| - 1}} \sum_{\substack{([u],[v]) \in \cM(H,J;\wt\gamma_-,\wt\delta) \times \\ \cM(H,J;\wt\delta,\wt\gamma_+)}} C(v) \circ C(u) = 0\,.$$
The Gromov compactification of $\cM(H,J;\wt\gamma_-,\wt\gamma_+)$ is obtained by adding two types of points, the usual Floer breaking, and bubbling, see \S\ref{sss:compactness_gluing}. A Maslov $2$ disk can bubble off only in the case $\gamma_- = \gamma_+ =: \gamma$, and it is attached to one of the ends points $q_i := \gamma(i)$, $i = 0,1$. In fact, such a bubble appears as follows: a sequence of Floer strips in $\wt \cM(H,J;\wt\gamma_-,\wt\gamma_+)$ degenerates into the $s$-independent strip $(s,t) \mapsto \gamma(t)$ and a Maslov two disk attached to it at one of its boundary components. The noncompact connected components of $\ol\cM(H,J;\wt\gamma_-,\wt\gamma_+)$ are thus subdivided into three types according to the number of ends corresponding to Floer breaking. The arguments used in proving that the boundary operator squares to zero when no bubbling is present, allow us to see that it is enough to show the vanishing of the sum
$$\sum C(v) \circ C(u)$$
where it is taken over all pairs $([u],[v])$ which appear as boundary points in those components of the compactified space $\ol\cM(H,J;\wt\gamma_-,\wt\gamma_+)$ whose other boundary point is a bubble.

We let $\wt\cM(J,q)$ be the space of $J$-holomorphic disks with boundary on $L$ passing through $q \in L$, and denote by $\cM(J,q)$ the quotient by the conformal automorphism group $\Aut(D^2,\partial D^2,1)$. Let $\Delta \subset \ol\cM(H,J;\wt\gamma_-,\wt\gamma_+)$ be a connected component whose boundary points are a broken Floer trajectory $([x],[y])$ and a holomorphic disk $[u] \in \cM(J_i,q_i)$ attached at one of the endpoints of $\gamma$. Recall that there is a bijection between isomorphisms $C(\wt\gamma_-) \simeq C(\wt\gamma_+)$ and orientations of $T_z\wt\cM(H,J;\wt\gamma_-,\wt\gamma_+)$, given by gluing, where $[z] \in \Delta$, and that moreover the isomorphism $C(y) \circ C(x)$ corresponds to the orientation $-\partial_z \wedge \text{inward}$. Similarly, using gluing, one can establish a bijection between isomorphisms $C(\wt\gamma_-) \simeq C(\wt\gamma_+)$ and orientations of $T_u\wt\cM(J_i,q_i)$. Namely, one has the following exact triple:
$$\xymatrix{0 \ar[r] & W^{1,p}(u)\sharp 0 \ar[r] \ar[d]^{D_u\sharp 0} & W^{1,p}(u) \sharp W^{1,p}(\wt\gamma_-)  \ar[r] \ar[d]^{D_u \sharp D_{\wt\gamma_-}} & W^{1,p}(\wt\gamma_-)  \ar[r] \ar[d]^{D_{\wt\gamma_-}} &  0 \\ 0 \ar[r] & L^p(u)  \ar[r] & L^p(u) \oplus L^p(\wt\gamma_-)  \ar[r] & L^p(\wt\gamma_-) \ar[r] & 0}$$
where $W^{1,p}(u)\sharp 0 = \{\xi \in W^{1,p}(u)\,|\,\xi(1) = 0\}$ and $D_u\sharp 0 = D_u|_{W^{1,p}(u)\sharp 0}$. This triple induces an isomorphism
$$\ddd(D_u\sharp 0) \otimes \ddd(D_{\wt\gamma_-}) \simeq \ddd(D_u\sharp D_{\wt\gamma_-})\,.$$
The latter operator can be identified with $D_{\wt\gamma_+}$, therefore using the fact that $\ker D_u \sharp 0 = T_u\wt\cM(J_i,q_i)$, we have obtained the desired bijection. By regularity, the operator $D_u$ is onto, and by genericity, so is $D_u \sharp 0$, and so its kernel is $2$-dimensional. The infinitesimal action of the isomorphism group induces a map $\C \to \ker D_{u}\sharp 0$, which is an isomorphism. Thus we see that the isomorphism $C(y) \circ C(x)$ corresponds to an orientation of $\C$. A computation shows that the isomorphism $\C \to \ker D_u\sharp 0$ induces the standard orientation on $\C$ if the bubble is attached at $\gamma(0)$, and the negative of the standard orientation on $\C$ if the bubble is attached at $\gamma(1)$. We see now that the sum
$$\sum C(y) \circ C(x)$$
over the boundary points of connected components of $\ol\cM$ whose other endpoints are bubbles equals the sum
$$\sum C(u)$$
running over those bubbles appearing as boundary points of these connected components, where $C(u)$ is the isomorphism $C(\wt\gamma_-) \simeq C(\wt\gamma_+)$ induced via gluing from the orientation on $u$ which is the standard orientation on $\C$ if $u$ is attached at $\gamma(0)$ and the negative of the standard orientation otherwise.

Fix a representative $\wh \gamma_- \in \wt\gamma_-$ and let $\{a(t)\}_{t\in[0,1]}$ be a parametrization of $\wh\gamma_-|_{\partial\dot D^2}$, up to the endpoints of $\gamma$, and define the following space:
$$\wt\cM(a) = \{(t,u)\,|\, t\in[0,1],\, u \in \wt\cM(J_t,a(t))\}\,.$$
This is a smooth $3$-dimensional manifold. Let also denote $\cM(a)$ its quotient by the automorphism group $\Aut(D^2,1)$; $\cM(a)$ is a smooth $1$-dimensional compact manifold with boundary. We claim that whenever $\Delta$ is a connected component of $\cM(a)$ with boundary points $(t,[u]),(t',[u'])$, then $C(u) + C(u')$ vanishes. A little combinatorial argument shows that this implies the vanishing of the above sum which only involves disks bubbling off at the ends of components whose other ends are broken Floer trajectories. We have therefore reduced the problem to the following claim: for any connected component $\Delta \subset \cM(a)$ the sum $C(u) + C(u')$ vanishes, where $u,u'$ represent the boundary points of $\Delta$.

So let indeed $\Delta$ be such a component and let $(t,[u]),(t',[u'])$ be its boundary points ($t,t' = 0,1$). Let $\{(t(\tau),[u(\tau)])\}_{\tau \in [0,1]}$ be a parametrization of $\Delta$ with $(t(0),[u(0)]) = (t,[u])$ and $(t(1),[u(1)]) = (t',[u'])$. Using a gluing argument as above, we get the isomorphism
$$\ddd(D_{u(\tau)}\sharp 0) \otimes \ddd(D_{\wt\gamma_-}) \simeq \ddd(D_{\wt\gamma_+})$$
which is continuous in $\tau$. This means that there is a bijection between isomorphisms $C(\wt\gamma_-) \simeq C(\wt\gamma_+)$ and orientations of the line bundle $\{\ddd(D_{u(\tau)}\sharp 0)\}_\tau$ over $[0,1]$. Using an obvious exact triple, we obtain the isomorphism
$$\ddd(D_{u(\tau)}\sharp 0) \otimes \ddd(\R) \simeq \ddd(T_{(t(\tau),u(\tau))}\wt\cM(a))$$
which is continuous in $\tau$. Therefore we have a bijection between the above isomorphisms and orientations of the component of $\wt\cM(a)$ above $\Delta$. We have the following isomorphism of short exact sequences for every $\tau$ for which $D_{u(\tau)} \sharp 0$ is surjective:
$$\xymatrix{0 \ar[r] & \C \ar[r] \ar[d] & T\wt\cM(a) \ar[r] \ar@{=}[d] & T\Delta \ar[r] \ar[d]^{\pi_*} & 0 \\ 0 \ar[r] & \ker D_{u(\tau)} \sharp 0 \ar[r] & T\wt\cM(a) \ar[r] & \R \ar[r] & 0}$$
where the left vertical arrow is the infinitesimal action of the automorphism group while the right vertical arrow is the differential of the map $\pi\fc (t,[u]) \mapsto t$. Orient $\C$, $\R$ with their standard orientations, and pick an orientation of $\Delta$. These induce an orientation on $T\wt\cM(a)$ and an orientation on $\ker D_{u(\tau)} \sharp 0$ via the short exact sequences, and it is easy to see that the left vertical arrow and the right vertical arrow are either both orientation-preserving or orientation-reversing.

Assume now that $\Delta$ connects two disks attached at the same point. Then $\pi_{*,(t,[u])}$ is orientation-preserving if and only if $\pi_{*,(t',[u'])}$ is orientation-reversing. This means that the isomorphism $\C \simeq \ker D_{u}\sharp 0$, which is the left vertical arrow, is orientation-preserving if and only if the isomorphism $\C \simeq \ker D_{u'}\sharp 0$ is orientation-reversing. This implies the following. Assume $u,u'$ are attached at $\gamma(0)$ and we have oriented both $\ker D_{u'}\sharp 0$ and $\ker D_{u} \sharp 0$ using the standard orientation on $\C$ and the above isomorphisms induced by the infinitesimal action. Then, since both left vertical arrows in the above diagram are orientation-preserving, we see that the orientations on $T\wt\cM(a)$ coming from the chosen orientations on $\ker D_{u'}\sharp 0$ and $\ker D_{u}\sharp 0$ and the standard orientation on $\R$ are opposite. This means that the isomorphisms $C(u),C(u')$ are opposite. A similar argument shows that when $\Delta$ connects two disks attached at $\gamma(1)$, or two disks attached at different points, we still have $C(u)+C(u') = 0$.

This finishes the proof that $\partial_{H,J}^2 = 0$ in the presence of bubbling.

\subsection{$\partial_\cD^2 = 0$}\label{ss:boundary_op_squares_zero_bubbling_QH}

The above compactification structure of the $1$-dimensional unparametrized pearly spaces can fail in case $q=q''$ and $N_L = 2$, and a bubble at $q$ can form. Here we prove that even if this happens, the quantum boundary operator still squares to zero.

The point is that the $1$-dimensional part $\cP^1(q,q)$ can no longer be compactified by adding Morse breaking alone (there is no holomorphic disk breaking or collision since we are in the minimal Maslov case), and has to be supplemented by adding disks of Maslov $2$ at $q$. Let us denote by $\ol\cP{}^1(q,q)$ the compact $1$-dimensional topological manifold with boundary obtained by adding the bubbles at $q$ and identifying two boundary components if they correspond to the same bubble. Pick a connected component $\Delta \subset \ol\cP{}^1(q,q)$ having two boundary points, both of which are pairs of the form $([u],[v]) \in \cP_0(q,q') \times \cP_1(q',q)$ or $([u],[v]) \in \cP_1(q,q') \times \cP_0(q',q)$, and such that there is a unique interior point represented by a bubble at $q$. Let $\partial \Delta = \{\delta,\delta'\}$ where $\delta = ([u],[v])$ and $\delta' = ([u'],[v'])$ with $u,v$ as described. We have to show that
$$C(v)\circ C(u) + C(v')\circ C(u') = 0$$
as a homomorphism $C(q,A) \to C(q,AB)$, where $B$ is the class of the bubble.

This is done as follows. Consider
$$\ev \fc C^\infty(D^2,S^1;M,L) \times S^1 \times S^1 \to L^2\,,\quad (w,\theta_1,\theta_2) \mapsto (w(\theta_1),w(\theta_2))\,.$$
For $(w,\theta_1,\theta_2) \in \ev^{-1}(\cU(q) \times \cS(q))$ we have the spaces
$$Y^{(w,\theta_1,\theta_2)} = \{(\xi,\Theta_1,\Theta_2) \in W^{1,p}(w) \times T_{\theta_1} S^1 \times T_{\theta_2} S^1\,|\, \xi(\theta_1) + \Theta_1 \in T_{w(\theta_1)}\cU(q)\,,\xi(\theta_2) + \Theta_2 \in T_{w(\theta_2)}\cS(q)\}$$
$$Y^{(w,\theta_1,\theta_2)}_0 = Y^{(w,\theta_1,\theta_2)} \cap (W^{1,p}(w) \times 0 \times 0)\,.$$
For $(w,\theta_1,\theta_2) \in \ev^{-1}(\cU(q) \times \cS(q))$ we have the natural Fredholm triple
$$\xymatrix{0 \ar[r] & Y^{(w,\theta_1,\theta_2)}_0 \ar[r] \ar[d]^{D_w|_{Y^{(w,\theta_1,\theta_2)}_0}} & Y^{(w,\theta_1,\theta_2)} \ar[r] \ar[d]^{D_w|_{Y^{(w,\theta_1,\theta_2)}}} & T_{(\theta_1,\theta_2)}(S^1 \times S^1) \ar[d] \ar[r] & 0 \\
0 \ar[r] & L^p(w) \ar@{=}[r] & L^p(w) \ar[r]& 0}$$
where the operators
$$D_w|_{Y^{(w,\theta_1,\theta_2)}_0} \quad \text{and} \quad D_w|_{Y^{(w,\theta_1,\theta_2)}}$$
are defined by pulling back $D_w \fc W^{1,p}(w) \to L^p(w)$ via the projections
$$Y^{(w,\theta_1,\theta_2)}_0 \to W^{1,p}(w) \quad \text{and} \quad Y^{(w,\theta_1,\theta_2)} \to W^{1,p}(w)\,.$$
This triple induces the canonical isomorphism
\begin{equation}\label{eqn:iso_detl_lines_proof_d_quantum_squares_zero}
\ddd(D_w|_{Y^{(w,\theta_1,\theta_2)}}) \simeq \ddd(D_w|_{Y^{(w,\theta_1,\theta_2)}_0}) \otimes \ddd(TS^1 \times TS^1)\,,
\end{equation}
which is continuous in $(w,\theta_1,\theta_2)$. Using a generalization of the construction in \S\ref{sss:boundary_op_Lagr_QH}, we obtain a canonical bijection between isomorphisms $C(q,A) \simeq C(q,AB)$ and orientations of $D_w|_{Y^{(w,\theta_1,\theta_2)}_0}$, which is continuous in $(w,\theta_1,\theta_2)$. Let us orient $S^1 \times S^1$ by the standard positive orientation. Then the isomorphism \eqref{eqn:iso_detl_lines_proof_d_quantum_squares_zero} yields a bijection between isomorphisms $C(q,A) \simeq C(q,AB)$ and orientations of $D_w|_{Y^{(w,\theta_1,\theta_2)}}$.

Let $\wt\cM(q,q) = \ev^{-1}(\cU(q) \times \cS(q)) \cap \wt\cM_2(J)$, where $\wt\cM_2(J)$ denotes the space of Maslov $2$ holomorphic disks with boundary on $L$. This is a smooth $4$-dimensional manifold. The group $\Aut(D^2)$ acts freely on it and we let $\cM(q,q)$ be the $1$-dimensional quotient. For $(w,\theta_1,\theta_2) \in \wt\cM(q,q)$ with $w$ in class $B$ we see that the operator $D_w|_{Y^{(w,\theta_1,\theta_2)}_0}$ is onto and has index $2$, the operator $D_w|_{Y^{(w,\theta_1,\theta_2)}}$ is onto, has index $4$ and its kernel equals $T_{(w,\theta_1,\theta_2)}\wt\cM(q,q)$. By the above, there is a bijection between isomorphisms $C(q,A) \simeq C(q,AB)$, orientations of $D_w|_{Y^{(w,\theta_1,\theta_2)}_0}$, and orientations of $T_{(w,\theta_1,\theta_2)}\wt\cM(q,q)$, that is orientations of the connected component of $\wt\cM(q,q)$ containing the bubble $z$ at $q$ in the form of the element $(z,1,1)$.

We therefore need to compute the orientations of $\wt\cM(q,q)$ corresponding to the isomorphisms $C(v)\circ C(u)$ and $C(v') \circ C(u)$. Let $w \in \wt\cP(q,q)$ and $w' \in \wt\cP(q,q)$ be obtained by gluing $u,v$ and $u',v'$, respectively. The disk $w$ defines an element of $\wt\cM(q,q)$ via $(w,-1,1)$ and we have an identification
$$\ker (D_w|_{Y_\Gamma^w}) = \ker(D_w|_{Y^{(w,-1,1)}_0})\,.$$
We know by \S\ref{sss:boundary_op_Lagr_QH} that the isomorphism $C(v) \circ C(u)$ corresponds to the orientation $-\partial_w \wedge \inward_\delta$ of $\ker (D_w|_{Y_\Gamma^w})$. Analogously we see that the isomorphism $C(v') \circ C(u')$ corresponds to the orientation $-\partial_{w'} \wedge \inward_{\delta'}$. These vectors, together with the standard orientation of $T_{(-1,1)}(S^1 \times S^1)$ by $\partial_{\theta_1} \wedge \partial_{\theta_2}$, induce the following orientations on $\wt\cM(q,q)$:
$$-\partial_w \wedge \inward_\delta \wedge \partial_{\theta_1} \wedge \partial_{\theta_2} \quad \text{and} \quad -\partial_{w'} \wedge \inward_{\delta'} \wedge \partial_{\theta_1} \wedge \partial_{\theta_2}\,,$$
and we have to show that these two orientations are opposite. We cannot compare them directly because the vector $\partial_w$ does not extend continuously pass the bubble.

We proceed as follows. We have the following exact sequence:
$$0 \to \Lie \Aut(D^2) \to T_{(w,-1,1)}\wt\cM(q,q) \to T_{[w,-1,1]}\cM(q,q) \to 0\,,$$
and similarly for $w'$. Orient $\Lie \Aut(D^2)$ by $\epsilon_1 \wedge \epsilon_2 \wedge \epsilon_3$ (see \S\ref{ss:general_technique_tori} for the definition of the vectors $\epsilon_i \in \Lie \Aut(D^2)$), and $T_{[w,-1,1]}\cM(q,q)$ by $\inward_\delta$. We see that the vector $\epsilon_2$ equals $\partial_w$ (by the definition of the latter, see \S\ref{sss:boundary_op_Lagr_QH}) while the projection $T\wt\cM(q,q) \to TS^1 \times TS^1$ maps $\epsilon_1$ to $\partial_{\theta_1} + \partial_{\theta_2}$ and $\epsilon_3$ to $\partial_{\theta_2}$. Therefore we obtain the induced orientation
$$\epsilon_1 \wedge \epsilon_2 \wedge \epsilon_3 \wedge \inward_\delta\,,$$
which equals the orientation
$$\partial_{\theta_1} \wedge \partial_w \wedge \partial_{\theta_2} \wedge \inward_\delta = - \partial_w \wedge \inward_\delta \wedge \partial_{\theta_1} \wedge \partial_{\theta_2}\,.$$
Since this orientation equals $\epsilon_1 \wedge \epsilon_2 \wedge \epsilon_3 \wedge \inward_\delta$, it extends continuously past the bubble, and equals
$$-\epsilon_1 \wedge \epsilon_2 \wedge \epsilon_3 \wedge \inward_{\delta'} = \partial_{w'} \wedge \inward_{\delta'} \wedge \partial_{\theta_1} \wedge \partial_{\theta_2}$$
at $(w',-1,1)$, as a similar computation shows. This implies that the orientations of $\wt\cM(q,q)$ induced by the isomorphisms $C(v) \circ C(u)$ and $C(v') \circ C(u')$ are opposite, which means that these isomorphisms themselves are opposite and the claim is proved.

\section{Quotient complexes}\label{s:quotient_cxs}

The above Floer and quantum complexes distinguish homotopy classes of cappings, and as such, at times they are too large. In applications it is often more convenient to work with smaller quotient complexes. In this section we describe how to construct such quotients in case $L$ is relatively $\Pin^\pm$. We also describe the familiar algebraic structures such as the module structure over a Novikov ring, which in this context turns out to be a nontrivial issue, due to the fact that the canonical complexes take into account the action of the fundamental group on the homotopy groups.

In \S\ref{ss:relative_Pin_structs_coh_ors_disks} we define relative Pin structures via \v Cech cochains and show how such a structure allows one to construct a coherent system of orientations on operators of the form $D_u \sharp 0$ where $u \fc (D^2,S^1) \to (M,L)$ is a smooth map with even Maslov index. \S\ref{ss:quotient_cxs_Ham_HF_QH_of_M} and \S\ref{ss:quotient_cxs_Lagr_HF_QH} deal with the construction of quotient complexes in the closed and open (Lagrangian) case, respectively.

\subsection{Relative Pin structures and coherent orientations for disks}\label{ss:relative_Pin_structs_coh_ors_disks}

We say that $L$ is \tb{relatively} $\Pin^\pm$ if $w^\pm \equiv w^\pm(L) \in \im(H^2(M;\Z_2) \to H^2(L;\Z_2))$, where $w^+(L) = w_2(TL)$ and $w^-(L) = w_2(TL) + w_1^2(TL)$. Not being relatively $\Pin^\pm$ is an obstructing for the existence of a relative $\Pin^\pm$-structure on $L$.
\begin{rem}
Assume $L$ is relatively $\Pin^\pm$ and let $w$ be such that $w|_L = w^\pm$. Since $H^1(S^2;\Z_2) = 0$, we see that $w_2(TL) \circ \partial = w^\pm \circ \partial$ as a map $\pi_3(M,L) \to \partial_2(L)$. It then factors as
$$\pi_3(M,L) \to \pi_2(L) \to \pi_2(M) \xrightarrow{w} \Z_2\,,$$
implying that $w_2(TL) \circ \partial = 0$, which means that being relatively $\Pin^\pm$ implies assumption \tb{(O)}. It is however a strictly weaker assumption, as the example of $\R P^5 \subset \C P^5$ shows: the map
$$H^2(\C P^5 ;\Z_2) \to H^2(\R P^5;\Z_2)$$
vanishes, the class $w_2(\R P^5) = w_2(\R P^5) + w_1^2(\R P^5)$ is nonzero, which means that $\R P^5$ is not relatively $\Pin^\pm$, however $\pi_2(\R P^5) = 0$, which means that assumption \tb{(O)} is satisfied in this case. \footnote{We thank J.-Y. Welschinger for this example.}

\end{rem}

Let us now describe the notion of a relative $\Pin^\pm$ structure on $L$ and a how choice of such a structure yields a system of coherent orientations for Cauchy--Riemann operators coming from disks.

We start with some generalities. This material is essentially contained in \cite{Wehrheim_Woodward_Orientations_pseudoholo_quilts}. Consider Lie groups $G,H$ and let $\phi \fc G \to H$ be a surjective Lie group homomorphism with finite abelian kernel $A = \ker \phi$. It is well-known that principal $H$-bundles on a smooth manifold $X$ are classified by the nonabelian cohomology $H^1(X;H)$. Let us recall the definition of $H^1(X;H)$. For an open cover $\cU = (U_i)_{i \in I}$ of $X$ we have the \v Cech cochain groups
$$C^k(\cU;H) = \prod_{(i_0,\dots,i_k)\in I^{k+1}}C^\infty(U_{i_0\dots i_k},H)\,,$$
where $k \geq 0$ and $U_{i_0\dots i_k} = U_{i_0}\cap \dots \cap U_{i_k}$. Define $\delta^k \fc C^k(\cU;H) \to C^{k+1}(\cU;H)$ by
$$(\delta c)_{i_0\dots i_{k+1}} = \prod_{j=0}^{k+1}\big(c_{i_0\dots\wh{i_j}\dots i_{k+1}}|_{U_{i_0\dots i_{k+1}}}\big)^{(-1)^j}\,.$$
The set $Z^1(\cU;H) = \ker \delta^1$ is the set of $1$-cocycles, that is $f = (f_{ij}) \in Z^1(\cU;H)$ if and only if for all $i,j,k \in I$ we have $f_{ik} = f_{ij}f_{jk}$. The group $C^0(\cU;H)$ acts on $Z^1(\cU;H)$ on the left as follows: for $c = (c_i) \in C^0(\cU;H)$ and $f = (f_{ij})$ we have $(c\cdot f)_{ij}=c_i f_{ij} c_j^{-1}$. We let
$$H^1(\cU;H) = Z^1(\cU;H)/C^0(\cU;H)\,.$$
Given another cover $\cV = (V_j)_{j \in J}$ a refinement map $\tau \fc \cV \to \cU$ by definition is a map $\tau \fc J \to I$ such that $V_j \subset U_{\tau(j)}$ for every $j \in J$. A refinement map induces homomorphisms $C^k(\cU;H) \to C^k(\cV;H)$ commuting with $\delta$ and therefore a well-defined map $H^1(\cU;H) \to H^1(\cV;H)$, which can be shown to be injective. Taking the direct limit over the set of covers directed by the relation of refinement, we obtain the first nonabelian cohomology $H^1(X;H)$. It can be shown that if $\cU$ is a good cover, then the canonical map $H^1(\cU;H) \to H^1(X;H)$ is a bijection.

There is a map from the set of isomorphism classes of principal $H$-bundles over $X$ to $H^1(X;H)$, defined by taking a sufficiently fine cover $\cU$, trivializing the bundle and taking the transition maps, which satisfy the cocycle relation, meaning they define an element in $Z^1(\cU;H)$. The corresponding class in $H^1(X;H)$ is well-defined, that is it only depends on the isomorphism class of the bundle. This map is a bijection. The inverse is obtained by gluing according to a cocycle.

Assume now that $Q$ is a principal $H$-bundle over $X$ and let $\cU$ be a good cover. Let $h = (h_{ij}) \in Z^1(\cU;H)$ be the transition cocycle corresponding to a trivialization of $Q$. Since $\cU$ is a good cover, every map $h_{ij} \fc U_{ij} \to H$ lifts to a map $g_{ij} \fc U_{ij} \to G$. The homomorphism $\phi$ induces in an obvious way homomorphisms $C^k(\cU;G) \to C^k(\cU;H)$ commuting with the differentials. Since clearly $\phi(g) = h$, we see that $\phi\delta g = \delta \phi g = \delta h = 0$, meaning $\delta g \in C^2(\cU;A)$. Since clearly $\delta(\delta g)= 0$, we see that $\delta g$ is in fact a cocycle in $Z^2(\cU;A)$. The corresponding class $[\delta g] \in H^2(X;A)$ (this is the ordinary second cohomology with coefficients in $A$) is well-defined, that is it only depends on the isomorphism class of $Q$, and is the characteristic class of $Q$. It vanishes if and only if $g$ can be corrected to a cocycle in $Z^1(\cU;G)$, meaning in this case the bundle $Q$ is in fact covered by a $G$-bundle $P$ via a $\phi$-equivariant map. A cocycle $g \in Z^1(\cU;G)$ with $\phi(g) = h$ is called a $G$-\tb{trivialization} of $h$. The set of such trivializations quotiented out by the equivalence relation induced by multiplication by $C^0(\cU;A)$ is the set of \tb{$G$-structures} on the bundle $Q$ relative to the cover $\cU$. Since it is a good cover, we get the same notion if we allow arbitrary covers and take limits over refinement maps. It can be seen that the set of $G$-structures on $Q$ is a torsor over the group $H^1(X;A)$.

Now we consider relative $G$-structures. Let $f \fc X \to Y$ be a smooth map and let $\cU$, $\cV$ be good covers of $X$, $Y$, respectively, where $\cU$ is a refinement of $f^{-1}\cV$. Let $Q \to X$ be an $H$-bundle and assume $h = (h_{ij}) \in Z^1(\cU;H)$ is a transition cocycle for $Q$. A \tb{$G$-trivialization on $Q$ relative to $f$} is a pair $(g,b) \in C^1(\cU;G) \times Z^2(\cV;A)$ such that $\phi(g) = h$ and $f^*b = \delta g$, that is $g$ is a lift of $h$ to $G$ and $b$ is a cocycle on $Y$ which pulls back to a cocycle on $X$ which is a boundary and in fact is the boundary of $g$. Two relative trivializations $(g,b)$ and $(g',b')$ are called \tb{equivalent} if there is $(a,a') \in C^0(\cU;A) \times C^1(\cV;A)$ with
$$(\delta a\cdot g,\delta a' \cdot b) = (g',b')\,.$$
An equivalence class of relative $G$-trivializations is called a \tb{relative $G$-structure on $Q$} (relative to $f$). It is a torsor over the group $H^1(f;A)$, defined as follows. Let $Z^1(f;A)$ be the group consisting of pairs $(a,a') \in C^1(\cU;A) \times C^2(\cV;A)$ such that $\delta a' = 0$ and $f^*a' = \delta a$. Now quotient it out by the subgroup of relative coboundaries whose elements are pairs $(\delta c\cdot c',\delta c')$ for $(c,c') \in C^0(\cU;A) \times C^1(\cV;A)$.

A \tb{relative $\Pin^\pm$-structure on $L$} is a relative $\Pin^\pm$-structure on $TL$ relative to the embedding $L \hookrightarrow M$. To spell it out, let $\cV$ be a good cover of $M$ and $\cU$ a good cover of $L$ which is a refinement of $\cV|_L$. Trivialize $TL$ to get transition functions $h = (h_{ij}) \in Z^1(\cU;O(n))$ (here we are using an auxiliary Riemannian metric on $L$). A trivialization of $h$ relative to the embedding $L \to M$ is a lift of $h$ to $g = (g_{ij}) \in C^1(\cU;\Pin^\pm(n))$ and a cocycle $b \in Z^2(\cV;\Z_2)$ on $M$ such that $\delta g = b|_L$. A relative $\Pin^\pm$-structure is an equivalence class of such trivializations. The set of relative $\Pin^\pm$-structures is a torsor over the group $H^1(\iota;\Z_2) = H^2(M,L;\Z_2)$ where $\iota \fc L \to M$.

Note that if we have a commutative diagram of smooth manifolds and maps
$$\xymatrix{X' \ar[r]^{f'} \ar[d] & Y'\ar[d] \\ X \ar[r]^f & Y}$$
then a relative $G$-structure on an $H$-bundle $Q$ on $X$ relative to $f$ canonically induces a relative $G$-structure on the pullback $H$-bundle $Q' \to X'$ by the map $X' \to X$ relative to the map $f'$.

We note the following obvious fact.
\begin{lemma}
Given a vector bundle $V \to S^1$, the canonical map sending $\Pin^\pm$-structures on $V$ to relative $\Pin^\pm$-structures on $V$ relative to the embedding $S^1 \to D^2$ is a bijection. This bijection is equivariant with respect to the natural actions of $H^1(S^1;\Z_2)$ and $H^2(D^2,S^1;\Z_2)$, connected by the boundary morphism $H^1(S^1;\Z_2) \to H^2(D^2,S^1;\Z_2)$, which is an isomorphism. \qed
\end{lemma}

\begin{coroll}
Assume a relative $\Pin^\pm$-structure on $L$ is given. Then for any disk $u \fc (D^2,S^1) \to (M,L)$ the bundle $F_u = (u|_{S^1})^*TL$ acquires a canonical $\Pin^\pm$-structure.
\end{coroll}
\begin{prf}
The map $u$ induces a relative $\Pin^\pm$ structure on $F_u$, which by the lemma corresponds to a unique $\Pin^\pm$-structure. \qed
\end{prf}

We will now prove the following result.
\begin{prop}\label{prop:canonical_oris_disks_Pin_struct}
A relative $\Pin^\pm$-structure on $L$ determines a system of orientations of the operators $D_u \sharp 0$ over the space of smooth disks $u$ with even Maslov number, which is coherent with respect to boundary gluing.
\end{prop}
\begin{prf}
Let $(E^0,F^0) \to (D^2,S^1)$ be a Hermitian bundle pair and assume $F^0$ has even Maslov number. Let $E^1 \to \C P^1$ be a Hermitian bundle with Chern number $-\mu(F^0)/2$. Then we can glue $E^0$ and $E^1$ at $0 \in D^2$ (see \S\ref{ss:boundary_gluing}) to obtain a Hermitian bundle pair $(E,F)$ with zero Maslov number. We have
$$\ddd(D_{E,F}) \simeq \ddd(D_{E^0,F^0}) \otimes \ddd(D_{E^1})\,,$$
where we omitted the factor $\ddd(E^0_0)$ since it's canonically oriented. The operator $D_{E^1}$, being a Cauchy--Riemann operator on a closed Riemann surface, is canonically oriented, therefore we get a canonical isomorphism
$$\ddd(D_{E^0,F^0}) = \ddd(D_{E,F})\,.$$
Similarly we get a canonical isomorphism
$$\ddd(D_{E^0,F^0}\sharp 0) = \ddd(D_{E,F} \sharp 0)\,.$$
It remains to establish a bijection between $\Pin^\pm$-structures on $F$ and orientations of the latter operator. Fix a unitary trivialization of $E$ and denote by $F$ the resulting Lagrangian loop in $\C^n$. The set of $\Pin^\pm$-structures on $F$ has two points, and varying $F$, we obtain a double cover of the space of contractible loops in the Lagrangian Grassmannian of $\C^n$. On the other hand, the set of orientations of the operator $D_{\C^n,F}\sharp 0$ also has two points and thus defines another double cover over the same space. The two covers have the same $w_1$, therefore they are isomorphic. It remains to choose an isomorphism at a point in this space. We do this for the constant loop. The constant loop $F \to S^1$ produces the operator $D_{\C^n,F} \sharp 0$, which can easily be seen to be an isomorphism. Since $F$ is the constant loop, of the two $\Pin^\pm$-structures over it one is the trivial structure (corresponding to transition functions at the identity in $\Pin^\pm(n)$). We associate the trivial structure to the canonical positive orientation, and the other structure to the negative orientation. Now note that this is independent of the chosen unitary trivialization of $E$.

The case of smooth maps is now obtained by noting that $(E_u,F_u)$ is a Hermitian bundle pair with even Maslov, and applying the above construction.

The coherence follows from the fact that boundary gluing two trivial bundle pairs yields a trivial bundle pair; the boundary gluing isomorphism maps the orientations of the corresponding isomorphism operators via multiplication in $\Z_2 = \{\pm 1\}$, and it does the same with the $\Pin^\pm$-structures. \qed
\end{prf}

\subsection{Hamiltonian Floer homology and quantum homology of $M$}\label{ss:quotient_cxs_Ham_HF_QH_of_M}

We first describe the quotient complexes for the periodic orbit Floer homology $HF_*(M)$ and the corresponding quantum homology $QH_*(M)$ since in this case there is a canonical way of doing it.

Let us consider the Floer complex corresponding to a regular time-periodic Floer datum $(H,J)$:
$$CF_*(H) = \bigoplus_{\wt x \in \Crit \cA_H} C(\wt x)\,.$$
In order to be able to define a quotient complex of $CF_*(H)$, we need to understand the necessary identifications that go into forming such a quotient. Of course, we wish only to identify the spaces $C(\wt x)$, $C(\wt x')$ where $\wt x = [x,\wh x]$, $\wt x' = [x,\wh x']$ share the same orbit and the only difference is in the capping. Let $q = x(0)$ and let $A \in \pi_2(M,q)$ be such that $[\wh x' \sharp - \wh x] = A$. We symbolize this by writing $\wt x' = A\cdot \wt x$. Let us see what it means to identify $C(\wt x)$ and $C(\wt x')$.

Gluing $D_A$ and $D_{\wt x}$ at a point close to $q$ produces an exact triple
$$0 \to D_A \sharp D_{\wt x} \to D_A \oplus D_{\wt x} \to 0_{T_qM} \to 0\,,$$
where the penultimate arrow is the difference of the evaluation maps to $T_qM$. This, together with the direct sum isomorphism, yields the isomorphism
$$\ddd(D_A) \otimes \ddd(D_{\wt x}) \simeq \ddd(D_A \sharp D_{\wt x}) \otimes \ddd(T_qM)\,.$$
Recall (Lemma \ref{lem:D_A_has_canonical_orientation}) that the operators $D_A$ have a canonical orientation. This, together with the canonical orientation of $T_qM$, yields the isomorphism
$$\ddd(D_{\wt x}) \simeq \ddd(D_A \sharp D_{\wt x}) \simeq \ddd(D_{\wt x'})\,,$$
where the second isomorphism comes from deformation. This means that for any $\wt x,\wt x' \in p^{-1}(x)$ we have an isomorphism
\begin{equation}\label{eqn:can_iso_C_x_C_x_prime_Ham_HF}
C(\wt x) = C(\wt x')
\end{equation}
which is independent of any choices.

Next, assume that $x_\pm$ are two periodic orbits of $H$ and $\wt x_\pm,\wt x_\pm'\in p^{-1}(x_\pm)$. Assume that $|\wt x_-| = |\wt x_+| + 1$, $|\wt x_-'| = |\wt x_+'| + 1$ and that $u \fc \R \times S^1 \to M$ is a Floer cylinder so that $u \in \wt \cM(H,J;\wt x_-,\wt x_+)$ and also $u \in \wt \cM(H,J; \wt x_-', \wt x_+')$. We have the diagram
\begin{equation}\label{dia:can_isos_commute_C_u_Ham_HF}
\xymatrix{C(\wt x_-) \ar@{=}[d] \ar[r]^{C(u)} & C(\wt x_+) \ar@{=}[d] \\ C(\wt x_-') \ar[r]^{C(u)} & C(\wt x_+')}
\end{equation}
where the vertical equality signs denote the above canonical isomorphisms. We claim that this diagram commutes. This follows from the associativity of the direct sum and gluing isomorphisms. The nontrivial point here is that the gluing of the difference classes $A_\pm = [\wh x_\pm' \sharp - \wh x_\pm]$ happens at different points of $M$, namely $x_\pm(0)$, and that $A_+$ is the result of transfer of $A_-$ to $x_+(0)$ along the curve $u(\cdot, 0)$. The isomorphism $\ddd(D_{A_-}) \simeq \ddd(D_{A_+})$ induced by this transfer coincides with the obvious isomorphism coming from the canonical orientations of $\ddd(D_{A_\pm})$, and it is what makes the diagram commute.

Now we can describe the construction of a quotient complex. Recall that the second homotopy groups of $M$ assemble into a local system of abelian groups (see \S\ref{ss:arbitrary_rings_loc_coeffs} for the definition of a local system) over $M$ with the group $\pi_2(M,q)$ being attached to $q \in M$, and with homotopy classes of paths between $q,q'$ inducing group isomorphisms $\pi_2(M,q) \simeq \pi_2(M,q')$ (see \cite{Hatcher_AG}). We denote this local system by $\pi_2(M)$. The datum that goes into the definition of a quotient complex in this case is a local subsystem $G$ of $\pi_2(M)$, that is a subgroup $G_q < \pi_2(M,q)$ for each $q \in M$ such that the isomorphisms coming from paths in $M$ preserve these subgroups. Fix such a subsystem $G$. Note that $\pi_2(M)$, being a local system, is itself a groupoid, and as such it acts on $\wt \Omega$ over $\Omega$, in the following sense: for every $q \in M$ let $\wt \Omega_q = \{\wt x = [x,\wh x] \in \wt \Omega\,|\, x(0) = q\}$; this is a covering space of $\Omega_q = \{x \in \Omega\,|\, x(0) = q\}$, and $\pi_2(M,q)$ acts on $\wt\Omega_q$ by attaching spheres, and all of these transformations happen over $\Omega_q$. Since $G$ is a local subsystem, it too acts on $\wt\Omega/\Omega$. The quotient of this action is a covering space of $\Omega$:
$$\wt\Omega/G$$
where two cappings of the same orbit through $q$ are identified if their difference lies in $G_q$. We have the quotient map $\wt\Omega \to \wt\Omega/G$ and we denote the image of $\wt x$ by this map via $[\wt x]_G$. Note that the action of $G$ on $\wt \Omega$ restricts in an obvious manner to an action on $\Crit \cA_H$, and we let $\Crit \cA_H/G$ be the quotient. For $[\wt x]_G \in \Crit \cA_H/G$ we let
$$C([\wt x]_G)$$
be the limit of the direct system of modules $(C(\wt x))_{\wt x \in [\wt x]_G}$ connected by the isomorphisms \eqref{eqn:can_iso_C_x_C_x_prime_Ham_HF}. The \tb{quotient complex} then is
$$CF_*^G(H) = \bigoplus_{[\wt x]_G \in \Crit \cA_H/G} C([\wt x]_G)$$
as a module. The commutativity of the diagram \eqref{dia:can_isos_commute_C_u_Ham_HF} ensures that the boundary operator $\partial_{H,J}$ on $CF_*(H)$ descends to a boundary operator $\partial_{H,J}^G$ so that the quotient map
$$(CF_*(H,),\partial_{H,J}) \to (CF_*^G(H),\partial_{H,J}^G)$$
is a chain map.

The local system $\pi_2(M)$ has a natural subsystem $\pi_2^0(M)$ consisting of spheres of zero area, and therefore of zero Chern number, by monotonicity. Since we are ultimately interested in spectral invariants, we only use subsystems $G$ contained in $\pi_2^0(M)$. If $G$ is such a subsystem, the action functional $\cA_H \fc \wt \Omega \to \R$ descends to a functional
$$\cA_H \fc \wt\Omega/G \to \R\,.$$
Moreover, the $\Z$-grading of $CF_*(H)$ descends to a $\Z$-grading on $CF_*^G(H)$ in this case as well.

There are two typical local subsystems of $\pi_2^0(M)$. One is $\pi_2^0(M)$ itself. Note that this subsystem is the kernel of the local system morphism $c_1\fc \pi_2(M) \to \Z$. Therefore in case $c_1$ does not vanish on $\pi_2(M)$, the quotient space $\wt\Omega/\pi_2^0(M)$ inherits an action of the trivial system $\Z$. This is the familiar Novikov action. It also induces an action of $\Z$ on the corresponding Floer complex $CF^{\pi_2^0(M)}_*(H)$. Letting $t$ be the positive generator of $\Z$, we see that the complex $CF^{\pi_2^0(M)}_*(H)$ then becomes a module over the group ring $\Z[t,t^{-1}]$, which is the familiar Novikov module structure.

The other subsystem of $\pi_2^0(M)$ is the kernel of the Hurewicz morphism $G:=\ker(\pi_2(M) \to H_2(M;\Z))$. In this case the quotient space $\wt\Omega/G$ inherits an action of $H_2^S(M) = \im(\pi_2(M) \to H_2(M;\Z))$, and therefore the quotient Floer complex $CF_*^G(H)$ becomes a module over the group ring $\Z[H_2^S(M)]$, also familiar in Floer homology.

It is similarly checked that the continuation morphisms respect the identifications by $G$, which means that we have well-defined abstract Floer homology
$$HF_*^G(M)\,.$$

Next we consider the effect of the quotient construction on products. Let $(H_i,J_i)$, $i=0,1,2$ be regular time-periodic Floer data, $(K,I)$ a regular compatible perturbation datum on the thrice-punctured sphere, and consider the resulting product map on chain level:
$$*=*_{K,I} \fc CF_*(H_0) \otimes CF_*(H_1) \to CF_*(H_2)\,.$$
We claim that this product descend to a well-defined product on the quotient complexes, to wit:
$$* \fc CF_*^G(H_0) \otimes CF_*^G(H_1) \to CF_*^G(H_2)$$
is well-defined and intertwines the quotient maps $CF_* \to CF_*^G$. To show this, it is enough to show the following. Let $\wt x_i, \wt x_i' \in \Crit \cA_{H_i}$ be such that $|\wt x_0| + |\wt x_1| - |\wt x_2| = 2n$, $|\wt x_0'| + |\wt x_1'| - |\wt x_2'| = 2n$, and let $u \in \cM(K,I;\{\wt x_i\}_i)$ and $u \in \cM(K,I;\{\wt x_i'\}_i)$. Then the diagram
$$\xymatrix{C(\wt x_0) \otimes C(\wt x_1) \ar@{=}[d] \ar[r]^-{C(u)} & C(\wt x_2) \ar@{=}[d] \\ C(\wt x_0') \otimes C(\wt x_1') \ar[r]^-{C(u)} & C(\wt x_2')}$$
commutes. The commutativity of this diagram is proved similarly to that of the diagram \eqref{dia:can_isos_commute_C_u_Ham_HF} for the boundary operator. Its commutativity is ensured by the fact that the operators $D_A$, $A \in \pi_2(M)$, are all canonically oriented, and that the manifold $M$ is oriented. Note as well that it is in defining the product operation on the quotient complexes that we really use the group structure on $G$: indeed, assuming $A_i = [\wh x_i' \sharp - \wh x_i]$, we see that for $\wt x_i$ to be $G$-equivalent to $\wt x_i'$, for all $i$, we need the product of $A_0$ and $A_1$ to lie in $G$, after transferring both of them to $x_2(0)$ along paths dictated by $u$, and we need this product to be equal $A_2$.

The reader will have no trouble checking that if $G = \pi_2^0(M)$ and $c_1|_{\pi_2(M)} \neq 0$, then the resulting homology $HF_*^G(M)$ becomes an algebra over the Novikov ring $\Z[t,t^{-1}]$, and if $G = \ker (\pi_2(M) \to H_2(M;\Z))$, then $HF_*^G(M)$ has the structure of an algebra over the ring $\Z[H_2^S(M)]$.

Almost identical arguments apply to quantum homology. Recall that the quantum complex for a datum $\cD = (f,\rho,J)$ is
$$QC_*(\cD) = \bigoplus_{q \in \Crit f}\bigoplus_{A \in \pi_2(M,q)}C(q,A)\,.$$
We have canonical identifications
$$C(q,A) = C(q,A')$$
for any $A,A' \in \pi_2(M,q)$. In fact, we have a canonical isomorphism
$$\ddd(D_A) \otimes \ddd(T\cS(q)) \simeq \ddd(D_A \sharp T\cS(q)) \otimes \ddd(T_qM)\,,$$
which shows that if we use the canonical orientations of $D_A$ and of $T_qM$, we have in fact an identification
$$C(q,0) = C(q,A)\,.$$
Thus a quotient complex corresponding to a local subsystem $G < \pi_2(M)$ can be formed, and we denote it $QC_*^G(\cD)$. In case $G < \pi_2^0(M)$, the quotient complex inherits a $\Z$-grading. The quotient map $QC_*(\cD) \to QC_*^G(\cD)$ is a chain map.

The chain-level product operation descends to quotient complexes as well, as do the units.

We note also that PSS maps $CF_*(H) \leftrightarrows QC_*(\cD)$ descend to well-defined maps on the quotient complexes, and they induce isomorphisms on homology. Therefore we have a well-defined abstract quantum homology $QH_*^G(M)$, which is canonically isomorphic to $HF_*^G(M)$ by the abstract PSS map. This is an isomorphism of supercommutative unital rings. In special cases of subsystems of $\pi_2^0(M)$ we get algebras over the corresponding Novikov rings.

\subsection{Lagrangian Floer and quantum homology}\label{ss:quotient_cxs_Lagr_HF_QH}

The Lagrangian case is significantly more involved due to the lack of canonical orientations of the operators $\ddd(D_A)$, $A\in \pi_2(M,L)$, and of the Lagrangian $L$ itself.

We first cover some preliminaries about the local systems involved. There is a natural local system of groups on $L$ given by $q \mapsto \pi_2(M,L,q)$, with natural isomorphisms $\pi_2(M,L,q) \simeq \pi_2(M,L,q')$ corresponding to homotopy classes of paths from $q$ to $q'$ \cite{Hatcher_AG}. We denote this local system by $\pi_2(M,L)$. In particular, looking at loops at $q$, we obtain an action of $\pi_1(L,q)$ on $\pi_2(M,L,q)$ by automorphisms. This action has the property that the boundary operator $\partial \fc \pi_2(M,L,q) \to \pi_1(L,q)$ intertwines it with the action of $\pi_1(M,q)$ on itself by conjugation. Moreover, if $A,B \in \pi_2(M,L,q)$ then $\partial(A)\cdot B = ABA^{-1}$. This local system acts, as a groupoid, on the space $\wt\Omega_L$, as follows. Let $q \in L$ and let $\Omega_{L,q} = \{\gamma \in \Omega_L \,|\, \gamma(0) = q\}$ and $\wt\Omega_{L,q} = p^{-1}(\Omega_{L,q})$ where $p \fc \wt\Omega_L \to \Omega_L$ is the projection. Then $p \fc \wt\Omega_{L,q} \to \Omega_{L,q}$ is a trivial covering space on which $\pi_2(M,L,q)$ acts simply transitively by attaching disks.

Let now $G < \pi_2(M,L)$ be a subsystem, that is $G_q$ is a subgroup of $\pi_2(M,L,q)$ and these subgroups are preserved by the isomorphisms induced by paths on $L$. In particular $G_q$ must be preserved by the action of elements of the form $\partial(A)$ for $A \in \pi_2(M,L,q)$, which by the above means that $G_q$ is normal in $\pi_2(M,L,q)$. This $G$ then acts on $\wt\Omega_L$ and we let
$$\wt \Omega_L/G$$
be the quotient. This is a covering space of $\Omega_L$ and it has an inherited action of the quotient system $\pi_2(M,L)/G$. If $H$ is a Hamiltonian on $M$, the action functional $\cA_{H:L}$ will descend to $\wt\Omega_L/G$ provided $G$ consists of disks of zero area, which is equivalent to saying that it is a subsystem of the local system $\pi_2^0(M,L)$ which is just the kernel of the local system morphism $\mu \fc \pi_2(M,L) \to \Z$ which is the Maslov index, or equivalently of the area morphism $\omega$, due to monotonicity. Henceforth we will only consider subsystems of $\pi_2^0(M,L)$.

Fix therefore a subsystem $G < \pi_2^0(M,L)$ and consider the quotient covering space $\wt\Omega_L/G \to \Omega_L$. Fix a regular Floer datum $(H,J)$. The set $\Crit \cA_{H:L}$ inherits an action of $\pi_2(M,L)$ in an obvious way, and therefore we can consider the quotient $\Crit \cA_{H:L}/G$, which naturally maps onto the set of Hamiltonian orbits of $H$ with boundary on $L$. We let $[\wt\gamma]_G$ be the image in $\Crit \cA_{H:L}/G$ of a point $\wt\gamma \in \Crit\cA_{H:L}$. Recall the Floer complex
$$CF_*(H:L) = \bigoplus_{\wt\gamma \in \Crit \cA_{H:L}}C(\wt\gamma)\,.$$
We wish to construct a quotient complex
$$CF_*^G(H:L) = \bigoplus_{[\wt\gamma]_G \in \Crit \cA_{H:L}/G}C([\wt\gamma]_G)$$
similarly to the periodic orbit case. In order to do so we need to construct identifications $C(\wt\gamma) = C(\wt\gamma')$ if $[\wt\gamma]_G = [\wt\gamma']_G$. Let us see what goes into such an identification.

Let $A = [\wh\gamma' \sharp - \wh\gamma] \in \pi_2(M,L,q)$ where $q = \gamma(0)$. Recall that $D_A \sharp 0$ denotes the restriction of $D_A$ to the subspace $\{\xi\in W^{1,p}(A)\,|\, \xi(1) = 0\}$. We have the operator $D_A\sharp D_{\wt\gamma}$ obtained by boundary gluing. Recall that this operator is just the restriction of $D_A \oplus D_{\wt\gamma}$ to the subspace of $W^{1,p}(D) \oplus W^{1,p}(\wt\gamma)$ consisting of pairs of sections agreeing at the point $q$. Thus the incidence condition, in the sense of \eqref{eqn:iso_abstract_deformation_incidence_condition_W}, is given by the diagonal $\Delta_{T_qL} \subset T_qL \oplus T_qL$. On the other hand, the family $D_A \sharp 0\oplus D_{\wt\gamma}$ is the restriction of $D_A \oplus D_{\wt\gamma}$ to the subspace where the incidence condition is $0 \oplus T_qL \subset T_qL \oplus T_qL$. The isomorphism \eqref{eqn:iso_abstract_deformation_incidence_condition_W} yields
$$\ddd(D_A \sharp 0 \oplus D_{\wt\gamma}) \simeq \ddd(D_A \sharp D_{\wt\gamma})\,,$$
which combined with the direct sum and deformation isomorphisms, gives us finally
$$\ddd(D_A \sharp 0) \otimes \ddd(D_{\wt\gamma}) \simeq \ddd(D_A \sharp D_{\wt\gamma}) \simeq \ddd(D_{\wt\gamma'})\,.$$
This means that isomorphisms $C(\wt\gamma) \simeq C(\wt\gamma')$ are in a natural bijection with orientations of the operator $D_A \sharp 0$. In order to form a quotient complex as above, we therefore need to choose an orientation of this operator. Recall Proposition \ref{prop:canonical_oris_disks_Pin_struct} which states that a choice of a relative $\Pin^\pm$-structure on $L$ determines a system of orientations of the families $D_A\sharp 0$ for $A$ varying in $\pi_2^0(M,L,q)$ for all $q$. Let us therefore assume that $L$ is relatively $\Pin^+$ or relatively $\Pin^-$ and fix a relative $\Pin^\pm$-structure on it, and endow all the lines $\ddd(D_A \sharp 0)$ with the corresponding orientations.

\begin{rem}
Even though more precise conditions on the existence of such a system of orientations required for the construction of quotient complexes can be formulated, in applications it is enough to limit oneself to the relatively Pin case, which is what we do here.
\end{rem}

Therefore we have isomorphisms $C(\wt\gamma) \simeq C(A\cdot\wt\gamma)$ for all $\wt\gamma \in \Crit \cA_{H:L}$ and so we can define
$$C([\wt\gamma]_G)$$
as the limit of the direct system of modules $(C(\wt\delta))_{\delta \in [\wt\gamma]_G}$ and above isomorphisms. The quotient module then is
$$CF_*^G(H:L) = \bigoplus_{[\wt\gamma]_G \in \Crit \cA_{H:L}/G}C([\wt\gamma]_G)\,.$$
For the boundary operator to descend it is enough to require that the diagram
$$\xymatrix{C(\wt\gamma_-) \ar[r]^{C(u)} \ar[d] & C(\wt\gamma_+) \ar[d] \\ C(\wt\gamma_-') \ar[r]^{C(u)} & C(\wt\gamma_+')}$$
commute for all $\wt\gamma_\pm,\wt\gamma_\pm' \in \Crit \cA_{H:L}$ with $|\wt\gamma_-| = |\wt\gamma_+| + 1$, $|\wt\gamma_-'| = |\wt\gamma_+'| + 1$, and $u \in \wt\cM(H,J;\wt\gamma_-,\wt\gamma_+)$, $u \in \wt\cM(H,J;\wt\gamma_-',\wt\gamma_+')$. The coherence of the chosen system of orientations with respect to boundary gluing ensures the commutativity of the diagram. 

Therefore the boundary operator $\partial_{H,J}$ descends to a well-defined boundary operator $\partial_{H,J}^G$ on $CF_*^G(H:L)$ so that the quotient map
$$(CF_*(H:L),\partial_{H,J}) \to (CF_*^G(H:L),\partial_{H,J}^G)$$
is a chain map. In any case we have a well-defined Floer homology $HF_*^G(H,J:L)$.

Next, the coherence of the system of orientations implies that we have well-defined continuation maps. Thus we have a well-defined abstract Floer homology $HF_*^G(L)$.

Lastly, we wish to endow the quotient complexes with a product structure coming from the product operation $\star$. Therefore fix regular Floer data $(H^i,J^i)$ associated to the punctures of the thrice-punctured disk used to define $\star$, and a regular compatible perturbation datum $(K,I)$. Again, due to coherence, diagrams of the following kind:
$$\xymatrix{C(\wt\gamma_0)  \otimes C(\wt\gamma_1) \ar[d] \ar[r]^-{C(u)} & C(\wt\gamma_2) \ar[d] \\ C(\wt\gamma_0') \otimes C(\wt\gamma_1') \ar[r]^-{C(u)}& C(\wt\gamma_2')}$$
commute, where $\wt\gamma_i,\wt\gamma_i' \in \Crit \cA_{H^i:L}$ are such that $A_i := [\wh\gamma_i'\sharp - \wh\gamma_i] \in G_{\gamma_i(0)}$, $|\wt\gamma_0| + \wt\gamma_1| - |\wt\gamma_2| = n$, $|\wt\gamma_0'| + \wt\gamma_1'| - |\wt\gamma_2'| = n$, and $u \in \cM(K,I;\{\wt\gamma_i\}_i)$, $u \in \cM(K,I;\{\wt\gamma_i'\}_i)$.

Transferring $A_0$ and $A_1$ to $\gamma_2(0)$ along the boundary of the disk $\wh\gamma_0\sharp \wh\gamma_1 \sharp u$ we see that $A_2$ must be equal to the product of these, which explains why we required $G_q$ to be a subgroup of $\pi_2(M,L,q)$ for every $q$.

We therefore obtain a well-defined product structure
$$\star \fc CF_*^G(H^0:L) \otimes CF_*^G(H^1:L) \to CF_*^G(H^2:L)\,.$$
It can be checked that this structure is intertwined by the quotient maps, and that the result is a well-defined product structure on the abstract Floer homology $HF_*^G(L)$, which becomes an associative unital ring.

Lastly we wish to comment on the usual structure of Floer homology as a module over Novikov rings. Firstly, let us point out that since $\pi_2(M,L)$ is a local system, that is a groupoid rather than a group, it does not make sense to speak of its group ring. One can define $H_2^D:= \im (\pi_2(M,L) \to H_2(M,L;\Z))$ and look at the corresponding group ring $\Z[H_2^D]$, which is the usual one appearing in Floer theory. We note that if one picks the local system $G = \ker (\pi_2(M,L) \to H_2(M,L;\Z))$, then the quotient space $\wt\Omega_L/G \to \Omega_L$ is a normal covering space with deck group canonically isomorphic to $H_2^D$, so in particular $H_2^D$ acts on $\Crit \cA_{H:L}/G$. But this is not enough to make $H_2^D$ act on the complex as this involves orienting operators $D_A\sharp 0$ for $A$ with nonzero Maslov, as explained above. This is possible if and only if $\mu(A)$ is even. This means that if $\mu(A)$ is odd, which of course only happens if $L$ is nonorientable, we simply cannot orient these operators in a coherent manner. Therefore the Floer complex $CF_*^G(H:L)$ is not, in general, a module over the Novikov ring $\Z[H_2^D]$. We illustrate this point in the example of $\R P^n \subset \C P^n$ for $n$ even, see \S\ref{ss:RPn}. Note however that since we do have a coherent system of orientations for classes $A$ with even Maslov, we see that the Floer complex $CF_*^G(H:L)$ is a module over the subring $\Z[H_2^{D,\text{even}}]$ where $H_2^{D,\text{even}}$ consists of classes with even Maslov number. We emphasize that this module structure depends on the chosen relative Pin structure.

We also consider the particular case of a class $A \in \pi_2(M,L,q)$ lying in the image of the morphism $\pi_2(M,q) \to \pi_2(M,L,q)$. In this case the operator $D_{A,q}$ is canonically oriented, because it can be considered as the result of gluing the operator $D_A$ where we view $A$ as a sphere attached at $q$, and the operator $D_0 \sharp 0$ where $0$ is the trivial class in $\pi_2(M,L)$. Since $\ddd(D_0 \sharp 0) = \ddd(0) \equiv \R$, we see that this operator is canonically oriented. This orientation also agrees with the one induced by any relative Pin structure, as can be seen from the multiplicative property. The conclusion is that convenient local systems $G < \pi_2^0(M,L)$ will contain the image of $\pi_2^0(M) \to \pi_2^0(M,L)$, because one has these canonical orientations. For the Novikov module structure, this has the following implication. Since $D_A \sharp 0$ is always canonically oriented for spherical classes $A\in \pi_2(M,L)$, the quotient Floer complex $CF_*^G(H:L)$ with $G = \ker (\pi_2(M,L) \to H_2(M,L;\Z))$ inherits the structure of a module over the group ring $\Z[H_2^S]$.

All of the above can be carried over to quantum homology verbatim. We consequently have quotient complexes $QC_*^G(\cD:L)$ with homology $QH_*^G(\cD:L)$. The PSS morphisms are well-defined isomorphisms $CF_*^G \simeq QH_*^G$, which allows us to define the abstract quantum homology $QH_*^G(L)$, which inherits the structure of a unital associative ring. In case $G$ is the kernel of the Hurewicz morphism, $QH_*^G(L)$ carries the structure of an algebra over the Novikov ring $\Z[H_2^S]$, and over the ring $\Z[H_2^{D,\text{even}}]$, and the two are compatible.

\subsection{Quantum module structure}\label{ss:quotient_cxs_quantum_module_struct}

Finally we wish to combine the two previous constructions and define a quantum module structure of Lagrangian Floer homology over the quantum homology of $M$. The necessary ingredients here are a local subsystem $G < \pi_2^0(M,L)$ containing $\im(\pi_2^0(M) \to \pi_2^0(M,L))$, and a relative Pin structure for $L$. The read will have no trouble checking that the quantum module operation
$$\bullet \fc CF_*(H) \otimes CF_*(H_0:L) \to CF_*(H_1:L)$$
descends to a well-defined operation
$$\bullet \fc CF_*^{G'}(H) \otimes CF_*^G(H_0:L) \to CF_*^G(H_1:L)$$
where $G' < \pi_2^0(M)$ is any local system such that the morphism $\pi_2(M,q) \to \pi_2(M,L,q)$ maps $G'_q$ into $G_q$ for $q \in L$.

\section{Examples and computations}\label{s:examples_computations}

Here we compute the canonical quantum complex for a number of Lagrangians.

\subsection{$\R P^n$ in $\C P^n$, $n \geq 2$}\label{ss:RPn}

We consider $M = \C P^n$ with the Fubini-Study form $\omega$ and $L = \R P^n$, which is monotone with minimal Maslov number $N_L = n+1$. Here we do not need to find the holomorphic disks in order to compute the homology, as we will see.

Choose a Morse function $f$ on $L$ with a unique critical point $q_i$ of index $i$, where $i = 0,\dots, n$. Assume that $J$ is an almost complex structure on $M$ such that $\cD = (f,\rho,J)$ is a regular quantum datum, where $\rho$ is a Riemannian metric for which there are precisely two gradient trajectories of $f$ from $q_{i+1}$ to $q_i$ for every $i < n$. The complex as a module is
$$QC_*(\cD:L) = \bigoplus_{i = 0}^n \bigoplus_{A \in \pi_2(M,L,q_i)}C(q_i,A)\,.$$
A computation shows that $\pi_2(M,L) \simeq \Z$, therefore the complex has rank $1$ in every degree. We see that $C(q_i,0)$ is generated by the orientations of $T\cS(q_i)$. As we know, the module $C(q_n,0)$ has a canonical generator, which we denote by $q_n$ by abuse of notation, corresponding to the positive orientation of $\ddd(T\cS(q_n)) \equiv \ddd(\R)$. We also let $q_i$, $i < n$, denote a generator of $C(q_i,0)$, by abuse of notation. By degree reasons, the boundary operator $\partial \fc C(q_{i+1},0) \to C(q_i,0)$ for $i< n$ coincides with the Morse boundary operator. We have therefore $\partial q_n = 0$. The Morse homology of $\R P^n$ in degree $n-1$ is $\Z_2$ if $n$ is even and $0$ if $n$ is odd; this forces the boundary operator $\partial \fc C(q_{n-1},0) = \Z\langle q_{n-1}\rangle \to C(q_{n-2},0) = \Z\langle q_{n-2}\rangle$ to send $q_{n-1} \mapsto \pm 2 q_{n-2}$, where the sign depends on the choice of generators. Since $0 = \partial^2q_{n-1} = \pm 2\partial q_{n-2}$, we see that $\partial q_{n-2} = 0$, that is the complex in degrees $n,n-1,n-2$ has the form
$$\dots \to \Z \langle q_{n} \rangle \xrightarrow{0} \Z\langle q_{n-1}\rangle  \xrightarrow{\pm 2} \Z \langle q_{n-2}\rangle \xrightarrow{0} \dots$$
Its homology therefore is $QH_{n-1}(\cD:L) \simeq 0$, $QH_{n-2}(\cD:L) \simeq \Z_2$. Since the quantum homology of $\C P^n$ is isomorphic to $\Z$ in every even degree, and all the homogeneous generators are invertible with respect to the quantum product, using the quantum module structure of $QH_*(\cD:L)$ over $QH_*(M)$, we see that the former is $2$-periodic with respect to degree, therefore we obtain finally
$$QH_*(L) \simeq \left \{ \begin{array}{ll}\Z_2\,, & n-*\text{ even} \\ 0\,, & n-*\text{ odd} \end{array}\right.\,.$$

From this we can compute the quantum homology of $\R P^n$ over an arbitrary ring $R$. Namely, the complex in degrees $n,n-1,n-2$ is
$$\dots \to R \langle q_n \rangle \xrightarrow{0} R \langle q_{n-1} \rangle \xrightarrow{2} R \langle q_{n-2} \rangle \xrightarrow{0} \dots $$
from which we see that $QH_{n-1}(L) \simeq \ker 2$ and $QH_{n-2} \simeq R /(2)$, where we view $2 \fc R \to R$ as the map corresponding to multiplication by $2$ and $(2) \subset R$ is the corresponding principal ideal. The rest of $QH_*(L)$ is obtained from the $2$-periodicity with respect to degree.

For instance if $2$ is invertible in $R$, we see that $QH_*(L) = 0$. In the other extreme where $R$ has characteristic $2$, we see that $QH_*(L) \simeq R$ in every degree.

When $n$ is even, $\R P^n$ is nonorientable. We see that in this case the quantum homology over $\Z$ is $2$-periodic in degree. It follows that the usual Novikov ring $\Z[t,t^{-1}]$ where $t$ is the Novikov parameter of degree $|t| = -N_L = -n-1$ does not act on this module, because the degree of $t$ is odd. This illustrates the point raised at the end of \S\ref{ss:quotient_cxs_Lagr_HF_QH} that in general the Floer complexes and the Floer homology of nonorientable Lagrangians are not modules over the full Novikov ring. We see, however, that the subring generated by the even powers of $t$ does in fact act on the complex and on the homology.

\subsection{A general technique for tori}\label{ss:general_technique_tori}

In order to efficiently compute the quantum complex, it makes sense to choose a basis for the underlying module in a suitable way. Here we present such a way in case $L$ is a Lagrangian torus.

We identify $L$ with the standard Euclidean torus $\T^n$ by means of a diffeomorphism. We note that $\T^n$ carries a canonical trivial $\Pin^\pm$-structure, which is determined by requiring the transition maps to equal the identity element of $\Pin^\pm(n)$. Proposition \ref{prop:canonical_oris_disks_Pin_struct} yields a coherent system of orientations for the operator families $D_A\sharp 0$ for $A \in \pi_2(M,L,q)$ and any $q$. Recall that the module $C(q,A)$ is generated by the orientations of $D_A\sharp T\cS(q)$. The canonical isomorphism
$$\ddd(D_A \sharp T\cS(q)) \simeq \ddd(D_A \sharp 0) \otimes \ddd(T\cS(q))$$
together with the orientation of $D_A \sharp 0$ gives us a canonical identification $C(q,A) = C(q,0)$. In order to remember $A$, we write this as
\begin{equation}\label{eqn:C_q_A_isomorphic_e_to_A_otimes_C_q_0}
C(q,A) = e^A \otimes C(q,0)\,.
\end{equation}
If we now fix an orientation for every $\cS(q)$, it gives us a generator of $C(q,0)$, which we denote by $q$ by abuse of notation. Therefore the complex becomes
$$QC_*(\cD:L) = \bigoplus_{q \in \Crit f} \Z[\pi_2(M,L,q)]\otimes q\,,$$
where we write elements of the group ring $\Z[\pi_2(M,L,q)]$ as sums $\sum_i c_i e^{A_i}$.

Assume now that the local system $\pi_2(M,L)$ is trivial, which means that the natural action of $\pi_1(L,q)$ on $\pi_2(M,L,q)$ is trivial for every $q$. The complex therefore can be written as
$$QC_*(\cD:L) = \Z[\pi_2(M,L)]\otimes \bigoplus_{q \in \Crit f} \Z\cdot q\,.$$
In this case, using the fact that $N_L$ is even, it is not hard to show that $\partial_Q$ is linear over the group ring $\Z[\pi_2(M,L)]$, in the sense that for any $q,q' \in \Crit f$ and any $A, B \in \pi_2(M,L)$ and any $u \in \wt\cP(q,q')$ with $u$ representing the class $B \in \pi_2(M,L,q)$, and where $|q| - |q'| + \mu(u) - 1 = 0$, we have the commutative diagram
$$\xymatrix{C(q,0) \ar[r]^{C(u)} \ar[d]^{e^A} & C(q',B) \ar[d]^{e^A} \\ C(q,A) \ar[r]^{C(u)} & C(q',AB)}$$
This means the following in terms of computing the quantum complex: it suffices to compute the isomorphisms $C(u) \fc C(q,0) \to C(q',B)$, and then the isomorphisms $C(u) \fc C(q,A) \to C(q',AB)$ are given by tensoring with $e^A$, that is using the isomorphism \eqref{eqn:C_q_A_isomorphic_e_to_A_otimes_C_q_0} above.

Another useful piece of information is as follows. Assume $u$ is a $J$-holomorphic disk of Maslov index $2$ and assume that the evaluation maps $\ev_\theta \fc \wt\cM(J) \to L$, $\ev_\theta(v) = v(\theta)$, for $\theta \in S^1$ have surjective differentials at $u$, that is $\ev_{\theta*,u} \fc \ker D_u \to T_{u(\theta)}L$ is onto for every $\theta$. The operator $D_u \sharp 0$ is surjective and has index $2$, therefore its kernel is $2$-dimensional. The infinitesimal action of the conformal automorphism group of $D^2$ preserving $1$ is an isomorphism $\C = \Lie \Aut(D^2,1) \to \ker D_u \sharp 0$, therefore using the canonical orientation of $\C$ we obtain an orientation of $D_u \sharp 0$. On the other hand $D_u \sharp 0$ is oriented by the canonical $\Pin^\pm$-structure on $L$.

\begin{lemma}
The canonical orientation on $D_u \sharp 0$ coming from the $\Pin$-structure coincides with the orientation of $\ker (D_u \sharp 0)$ coming from its identification with $\C$.
\end{lemma}
\begin{prf}
We note that the bundle pair $(E_u = u^*TM, F_u = (u|_{S^1})^*TL)$ splits into the direct sum of bundle pairs $(\im u_*,\im(u|_{S^1})_*)$ and a complement $(E',F')$. The latter bundle pair can be identified with $(\C^{n-1},\R^{n-1})$ and the corresponding operator on it is surjective with kernel isomorphic to $F'_1$. It follows that the restricted operator $D_{E',F'}\sharp 0$ is an isomorphism and thus it's oriented by the canonical positive orientation.

The remaining bundle pair $(\im u_*,\im(u|_{S^1})_*)$ has Maslov index $2$. It can be seen that the corresponding operator can be represented as the gluing of the Dolbeault operator on $\cO(1) \to \C P^1$ and the Dolbeault operator on the standard pair $(\C,\R)$. The orientation of $D_u\sharp 0$ corresponding to the canonical trivial $\Pin$-structure comes from the complex orientation of the kernel of $\ol\partial$ on $\cO(1) \to \C P^1$ restricted to the subspace of sections vanishing at $0 \in \C P^1$. It can be seen that this complex orientation coincides with the orientation induced by identifying $\ker (D_u \sharp 0) = \C$. The lemma is proved. \qed
\end{prf}

This result allows us to compute the isomorphism $C(u)$ coming from such a disk. Let therefore $q,q' \in \Crit f$ with $|q| = |q'| - 1$ and let $u \in \wt\cP(q,q')$ be a Maslov $2$ disk such that $(u(-1),u(1)) \in \cU(q) \times \cS(q')$ and assume that the evaluation maps from the space of $J$-disks to $L$ at points of $S^1$ have surjective differentials at $u$. The space $\ker D_u|_{X_\Gamma} = \{\xi \in \ker D_u\,|\, \xi(1) \in T_{u(1)}\cS(q')\}$ enters into the following two exact sequences:
$$0 \to \ker D_u|_{Y_\Gamma} \to \ker D_u|_{X_\Gamma} \to T_{u(-1)}L/T_{u(-1)}\cU(q) \to 0\,.$$
$$0 \to \ker D_u \sharp 0 \to \ker D_u|_{X_\Gamma} \to T_{u(1)}\cS(q') \to 0\,.$$
The first sequence yields
$$\ddd(D_u|_{X_\Gamma}) \simeq \ddd(D_u|_{Y_\Gamma}) \otimes \ddd(TL/T\cU(q)) = \ddd(D_u|_{Y_\Gamma}) \otimes \ddd(T\cS(q))$$
while the second one gives
$$\ddd(D_u|_{X_\Gamma}) \simeq \ddd(D_u \sharp 0) \otimes \ddd(T\cS(q'))$$
From the definitions in \S\ref{ss:Lagr_QH} it follows that the isomorphism $C(u) \fc C(q,0) \to C(q',A)$, where $A = [u]$, is obtained as follows: a generator of $C(q,0)$, that is an orientation of $T\cS(q)$,  together with the canonical orientation of $D_u|_{Y_\Gamma}$ by the infinitesimal action of $\R$ determine an orientation of $D_u|_{X_\Gamma}$ by the first sequence. The second sequence gives an orientation of $\ddd(D_u\sharp 0) \otimes \ddd(T\cS(q'))$, or equivalently an element of $C(q',A)$. Now using the canonical orientation of $D_u \sharp 0$ we obtain an orientation of $T\cS(q')$.

We now define a basis of the Lie algebra of the conformal automorphism group $\Aut(D^2)$: $\epsilon_1$ is the infinitesimal vector of the elliptic counterclockwise rotation about $0 \in D^2$; $\epsilon_2$ is the infinitesimal vector of the hyperbolic translation from $-1$ to $1$; $\epsilon_3$ is the infinitesimal vector of the parabolic rotation around $1$ which evaluates to $-i$ at $0 \in D^2$, so that it induces a counterclockwise translation on the boundary away from $1$. Note that $\epsilon_2,\epsilon_3$ form a \emph{negative} basis of $\Lie \Aut (D^2,1) = \C$. When we have a holomorphic disk $u$ with boundary on $L$, by abuse of notation we will denote the vector fields along $u$ corresponding to the infinitesimal actions of these vectors by the same letters.

We will apply this in case $L$ is a two-dimensional torus. There are two possibilities. The first one is when $q$ has index $1$ and $q'$ has index $2$, meaning it's a maximum. Let us choose an orientation of $\cS(q)$, which consists of a tangent vector field $v$ along $\cS(q)$. Assume $\epsilon_3(-1) = \epsilon v$. Then the orientation of $\ker D_u|_{X_\Gamma}$ induced by the first sequence is $\epsilon_2 \wedge \epsilon\epsilon_3$. We see that in this case $D_u|_{X_\Gamma} = D_u \sharp 0$ and therefore we just obtain $-\epsilon$ times the standard orientation. We see that in this case the isomorphism $C(u) \fc C(q,0) \to C(q',A)$, where $A = [u]$, sends $v$ to $-\epsilon e^A\otimes q'$, where $q'$ denotes the canonical positive orientation of $T\cS(q') = 0$.

The other possibility is that $q$ is a minimum and $q'$ has index $1$. In this case we choose an orientation of $TL = T\cS(q)$ by a pair of vectors $v_1,v_2$, for instance we can take the vectors corresponding to the chosen coordinates on $L$. We also choose an orientation of $\cS(q')$ by a nonvanishing tangent vector field $v$ along $\cS(q)$. In $\ker D_u|_{X_\Gamma}$ let $\eta$ be a vector such that $\epsilon_3(-1),\eta(-1)$ is a positive basis at $q$. Therefore the first sequence gives $\ker D_u|_{X_\Gamma}$ the orientation
$$\epsilon_2 \wedge \epsilon_3 \wedge \eta\,.$$
Now $\eta$ evaluates to $\epsilon v$ at $+1$. It follows from the second sequence that the induced orientation on $T\cS(q')$ is given by $-\epsilon v$. Thus $C(u) \fc C(q,0) \to C(q',A)$ maps $v_1\wedge v_2$ to $-\epsilon e^A \otimes v$.

We obtain from these considerations the following result.
\begin{lemma}\label{lem:computing_Lagr_QH_tori}
Let $L$ be a Lagrangian two-torus with $N_L = 2$ such that the local system $\pi_2(M,L)$ is trivial, and such that the evaluation maps from the space of Maslov $2$ disks at the points of $S^1$ all have surjective differentials. Let $f$ be a Morse function on $L$ with vanishing Morse boundary operator. Assume the critical points of $f$ are the maximum $q_2$, saddles $x,y$, and the minimum $q_0$. Orient the stable manifolds of $x,y$ somehow and orient $\cS(q_0)$ using the orientation coming from the identification $L \simeq \T^2$.
\begin{itemize}
\item Every unparametrized holomorphic disk $[u] \in \cM_1(J;q_2,A)$ yields contributions to the matrix elements of $\partial_\cD$, $C(x,0) \to C(q_2,A)$ and $C(y,0) \to C(q_2,A)$, as follows. The contribution to the matrix element $C(x,0) \to C(q_2,A) = e^A \otimes C(q_2,0)$ is $-1$ times the intersection number of $\cU(x)$ and the oriented parametrized curve $u|_{S^1}$, where $\cU(x)$ is cooriented by the chosen orientation of $\cS(x)$ and $u|_{S^1}$ is oriented by the counterclockwise orientation of $S^1$. The contribution to $C(y,0) \to C(q_2,A)$ is calculated similarly.
\item every unparametrized holomorphic disk $[u] \in \cM_1(J;q_0,A)$ yields the contribution to the matrix element $C(q_0,0) \to C(x,A) = e^A \otimes C(x,0)$ given by $-i$ times the chosen orientation of $\cS(x)$, where $i$ is the intersection number of the oriented curve $\cS(x)$ and the parametrized curve $u|_{S^1}$, cooriented by a vector $\eta \in \ker D_u|_{X_\Gamma}$ subject to the condition that $\epsilon_3(-1),\eta(-1)$ form an oriented basis at $q_0$. The same is true verbatim for the contribution to the matrix element $C(q_0,0) \to C(y,A) = e^A \otimes C(y,0)$. \qed
\end{itemize}
\end{lemma}

In case $L$ is a circle, the condition on the differentials of the evaluation maps is automatic, therefore we obtain
\begin{lemma}
Assume $L$ is a circle with $N_L = 2$, such that $\pi_2(M,L)$ is a trivial local system. Choose an orientation of $L$. Let $f$ be a Morse function on $L$ with maximum $q_1$ and minimum $q_0$. The contribution of $[u] \in \cM_1(J;q_1,A)$ to the matrix element $C(q_0,0) \to C(q_1,A)$ is given by $-1$ times the intersection number of the oriented parametrized curve $u|_{S^1}$ and the point $q_0$ cooriented by the orientation of $\cS(q_0)$. \qed
\end{lemma}

\subsection{$\R P^1$ in $\C P^1$}

Let $\cD = (f,\rho,J_0)$ be a quantum datum for $L = \R P^1$ where $f$ is a Morse function with a unique maximum $q_1$ and a unique minimum $q_0$, and $J_0$ is the standard complex structure on $\C P^1$. The local system $\pi_2(M,L)$ is trivial and the relative Hurewicz morphism $\pi_2(M,L) \to H_2(M,L;\Z)$ is an isomorphism. The latter group is generated by two classes $A,B$, defined as follows. Fix an orientation on $L$. The class $A$ is the class of an embedded disk representing a contraction of an orientation-preserving diffeomorphism $S^1 \to L$, while $B$ is represented by an embedded disk realizing the contraction of an orientation-reversing diffeomorphism $S^1 \to L$. For any $q \in L$ the spaces $\cM_1(J_0;q, A)$ and $\cM_1(J_0;q, B)$ each contain one point. Choose a basis of $T\cS(q_0)$ giving $\cS(q_0)$ the chosen orientation of $L$. Then Lemma \ref{lem:computing_Lagr_QH_tori} gives us the following formula for the boundary operator:
$$\partial_\cD(q_1) = 0\quad \text{and} \quad \partial_\cD(q_0) = (-e^A + e^B)q_1\,.$$

The homology of this canonical complex can be computed. To this end we identify the complex in every degree with a direct sum of a countable number of copies of $\Z$, meaning we can identify each $QC_i(\cD:L)$ with the space of functions $\Z \to \Z$ having finite support. The boundary operator acting, for instance, from $QC_0$ to $QC_{-1}$, can be identified with the following operator on functions. Let $\delta_j$ denote the function taking the value $1$ on $j$ and $0$ otherwise. Then $\partial(\delta_j) = \delta_{j+1} - \delta_j$. We see that $QH$ vanishes in even degrees. To compute $QH_{-1}$, for example, define the sum function $f \mapsto \sum_j f(j)$. We see that it is onto $\Z$ and that the kernel consists of functions having sum $0$, which is precisely the set of boundaries, as follows from the description of $\partial$ we just had. Thus the homology is isomorphic to $\Z$ in every odd degree. Thus

$$QH_*(L) \simeq \left\{ \begin{array}{ll} \Z\,, & *\text{ is odd} \\ 0\,, & *\text{ is even}\end{array}\right.\,.$$

There are two Spin structures on $L$. These allow us to form a quotient of $QC_*(\cD:L)$ by identifying all submodules in the same degree. The quotient complex then is isomorphic to $\Z$ in each degree. The induced boundary operator depends on the Spin structure: for one Spin structure it vanishes, meaning that $QH_*(L) \simeq \Z$ in every degree, while for the other Spin structure it is the multiplication by $2$ when going from an even to an odd degree, which means that $QH_*(L) \simeq \Z_2$ in odd degrees and $QH_*(L) = 0$ in even degrees.

\subsection{Some Lagrangian tori}

We compute the canonical quantum complex of three Lagrangian tori: the Clifford and the Che\-ka\-nov torus in $\C P^2$ and the exotic torus in $S^2 \times S^2$, based on our general technique of \S \ref{ss:general_technique_tori}. The only thing we need to know is the parametrizations of the boundary circles of $J_0$-holomorphic disks of Maslov $2$. It turns out that the conditions of Lemma \ref{lem:computing_Lagr_QH_tori} are satisfied for all the three tori and for $J_0$ being the standard complex structure.

Let us therefore list the disks of Maslov $2$ with boundary on the tori. Note that the local systems $\pi_2(M,L)$ are trivial and the relative Hurewicz morphism $\pi_2(M,L) \to H_2(M,L;\Z)$ is an isomorphism in all three cases.

\subsubsection{The Clifford torus in $\C P^2$}\label{sss:Cliff_torus}

The Clifford torus is the Lagrangian
$$L = \{[z_0:z_1:z_2] \,|\, |z_0| = |z_1| = |z_2|\} \subset \C P^2\,,$$
and we can identify
$$S^1 \times S^1 \simeq L \quad \text{via} \quad (e^{i\theta_0},e^{i\theta_1}) \mapsto [e^{i\theta_0}:e^{i\theta_1}:1]\,.$$
We endow $L$ with the trivial $\Pin^+$-structure \footnote{We can also choose the trivial $\Pin^-$-structure; for the present calculation it is immaterial.} corresponding to this identification. We have the $J_0$-holomorphic disks
$$u_0(z) = [z:1:1]\,,\quad u_1(z) = [1:z:1]\,, \quad u_2(z) = [1:1:z]\quad \text{for }z \in D^2\,,$$
and we let $A,B,C \in \pi_2(M,L)$ be the corresponding classes; these freely generate $\pi_2(M,L)$ (as an abelian group). We choose a Morse function $f$ with critical points $q_2,x,y,q_0$ such that $\cU(y)$ and $\cS(x)$ are both slight deformations of vertical curves $\theta_0 = \const$ while $\cS(y)$ and $\cU(y)$ are both slight deformations of horizontal curves $\theta_1 = \const$. The Morse boundary operator of $f$ vanishes. We orient $\cS(y)$ by $\partial_0 \equiv \partial_{\theta_0}$, $\cS(x)$ by $\partial_1 \equiv \partial_{\theta_1}$, and $\cS(q_0)$ by $\partial_0 \wedge \partial_1$.

The zero-dimensional part of the space $\cP(x,q_2)$ has two disks, one in each class $B,C$. By Lemma \ref{lem:computing_Lagr_QH_tori} we have to compute the intersection numbers of the parametrized boundaries of these disks with $\cU(x)$, cooriented by $\partial_1$. We can see that for the disk in $B$ this intersection number equals $1$ while for the disk in $C$ it equals $-1$, therefore by Lemma \ref{lem:computing_Lagr_QH_tori} we have
$$\partial x = (-e^B + e^C)q_2\,.$$

The zero-dimensional part of the space $\cP(y,q_2)$ has one disk in each of the classes $A,C$. The intersection number of the parametrized boundary of a disk in class $A$ with $\cU(y)$ cooriented by $\partial_0$ equals $1$ while the corresponding intersection number for a disk in class $C$ is $-1$, therefore by Lemma \ref{lem:computing_Lagr_QH_tori} we have
$$\partial y = (-e^A + e^C)q_2\,.$$

The zero-dimensional part of the space $\cP(q_0,x)$ has a disk in each of the classes $A,C$. Let us look at a disk $u \in \wt \cP(q_0,x)$ in class $A$. According to Lemma \ref{lem:computing_Lagr_QH_tori} we have to coorient its boundary by a vector $\eta$ such that $\epsilon_3(-1) = \partial_0,\eta(-1)$ form a positive basis. Therefore $\eta$ should be the vector $\partial_1$. It evaluates to the vector $\partial_1$ at $1$, therefore the intersection number of $u|_{S^1}$ cooriented by $\eta$ with $\cS(x)$ cooriented by $\partial_1$ is $1$. An analogous computation shows that the intersection number corresponding to a disk in class $C$ is $-1$.

The zero-dimensional part of $\cP(q_0,y)$ has a disk in each of the classes $B,C$, and the respective intersection numbers equal $-1, 1$, therefore by Lemma \ref{lem:computing_Lagr_QH_tori} we have
$$\partial q_0 = (-e^A + e^C)x + (e^B - e^C)y\,.$$

\subsubsection{The Chekanov torus in $\C P^2$}

We can view $\C P^2$ as the symplectic cut of the closed unit disk cotangent bundle $D^*\R P^2$ by the geodesic flow on the boundary relative to the round metric on $\R P^2$, with the symplectic form scaled by an appropriate factor. We let $C$ be the image of the unit cotangent bundle of $\R P^2$ by the quotient map; it is a smooth conic. If we fix a point $x \in \R P^2$, take a cotangent circle of radius $r$ at it, and let it flow with the geodesic flow, the resulting set $L \subset \C P^2$ is a monotone Lagrangian torus for a unique value of $r$. This is the Chekanov torus. We let $\alpha \in \pi_2(M,L)$ be the class of the cotangent disk at $x$, $\beta \in \pi_2(M,L)$ be the class of a disk contracting the geodesic circle whose intersection number with $C$ is $1$. Finally we let $h = [\C P^1] \in \pi_2(M,L)$ be the class of the line. The group $\pi_2(M,L)$ is freely generated by these three classes. It is known \cite{Chekanov_Schlenk_Notes_mon_twist_tori, Auroux_Mirror_symmetry_T_duality_compl_antican_divisor} that there are four classes containing Maslov $2$ disks: $\beta$ and $h - 2\beta + k\alpha$ for $k = -1,0,1$.

We put coordinates $(\theta_0,\theta_1)$ on $L$, where the boundary of $\alpha$ is given by $\theta_1 = \const$, while the boundary of $\beta$ is given by $\theta_0 = \const$. We choose a Morse function as in \S\ref{sss:Cliff_torus}. We orient the stable manifolds $\cS(x)$, $\cS(y)$, $\cS(q_0)$ by $\partial_1$, $\partial_0$, $\partial_0\wedge \partial_1$, respectively.

The zero-dimensional part of $\cP(x,q_2)$ contains seven disks, one in the class $\beta$, and two in each of the classes $h - 2\beta + k\alpha$. Performing a calculation of the corresponding intersection numbers as in \S\ref{sss:Cliff_torus}, we get
$$\partial x = (-e^\beta + 2(e^{h-2\beta-\alpha} +e^{h-2\beta} +e^{h-2\beta+\alpha}))q_2\,.$$

The zero-dimensional part of $\cP(y,q_0)$ has a disk in each of the classes $h - 2\beta \pm \alpha$, and the corresponding intersection numbers give us
$$\partial y = (e^{h-2\beta-\alpha} - e^{h-2\beta+\alpha})q_2\,.$$

The zero-dimensional part of $\cP(q_0,x)$ has a disk in each of the classes $h - 2\beta \pm \alpha$, while the zero-dimensional part of $\cP(q_0,y)$ has one disk in class $\beta$ and two disks in each of the classes $h - 2\beta +k\alpha$. We have:
$$\partial q_0 = (e^{h-2\beta-\alpha} - e^{h-2\beta+\alpha})x + (e^\beta - 2(e^{h-2\beta-\alpha} +e^{h-2\beta} +e^{h-2\beta+\alpha}))y\,.$$

\subsubsection{The exotic torus in $\C P^1 \times \C P^1$}

We identify $\C P^1 = S^2$, and view $S^2$ as the set of unit vectors in $\R^3$. The exotic torus is
$$L = \{(x,y) \in S^2 \times S^2\,|\, x\cdot y = -\tfrac 1 2\,, x_3 + y_3 = 0\}\,,$$
and put on it the coordinates $\theta_0, \theta_1$ where $\theta_0$ corresponds to the rotation around the $3$-axis in the counterclockwise direction while $\theta_1$ is the counterclockwise rotation of the pair $(x,y)$ around their sum $x + y$. We let the contraction of the curve $\{x_3 = -y_3 = \sqrt 3/2\}$ through the respective poles generate the class $\alpha$. The curve $\theta_0 = \const$ contracts via a disk whose class we denote by $\beta$. There are also the classes $A = [S^2 \times \pt]$ and $B = [\pt \times S^2]$ in $\pi_2(M,L)$. These four classes freely generate this abelian group.

It is known \cite{Chekanov_Schlenk_Notes_mon_twist_tori, Auroux_Mirror_symmetry_T_duality_compl_antican_divisor} that there are five classes in $\pi_2(M,L)$ containing Maslov $2$ disks, namely $\beta$, $A-\beta$, $B - \beta$, $A - \beta - \alpha$, and $B - \beta + \alpha$.

We choose a Morse function as above and orient its stable manifolds in the same manner. The zero-dimensional parts of the spaces $\cP(x,q_2)$, $\cP(q_0,y)$ each contain one disk in each one of the above five classes, while the spaces $\cP(y,q_2)$, $\cP(q_0,x)$ contain one disk in each of the classes $A - \beta - \alpha$, $B - \beta + \alpha$. We then have
$$\partial x = (-e^{\beta} + e^{A - \beta} + e^{B - \beta} + e^{A - \beta - \alpha} + e^{B - \beta + \alpha})q_2 \,, \quad \partial y = (e^{A - \beta - \alpha} - e^{B - \beta + \alpha})q_2\,,$$
$$\partial q_0 = (e^{A - \beta - \alpha} - e^{B - \beta + \alpha})x + (e^{\beta} - (e^{A - \beta} + e^{B - \beta} + e^{A - \beta - \alpha} + e^{B - \beta + \alpha}))y\,.$$

\bibliography{biblio}

\begin{thebibliography}{FOOO09b}

\bibitem[Abo13]{Abouzaid_Symp_H_Viterbo_thm}
Mohammed Abouzaid.
\newblock Symplectic cohomology and {V}iterbo's theorem.
\newblock 2013.
\newblock {\tt arXiv:1312.3354}.

\bibitem[Alb08]{Albers_Lagr_PSS_comparison_morphisms_HF}
Peter Albers.
\newblock A {L}agrangian {P}iunikhin--{S}alamon--{S}chwarz morphism and two
  comparison homomorphisms in {F}loer homology.
\newblock {\em Int. Math. Res. Not.}, (4):Art. ID rnm134, 56, 2008.

\bibitem[Aur07]{Auroux_Mirror_symmetry_T_duality_compl_antican_divisor}
Denis Auroux.
\newblock Mirror symmetry and {$T$}-duality in the complement of an
  anticanonical divisor.
\newblock {\em J. G\"okova Geom. Topol. GGT}, 1:51--91, 2007.

\bibitem[BC]{Paul_Octav_Private_comm_July_2015}
Paul Biran and Octav Cornea.
\newblock Private communication.
\newblock {J}uly 2015.

\bibitem[BC07]{Biran_Cornea_Quantum_structures_Lagr_submfds}
Paul Biran and Octav Cornea.
\newblock Quantum structures for {L}agrangian submanifolds.
\newblock 2007.
\newblock {\tt arXiv:0708.4221}.

\bibitem[BC09]{Biran_Cornea_Rigidity_uniruling}
Paul Biran and Octav Cornea.
\newblock Rigidity and uniruling for {L}agrangian submanifolds.
\newblock {\em Geom. Topol.}, 13(5):2881--2989, 2009.

\bibitem[BC12]{Biran_Cornea_Lagr_top_enumerative_invts}
Paul Biran and Octav Cornea.
\newblock Lagrangian topology and enumerative geometry.
\newblock {\em Geom. Topol.}, 16(2):963--1052, 2012.

\bibitem[Bir06]{Biran_Lag_non_intersections}
Paul Biran.
\newblock Lagrangian non-intersections.
\newblock {\em Geom. Funct. Anal.}, 16(2):279--326, 2006.

\bibitem[Cha12]{Charest_Source_spaces_perturbations_cluster_complexes}
Fran{\c c}ois Charest.
\newblock Source spaces and perturbations for cluster complexes.
\newblock 2012.
\newblock {\tt arXiv:1212.2923}.

\bibitem[CL05]{Cornea_Lalonde_Cluster_homology}
Octav Cornea and Fran{\c c}ois Lalonde.
\newblock Cluster homology.
\newblock 2005.
\newblock {\tt arXiv:math/0508345}.

\bibitem[CS10]{Chekanov_Schlenk_Notes_mon_twist_tori}
Yuri Chekanov and Felix Schlenk.
\newblock Notes on monotone {L}agrangian twist tori.
\newblock {\em Electron. Res. Announc. Math. Sci.}, 17:104--121, 2010.

\bibitem[Dju15]{Duretic_PSS_isos_spectral_invts_conormal_bundle}
Jovana Djureti\'c.
\newblock Piunikhin-salamon-schwarz isomorphisms and spectral invariants for
  conormal bundle.
\newblock 2015.
\newblock {\tt arXiv:1411.0852}.

\bibitem[EG91]{Eliashberg_Gromov_Convex_symp_mfds}
Yakov Eliashberg and Mikhael Gromov.
\newblock Convex symplectic manifolds.
\newblock In {\em Several complex variables and complex geometry, {P}art 2
  ({S}anta {C}ruz, {CA}, 1989)}, volume~52 of {\em Proc. Sympos. Pure Math.},
  pages 135--162. Amer. Math. Soc., Providence, RI, 1991.

\bibitem[FH93]{Floer_Hofer_Coherent_orientations}
A.~Floer and H.~Hofer.
\newblock Coherent orientations for periodic orbit problems in symplectic
  geometry.
\newblock {\em Math. Z.}, 212(1):13--38, 1993.

\bibitem[FHS95]{Floer_Hofer_Salamon_Transversality}
Andreas Floer, Helmut Hofer, and Dietmar Salamon.
\newblock Transversality in elliptic {M}orse theory for the symplectic action.
\newblock {\em Duke Math. J.}, 80(1):251--292, 1995.

\bibitem[Flo88a]{Floer_Morse_thry_Lagr_intersections}
Andreas Floer.
\newblock Morse theory for {L}agrangian intersections.
\newblock {\em J. Differential Geom.}, 28(3):513--547, 1988.

\bibitem[Flo88b]{Floer_unregularized_grad_flow_symp_action}
Andreas Floer.
\newblock The unregularized gradient flow of the symplectic action.
\newblock {\em Comm. Pure Appl. Math.}, 41(6):775--813, 1988.

\bibitem[Flo89]{Floer_Witten_cx_inf_dim_Morse_thry}
Andreas Floer.
\newblock Witten's complex and infinite-dimensional {M}orse theory.
\newblock {\em J. Differential Geom.}, 30(1):207--221, 1989.

\bibitem[FOOO09a]{FOOO_Lagr_intersection_Floer_thry_anomaly_obstr_I}
Kenji Fukaya, Yong-Geun Oh, Hiroshi Ohta, and Kaoru Ono.
\newblock {\em Lagrangian intersection {F}loer theory: anomaly and obstruction.
  {P}art {I}}, volume~46 of {\em AMS/IP Studies in Advanced Mathematics}.
\newblock American Mathematical Society, Providence, RI; International Press,
  Somerville, MA, 2009.

\bibitem[FOOO09b]{FOOO_Lagr_intersection_Floer_thry_anomaly_obstr_II}
Kenji Fukaya, Yong-Geun Oh, Hiroshi Ohta, and Kaoru Ono.
\newblock {\em Lagrangian intersection {F}loer theory: anomaly and obstruction.
  {P}art {II}}, volume~46 of {\em AMS/IP Studies in Advanced Mathematics}.
\newblock American Mathematical Society, Providence, RI; International Press,
  Somerville, MA, 2009.

\bibitem[Geo13]{Georgieva_Orientability_problem_open_GW}
Penka Georgieva.
\newblock The orientability problem in open {G}romov-{W}itten theory.
\newblock {\em Geom. Topol.}, 17(4):2485--2512, 2013.

\bibitem[Hat02]{Hatcher_AG}
Allen Hatcher.
\newblock {\em Algebraic topology}.
\newblock Cambridge University Press, Cambridge, 2002.

\bibitem[HL10]{Hu_Lalonde_Relative_Seidel_morphism_Albers_map}
Shengda Hu and Fran{\c{c}}ois Lalonde.
\newblock A relative {S}eidel morphism and the {A}lbers map.
\newblock {\em Trans. Amer. Math. Soc.}, 362(3):1135--1168, 2010.

\bibitem[HLL11]{HLL11}
Shengda Hu, Fran{\c{c}}ois Lalonde, and R{\'e}mi Leclercq.
\newblock Homological {L}agrangian monodromy.
\newblock {\em Geom. Topol.}, 15(3):1617--1650, 2011.

\bibitem[HS95]{Hofer_Salamon_HF_Nov_rings}
H.~Hofer and D.~A. Salamon.
\newblock Floer homology and {N}ovikov rings.
\newblock In {\em The {F}loer memorial volume}, volume 133 of {\em Progr.
  Math.}, pages 483--524. Birkh\"auser, Basel, 1995.

\bibitem[KM76]{Knudsen_Mumford_Projectivity_moduli_space_stable_curves_I_prelims_det_Div}
Finn~Faye Knudsen and David Mumford.
\newblock The projectivity of the moduli space of stable curves. {I}.
  {P}reliminaries on ``det'' and ``{D}iv''.
\newblock {\em Math. Scand.}, 39(1):19--55, 1976.

\bibitem[KM05]{Katic_Milinkovic_PSS_Lagr_intersections}
Jelena Kati{\'c} and Darko Milinkovi{\'c}.
\newblock Piunikhin--{S}alamon--{S}chwarz isomorphisms for {L}agrangian
  intersections.
\newblock {\em Differential Geom. Appl.}, 22(2):215--227, 2005.

\bibitem[KM09]{Katic_Milinkovic_coherent_oris_mixed_moduli_spaces}
Jelena Kati{\'c} and Darko Milinkovi{\'c}.
\newblock Coherent orientation of mixed moduli spaces in {M}orse-{F}loer
  theory.
\newblock {\em Bull. Braz. Math. Soc. (N.S.)}, 40(2):253--300, 2009.

\bibitem[KMS11]{Katic_Milinkovic_Simcevic_cohomology_rings_iso_HF_Morse}
Jelena Kati{\'c}, Darko Milinkovi{\'c}, and Tatjana Sim{\v{c}}evi{\'c}.
\newblock Isomorphism between {M}orse and {L}agrangian {F}loer cohomology
  rings.
\newblock {\em Rocky Mountain J. Math.}, 41(3):789--811, 2011.

\bibitem[Lec08]{Leclercq_spectral_invariants_Lagr_FH}
R{\'e}mi Leclercq.
\newblock Spectral invariants in {L}agrangian {F}loer theory.
\newblock {\em J. Mod. Dyn.}, 2(2):249--286, 2008.

\bibitem[LZ15]{Leclercq_Zapolsky_Spectral_invts_monotone_Lags}
R{\'e}mi Leclercq and Frol Zapolsky.
\newblock Spectral invariants for monotone {L}agrangians.
\newblock 2015.
\newblock {\tt arXiv:1505.07430}.

\bibitem[Oh93]{Oh_FH_Lagr_intersections_hol_disks_I}
Yong-Geun Oh.
\newblock Floer cohomology of {L}agrangian intersections and pseudo-holomorphic
  disks. {I}.
\newblock {\em Comm. Pure Appl. Math.}, 46(7):949--993, 1993.

\bibitem[Oh96a]{Oh_HF_spectral_sequences_Maslov_class}
Yong-Geun Oh.
\newblock Floer cohomology, spectral sequences, and the {M}aslov class of
  {L}agrangian embeddings.
\newblock {\em Internat. Math. Res. Notices}, (7):305--346, 1996.

\bibitem[Oh96b]{Oh_Relative_Floer_quantum_cohomology}
Yong-Geun Oh.
\newblock Relative {F}loer and quantum cohomology and the symplectic topology
  of {L}agrangian submanifolds.
\newblock In {\em Contact and symplectic geometry ({C}ambridge, 1994)},
  volume~8 of {\em Publ. Newton Inst.}, pages 201--267. Cambridge Univ. Press,
  Cambridge, 1996.

\bibitem[Oh97]{Oh_sympl_topology_action_fcnl_I}
Yong-Geun Oh.
\newblock Symplectic topology as the geometry of action functional. {I}.
  {R}elative {F}loer theory on the cotangent bundle.
\newblock {\em J. Differential Geom.}, 46(3):499--577, 1997.

\bibitem[Oh99]{Oh_sympl_topology_action_fcnl_II}
Yong-Geun Oh.
\newblock Symplectic topology as the geometry of action functional. {II}.
  {P}ants product and cohomological invariants.
\newblock {\em Comm. Anal. Geom.}, 7(1):1--54, 1999.

\bibitem[Oh05]{Oh_Construction_sp_invts_Ham_paths_closed_symp_mfds}
Yong-Geun Oh.
\newblock Construction of spectral invariants of {H}amiltonian paths on closed
  symplectic manifolds.
\newblock In {\em The breadth of symplectic and {P}oisson geometry}, volume 232
  of {\em Progr. Math.}, pages 525--570. Birkh\"auser Boston, Boston, MA, 2005.

\bibitem[PSS96]{PSS}
Sergey Piunikhin, Dietmar Salamon, and Matthias Schwarz.
\newblock Symplectic {F}loer-{D}onaldson theory and quantum cohomology.
\newblock In {\em Contact and symplectic geometry ({C}ambridge, 1994)},
  volume~8 of {\em Publ. Newton Inst.}, pages 171--200. Cambridge Univ. Press,
  Cambridge, 1996.

\bibitem[Sch]{Schwarz_PhD_thesis}
Matthias Schwarz.
\newblock {\em Cohomology operations from $S^1$-cobordisms in {F}loer
  {H}omology}.
\newblock PhD thesis, Ruhr-Universit\"at Bochum.

\bibitem[Sch93]{Schwarz_Morse_H_book}
Matthias Schwarz.
\newblock {\em Morse homology}, volume 111 of {\em Progress in Mathematics}.
\newblock Birkh\"auser Verlag, Basel, 1993.

\bibitem[Sei08]{Seidel_The_Book}
Paul Seidel.
\newblock {\em Fukaya categories and {P}icard-{L}efschetz theory}.
\newblock Zurich Lectures in Advanced Mathematics. European Mathematical
  Society (EMS), Z\"urich, 2008.

\bibitem[Sol06]{Solomon_Intersection_thry_moduli_space_holo_curves_Lag_boundary_conds}
Jake~P. Solomon.
\newblock Intersection theory on the moduli space of holomorphic curves with
  {L}agrangian boundary conditions.
\newblock 2006.
\newblock {\tt arXiv:math/0606429}.

\bibitem[Wel08]{Welschinger_Open_strings_Lag_conductors_Floer_fctr}
Jean-Yves Welschinger.
\newblock Open strings, {L}agrangian conductors and {F}loer functor.
\newblock 2008.
\newblock {\tt arXiv:0812.0276}.

\bibitem[WW15]{Wehrheim_Woodward_Orientations_pseudoholo_quilts}
Katrin Wehrheim and Chris Woodward.
\newblock Orientations for pseudoholomorphic quilts.
\newblock 2015.
\newblock {\tt arXiv:1503.07803}.

\bibitem[Zin13]{Zinger_Det_line_bundle_Fredholm_ops}
Aleksey Zinger.
\newblock The determinant line bundle for {F}redholm operators: Construction,
  properties, and classification.
\newblock 2013.
\newblock {\tt arXiv:1304.6368}.

\end{thebibliography}
\bibliographystyle{alpha}

\end{document}